\documentclass[11pt,twoside]{article}
\usepackage{latexsym}
\usepackage{amssymb,amsbsy,amsmath,amsfonts,amssymb,amscd}
\usepackage{mathrsfs}
\usepackage{epsfig, graphicx}
\usepackage{graphics,epsfig}
\usepackage{color}
\usepackage{enumerate}
\usepackage{epstopdf}
\usepackage{titlesec}
\usepackage{url} \usepackage{hyperref}
\usepackage[titletoc,toc,title]{appendix}
\usepackage{subfigure}
\usepackage[utf8]{inputenc}
        \usepackage[T1]{fontenc}
        \usepackage{lmodern}
        \usepackage{graphicx}
        \usepackage{caption}
         \usepackage{floatrow}%
\usepackage{tikz}
\usetikzlibrary{arrows,positioning, arrows}

\usepackage{amsthm}
\usepackage{amsfonts}
\usepackage{graphicx}
\usepackage{caption}
\usepackage{floatrow}
\usepackage{amsmath}
\usepackage{amssymb}
\usepackage{amsmath}
\usepackage{bm}
\usepackage{subeqnarray}
\usepackage[linesnumbered,ruled]{algorithm2e}
\usepackage[noend]{algpseudocode}
\usepackage{cases}
\theoremstyle{plain}
    \newtheorem{theorem}{Theorem}

    \newtheorem{proposition}{Proposition}

\theoremstyle{definition} 
    
    \newtheorem{remark}{Remark}

\newcommand{\half}{\frac{1}{2}}
\newcommand{\RR}{\mathbb{R}}
\newcommand{\rd}{\mathrm{d}}
\newcommand{\eps}{\mathrm{\varepsilon}}

\textheight 9in \textwidth 6.5in \numberwithin{equation}{section}

\title{To blow-up or not to blow-up for a granular kinetic equation}
\date{}

\author{Jos\'e A. Carrillo \thanks{Mathematical Institute, University of Oxford, Oxford OX2 6GG, UK. (carrillo@maths.ox.ac.uk)}, Ruiwen Shu \thanks{Department of Mathematics, University of Georgia, Athens, GA 30602. (ruiwen.shu@uga.edu)}, Li Wang\thanks{School of Mathematics, University of Minnesota, Twin cities, MN 55455, USA. (liwang@umn.edu)}, and Wuzhe Xu\thanks{Department of Mathematics and Statistics, University of Massachusetts Amherst, Amherst, MA 01003, USA. (wuzhexu@umass.edu)}}

\begin{document}

\maketitle

\begin{abstract}
    A simplified kinetic description of rapid granular media leads to a nonlocal Vlasov-type equation with a convolution integral operator that is of the same form as the continuity equations for aggregation-diffusion macroscopic dynamics. While the singular behavior of these nonlinear continuity equations is well studied in the literature, the extension to the corresponding granular kinetic equation is highly nontrivial. The main question is whether the singularity formed in velocity direction will be enhanced or mitigated by the shear in phase space due to free transport. We present a preliminary study through a meticulous numerical investigation and heuristic arguments. We have numerically developed a structure-preserving method with adaptive mesh refinement that can effectively capture potential blow-up behavior in the solution for granular kinetic equations. We have analytically constructed a finite-time blow-up infinite mass solution and discussed how this can provide insights into the finite mass scenario.
\end{abstract}

\section{Introduction}\label{sec:intro}
Granular flows are omnipresent in nature, from large scale mudslides to small scale table top experiments \cite{MeloUmbanhowarSwinney:1995,MeloUmbanhowarSwinney:1996}. Rapid granular flows or granular gases consist of a large number of small discrete grains, which interact by instantaneous hard-sphere like collisions \cite{brilliantov:2004}. The physical modeling of granular gases has been revolutionized in the past decades due to the maturity of kinetic theory. Unlike the ideal gas particles, the distinct feature of granular particles is the intrinsic inelasticity of the collisions between grains. As a result, there is a dissipation of energy, which brings a perturbed system quickly to rest. 

Deriving kinetic equations from statistical mechanics of particle systems undergoing inelastic collisions faces important issues such as the inelastic collapse \cite{mcnamara1992inelastic}, i.e. infinite many collisions in finite time. However, the kinetic description of granular gases \cite{JR85,JR852,Goldhirsch:2003} has been successful in computing transport coefficients for hydrodynamic descriptions used in experiments far from their supposed limits of validity \cite{TsimringAranson:1997,RBSS02,BMSS02,HM03,BSSP04,Bougieetal04,carrillo:2008granular}. We refer to the seminal book of Brilliantov and P\"oschel \cite{brilliantov:2004} and the recent review \cite{CHMRReview} for further physical and mathematical details and a comprehensive list of references in the subject. 

Let us consider particles moving in a one dimensioal setting. Denoting by $v,\,v_* \in \RR $ the velocities of two particles before their collision, and $v'$ and $v_*'$ the velocities after collision, we can write the postcollisional velocities in terms of the precollisional ones as 
\begin{align*}
    v' = v- \eps(v-v_*), \quad v_*' = v_* + \eps(v-v_*)
\end{align*}
for $0\leq\eps = \eps(|v-v_*|,\theta) \leq 1$. This immediately implies the energy dissipation: 
\begin{align*}
    |v'_*|^2 + |v'|^2 - |v|^2 - |v_*|^2 = -\eps(1-\eps) (v-v_*)^2 \leq 0\,,
\end{align*}
and therefore, $\eps$ controls the loss of kinetic energy and is referred as the restitution coefficient while $\theta$ is a parameter controlling the strength of inelasticity. Typical restitution coefficients are of the form $\eps(|v-v_*|,\theta)=\eps_0(1+\theta |v-v_*|^\beta)^{-1}$, with $\beta\ge -2$ and $0\leq \eps_0\leq 1$. 
Then following the same formal derivation of the Boltzmann equation from the particle dynamics, the kinetic description of the above inelastic collision takes the form \cite{MY93,DLK95,benedetto1997one,toscani2000one}
\begin{align} \label{000}
 \partial_t f + v \cdot \partial_x f = \int_{\RR} \int_{\RR_+} B(\theta) \left\{ |v'-v_*'| \frac{1}{J} f(v') f(v_*') - |v-v_*| f(v) f(v_*) \right\} \rd v_* \rd \theta \,,
\end{align}
where $f(t,x,v): (0, \infty) \times \RR \times \RR  \mapsto  \RR^+$ denotes the probability density function of grains and $J$ is the Jacobian of the transformation from $(v,v_*)$ to $(v', v_*')$, and $B$ is the rate function. 

As opposed to the Boltzmann equation whose equilibrium is the Maxwellian, the equilibrium of \eqref{000} is a Dirac mass located at the mean velocity of particles, which means that all particles are at rest in the comoving frame. This is a direct consequence of energy dissipation. If in addition to the friction, the granular material is put in interaction with a thermal bath, then a combined effect of friction and diffusion may lead to a non-Gaussian, non-Dirac equilibrium \cite{benedetto1998non, carrillo2003kinetic, gamba2004boltzmann}.

In the quasi-elastic regime, i.e., $\eps \simeq 1$, \eqref{000} can be formally approximated by \cite{MY93,DLK95,toscani2000one}
\begin{align} \label{eqn: gke}
\partial_t f + v \partial_x f &= \frac{\lambda }{2} \partial_v \left( f \int_\RR |v-w|^{\beta+1}(v-w) f(w) \rd w \right) 
 =: \frac{\lambda}{2} \partial_v ( (\partial_v W *_v f) f),,
\end{align}
where the collision now is described by a nonlocal interaction with kernel 
\begin{align} \label{kernel}
     W(|v|) = \frac{|v|^{\gamma}}{\gamma}\,,
\end{align}
where $\gamma = \beta + 3$. This is also a reminiscent of the {\it grazing collision limit} of Boltzmann operators.

The spatially homogeneous version of \eqref{eqn: gke}, that bares the name of aggregation equation (upon replacing $v$ by $x$), has been well studied in the literature \cite{Li:2004,bertozzi2009blow}. In particular, starting from a bounded initial data, there is a sharp threshold that distinguishes the finite time blow-up versus global in time solution in terms of the Osgood condition. However, by allowing a spatial dependence with a free transport dynamics, it generates substantial difficulties in understanding its solution behavior \cite{ACI15}. In essence, the lingering question is:

{ \hspace{1cm} Will the singularity formed in $v$-direction enhanced or mitigated by the shear? }

Here, we tackle this problem through a series of numerical experiments followed by a heuristic argument. As with many other physical systems, such as the Navier-Stokes equations for fluids and the Boltzmann equation for rarefied gas, a naive discretization would easily lead to unstable, physics-violated numerical solutions. In principle, this issue is addressed with the aid of numerical analysis, which provide a theoretical guarantee and practical insight to the numerical schemes. Despite that numerical analysis has its own theory and tools, it is undoubtedly that PDE analysis is the stepping stone for the development of the numerical analysis of PDEs. Said differently, one needs to have at least some well-posedness result or apriori estimate of solution for a PDE before one can simulate it numerically with confidence. Otherwise it is ambiguous to assert that the computed numerical solutions are physically relevant. 

However, when it comes to possible singular solutions, the theory is often lacking. One would instead rely on numerical solutions to give some guidance on the possible solutions \cite{huang2012asymptotics, luo2014potentially, luo2014toward}. This is the approach we are taking here. To this end, we will follow the splitting approach in \cite{agueh2016local} which  gives rise to a local in time weak solution, and our design principle of the numerical scheme is to preserve as much physical quantities as possible, including conservation and dissipation. Our main challenge to numerically solve \eqref{eqn: gke} is the possible singularity formation. For this problem, we will integrate a mesh refinement technique \cite{luo2014potentially} that rearranges a certain amount of grid points dynamically, and therefore ensures the high resolution within the possible blow-up region. A slight modification that guarantees mass conservation will also be added. 

Theoretically, for the case $\gamma=2$, we construct a family of self-similar analytic solutions with infinite total mass. The maximal time $T$ of existence for such solutions is always finite, and there is a critical threshold which determines whether the velocity support shrinks (in the super-critical case) or expands (in the sub-critical case) as $t\rightarrow T^-$. Heuristic arguments show that in the super-critical case such infinite-mass solution could be approximated by a finite-mass solution, leading to the finite-time blow-up of a finite-mass solution, while in the subcritical case a similar finite-mass solution will not have finite-time blow-up. This threshold for finite-mass solution is verified by numerical simulations. Further heuristic analysis suggests that this mechanism will not lead to a finite-time blow-up if $\gamma>2$. The case $\gamma=2$ has a special structure and global-in-time finite mass measure solutions have been constructed by optimal transport theory \cite{AC16}, where finite mass concentration is neither proved or disproved.

In section \ref{sec:numerical}, we first summarize the solution properties of the spatially homogeneous version of \eqref{eqn: gke}, as well as the physical quantities that we intend to preserve for the spatially inhomogeneous version. In the next subsection, we introduce the mesh refinement technique and then build it into a variational formulation---the Jordan-Kinderlehrer-Otto scheme. Section~\ref{sec:numerical_examples} is devoted to the numerical tests. We first verify the performance of our scheme validating our numerical solution with theoretical predictions for the spatially homogeneous case. We then carry out a series of systematic experiments to conjecture the solution behavior in the spatially inhomogeneous case\footnote{Matlab codes will be available on Github https://github.com/woodssss/Granular-kinetic-equation.}. We finally offer theoretical evidence in Section~\ref{sec: explicit} to substantiate our numerical findings. 


\section{Structure-Preserving Numerical Schemes}\label{sec:numerical}

We first discuss the main properties of solutions to the equations of interest. We focus next on the development of suitable numerical schemes keeping these properties at the discrete level.

\subsection{Solution properties}\label{sec: concentration}
Consider spatially homogeneous case of \eqref{eqn: gke}, i.e. ,
\begin{equation}\label{eqn: gke_hom}
    \partial_t f= \partial_v ( (\partial_v W *_v f) f)
\end{equation}
with radial kernel $W(v) = w(|v|)$. We first cite the following theorem \cite{bertozzi2009blow} regarding its sharp threshold dynamics.
\begin{theorem}
Starting from a bounded, compactly supported, nonnegative initial data, 
\begin{itemize}
    \item if $\frac{w'(r)}{r}$ is monotone decreasing, $w''(r)>0$ and
    \begin{equation} \label{osgood1}
        \int_0^1 \frac{1}{w'(r)} \rd r < +\infty\,,
    \end{equation}
    then the solution to \eqref{eqn: gke_hom} blows up in finite time. The blow-up time depends only on the initial data through its radius of support and total mass;
    \item if $\frac{w'(r)}{r}$ is monotone decreasing, $w''(r)>0$ and
    \begin{equation} \label{osgood2}
        \int_0^1 \frac{1}{w'(r)} \rd r = +\infty\,,
    \end{equation}
   then the solution to \eqref{eqn: gke_hom} stays bounded for all time and converges as $t\rightarrow \infty$ to a Dirac centered at the center of mass of the initial data. 
\end{itemize}
\end{theorem}
The conditions \eqref{osgood1} and \eqref{osgood2}  
stem from the Osgood condition for well-posedness of the ODE $\frac{dR}{dt} = -w'(R)$, with the former violating the condition and the latter satisfying it. This theorem, applied to the kernel \eqref{kernel}, implies that: infinite time blow-up for $\gamma> 2$ and finite time blow-up for $-1 < \gamma < 2$, see \cite{Li:2004} for further properties. The case of $\gamma=2$ leads to infinite time blow-up as it can be directly checked since it is a linear equation once conservation of mass and momentum are taken into account.

More particularly, when $\gamma = 3$ in \eqref{kernel}, with a smooth enough initial condition centered at origin with a compact support, the solution as $t\rightarrow \infty$  takes the asymptotic form $f(t,v) \sim \frac{1}{2} \delta_{\frac{1}{t}} + \frac{1}{2} \delta_{-\frac{1}{t}}$. Therefore, when $t\rightarrow \infty$, the solution tends to form two symmetric Dirac delta functions, whose centers are converging to $v=0$ at rate $\frac{1}{t}$.
More details can be found \cite{benedetto1997kinetic, toscani2000one, caglioti2002homogeneous, Li:2004}.

When $\gamma = 1$, we have explicit prediction for the blow-up time since  $\eqref{eqn: gke_hom}$ can be viewed as the derivative of the Burgers equation \cite{BCDP15}. Indeed,  $W = |v|$ and $\partial_v W = \text{sign}(v)$, $\partial^2_{v} W = 2 \delta(v)$. By setting $u(t,v) = - \int_{-\infty}^v f(t,y) \rd y$, we then have 
\[
\partial_t u + 2 u \partial_v u = 0
\]
with initial condition
\[
u_0(v) = - \int_{-\infty}^v f_0(y) \rd y\,.
\]
Then the blow-up time is equal to the shock formation time of Burgers equation, i.e.
\begin{equation*} 
T = - \frac{1}{ 2 \min_v u_0^{'} (v)}\,.
\end{equation*}

Next we cite a few apriori properties of the solutions to the spatially inhomogeneous case proposed in \cite{agueh2016local,ACI15}.
\begin{proposition}\label{prop:properties}
For $\gamma >1$, the solution to \eqref{eqn: gke} satisfies 
\begin{itemize}
    \item \text{Mass conservation}: \quad
    $ \displaystyle \frac{d}{ dt} \int_{\RR \times \RR} f(t,x,v) \rd x \rd v  =  0$.

    \item \text{Momentum conservation} : \quad
    $ \displaystyle \frac{d}{ dt} \int_{\RR \times \RR} v f(t,x,v) \rd x \rd v  =  0$.
    
    \item \text{Decrease of moments for $p \geq 2$} : \quad
    $ \displaystyle \frac{d}{ dt} \int_{\RR \times \RR} |v|^p f(t,x,v) \rd x \rd v  \leq  0$. 
    
    \item \text{Increase of internal energy} : \quad for a $C^1$ convex function $U$ with $U(0) = 0$, we have 
    \[
    \frac{d}{ dt} \int_{\RR \times \RR} U(f(t,x,v)) \rd x \rd v \geq 0\,.
    \]

\end{itemize}
\end{proposition}
In this proposition, the third part, when $p=2$, signifies the dissipation of energy, as expected in an inelastic collision. When $p = \infty$, it implies that for compactly supported initial data $f(0,x,\cdot)$,  the support of $f(t,x,\cdot)$ remains compact for all $t>0$. Moving to the fourth part, a notable selection is $U(f) = f \ln f$, representing the negative of the entropy.

\subsection{Numerical Method}

In this subsection, we present a variational semi-Lagrangian scheme for the granular kinetic equation  \eqref{eqn: gke}. As demonstrated in the theoretical paper \cite{agueh2016local}, we divide the equation into two stages: a transport step that is explicitly addressed, and a collision step that is treated implicitly.
More precisely, denote $f^n(x,v) \approx f(t_n, x, v)$, where $t_n = n \Delta t$ with $\Delta t$ being the time step, we have: 
\begin{equation}\label{eqn:transp}
\frac{f^{n+\half}-f^n}{\Delta t} + v \partial_x f^n = 0\,,
\end{equation}
and
\begin{equation}\label{eqn:collision}
\frac{f^{n+1}-f^{n+\half}}{\Delta t} = \partial_v((\partial_v W *_v f^{n+1}) f^{n+1})\,.
\end{equation}
Periodic boundary condition in $x$ and no flux boundary condition in $v$ will be used throughout the paper. 

\subsection{Semi-Lagrangian scheme for the transport step} \label{sec:SL}
For the transport step \eqref{eqn:transp}, we will use the semi-Lagrangian scheme, which combines the method of characteristic and interpolation. 
Specifically, let $[-L_x, L_x]$ be the computational domain and define 
\[
 x_i^0 = -L_x + i \Delta x/2\,, \qquad  1\leq i \leq N_x, \quad \Delta x= 2L_x/N_x\,,
\]
as the initial grid points. As time progresses, the grids $\{ x^n_i \}$ will be updated and become nonuniform due to mesh refinement (details of which will be provided in Section~\ref{sec:mr}). Given $\{  f^{n}(x^{n}_i) \}_{i=1}^{N_x}$, $\{  f^{n+1}(x^{n+1}_i) \}$ can be acquired by tracing back the characteristics:
\[
f^{n+1}(x^{n+1}_i, v) = f^{n}(x^{n+1}_i - v \Delta t, v)\,,
\]
where $ f^{n}(x^{n+1}_i - v \Delta t, v)$ is evaluated using piecewise cubic Hermite interpolation.  To guarantee mass conservation, we re-weight the interpolated value by the total mass, namely, 
\begin{equation}\label{eqn:semi_update}
    { f^{n+1}(x^{n+1}_i, v) } := \max \{{ f^{n}(x^{n+1}_i - v \Delta t, v) } , 10^{-10}\}\frac{M_n}{M_{n+1}}\,,
\end{equation}
where 
\begin{align*}
M_n = \sum_{i = 2}^{N_x-1} f^{n}(x^n_i, v) \frac{x^n_{i+1}-x^n_{i-1}}{2} + f^{n}(x^n_1, v) (L_x + \frac{x^n_2+x^n_1}{2})  + f^n (x^n_{N_x},v) (L_x- \frac{x^n_{N_x}+x^n_{N_x-1}}{2})\,,
\end{align*}
and let $y_i := x^{n+1}_i-v \Delta t $
\begin{align*} 
    \displaystyle M_{n+1} = \sum_{i = 2}^{N_x-1} f^{n+1}(y_i, v) \frac{y_{i+1}-y_{i-1}}{2} + f^{n+1}(y^n_1, v) (L_x + \frac{y_2+y_1}{2}) + f^{n+1} (y_{N_x},v) (L_x- \frac{y_{N_x}+y_{N_x-1}}{2})\,,
\end{align*}
where $f^{n+1}(y_i, v)$ is defined in \eqref{eqn:semi_update}.

\subsection{Regularized JKO scheme with adaptive mesh refinement for collision step}
This section centers around solving the collision step \eqref{eqn:collision},
for which we will adopt a regularized JKO scheme along with a mesh refinement technique. In the following, we first introduce the Fisher information regularized JKO scheme that was proposed in \cite{li2020fisher}, and extend it to the non-uniform mesh case here. Secondly, in order to increase the resolution at the blow-up region, we use a mesh refinement technique that was originally proposed in \cite{luo2014potentially}, which adaptively redistribute a  large portion of the grid points to the blow-up region.

\subsubsection{Regularized JKO scheme}
To start, rewrite \eqref{eqn: gke_hom} as 
\begin{equation}\label{eqn:hom_energy}
\partial_t f = \partial_v (f \partial_v \frac{\delta E(f)}{\delta f})\,,
\end{equation}
where the Energy functional is defined as 
\[
E(f) = \half \int_{\RR} \int_{\RR} W(v-w)f(v)f(w) \rd v \rd w \,,
\]
and $\frac{\delta E}{\delta f}$ is the functional derivative of $E$. 
Then $\eqref{eqn:hom_energy}$ can be interpreted as a gradient flow of the $E$ with respect to the Wasserstein metric, and therefore admits the following minimizing movement scheme, also called JKO scheme: 
\begin{equation}\label{eqn:JKO1}
f^{n+1} \in \underset{f \in \mathcal{P}_{a c}\left(\Omega_{v}\right)}{\operatorname{argmin}}\left\{\frac{1}{2} d_{\mathcal{W}}\left(f, f^{n+\half}\right)^{2}+\Delta t E\left(f \right)\right\}\,.
\end{equation}
Here $d_{\mathcal{W}}\left(f, f^n \right)$ is the Wasserstein distance between $f$ and $f^n$ and $\mathcal{P}_{a c}\left(\Omega_{v}\right)$ is the set of probability measures on $\Omega_v$ that are absolutely continuous with respect to Lebesgue measure. In this paper, we use the Benamou-Brenier’s dynamic formulation to numerically compute the
Wasserstein distance. In particular, the distance between two measures $f_0 \rd x$ and $f_1 \rd x$ is
\begin{equation}\label{eqn:waser_dis}
d_{\mathcal{W}}\left(f_{0}, f_{1}\right)^{2}=\min _{(f, m) \in \mathcal{C}_{1}} \int_{0}^{1} \int_{\Omega_{v}} \Phi(f(t, v),\|m(t, v)\|) \mathrm{d} v \mathrm{~d} t\,,
\end{equation}
where
\[
\Phi(f,\|m\|)=\left\{\begin{array}{cl}
\frac{\|m\|^{2}}{f} & \text { if } f>0\,, \\
0 & \text { if }(f, m)=(0,0)\,, \\
+\infty & \text { otherwise }\,,
\end{array}\right.
\]
and the constraint set $\mathcal{C}_{1}$ consists of 
\[
\begin{array}{l}
\partial_{t} f+\nabla_{v} \cdot m=0 \text { on } \Omega_{v} \times[0,1], \quad m \cdot \nu=0 \text { on } \partial \Omega_{v} \times[0,1], \\
f(\cdot, 0)=f_{0}, \quad f(\cdot, 1)=f_{1} \text { on } \Omega_{v}\,.
\end{array}
\]
Here $\nu$ is the outer normal direction of  $\Omega_{v}$. Substituting $\eqref{eqn:waser_dis}$ into $\eqref{eqn:JKO1}$, and write $m = fv$, we end up with a convex optimization with a linear constraint:
\begin{equation} \label{eqn:JKO_final}
    \left\{\begin{array}{l}
(f, m)= \displaystyle \arg \min_{f, m} \half \int_{0}^{1} \int_{\Omega_v} \Phi(f,\|m\|) \rd v \rd t+ \Delta t E(f(1, \cdot)) \\[3mm]
\text { s.t. } \quad \partial_{t} f +\nabla_v \cdot m=0, f(0, v)= f^{n}(v), ~~m \cdot \nu=0\,.
\end{array}\right.
\end{equation}

 In practice, to further increase the convexity of the objective function and preserve the positivity of the desired minimizer $f^*$, we add a Fisher information regularization term 
 $$
 \beta^{-2} \Delta t^{2} \int_{\Omega_v}\|\nabla \log f \|^{2} f \mathrm{d} v
 $$
to the objective function and reach the following regularized problem, originally proposed in \cite{li2020fisher}:
\begin{equation}\label{eqn:JKO_fisher}
f^{n+1}(v) \in \arg \inf _{m, f}  \int_{\Omega_v} \frac{\|m(t)\|^{2}}{f(t)}+\beta^{-2} \Delta t^{2}\|\nabla \log f(v)\|^{2} f(v)  \rd v+2 \Delta t E(f)\,,
\end{equation}
such that
\begin{equation}\label{eqn:JKO_fisher_const}
f- f^{n+\half}+\nabla_v \cdot m=0,~~ m \cdot \nu=0\,.
\end{equation}
Here the artificial time in the dynamic formulation \eqref{eqn:JKO_final} is replaced by one step finite difference in time. This is shown in \cite{li2020fisher} that significantly reduces the computational cost without violating the first order accuracy in time.

\subsubsection{Non-uniform velocity discretization}
Different from \cite{li2020fisher}, we will use a non-uniform discretization in $v$ so as to facilitate the later mesh refinement proposal. To be more specific, denote the 
computational domain in $v$ to be $v \in [-L,L]$ and the number of grid points as $N_v$. Let $\{v_i\}_{i=1}^{N_v}$ be the grid points, and define
\[
\Delta v_i := v_i -v_{i-1}, \quad
h_i = \half (\Delta v_i + \Delta v_{i+1})\,.
\]

 To discretize $\eqref{eqn:JKO_fisher}$ and $\eqref{eqn:JKO_fisher_const}$, we evaluate $f$ at the grid points and $m$ at the half grid points, i.e., i.e. $m_{i-\half} = m(\frac{v_i+v_{i-1}}{2})$. Then the fully discretized version of $\eqref{eqn:JKO_fisher}$ and $\eqref{eqn:JKO_fisher_const}$ reads
\begin{equation}\label{eqn:non-uniform JKO}
f^{n+1} \in \arg \min 
F (f,m)
+
 \Delta t E (f,m) \,,
\end{equation}
where
\[
F (f,m) = \sum_{j=2}^{N_{v}}\left[\frac{2 m_{j-\half}^{2}}{f_{j}+f_{j-1}}+\frac{\beta^{-2} \Delta t^{2}}{\Delta v_i^{2}}\left(\log f_{j}-\log f_{j-1}\right)^{2} \frac{f_{j}+f_{j-1}}{2}\right] \Delta v_j\,,
\]
and 
\[
E (f) = \sum_{i, l=1}^{N_v} W_{i, l} f_{i} f_{l} h_i h_l \,, \quad W_{i,l}=W(v_i - v_l)\,.
\]
The constraint function reads
\begin{equation}\label{eqn:linear_constraint}
f^{n+1}_j -f^n_j + \frac{m_{j+\half}-m_{j-\half}}{h_j} = 0.
\end{equation}

Denote $u = [f, m]$, then the constraint $\eqref{eqn:linear_constraint}$ can be reformulated as a linear system
$A u = b$.
By introducing the indicator function,  $\eqref{eqn:non-uniform JKO}$ together with $\eqref{eqn:linear_constraint}$ can be reformulated into an unconstrained optimization problem: 
\begin{equation}\label{eqn:uncon_opt}
\min _{u} J(u)+\chi(u), \quad \chi(u)=\left\{\begin{array}{cc}
0 & \mathrm{~A} u=b\,, \\
+\infty & \text { otherwise }\,,
\end{array}\right.
\end{equation}
where
\begin{equation} \label{JJ}
    J(u) := F (u)
+
 \Delta t  E (u)\,.
\end{equation}

As discussed in \cite{li2020fisher}, thanks to the convexity of the $J(u)$, we can invoke the sequential quadratic programming to solve $\eqref{eqn:uncon_opt}$:
\begin{equation*}
\left\{\begin{array}{l}
z^{(l+1)} \in \arg \min _{z} \frac{1}{2}\left(z-u^{(l)}\right)^{T} \mathsf{H}^{(l)}\left(z-u^{(l)}\right)+\nabla J\left(u^{(l)}\right)^{T}\left(z-u^{(l)}\right)+\chi(z), \\
u^{(l+1)}=u^{(l)}+t_{l}\left(z^{(l+1)}-u^{(l)}\right)
\end{array}\right.\,.
\end{equation*}
Here $\mathsf{H}^{(l)}$ is an approximation of $\nabla^2 J(u^{(l)})$, for our particular form of $J$ in \eqref{JJ}, it has the form
\begin{equation*}
(\mathsf{H}^{(l)})_{i,j} =
    \begin{cases}
       (\nabla^2 F(u^{(l)}))_{i,i}\,, \hfill  \quad \text{if $i=j$}, \\
       0, \hfill  \quad \text{else.}
    \end{cases}
\end{equation*}

We now summarize our one step regularized JKO scheme in Algorithm~\ref{alg:non_uni_JKO}, which is essentially the same as in \cite{li2020fisher}. Note that when the exit flag as defined in Algorithm~\ref{alg:non_uni_JKO} becomes $e_f = 1$, it signifies that the minimizer cannot be reached even after a sufficiently large number of iterations. This condition can be regarded as a numerical indicator of blow-up, and we will demonstrate its significance in the forthcoming numerical examples section.

\begin{algorithm}[h!] 
\SetKwInOut{Input}{Input}
\SetKwInOut{Output}{Output}
\caption{One step regularized JKO scheme }\label{alg:non_uni_JKO}
\Input{Grid points $\{ v_j \}_j$, function value $\{  f^n(v_j) =:f_j^n \}_j$, time step $\Delta t$, optimization step $\lambda$, max iteration number $I_{max}$ and a stopping criteria $\epsilon$}
\Output{$\{ f^{n+1}_j \}_j$ and exit flag $e_f$}
1. Compute $\Delta v_i := v_i -v_{i-1}$ and $h_i = \frac{\Delta v_i + \Delta v_{i+1}}{2}$.\\
2. Let $l = 0$, $\zeta = 1$, $e_f=0$ \;
\While{$l<I_{max}$}{
  \eIf{$\zeta > \epsilon$}{
    $z^{(l+1)} \in \arg \min _{z} \frac{1}{2} (z-u^{(l)} )^{T} \mathsf{H}^{(l)} (z-u^{(l)} )+\nabla J (u^{(l)} )^{T} (z-u^{(l)} )+\chi(z)$ \;
    $u^{(l+1)}=u^{(l)}+\lambda (z^{(l+1)}-u^{(l)} )$ \;
    $\zeta = \frac{\| u^{(l+1)} - u^{l}\|}{\|  u^{l} \|}$ \;
  }{Break \;}
     $l = l+1$ \;
}
3. $e_f = 1$ \textbf{if} $l=I_{max}$\;
\end{algorithm}

\subsubsection{Adaptive mesh refinement}\label{sec:mr}
To investigate the possible blow-up in the solution, we need to keep increasing the resolution in regions where the solution is concentrated. Here we follow the approach developed in \cite{luo2014potentially}. The main concept is to create a dynamic mapping between the original domain and the interval $[0,1]^d$. This mapping ensures that the grid points in $[0,1]^d$ remain uniformly spaced and finite. However, when transferring the grid back to the original domain, the points are concentrated more heavily in the vicinity of singularities. 

To better illustrate the idea, let us consider $\eqref{eqn: gke_hom}$ with $v \in [-L, L]$ and initial condition $f(0,v) = f_0(v)$. Choose two paramters $\delta_0 $ and $\delta  $, both in $(0,1)$, and let
\begin{align} \label{S0}
   S_0: = \{v : |f^0(v)| \geq \delta_0 \| f^0 \|_{\infty} \}\,, 
\end{align}
 which indicates the blow-up region. Now we intend to find a one-to-one map $v = \mu(s)$ between $v \in [-L,L]$ and an auxiliary variable $s \in [-1,1]$. Specifically, we seek a mapping such that if we distribute $N$ uniform grid points along $s$, then the resulting grid in $v$ contains $\delta N$ points within $S_0$. By doing so, we can densely pack grid points in the neighbor of the concentration region, i.e., $S_0$, while maintaining a sparser grid elsewhere. Note that the choice of the mapping function $\mu$ highly depends on the shape of $S_0$ in $\eqref{S0}$, and we will focus on two specific cases. 

The first case is when the concentration is symmetrically centered at a single point. Without loss of generality, let us assume it is centered at $v=0$. In this case, we only need to find $\mu(s): [0,1] \rightarrow [0,L]$  and subsequently perform an odd extension to extend the mapping to the entire domain $s \in [-1,1]$. In particular, 
define
\begin{equation}\label{eqn:Sr}
    S = \{ v_j :f(v_j)>\delta_0 (\max_i{f(v_i)}) , ~ v_j >0\}\,, \quad 
r = \min \{v_j: v_j \in S \}\,,
\end{equation}
then we seek a continuous and monotone increasing  function $v = \mu(s)$, such that 
\[
\mu(0) = 0, ~~\mu(\delta) = r, ~~\mu(1) = L\,.
\]
To this end, we use a straight line to connect $s=0$ to $s=\delta$, and a concave up curve to connect $s=\delta$ to $s=1$ if $\frac{r}{\delta}<L$; or a concave down curve otherwise. More specifically, $\mu$ takes the form
\begin{align} \label{mu1}
    \mu(s) = \left\{ \begin{array}{cc} \frac{r}{\delta} s & s \in [0, \delta] \,;\\ [2mm]
    a_1s^5 + b_1s & s \in [\delta, 1]\,, \quad\text{if}~ r/\delta < L\,; \\ [2mm]
    \frac{a_2}{1+e^{-b_2s}}  & s \in [\delta, 1]\,,  \quad \text{if}~ r/\delta \geq L\,, 
    \end{array} \right.
\end{align}
where $a_{i}$ and $b_i$, $i=1,2$ are two constants depending on $r$ and $\delta$ to ensure continuity of $\mu(s)$ at the interface $s=\delta$. 

In the second scenario, the concentrations are symmetrically positioned at two points. Once more, without loss of generality, we assume that the centers of these two concentrations are at $v=0$. Similar to the previous case, we focus on the mapping in the positive part of the domain and perform an odd extension to cover the entire domain. Define 
\[
S = \{ v_j :f(v_j)>\delta_0 (\max_i{f(v_i)}), ~ v_j>0 \}, \quad
r_1 = \min \{v_j: v_j \in S \},~
r_2 = \max \{v_j: v_j \in S \}\,.
\]
Then we aim to find a continuous and monotone increasing function $\mu(s)$ such that 
\begin{align} \label{rrr}
 \mu(0) = 0, ~~\mu(\half - \frac{\delta}{2}) = r_1,~~ \mu(\half + \frac{\delta}{2}) = r_2, ~~
\mu(1) = L\,.   
\end{align}
As with the previous case, we use a straight line to connect $s = \half - \frac{\delta}{2}$ to $  s= \half + \frac{\delta}{2}$, and use either $\mu(s) = a s^5 + b s$ or $\mu(s) = \frac{a}{1+e^{-bs}}$ to connect  $s=0$ to 
$s= \half - \frac{\delta}{2}$, and $s = \half + \frac{\delta}{2}$ to 
$s=1$, depending on the concavity.  To be more precise, the expression for $\mu(s)$ takes the form:
\begin{align} \label{mu2}
    \mu(s) = \left\{ \begin{array}{cc}
    \frac{r}{\delta} s & s \in [\frac{1}{2} - \frac{\delta}{2}, \frac{1}{2} + \frac{\delta}{2}]\,, \\[2mm]
    a_1 s^5 + b_1 s & s \in [ 0,  \frac{1}{2} - \frac{\delta}{2}]   \quad \text{if}~ \frac{r2-r1}{\delta} \geq L \,, \\[2mm] 
    \frac{a_2}{1+e^{-b_2s}}  & s \in [ 0,  \frac{1}{2} + \frac{\delta}{2}]   \quad \text{if}~ \frac{r2-r1}{\delta} < L\,,  \\[2mm]
    a_3 s^5 + b_3 s & s \in [  \frac{1}{2} - \frac{\delta}{2}, 1]   \quad \text{if}~ \frac{r2-r1}{\delta} < L \,, \\[2mm] 
    \frac{a_4}{1+e^{-b_4s}}  & s \in [\frac{1}{2} + \frac{\delta}{2}, 1]   \quad \text{if}~ \frac{r2-r1}{\delta} \geq L\,, 
    \end{array} \right.
\end{align}
where $a_{i}$ and $b_i$, $i=1,2, 3, 4$ are  constants depending on $r_1$, $r_2$ and $\delta$ to ensure continuity of $\mu(s)$ at the interfaces.

\begin{figure}[h!]
\centering
{\includegraphics[width=0.32\textwidth]{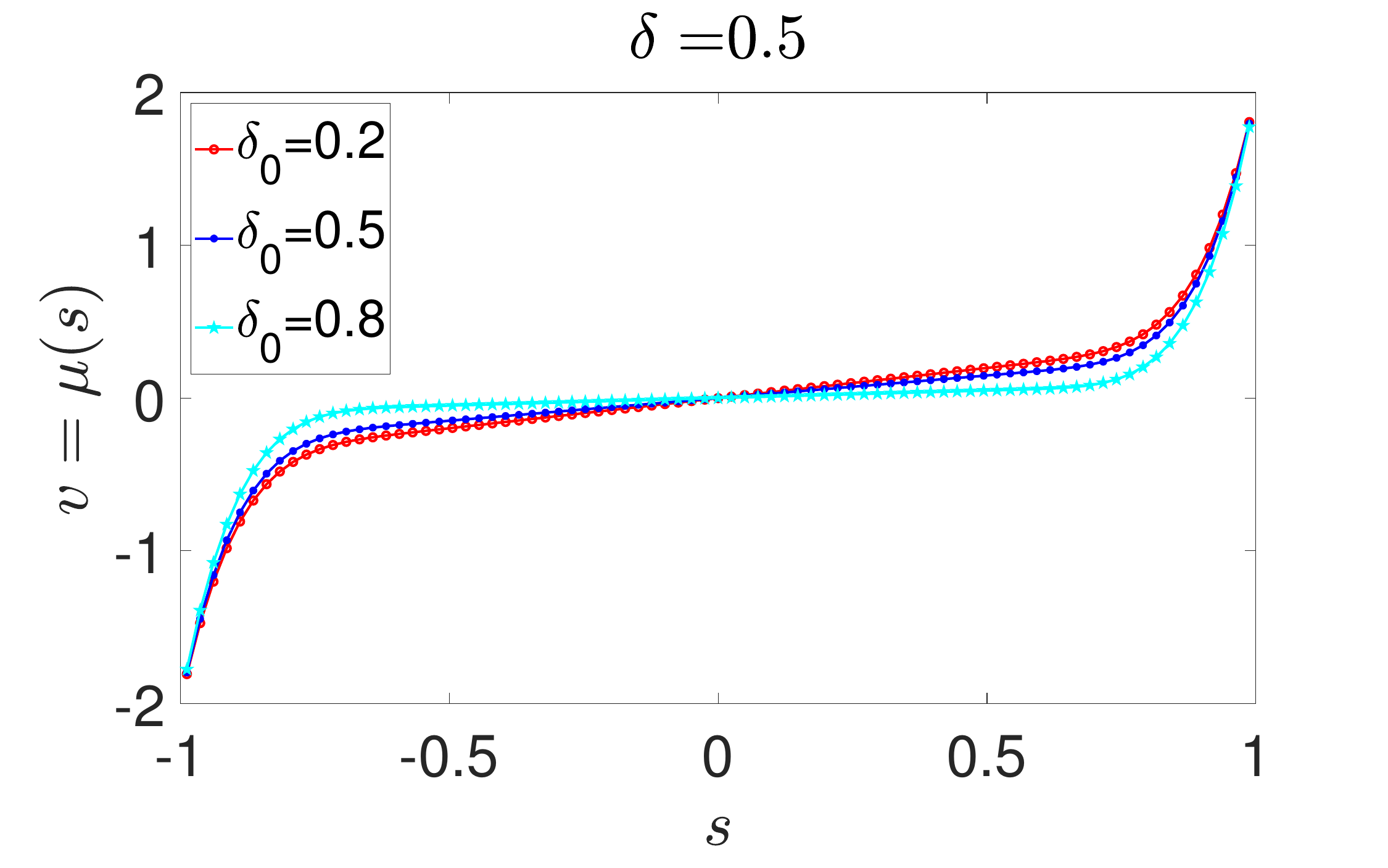}}
{\includegraphics[width=0.32\textwidth]{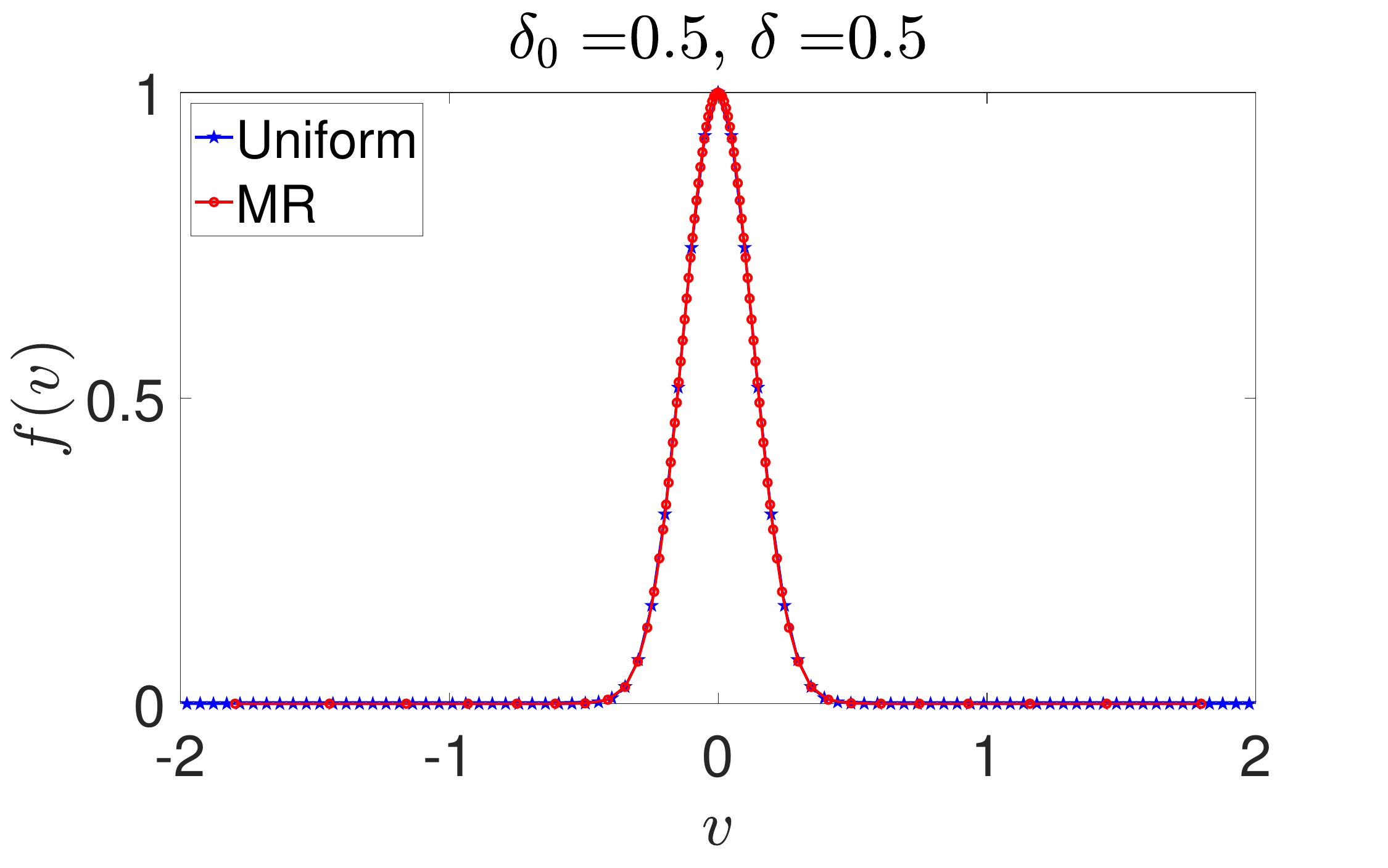}}
{\includegraphics[width=0.32\textwidth]{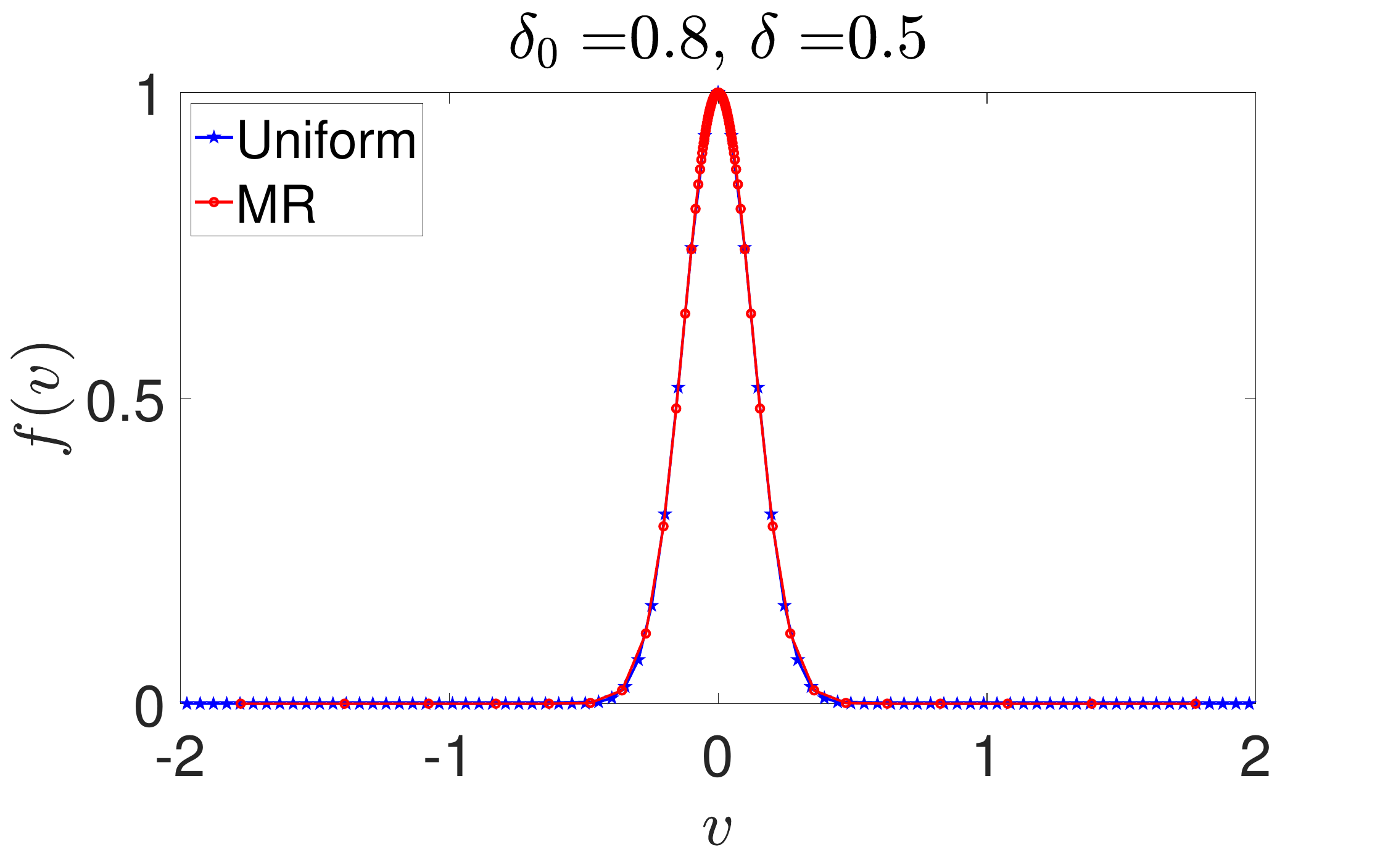}}
 \caption{Demonstration of Mesh refinement method for $f(v) = e^{-30v^2}$ for $v \in [-2, 2]$. The left figure shows the mapping function $\mu(s)$ with fixed $\delta=0.5$ and various $\delta_0=0.2, 0.5, 0.8$. The middle and right one compare function $f(v)$ in uniform grid and after mesh refinement (MR). }
 \label{fig:mr_demon}
\end{figure}

To visualize the mapping, we examine a simple function $f(v) = e^{-20v^2}$ with $v \in [-2, 2]$. We plot this function using both a uniform mesh and a non-uniform mesh. As illustrated in Fig~\ref{fig:mr_demon}, the mesh refinement approach results in higher resolution near the function's peak.

\begin{figure}[h!]
\centering
{\includegraphics[width=0.32\textwidth]{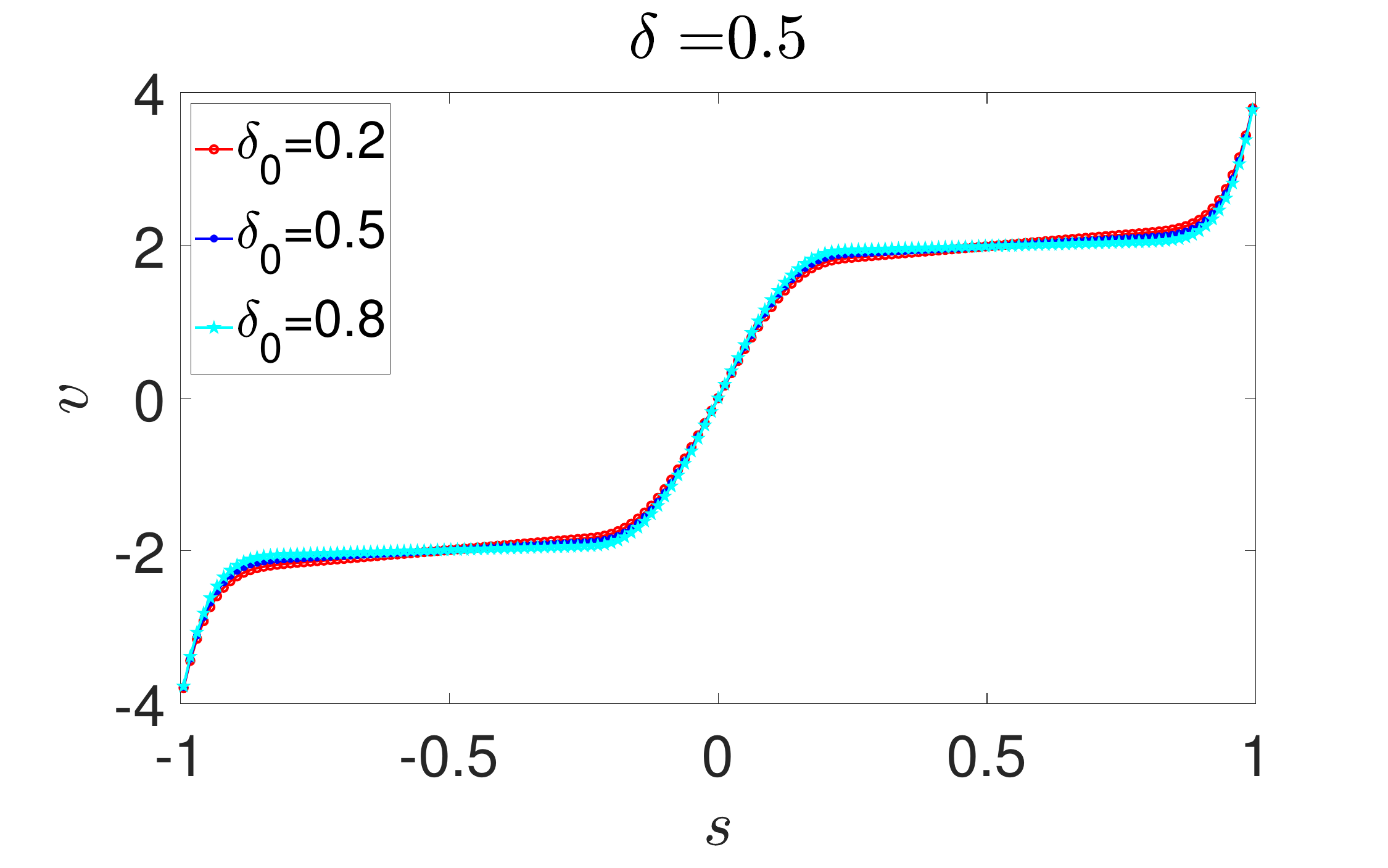}}
{\includegraphics[width=0.32\textwidth]{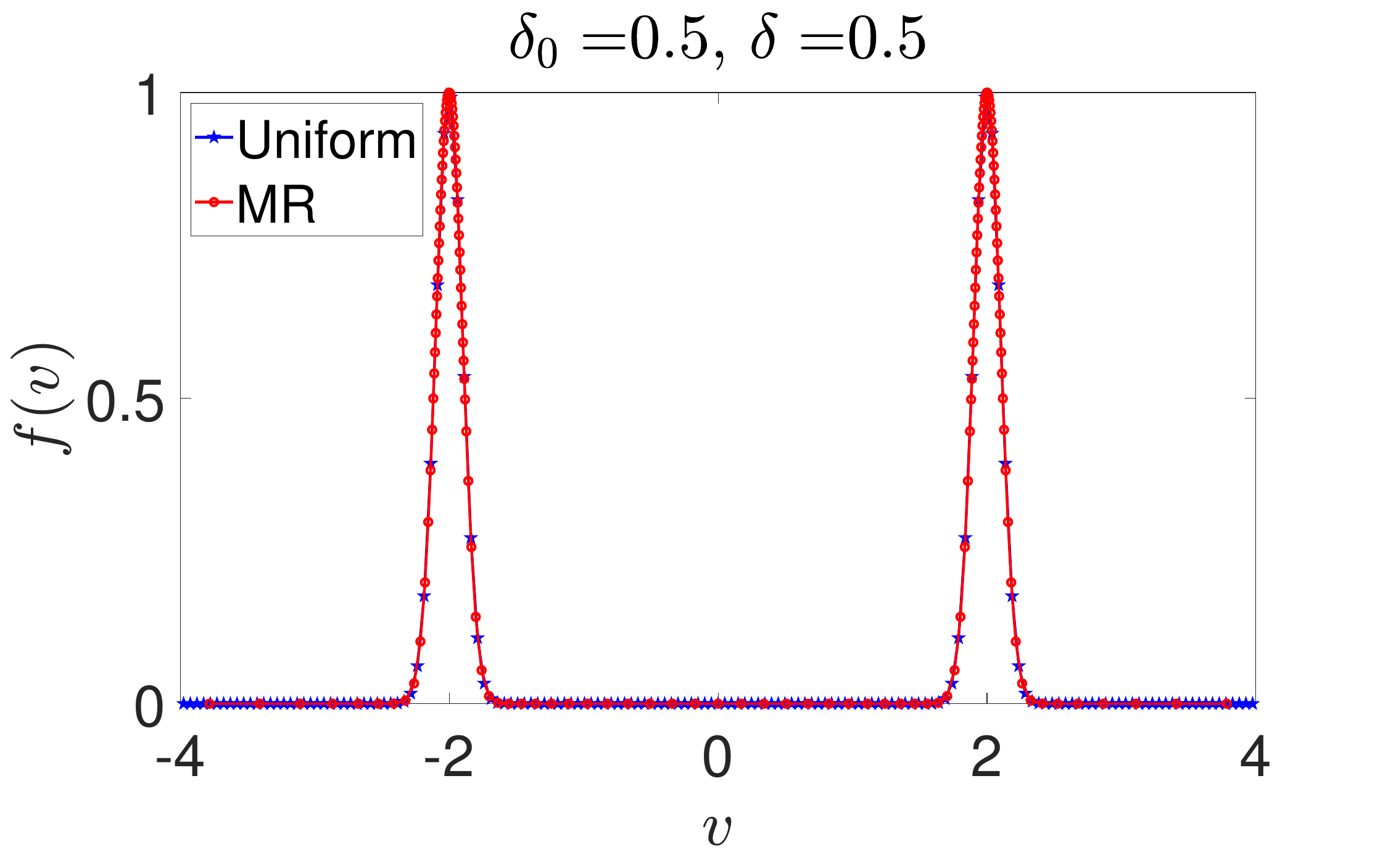}}
{\includegraphics[width=0.32\textwidth]{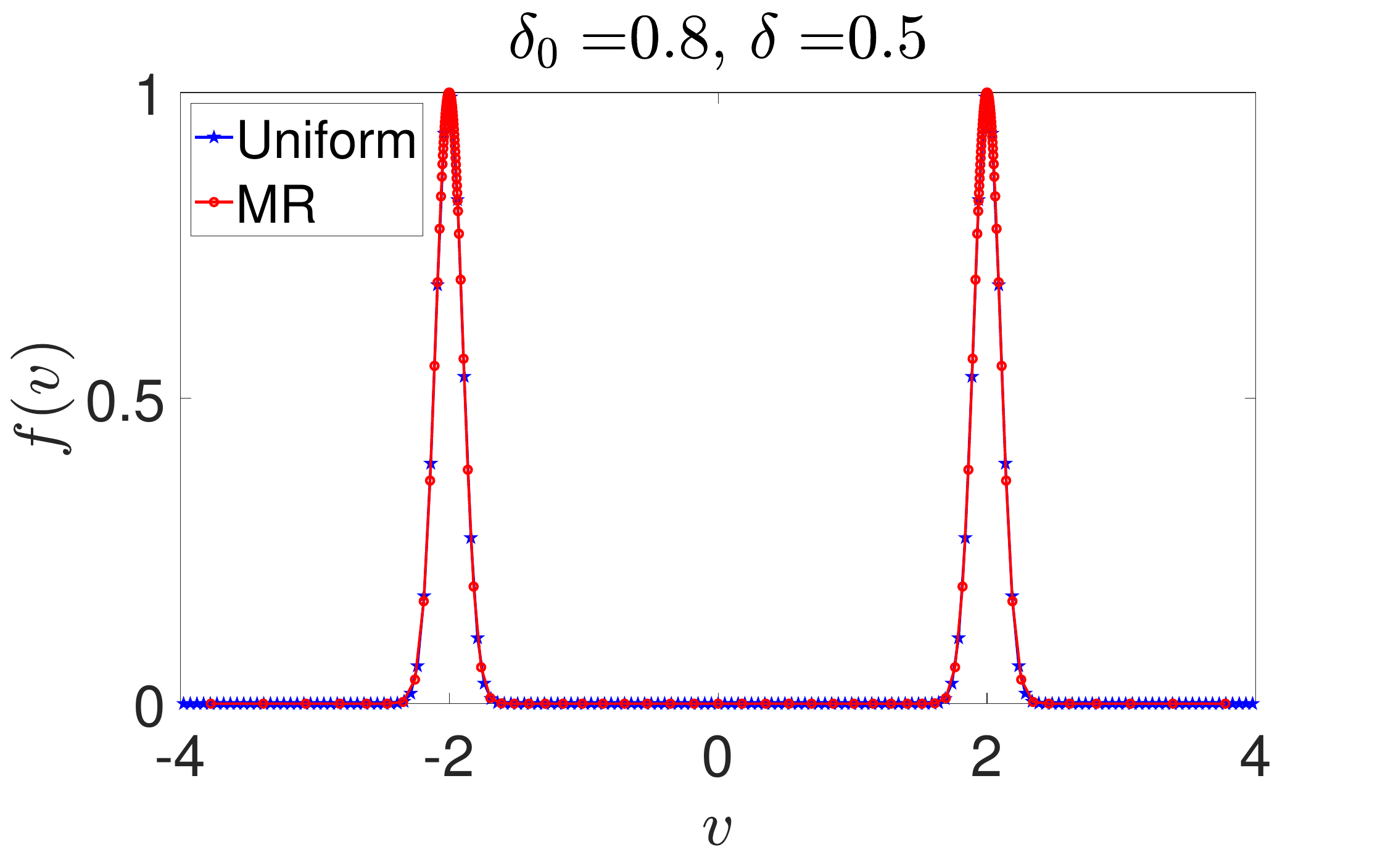}}
 \caption{Demonstration of Mesh refinement method for $f(v) = e^{-50(v-2)^2 - 50(v+2)^2}$ for $v \in [-4, 4]$. The left figure shows the mapping function $\mu(s)$ with fixed $\delta=0.5$ and various $\delta_0=0.2, 0.5, 0.8$. The middle and right one compare function $f(v)$ in uniform grid and after mesh refinement(MR).}
 \label{fig:mr_demon_mu2}
\end{figure}

Similarly, for the second case, we plot the function  $f(v) = e^{-50(v-2)^2 - 50(v+2)^2}$ for $v \in [-4, 4]$ in Fig.~\ref{fig:mr_demon_mu2}. Once more, the mesh refinement approach provides higher resolution near the function's peaks.

In practice, given a mesh ${ v^n_j }$ and corresponding function values ${ f^n(v^n_j) }$, we can construct an appropriate mapping $\mu^n(s)$ as described above to generate a new set of mesh grid ${ v^{n+1}_j }$ with higher resolution in the desired region. Note that our mesh refinement method differs somewhat from the approach proposed in \cite{luo2014toward}, where the author employed a Gaussian function as the mapping function, along with three adjustable parameters. Our proposed mesh refinement technique is more general and comprises an arbitrary number of components. This enables easy manipulation of the mapping's shape, including its monotonicity and concavity. Our approach is not only more flexible but also capable of handling multiple blow-ups concurrently.

When interpolation is applied, a rescaling is imposed to ensure exact mass conservation. Specifically, after obtaining $f(v_j^{n+1})$, the following procedure is carried out:
\begin{equation}\label{eqn:mass_conserv_v}
    { f(v_j^{n+1}) } := \max \{{ f(v_j^{n+1}) , 10^{-10}}\}\frac{M_n}{M_{n+1}}\,,
\end{equation}
where
\begin{align*}
M_n &= \sum_{j = 2}^{N_v-1} f(v^{n}_j) \frac{v^{n}_{j+1}-v^{n}_{j-1}}{2} + f(v^{n}_j) (L_v + \frac{v^{n}_2+v^{n}_1}{2})  + f (v^{n}_{N_v}) (L_v- \frac{v^{n}_{N_v}+v^{n}_{N_v-1}}{2})\,,
\\ M_{n+1} &= \sum_{j = 2}^{N_v-1} f(v^{n+1}_j) \frac{v^{n+1}_{j+1}-v^{n+1}_{j-1}}{2} + f(v^{n+1}_j) (L_v + \frac{v^{n+1}_2+v^{n+1}_1}{2})  + f (v^{n+1}_{N_v}) (L_v- \frac{v^{n+1}_{N_v}+v^{n+1}_{N_v-1}}{2})\,.
\end{align*}

\subsubsection{Mesh refinement JKO scheme}
Combining both the mesh refinement mechanism and the nonuniform JKO scheme results in a mesh refinement JKO scheme. First, we outline the following algorithm for one step of the JKO scheme:

\begin{algorithm}[H]
\SetAlgoLined
\SetKwInOut{Input}{Input}
\SetKwInOut{Output}{Output}
\Input{Parameters $\delta_0, \delta \in (0,1)$ and a uniform mesh $\{s_j\}_j$, grid points  $\{  v^n_j \}_j$, function value $\{f^n(v^n_j)\}_j$, time step size $\Delta t$ and threshold value $\epsilon_v$.}
\Output{updated grid 
$\{  v^{n+1}_j \}$, function value $f^{n+1}(v^{n+1}_j)$ and blow-up indicator $e_f$.}
1. Acquire $f^{n+1}(v^{n}_j)$ and $e_f$ from Algorithm~\ref{alg:non_uni_JKO}. \\
2. Determine $r$ according to \eqref{eqn:Sr} and construct $\mu$ via \eqref{mu1}; or determine $r_1$ and $r_2$ according to \eqref{rrr} and construction $\mu$ via \eqref{mu2}. Update mesh by setting $v_j^{n+1} = \mu(s_j)$. \\
3. Compute $f^{n+1}(v^{n+1}_j)$ using the interpolation and rescaling \eqref{eqn:mass_conserv_v}. 
  \caption{One step mesh refinement JKO scheme}
\label{alg:one_step_MR_JKO}
\end{algorithm}

\

For the evolution problem, we execute Algorithm~\ref{alg:one_step_MR_JKO} with two time stepping strategies. One is to use fixed time step and the other one is the adaptive time stepping. The latter choice is more expensive but useful in predicting the more accurate finite blow-up time. The algorithm is summarized in Algoirthm~\ref{alg:MR_JKO_ver2}.

\begin{algorithm}[h!]
\SetAlgoLined
\SetKwInOut{Input}{Input}
\SetKwInOut{Output}{Output}
\Input{Parameters $\delta_0, \delta \in (0,1)$ and a uniform mesh $\{s_j\}_j$, initial grid points  $\{  v^0_j \}$, initial condition $\{f^0(v^0_j)\}_j$, initial time step $\Delta t_0$, final time $T$ and threshold values $\epsilon_t$ and $\epsilon_v$.}
\Output{Numerical solution of equation $\eqref{eqn: gke_hom}$ with high resolution at singularity and numerical blow-up time $T_b$.}
Let $t=0$, $n=0$, $\Delta t = \Delta t_0$, $\Delta v_{min} = \min_j\{v^0_{j+1} - v^0_{j}\}$, $e_f = 0$\;
\While{$t < T$ and $\textstyle \begin{cases} 
\Delta t > \epsilon_t & \text{with adapt time step}, \\ 
\min \Delta v > \epsilon_v & \text{with fixed time step}
\end{cases}$
}{
1. Run Algorithm~\ref{alg:one_step_MR_JKO} to get $\{ v^{n+1}_j \}$ and $\{ f^{n+1}(v^n_j) \}$ and update $e_f$\;
2. If using adapt time step  \\
\While{$e_f$ = 1}{$\Delta t = \Delta t/ 2$ \;
Run Algorithm~\ref{alg:one_step_MR_JKO} to get $\{ v^{n+1}_j \}$ and $\{ f^{n+1}(v^n_j) \}$ and update $e_f$\;}
3. $t = t+\Delta t$, $n=n+1$.
}
$T_b(\epsilon) = t$.
  \caption{Mesh refinement JKO scheme with/without adaptive time stepping}
 \label{alg:MR_JKO_ver2}
\end{algorithm}

\begin{remark} 
In Algorithm~\ref{alg:MR_JKO_ver2}, there are two indicators for blow-up. One is when the mesh size $\Delta v_{\text{min}}$ falls below a threshold $\epsilon_v$, and the other is when optimization fails to converge within the maximum number of iterations. The latter condition is particularly valuable in the finite time blow-up case when a more accurate prediction of the blow-up time is desired. In such instances, as time approaches the blow-up time, a large time step may easily surpass the blow-up time, causing the optimization step to diverge. In response, the time step can be decreased, and this step can be recomputed iteratively, progressively reducing the time step until the optimization converges.
See Fig.~\ref{fig:gk1_adpt_tau} and the discussion in Section~\ref{sec31} for more details. 
\end{remark}

\subsection{Combining the two steps}
We now summarize the final algorithm in Algorithm~\ref{alg:SPMCSL-JKO} for calculating the spatially inhomogeneous granular kinetic equation. Note that for the time step, we can use either a fixed time step $\Delta t_0$ or adapt the time step based on the mesh size $\Delta t = \min \{ \Delta t_0, \min { 0.9 \frac{\Delta v_{\text{min}}}{L_v}, 0.9 \frac{\Delta x_{\text{min}}}{L_x} } \}$. In most examples, a fixed time step will be employed unless the accurate prediction of the finite blow-up time necessitates the use of an adaptive time step.

\begin{algorithm}[h!]
\SetAlgoLined
\SetKwInOut{Input}{Input}
\SetKwInOut{Output}{Output}
\Input{Parameters $\delta_0, \delta \in (0,1)$ and a uniform mesh $\{s_j\}_j$, initial time step $\Delta t_0$, time stepping strategy, final time $T$, stopping threshold $\epsilon_{x,v}$, initial mesh gird $\{  x_i^0 \}$, $\{  v_j^0 \}$ and 
initial condition $f(0,x,v) = f_0(x,v)$. }
\Output{Numerical solution of equation $\eqref{eqn: gke}$ with high resolution at singularity,  numerical blow-up time $T_b$ and blow-up indicators $B_x$, $B_v$.}
1. Compute  $\Delta v_{min} = \min_j\{v^0_{j+1} - v^0_{j}\}$, $\Delta x_{min} = \min_i{x^0_{i+1} - x^0_{i}}$, and set $B_x=0$ and $B_v=0$, $\Delta t = \Delta t_0$ \;
\While{$t < T$, $\Delta v_{min} > \epsilon_{x,v}$ and $\Delta x_{min} > \epsilon_{x,v}$, $e_f=0$}{
{1. Either use fixed time step $\Delta t = \Delta t_0$ or adjust the time step according to: $\Delta t = \min \{ \Delta t_0, \min \{ 0.9 \frac{\Delta v_{min}}{L_v}, 0.9 \frac{\Delta x_{min}}{L_x} \} \}$ }\;
2. Transport step \\
\For{$j = 1 \cdots N_v$}{
 Compute $f^*(x_i^n, v^n_j)$ ~ for $i = 1 \cdots N_x$ by \eqref{eqn:semi_update}  \;
}
3. Collision step \\
\For{$i = 1 \cdots N_x$}{
Compute $f^{n+1}(x^n_i, v^n_j)$ by Algorithm~\ref{alg:non_uni_JKO} given $f^*(x^n_i, v^n_j)$ and update $e_f$ \;
}
4. Mesh refinement\\
For one bump case, determine\\
$S = \{ (x_i, v_j): f^*(x_i, v_j) > \delta_0 (\max_{i,j} f^*(x_i, v_j)), x_i>0, v_j>0\}$, then $r_x = \min \{ x_i: (x_i, \cdot) \in S \}$ and $r_v = \min \{ v_j: (\cdot, v_j) \in S \}$, 
then find $\mu_x(s)$ via $r_x$ and \eqref{mu1}, find $\mu_v(s)$ via $r_v$ and \eqref{mu1}. \\
For two bumps case, determine $S = \{ (x_i, v_j): f^*(x_i, v_j) > \delta_0 (\max_{i,j} f^*(x_i, v_j)), x_i>0, v_j>0\}$, then $r_{x1} = \min \{ x_i: (x_i, ) \in S \}$, $r_{x2} = \max \{ x_i: (x_i, \cdot) \in S \}$, $r_{v1} = \min \{ v_j: (\cdot, v_j) \in S \}$ and $r_{v2} = \max \{ v_j: (\cdot, v_j) \in S \}$. Then find $\mu_x(s)$ via $r_{x1}, r_{x2}$ and \eqref{mu2}, find $\mu_v(s)$ via $r_{v1}, r_{v2}$ and \eqref{mu2}. \\
5. Interpolation \\
 Use `pchip' interpolation to get function value at the new gird points $\{  f^{n+1}(x^{n+1}_i, v^{n+1}_j) \}$ and rescale the total mass dimensionwise according to \eqref{eqn:semi_update} and \eqref{eqn:mass_conserv_v} respectively.\; 
5. Compute $\Delta x_{min} = \min_i(x^{n+1}_{i+1} - x^{n+1}_{i})$, $\Delta v_{min} = \min_j(v^{n+1}_{j+1} - v^{n+1}_{j})$ \;
6. $t = t + \Delta t$. }
$T_b = t$; set $B_x=1$ if $\Delta x_{min}<\epsilon_{x,v}$; set $B_v=1$ if $\Delta v_{min}<\epsilon_{x,v}$ or $e_f = 1$. 
  \caption{Semi-Lagrangian JKO scheme}
 \label{alg:SPMCSL-JKO}
\end{algorithm}

\section{Numerical examples}\label{sec:numerical_examples}
This section is dedicated to presenting the numerical findings regarding the (non)blow-up behavior of the granular kinetic equation. Specifically, there are several solution behaviors:
\begin{itemize}
    \item[1)]no blow-up occurs;
    \item[2)]blow-up initially forms in the spatial $x$ direction;
    \item[3)] blow-up initially forms in the velocity $v$ direction; 
    \item[4)] blow-up occurs simultaneously in both the spatial  and velocity directions.
\end{itemize}
The last case is theoretically possible but exceptionally rare and challenging to construct. Therefore, our focus will be solely on the first three cases.

We first validate the capability of our proposed numerical solver (Algorithm~\ref{alg:MR_JKO_ver2}) for the spatially homogeneous case, where analytical results are well-established. We then disclose the numerical outcomes derived from the implementation of our Algorithm~\ref{alg:SPMCSL-JKO}. Based on these results, we formulate a conjecture regarding the blowup behavior for the spatially inhomogeneous granular kinetic equation.

\subsection{Finite time blow-up verification: homegeneous problem} \label{sec31}
We begin by examining the spatially homogeneous case and establishing an appropriate blow-up criterion ($\epsilon_{x,v}$ in Algorithm ~\ref{alg:SPMCSL-JKO}). To this end, we consider $\eqref{eqn: gke}$ with the kernel $W(v) = |v|$ and an initial condition $f(0,v) = g(v)$. The analytical blow-up time for this scenario is given by:
\begin{equation*} 
    T = \frac{1}{2 \max_v g(v)}.
\end{equation*}

\begin{figure}[h!]
\centering
{\includegraphics[width=0.4\textwidth]{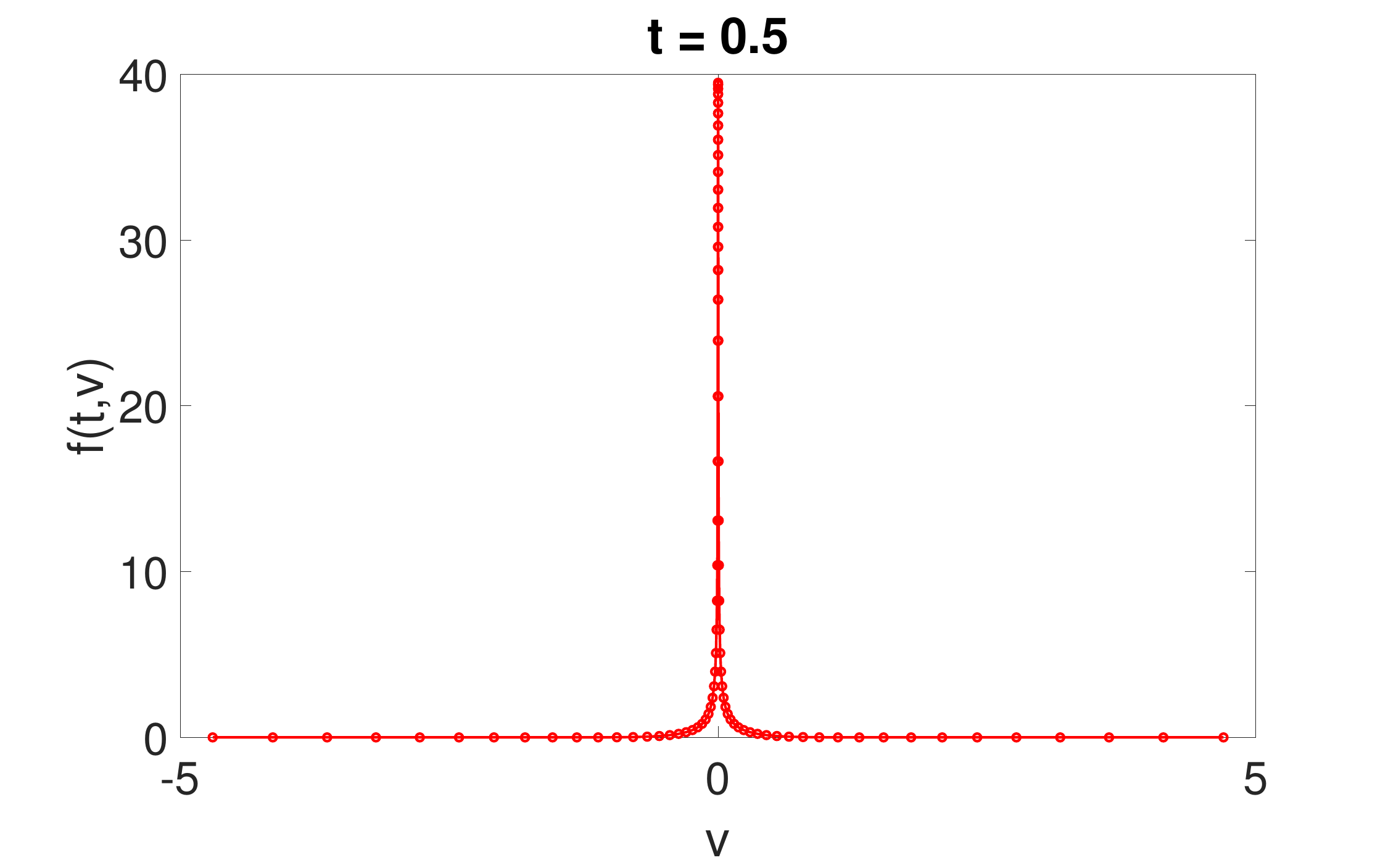}}
{\includegraphics[width=0.4\textwidth]{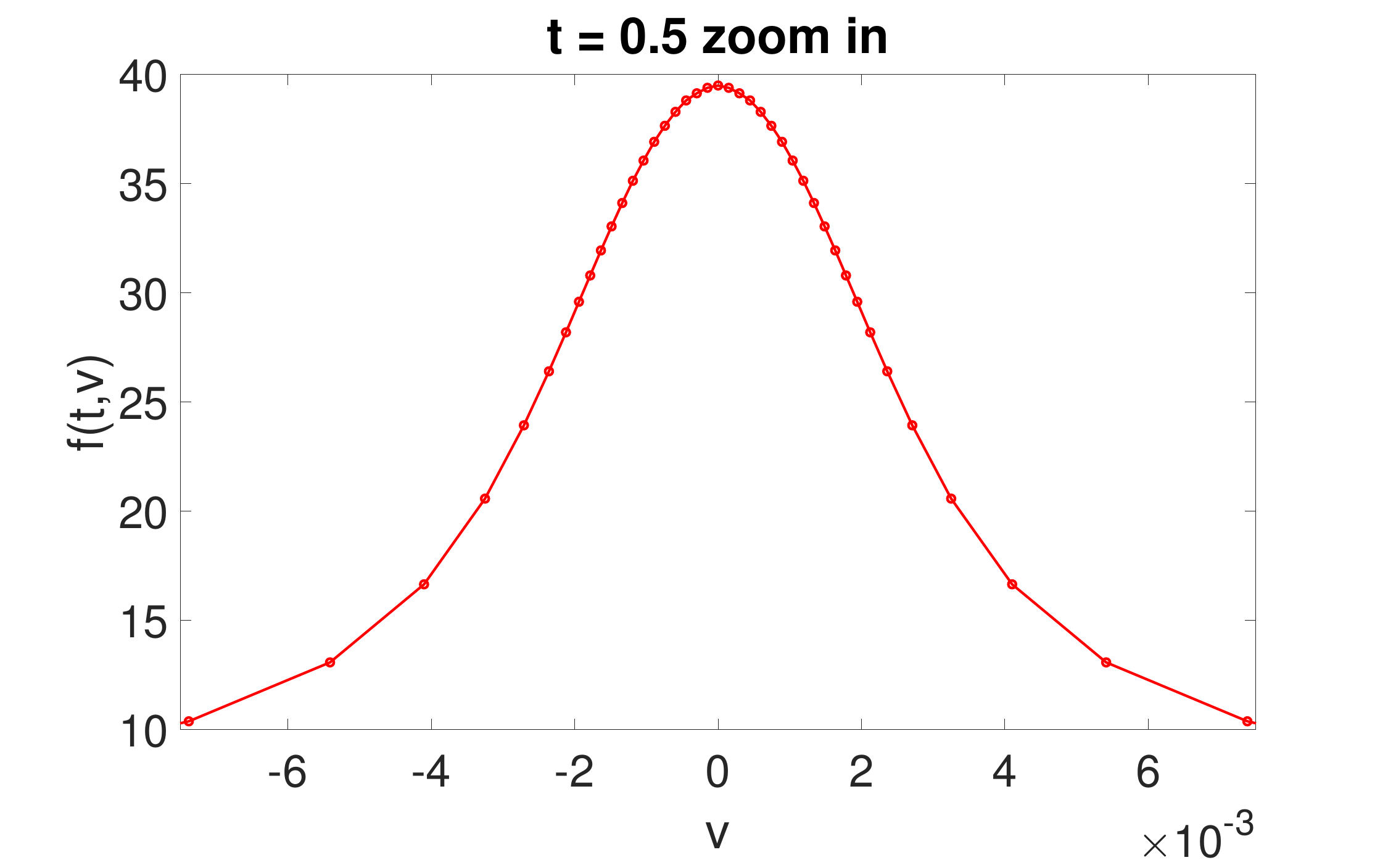}}
 \caption{Numerical solution with initial condition $g_1(v)= e^{-2v^2}$, $\delta=0.2$, $\delta_0 = 0.5$ and fixed time step $\Delta t=0.01$ with the kernel $W(v) = |v|$.}
 \label{fig:gk1_final}
\end{figure}

To evaluate our numerical approach, we consider three one-bump initial conditions: $g_1(v) = e^{-2v^2}$, $g_2(v) = 2e^{-2v^2}$, and $g_3(v) = 4e^{-2v^2}$ and one two-bumps initial condition $g_4(v) = e^{-10(v-1.5)^2} +  e^{-10(v+1.5)^2}$, with analytic blow-up times $T_1=0.5$, $T_2 = 0.25$, $T_3=0.125$ and $T_4 = 0.5$, respectively. To fully investigate the blow-up behavior, we employ two time-stepping strategies. One is with a fixed time step $\Delta t$. The numerical results with the initial condition $g_1$ at $t=0.5$  are depicted in Figure~\ref{fig:gk1_final}, where $\Delta t=0.01$ and the solution  $f(t=0.5, v)$ nearly converges to a Dirac delta function centered at $v=0$.

\begin{figure}[h!]
\centering
{\includegraphics[width=0.45\textwidth]{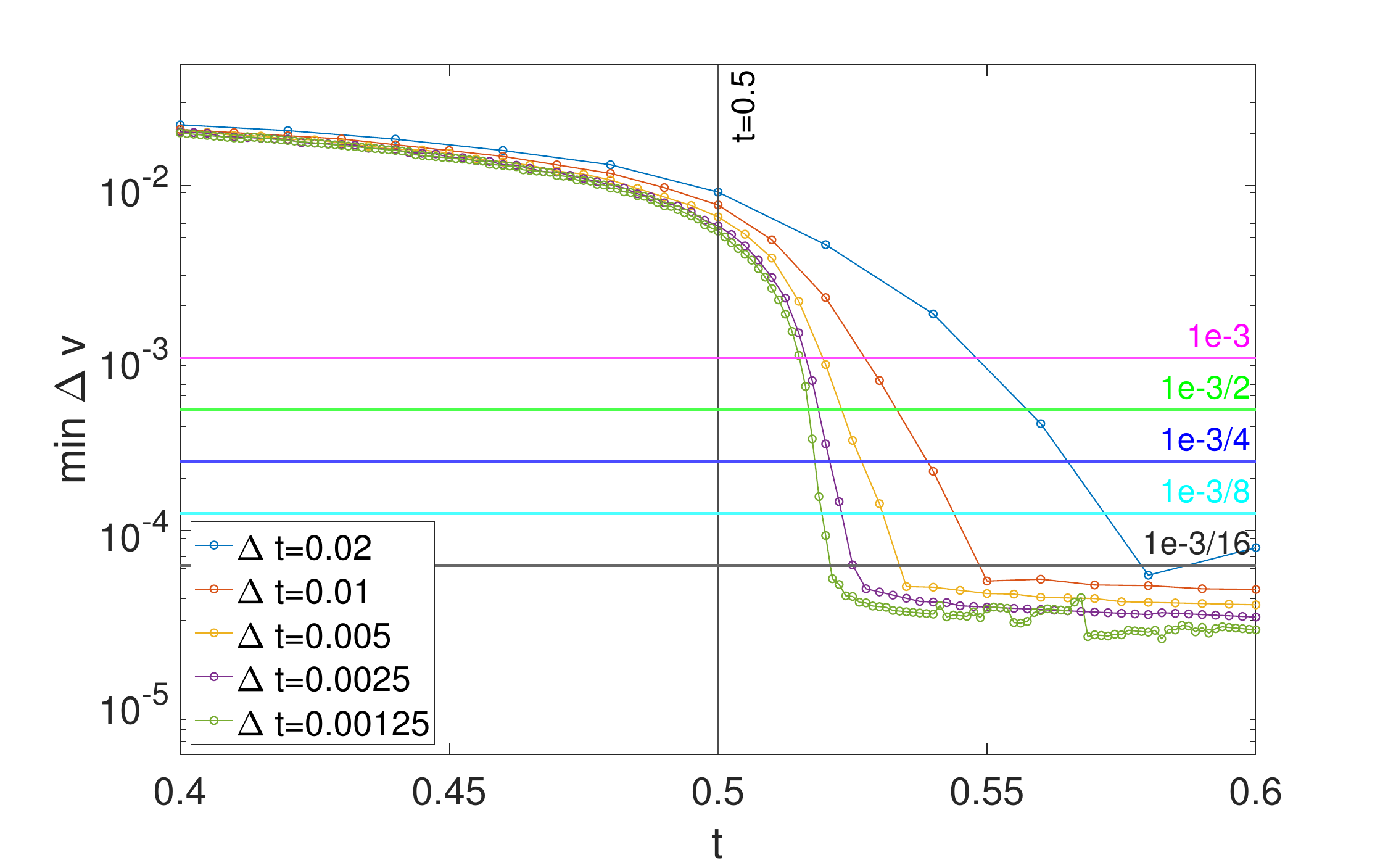}}
{\includegraphics[width=0.45\textwidth]{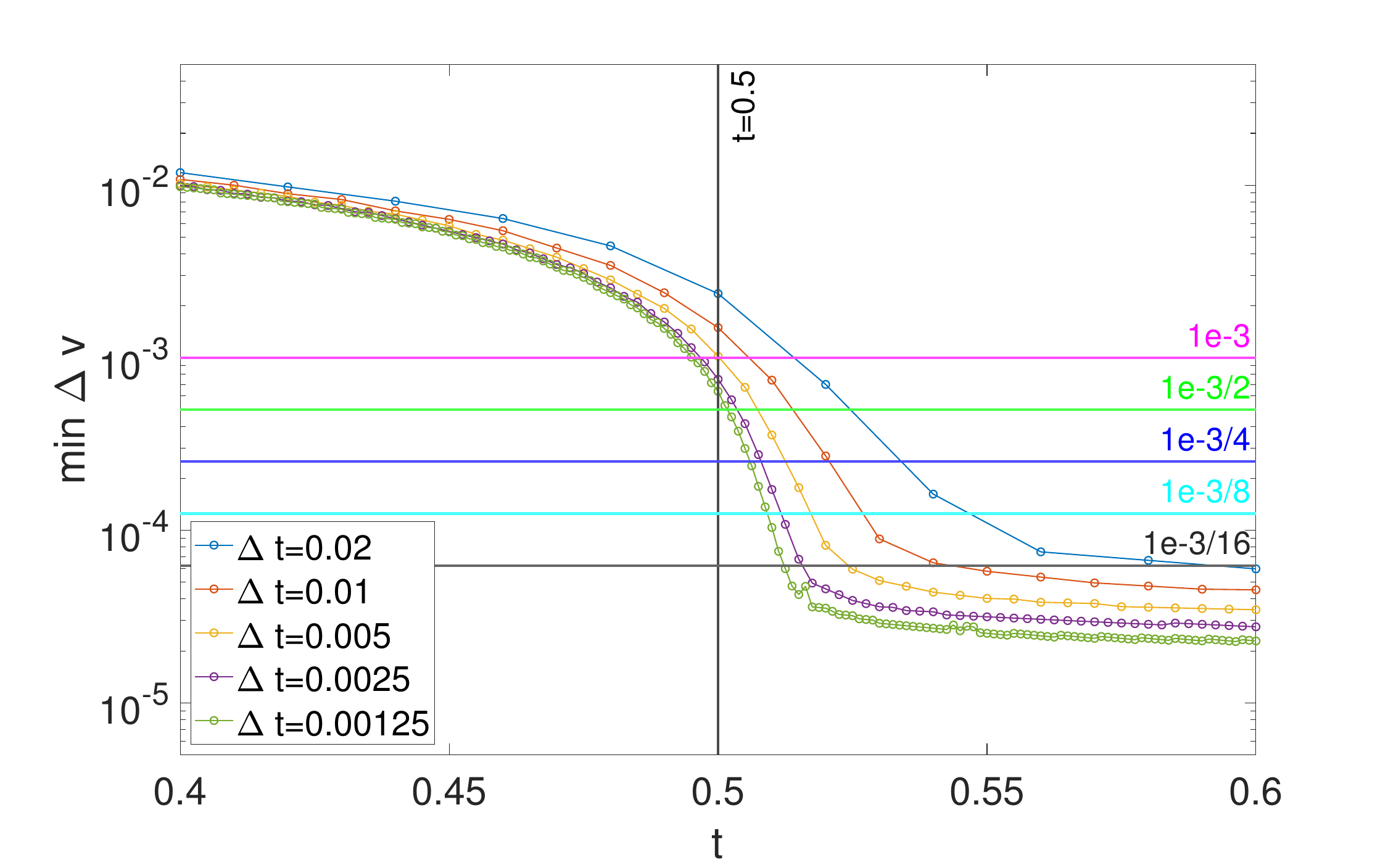}}

{\includegraphics[width=0.45\textwidth]{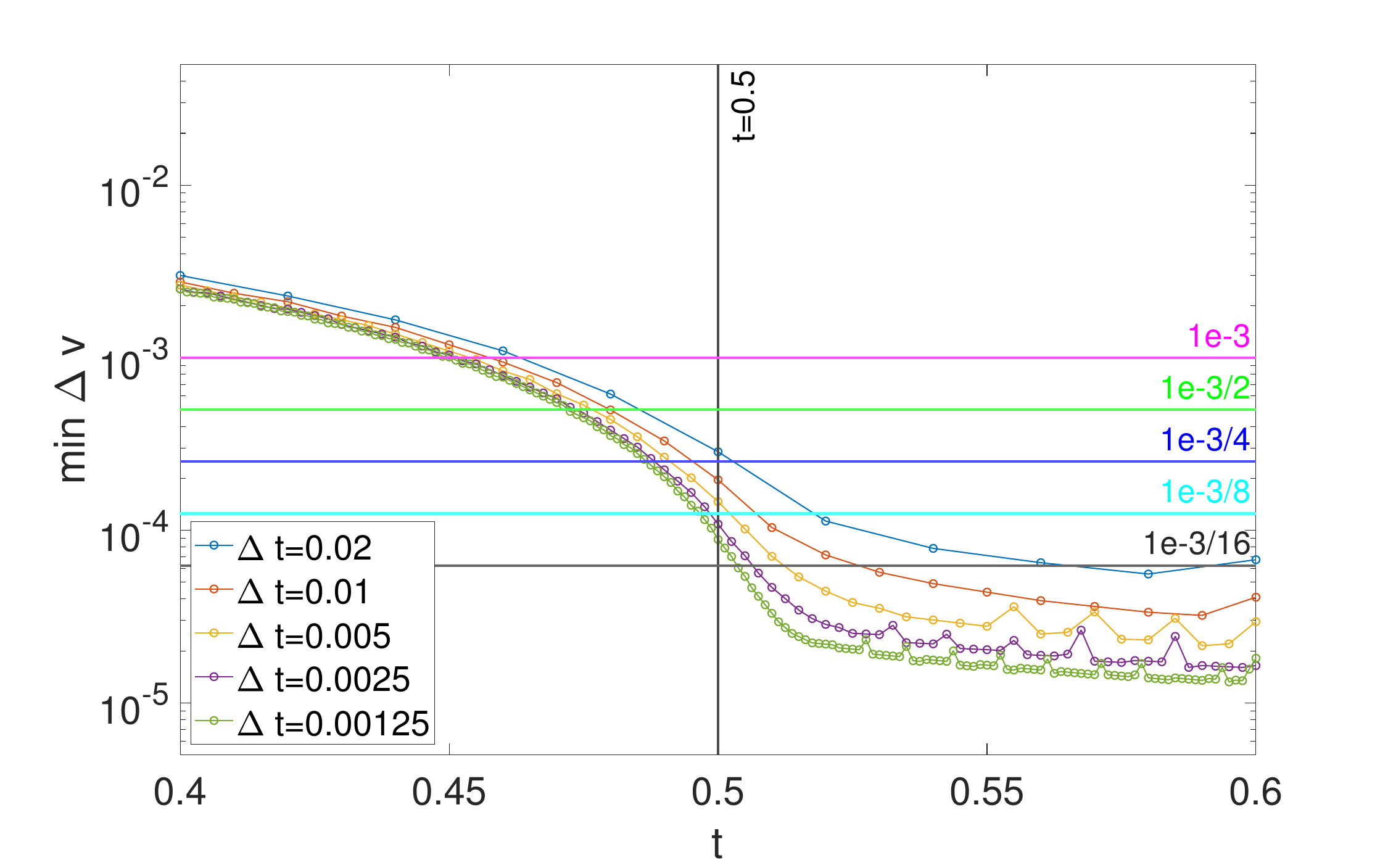}}
{\includegraphics[width=0.45\textwidth]{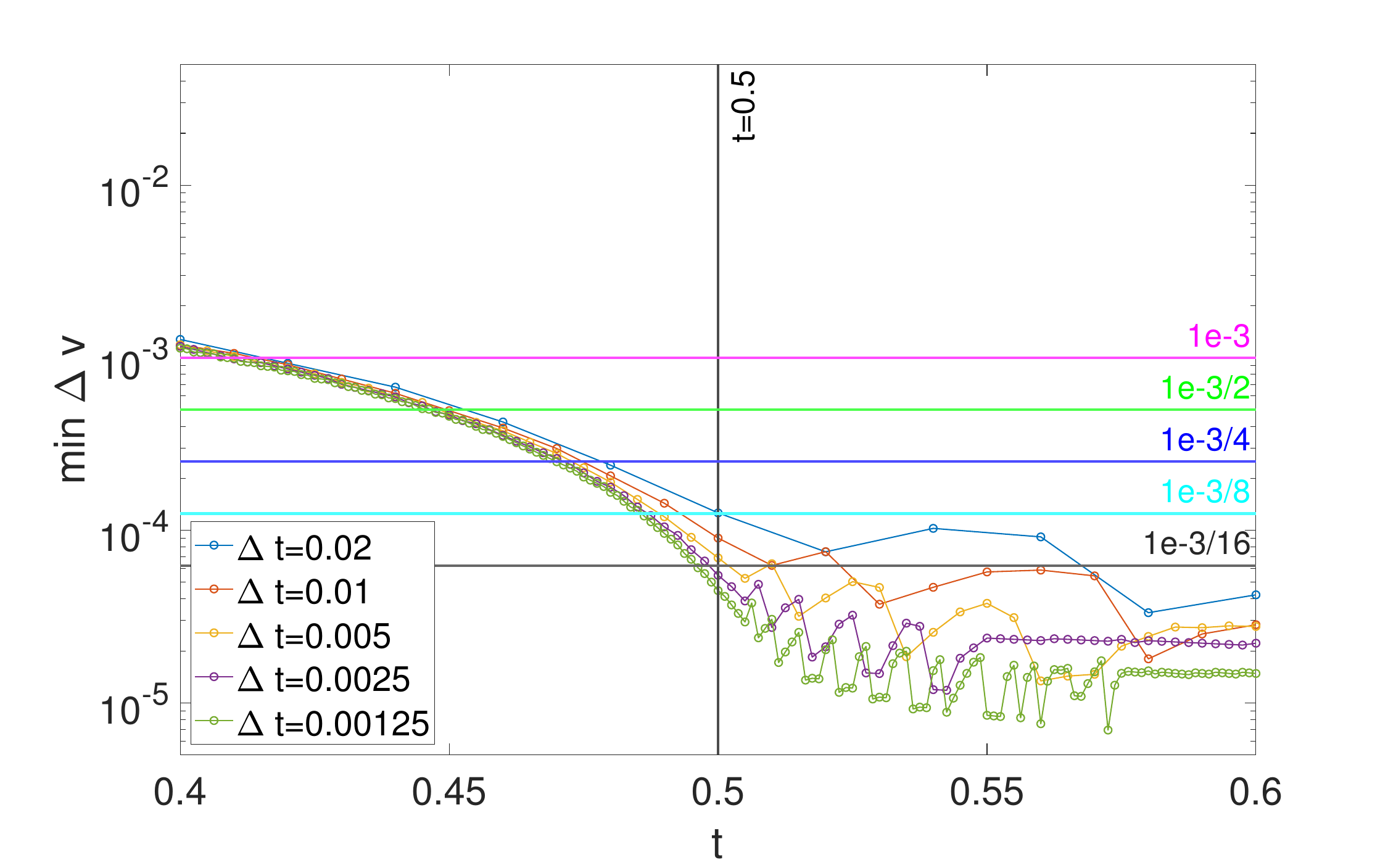}}

 \caption{Time versus minimum $\Delta v$ for initial condition $g_1(v)$ with fixed $\Delta t$ with the kernel $W(v) = |v|$. Here $N_v=121$ and $\delta=0.5$, the images from top left to bottom right correspond to $\delta_0 = 0.05, 0.1, 0.5, 0.8$ respectively.}
 \label{fig:gk1_2}
\end{figure}

More detailed results are provided in the Table~\ref{tab:bu},~\ref{tab:bu_tb2} and Figure~\ref{fig:gk1_2}. From Table~\ref{tab:bu},~\ref{tab:bu_tb2}, we observe that as we reduce the stopping criterion (as indicated by the columns), the numerical blow-up time increases. Conversely, when we decrease the time step size, the numerical blow-up time decreases. Notably, as we reduce both the stopping criterion and the time step size simultaneously, the numerical blow-up time approaches convergence with the analytical blow-up time. In Figure~\ref{fig:gk1_2}, we observed that the minimum $\Delta v$ decreased very rapidly near the analytic numerical blow-up time. In this particular test, we have chosen to omit the exit flag as a stopping criterion, allowing the algorithm to continue running even if the minimizer is not reached, and this is why we observe oscillations at the tail of the results when the simulation time extends beyond the analytical blow-up time.

\begin{figure}[h!]
\centering
{\includegraphics[width=0.45\textwidth]{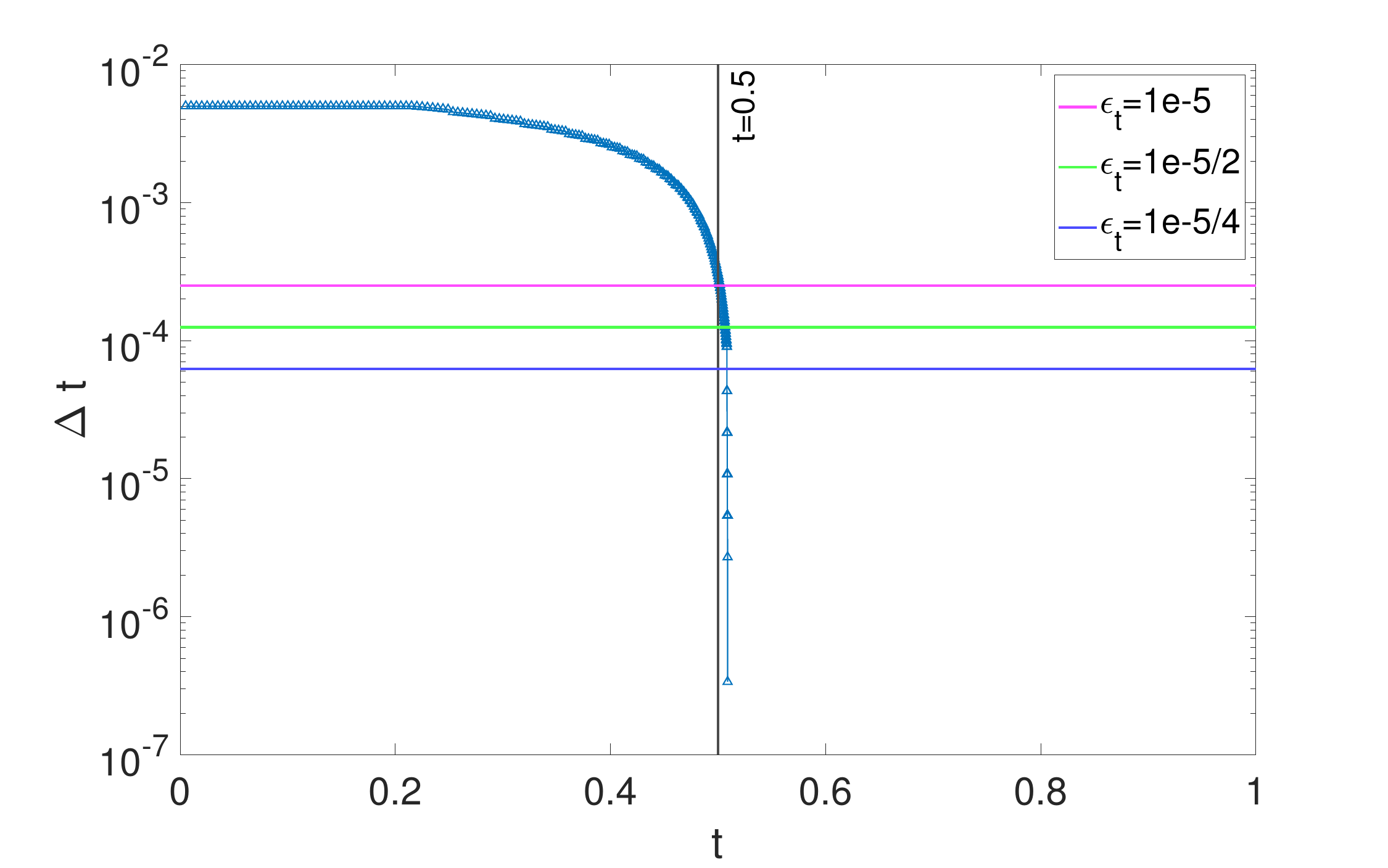}}
{\includegraphics[width=0.45\textwidth]{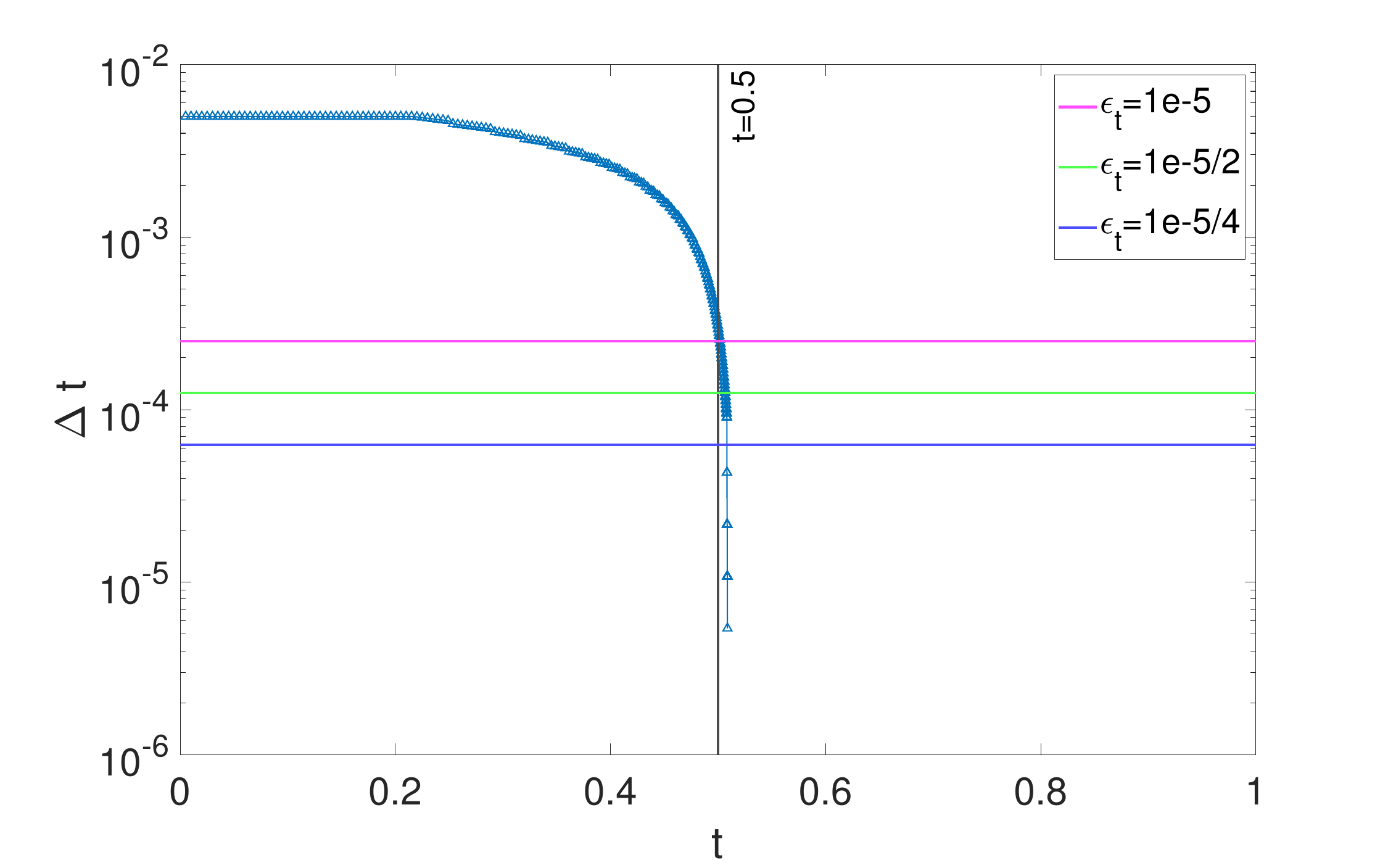}}
 \caption{Time versus adaptive time step size $\Delta t$ with the kernel $W(v) = |v|$. Here, we consider the initial condition $g_1(v)$ with $N_v=121$ and $\delta=0.5$. On the left is $\delta_0=0.05$, and on the right is $\delta_0=0.5$.}
 \label{fig:gk1_adpt_tau}
\end{figure}

\begin{table}[h!]
\centering
 \begin{tabular}{|c c | c c c c c |} 
 \hline
    & & $\Delta t = 0.02$ & $\Delta t = 0.01$ & $\Delta t = 0.005$ & $\Delta t = 0.0025$ & $\Delta t = 0.00125$ \\ [0.5ex] 
 \hline
 $f(0,v) = g_1(v)$ & $\epsilon_v=1e-3$ & 0.48 & 0.48 & 0.47 & 0.4675 & 0.4512 \\ 
 & $\epsilon_v=1e-3/2$ & 0.5 & 0.49 & 0.49 & 0.4875 & 0.4725\\
 & $\epsilon_v=1e-3/4$ & 0.52 & 0.51 & 0.5 & 0.5 & 0.4875 \\
 & $\epsilon_v=1e-3/8$ & 0.54 & 0.52 & 0.51 & 0.5075 & 0.4975 \\
 & $\epsilon_v=1e-3/16$ & 0.58  & 0.53  & 0.515  & 0.5075 & 0.50375 \\
 \hline
 $f(0,v) = g_2(v)$ & $\epsilon_v=1e-3$ & 0.24 & 0.24 & 0.235 & 0.235 & 0.2275 \\ 
 & $\epsilon_v=1e-3/2$ & 0.26 & 0.25 & 0.245 & 0.2425 & 0.2375\\
 & $\epsilon_v=1e-3/4$ & 0.26 & 0.26 & 0.25 & 0.2475 & 0.245\\
 & $\epsilon_v=1e-3/8$ & 0.28 & 0.26 & 0.255 & 0.2525 & 0.25\\
 & $\epsilon_v=1e-3/16$ & 0.28  & 0.26 & 0.26 & 0.255  & 0.2525 \\
 \hline
 $f(0,v) = g_3(v)$ & $\epsilon_v=1e-3$ & 0.12 & 0.12 & 0.12 & 0.1175 & 0.1137 \\
 & $\epsilon_v=1e-3/2$ & 0.12 & 0.13 & 0.125 & 0.1225 & 0.1187\\
 & $\epsilon_v=1e-3/4$ & 0.14 & 0.13 & 0.125 & 0.125 & 0.1225\\
 & $\epsilon_v=1e-3/8$ & 0.14 & 0.13 & 0.13 & 0.1275 &0.1265\\
 & $\epsilon_v=1e-3/16$ &  0.14 & 0.13 & 0.13 & 0.1275 & 0.12625 \\
 \hline
 \end{tabular}
 \caption{Numerical blow-up time with various $\Delta t$, $\epsilon_v$ and initial conditions with the kernel $W(v) = |v|$. Here we use $N_v=121$, $\delta=0.5$ and $\delta_0=0.5$.}\label{tab:bu}
\end{table}

The second strategy is to employ adaptive step sizes, with the initial step size set at $\Delta t_0 = 0.01$, and subsequently halving it until the minimization process converges. As depicted in Figure~\ref{fig:gk1_adpt_tau}, one can observe that the time step size $\Delta t$ decays significantly as it approaches the analytic blow-up time. Based on these observations, we believe that our proposed numerical solver is capable of reproducing the blow-up solution and accurately capturing the analytic blow-up time.

\begin{table}[h!]
\centering
 \begin{tabular}{|c | c c c c |} 
 \hline
     & $\Delta t = 0.01$ & $\Delta t = 0.005$ & $\Delta t = 0.0025$ & $\Delta t = 0.00125$ \\ [0.5ex] 
 \hline
   $\epsilon_v=$1e-3  & 0.49 & 0.48 & 0.4775 & 0.4775 \\ 
  $\epsilon_v=1e-3/2$  & 0.51 & 0.505 &0.50125 & 0.50125\\
  $\epsilon_v=1e-3/4$  & 0.53 & 0.52 & 0.52 &0.51875\\
  $\epsilon_v=1e-3/8$  & 0.56 &0.54 & 0.5325 & 0.52875 \\
 \hline
 \end{tabular}
 \caption{{Numerical blow-up time with various $\Delta t$ and $\epsilon_v$ for initial condition $g_4(v)$ with $N_v=121$, $\delta=0.5$ and $\delta_0=0.5$}.}\label{tab:bu_tb2}
\end{table}

\begin{figure}[h!]
\centering
{\includegraphics[width=0.45\textwidth]{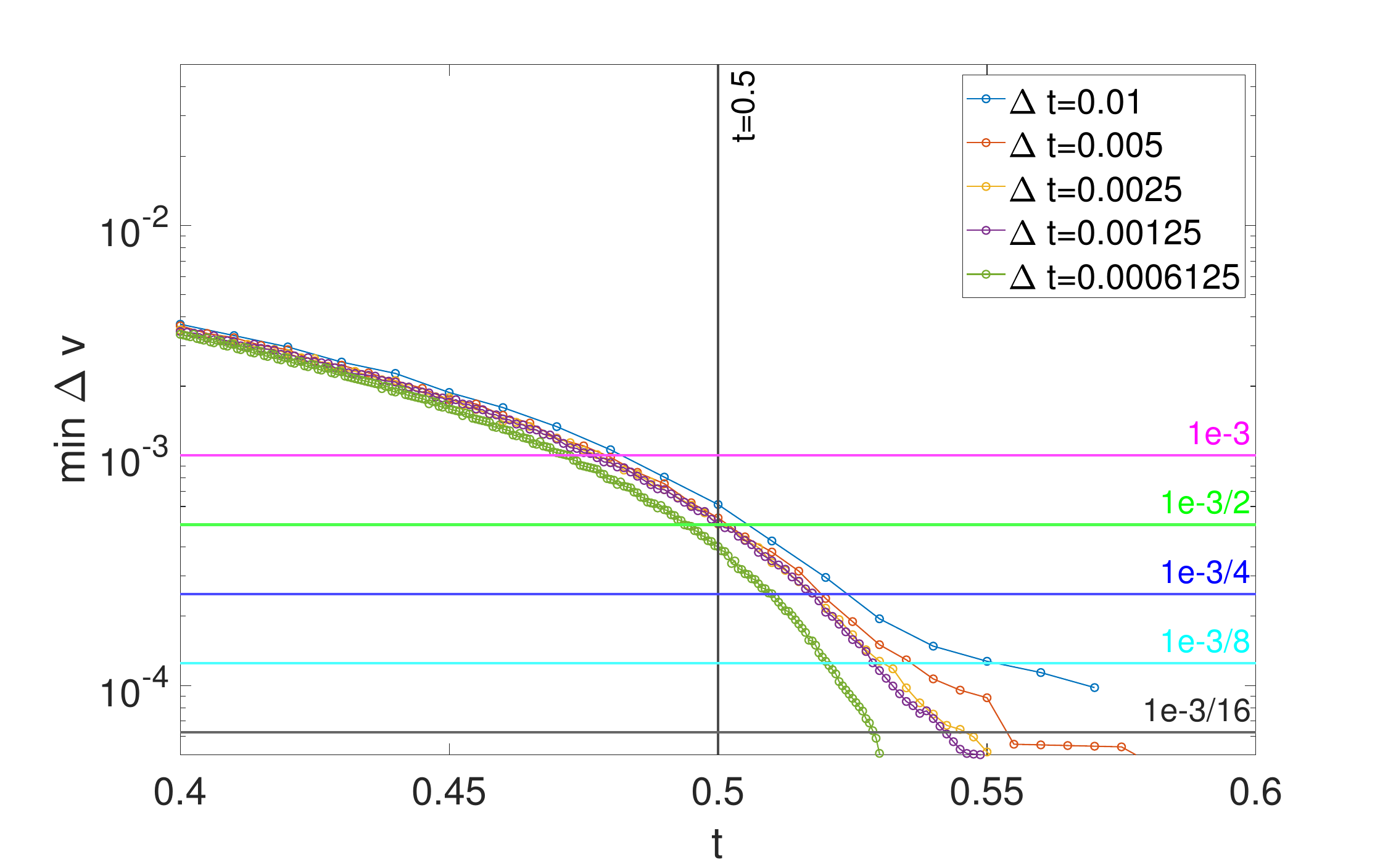}}
{\includegraphics[width=0.45\textwidth]{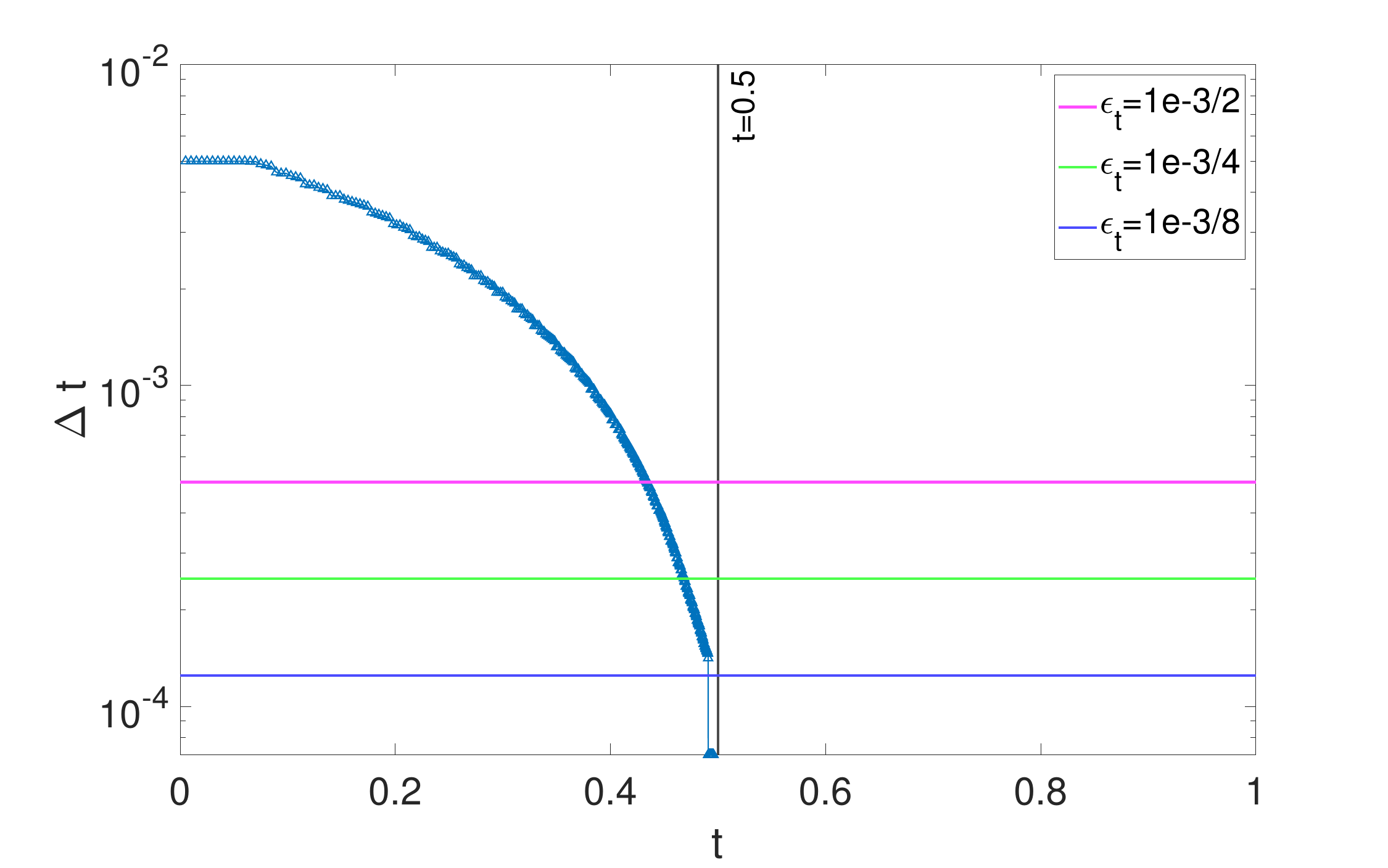}}
 \caption{Plot with initial condition $g_4(v)$ and kernel $W(v) = |v|$, $N_v=121$, $\delta=0.5$ and $\delta_0=0.5$. On the left is time versus minimum $\Delta v$ and the right is time versus adaptive time step size $\Delta t$.}
 \label{fig:gk1_tb}
\end{figure}

\subsection{Infinite time blow-up verification: homogeneous problem}\label{sec:ibu}
\begin{figure}[h!]
\centering
{\includegraphics[width=0.3\textwidth]{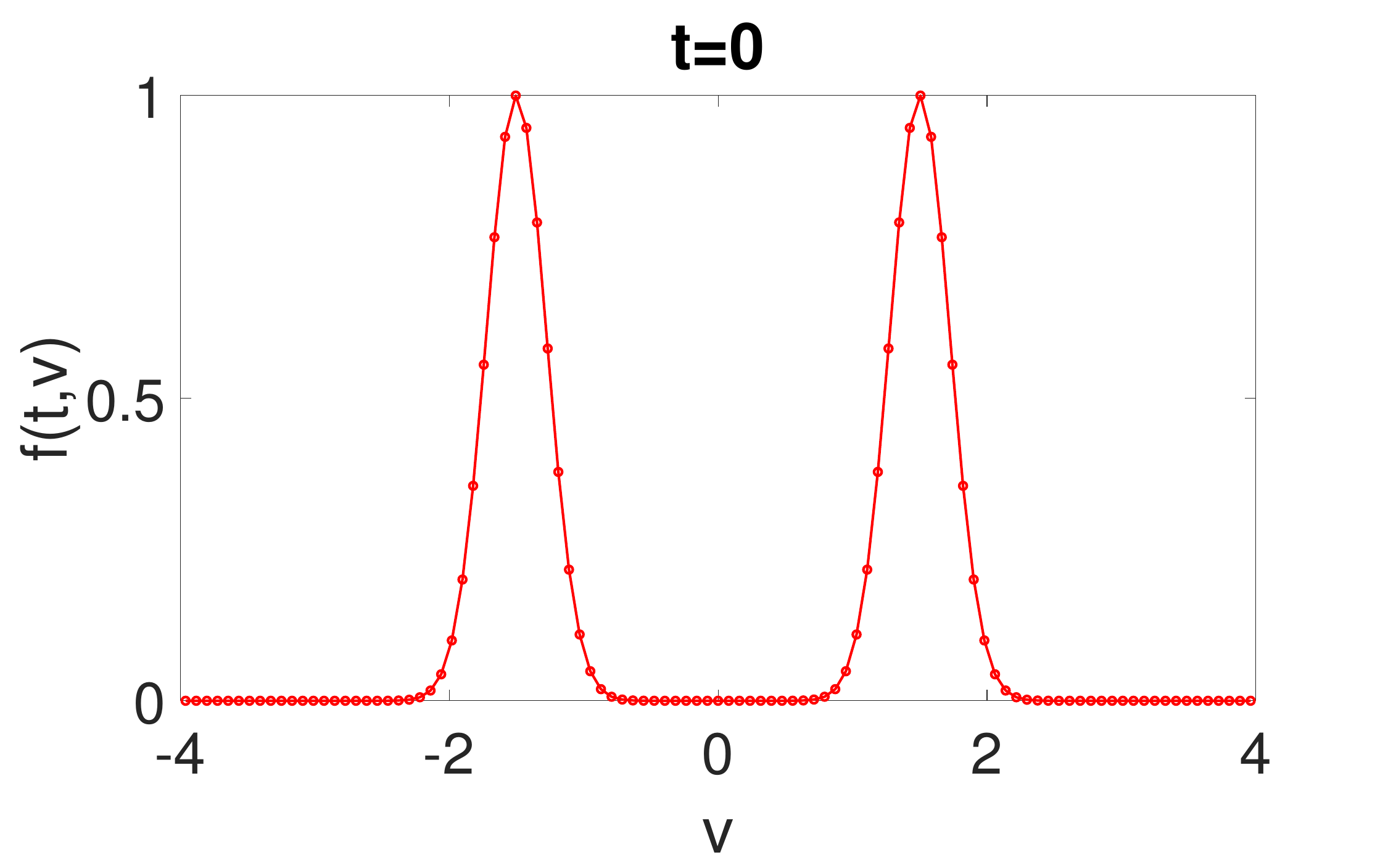}}
{\includegraphics[width=0.3\textwidth]{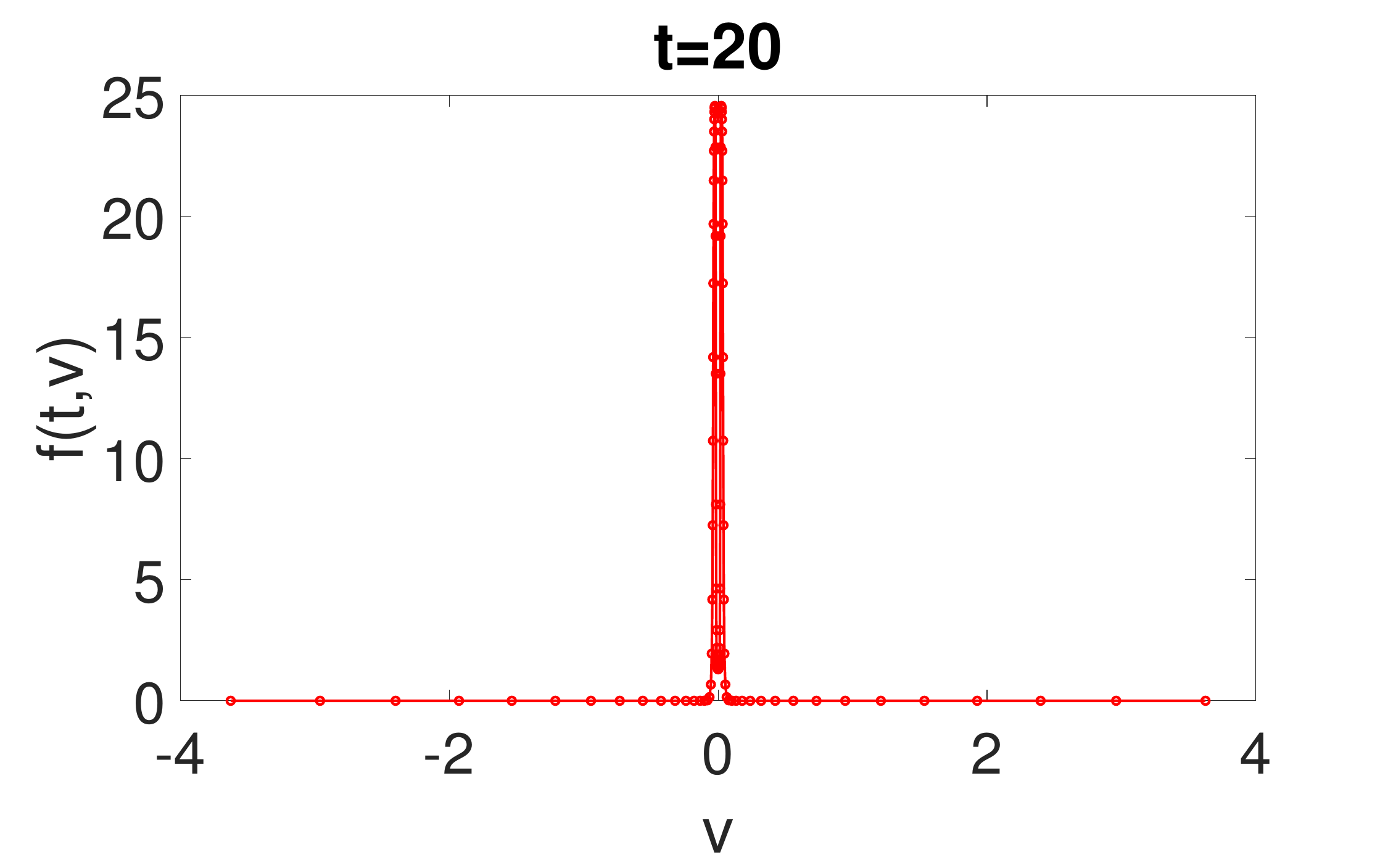}}
{\includegraphics[width=0.3\textwidth]{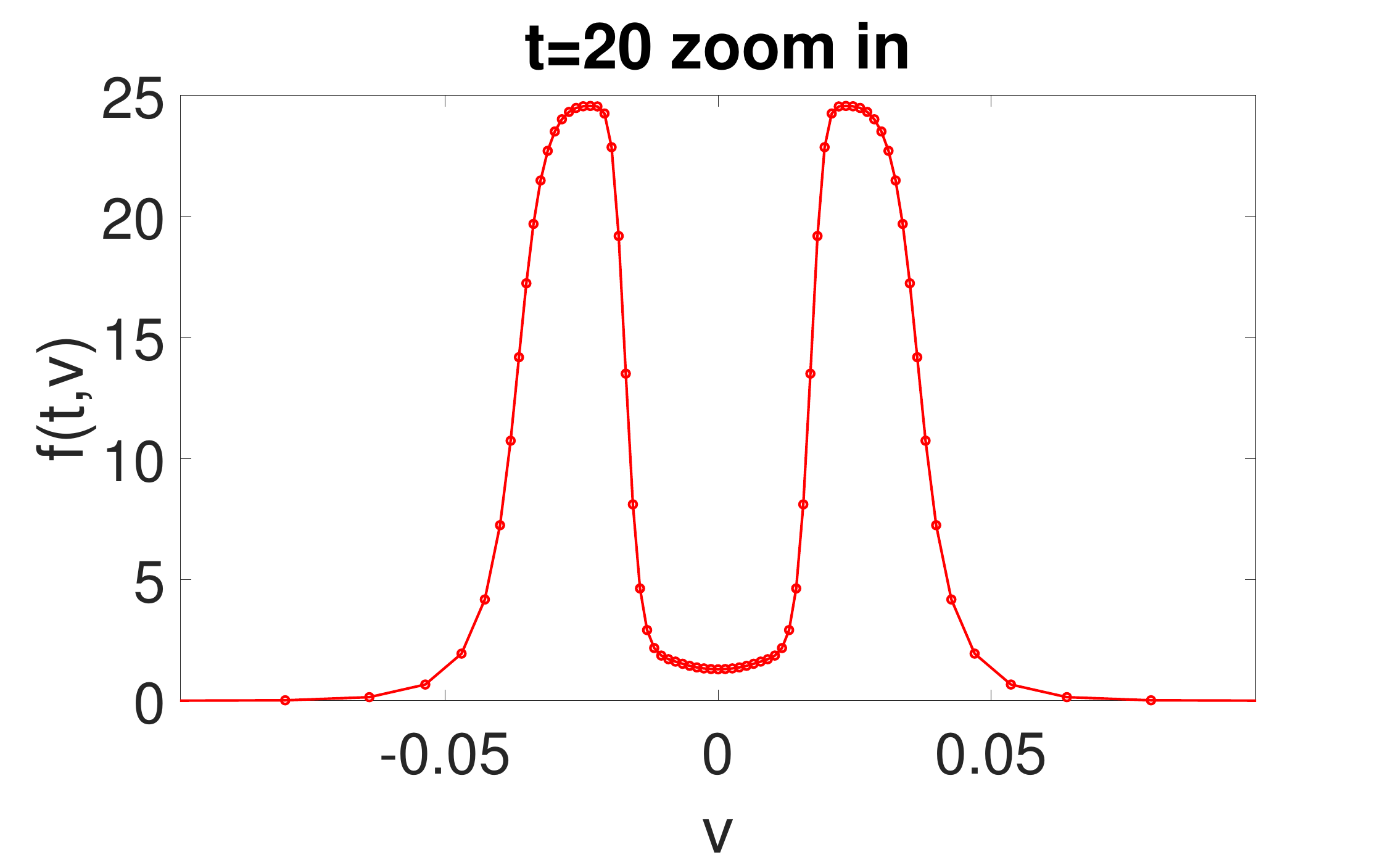}}
 \caption{Numerical solution with $\delta_0=0.8$, $\delta = 0.5$, $\Delta t=0.01$ and $N_v=101$ with the kernel $W(v) = |v|^3$. The left plots the initial condition, the middle depicts $f(20,v)$, and the right plot provides a zoom-in of $f(20,v)$ at $v=0$. }
 \label{fig:gk3}
\end{figure}

\begin{figure}[h!]
\centering
{\includegraphics[width=0.45\textwidth]{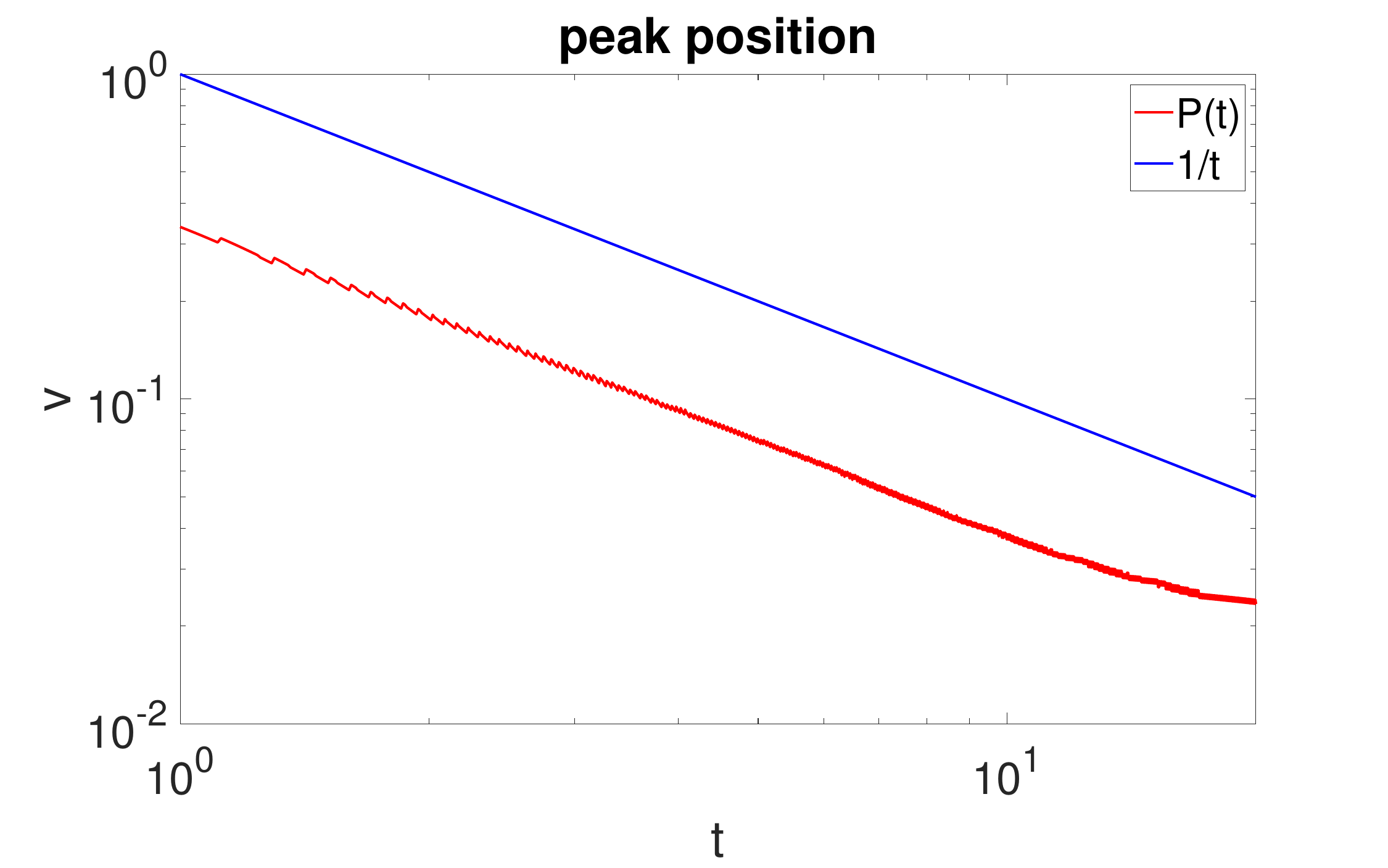}}
{\includegraphics[width=0.45\textwidth]{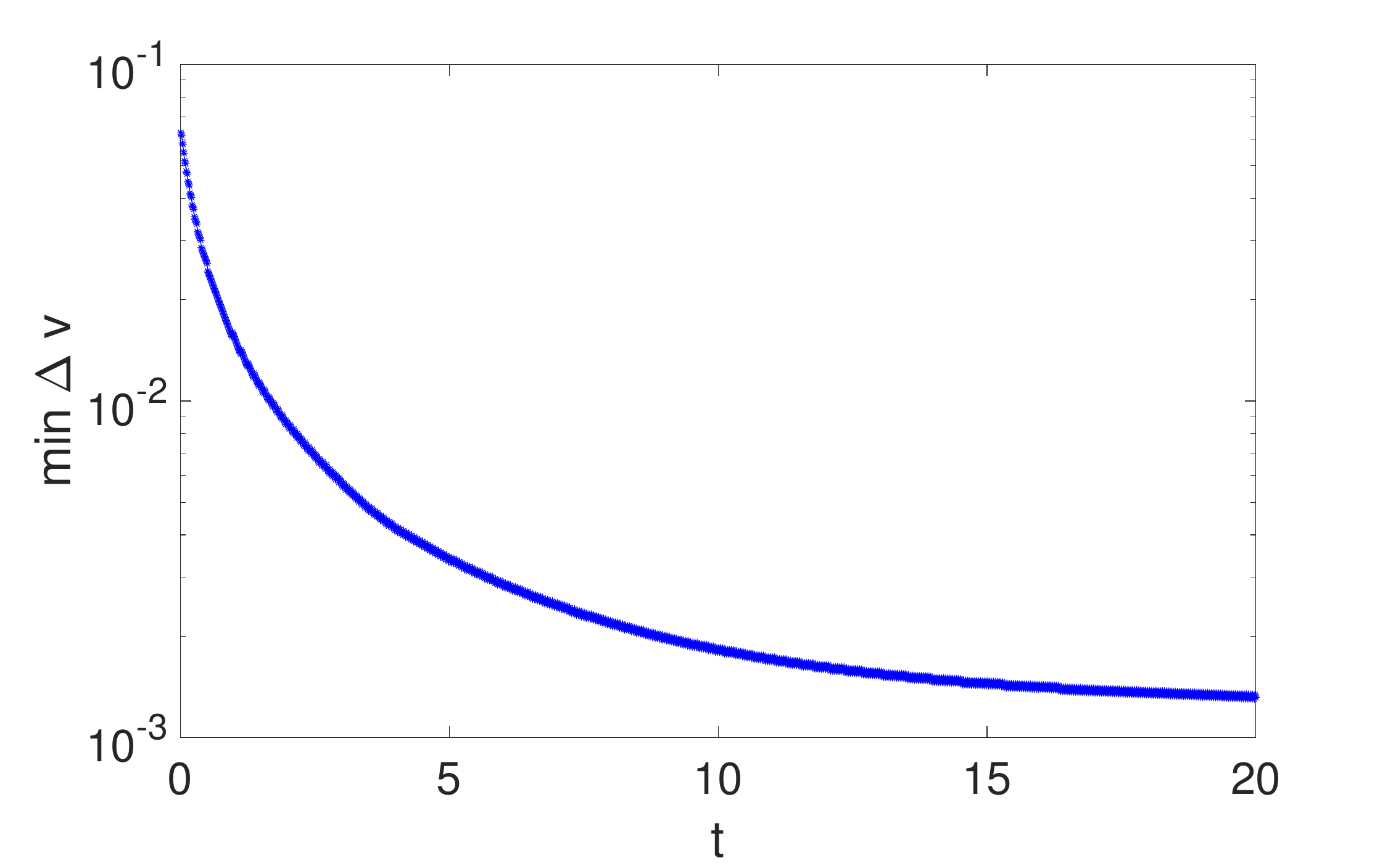}}
 \caption{Numerical solution with $\delta_0=0.8$, $\delta = 0.5$, $\Delta t=0.01$ and $N_v=101$ with the kernel $W(v) = |v|^3$. The left plot is the position of maximum value of $f(t,v)$, i.e., $P(t) = \max_v f(t,v)$ and the right plot is the minimum $\Delta v$ along time. }
 \label{fig:gk3_2}
\end{figure}

We now consider the spatially homogeneous case of $\eqref{eqn: gke}$ with kernel $W(v) = |v|^3$, whose solution converges asymptotic approaches  $f(t, v) = \frac{M}{2} \delta_{\frac{1}{t}} + \frac{M}{2} \delta_{-\frac{1}{t}}$ ($M$ is the total mass) as time approaches infinity. To demonstrate that our method can capture the desired blow-up behavior, we consider an initial condition with two symmetric bumps:
\begin{align} \label{in_two}
f(0,v) = e^{-10(x-1.5)^2} + e^{-10(x+1.5)^2} \,.
\end{align}
Running algorithm~\ref{alg:MR_JKO_ver2} with $\delta_0=0.8$, $\delta = 0.5$, $\Delta t=0.01$ and $N_v=101$, we observed in Figure~\ref{fig:gk3} that the two initially symmetric bumps shrinking to two symmetric Dirac delta bumps and approaching $v=0$ from two sides. Additionally, we examined the rate at which the two peaks converge to the origin, and as shown in Figure~\ref{fig:gk3_2}, they converge with a rate of $O(\frac{1}{t})$, which is consistent with the analytical result. In contrast to the finite time blow-up cases (see Figure~\ref{fig:gk1_2}), the minimum of $\Delta v$ decreases very slowly. 

Comparing Figure~\ref{fig:gk1_2} with Figure~\ref{fig:gk3_2}, it is important to point out that the minimum of $\Delta v$ exhibits distinct patterns. In cases with finite-time blow-up, the minimum of $\Delta v$ experiences a rapid decline as time approaches or slightly surpasses the analytic blow-up time. Conversely, in scenarios without finite-time blow-up, the minimum of $\Delta v$ decreases at a considerably slower rate. This observation suggests that one can confirm the presence of finite-time blow-up by monitoring the rate of decay of the minimum of $\Delta v$, serving as an additional criterion for detecting finite-time blow-up behavior. In practice, when employing a fixed time stepping strategy as in Algorithm \ref{alg:MR_JKO_ver2}, we suggest setting the threshold for detecting blow-up in the single bump scenario at $\epsilon_v = 1e-3/16$ (as shown in Figure~\ref{fig:gk1_2},Table~\ref{tab:bu}) and setting $\epsilon_v = 1e-3/8$ for two bumps case (see Figure~\ref{fig:gk1_tb},Table~\ref{tab:bu_tb2}). Furthermore, in situations aiming at accurately capturing the analytical blow-up time, we implement the adaptive time step strategy with a threshold set as in Algorithm \ref{alg:MR_JKO_ver2} at $\epsilon_t = 5e-6$ (See Figure~\ref{fig:gk1_adpt_tau}).

\subsection{Spatially inhomogeneous case with $\gamma=3$}\label{sec:sp_inhom_gamma_3}
As previously discussed in Section~\ref{sec:intro} and numerically verified in Section~\ref{sec:ibu}, the spatially homogeneous case of $\eqref{eqn: gke}$ with kernel $W(v)=|v|^3$ exhibits an infinite-time blow-up solution. In this section, we aim to numerically investigate how the spatial dependence will affect the behavior of the solution solving \eqref{eqn: gke}.

\begin{figure}[h!]
\centering
{\includegraphics[width=0.45\textwidth]{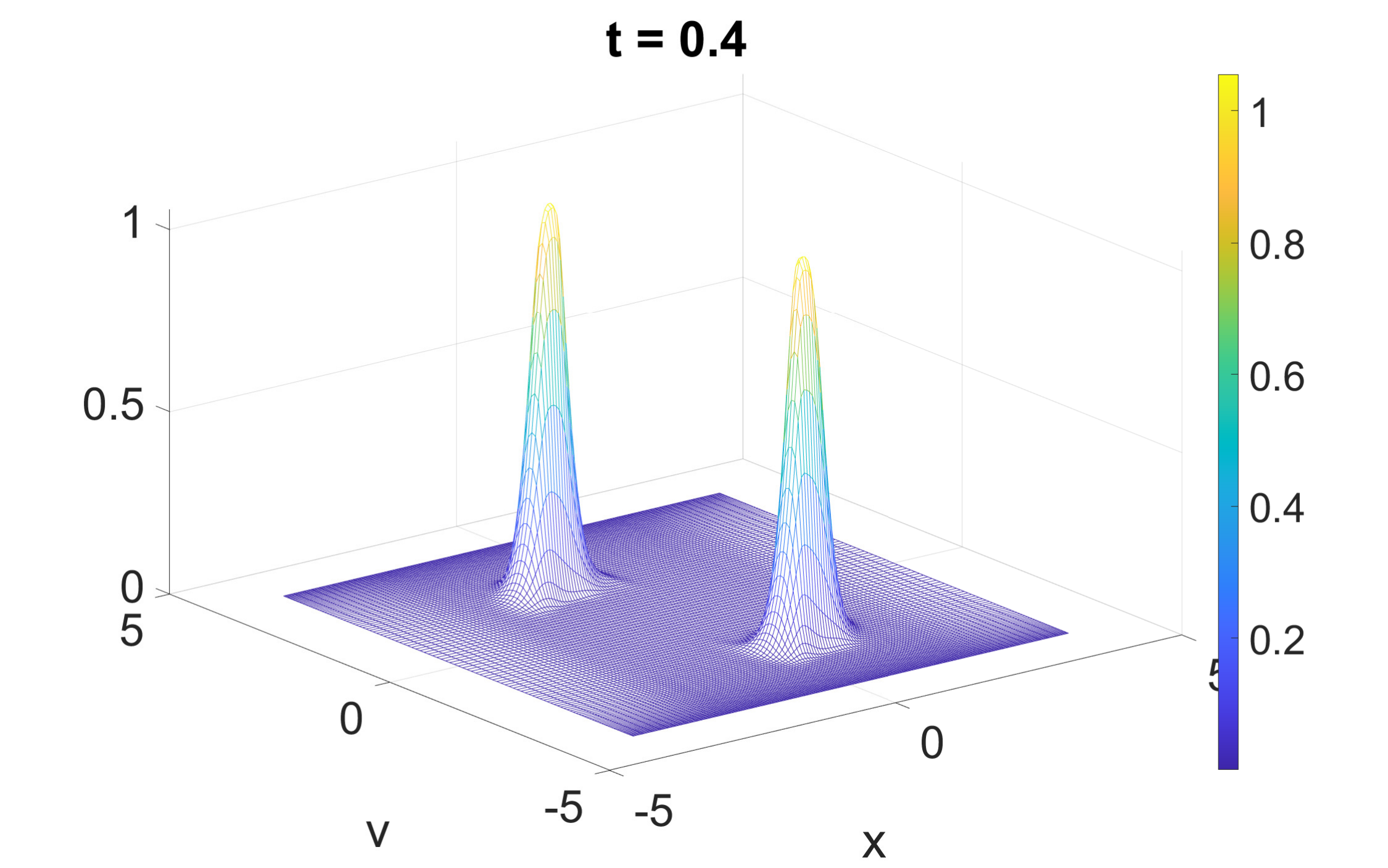}}
{\includegraphics[width=0.45\textwidth]{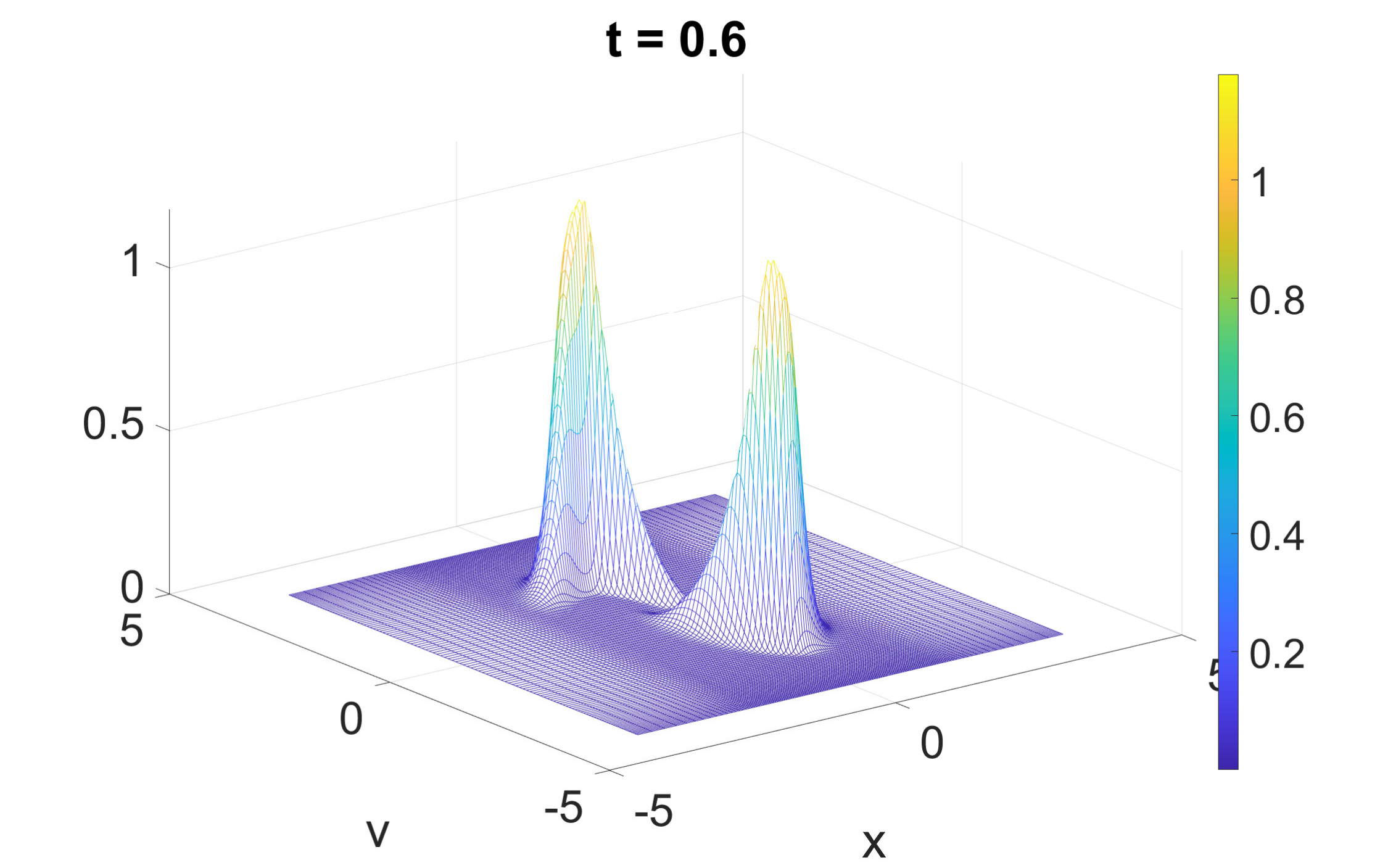}}
{\includegraphics[width=0.45\textwidth]{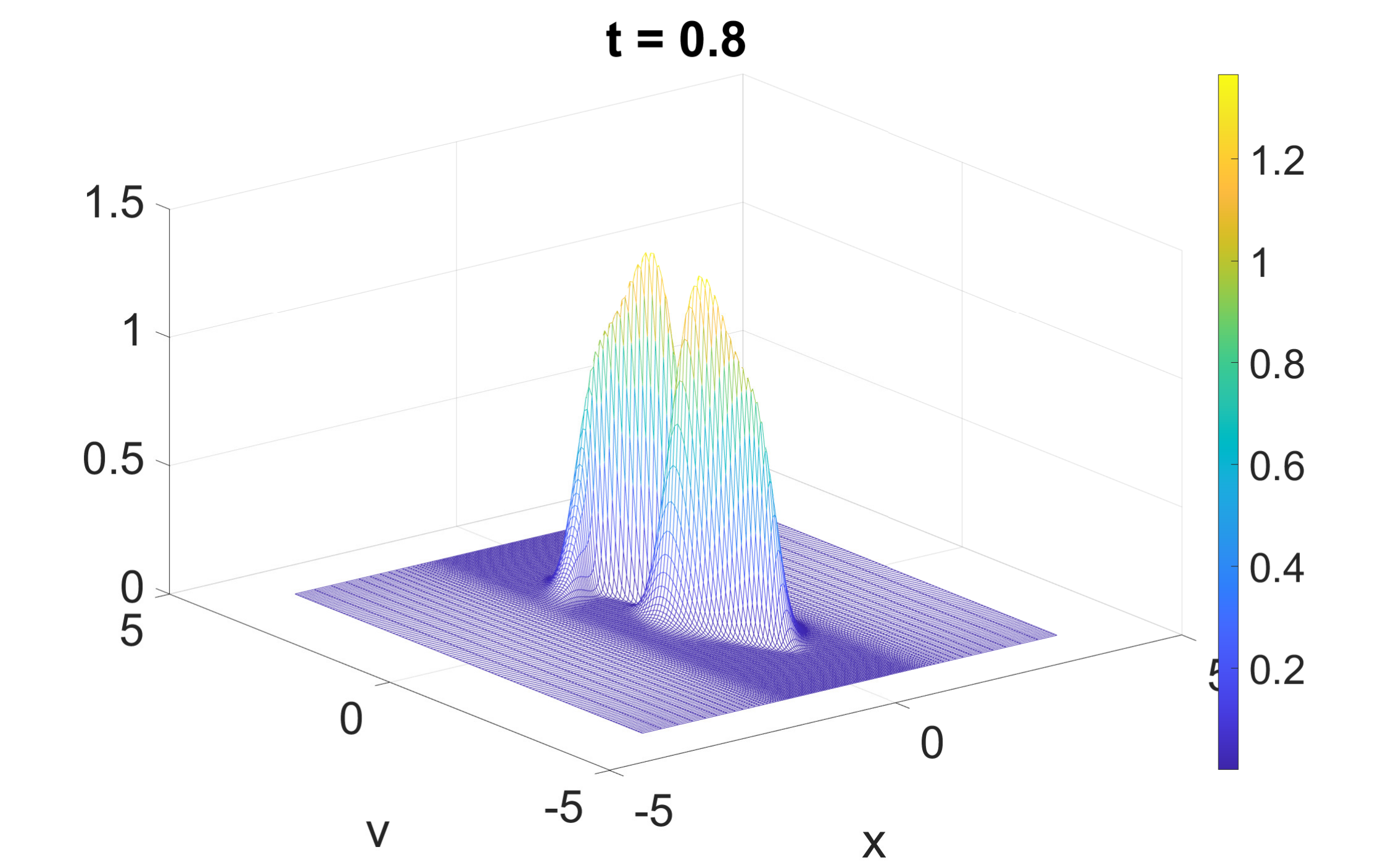}}
{\includegraphics[width=0.45\textwidth]{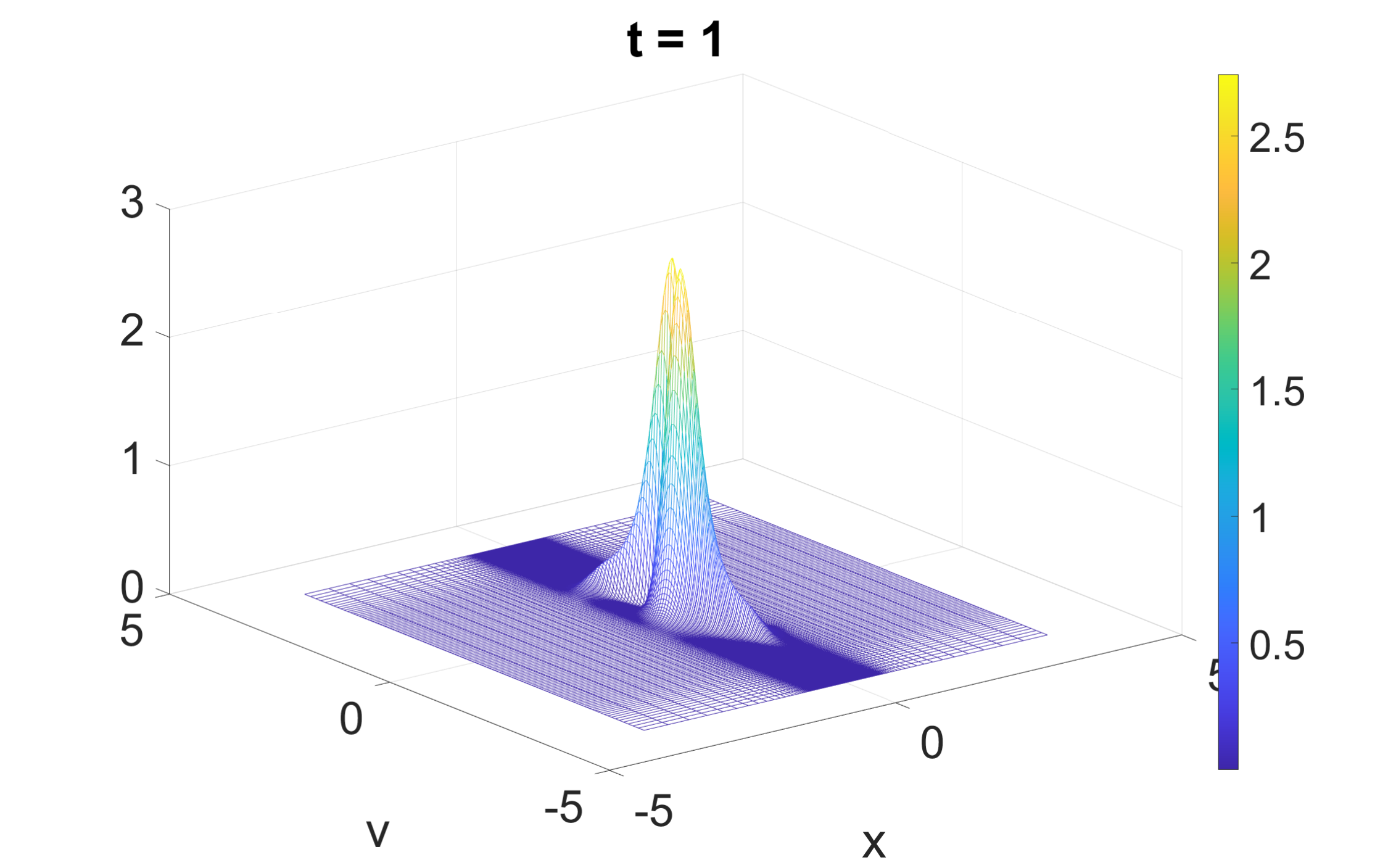}}
{\includegraphics[width=0.45\textwidth]{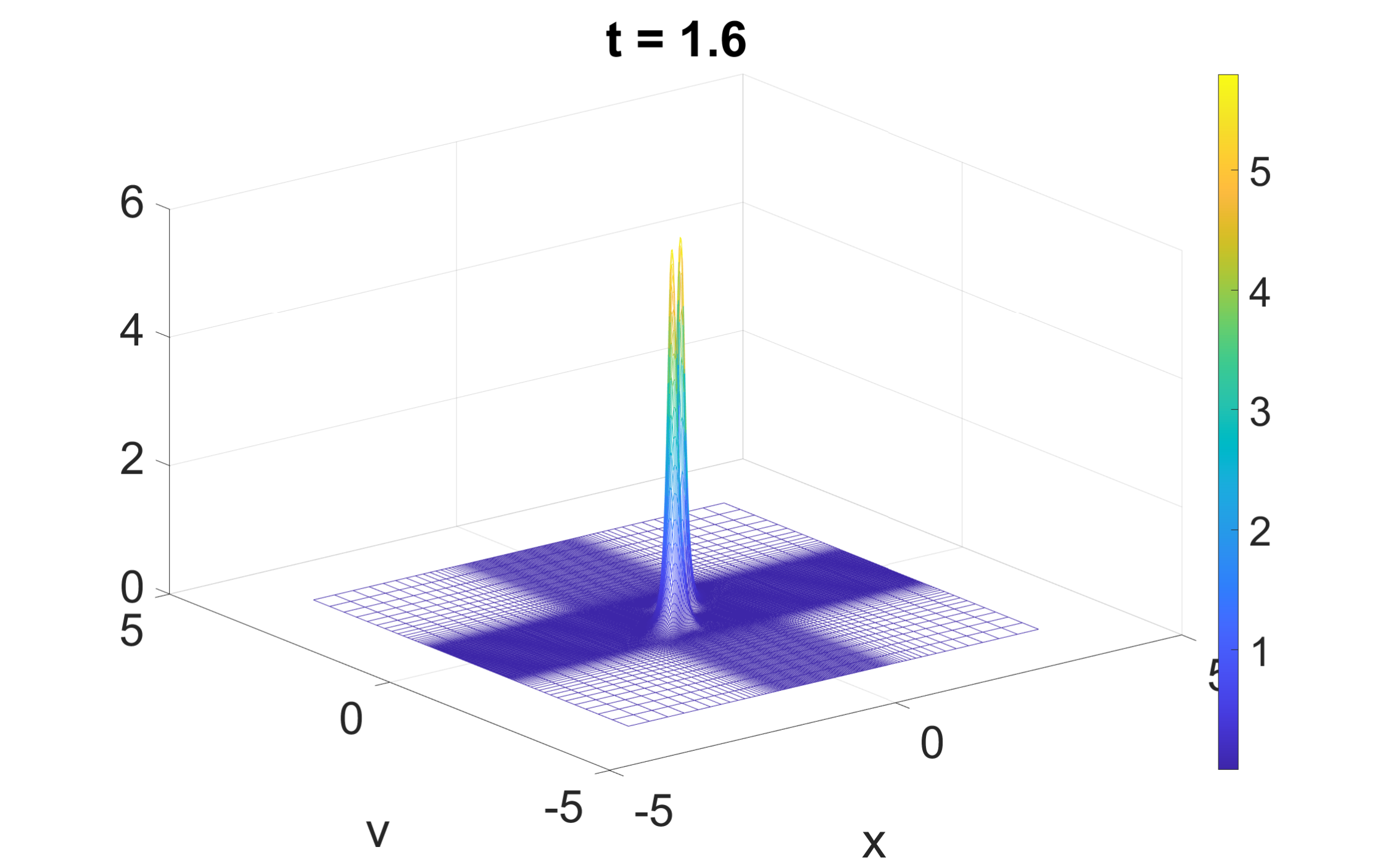}}
{\includegraphics[width=0.45\textwidth]{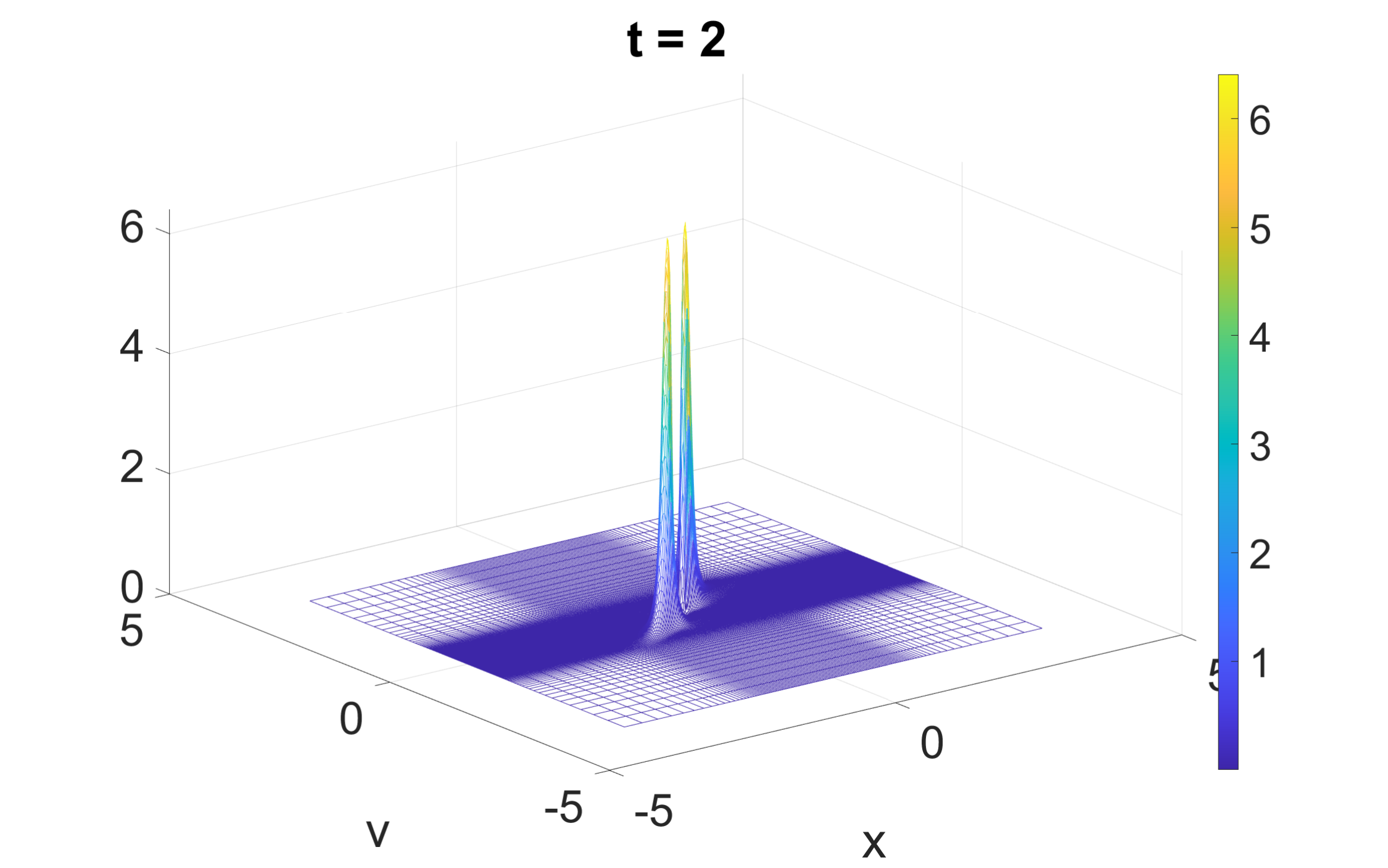}}
 \caption{Numerical solution of \eqref{eqn: gke} with kernel $W(v)=|v|^3$ for $N_x=121$, $N_v=121$, $\lambda = 4$. From top left to bottom right are time snapshots at $t=0.4, 0.6, 0.8, 1, 1.6, 2$.}
 \label{fig:gk3_ih_2_C2_2}
\end{figure}

To achieve this, we consider  the following initial condition, which consists of two symmetrically located bumps with respect to the origin:
\begin{equation}\label{eqn:f0}
    f_0(x,v) = e^{-a(x+c)^2}e^{-b(v-d)^2} + e^{-a(x-c)^2}e^{-b(v+d)^2}\,.
\end{equation}
Note that when the spatial dependence is disregarded, \eqref{eqn:f0} reduces to the two bump case discussed in Section~\ref{sec:ibu}. Consequently, the two bumps, situated with opposite 
spatial locations, will both 
move toward the origin while aggregating their velocity, as described in Section~\ref{sec:ibu}. It becomes intriguing to examine which direction exhibits blow-up first, if indeed there is a blow-up at all.

In practice, we set $L_x=4$, $L_v=4$, $a=b=6$, $c=1.5$, $d=2$ and then execute Algorithm~\ref{alg:SPMCSL-JKO} with $\delta_0 = 0.02$, $\delta = 0.5$, and various values of the strength of inelasticity $\lambda$ solving \eqref{eqn: gke}. Note that the values of $\delta_0$ and $\delta$ in this case differ from those used in the spatially homogeneous scenario. This difference is necessary because the small value of $\delta_0$ helps prevent oscillations caused by low resolution in regions with small values of $f$. A fixed time step of $\Delta t=0.05$ will be used throughout the evolution.

\begin{figure}[h!]
\centering
{\includegraphics[width=0.45\textwidth]{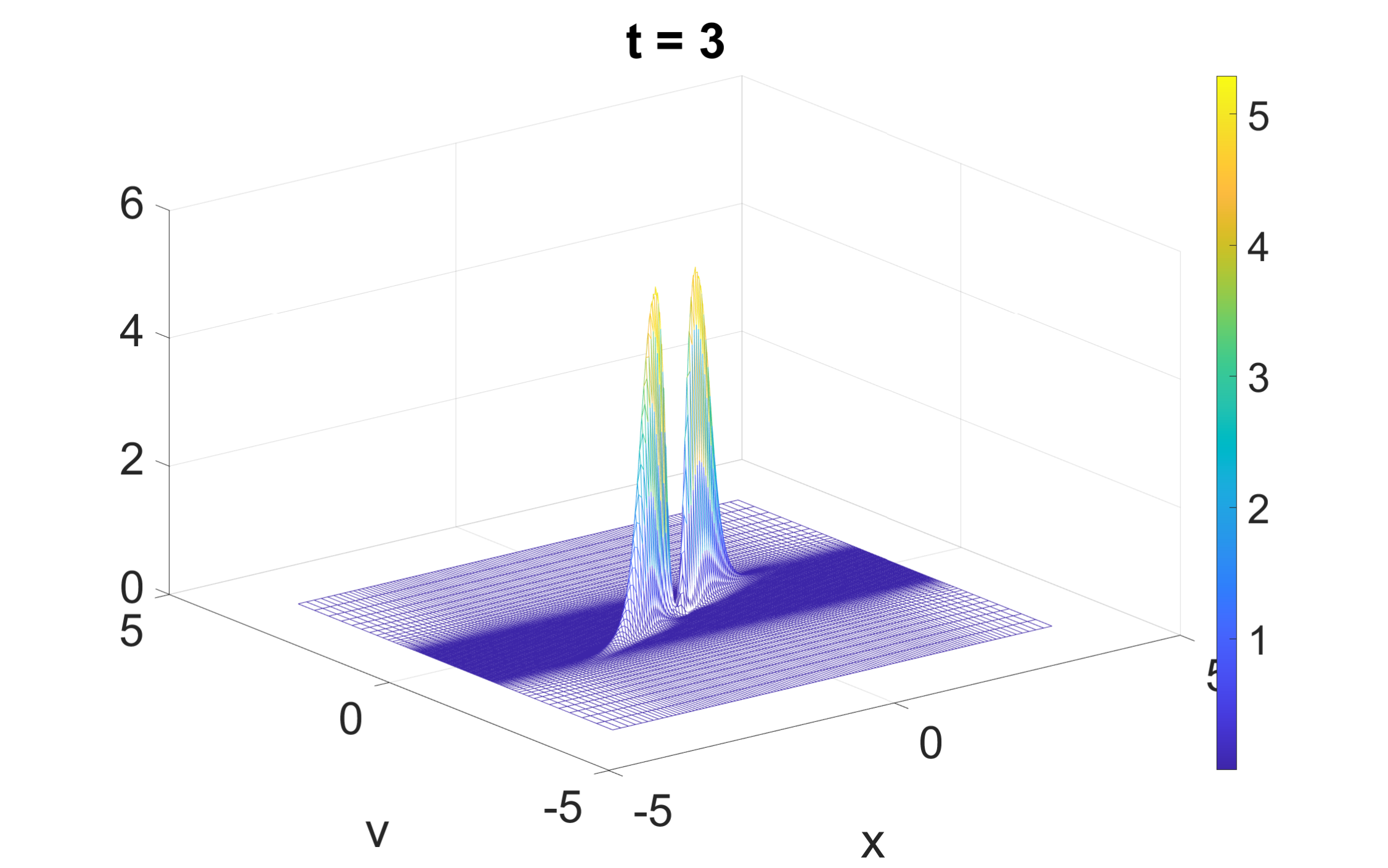}}
{\includegraphics[width=0.45\textwidth]{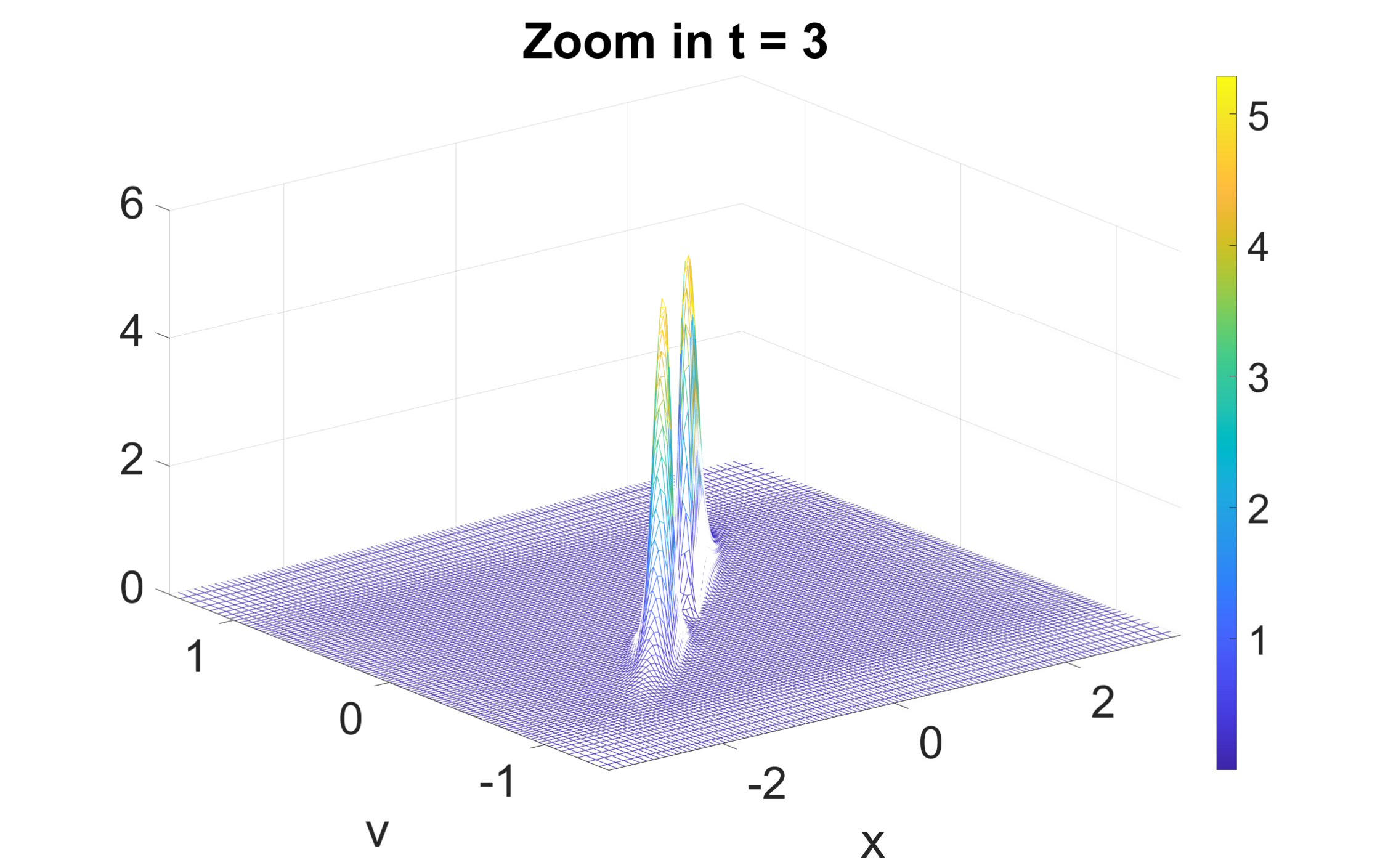}}
{\includegraphics[width=0.45\textwidth]{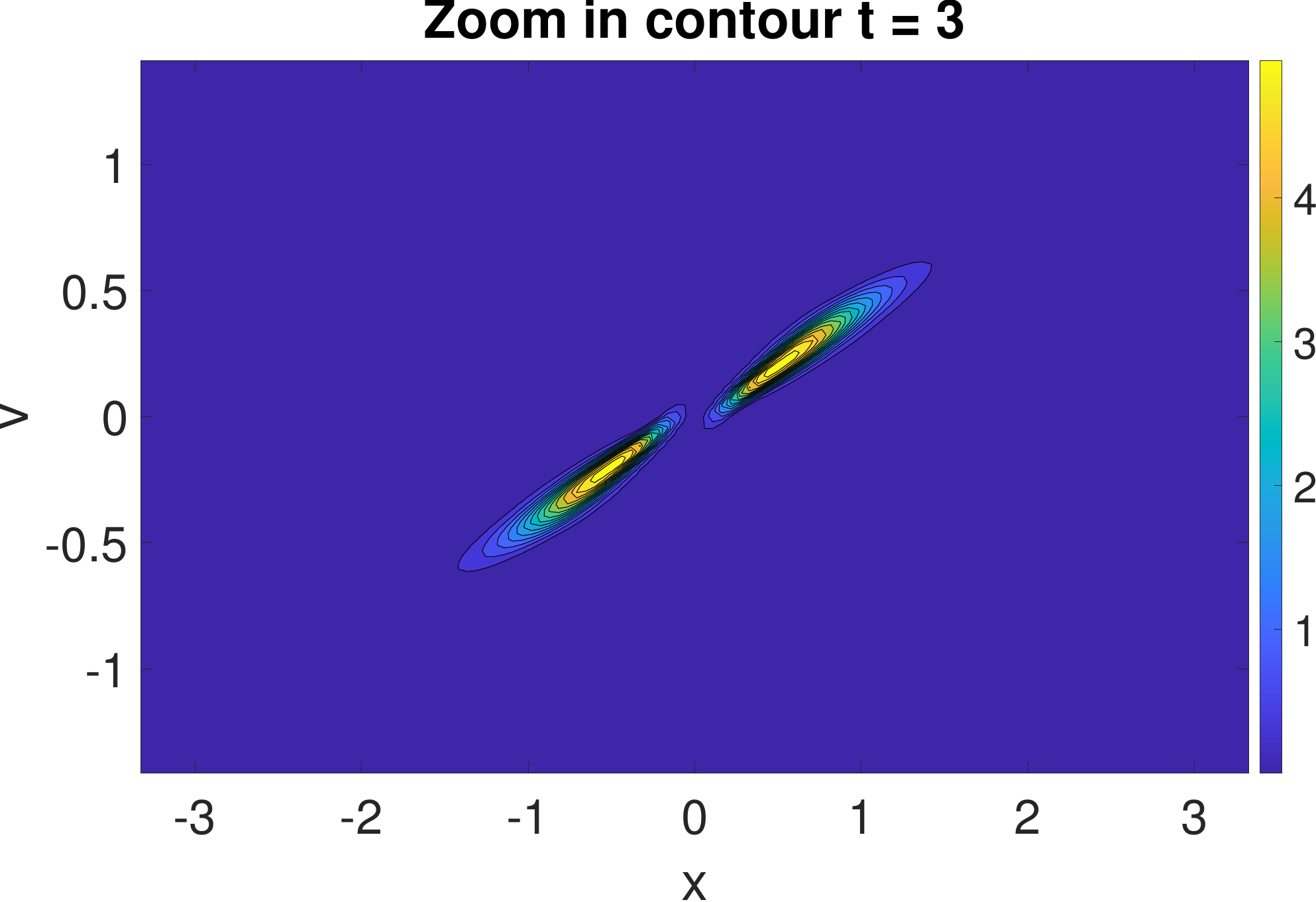}}
{\includegraphics[width=0.45\textwidth]{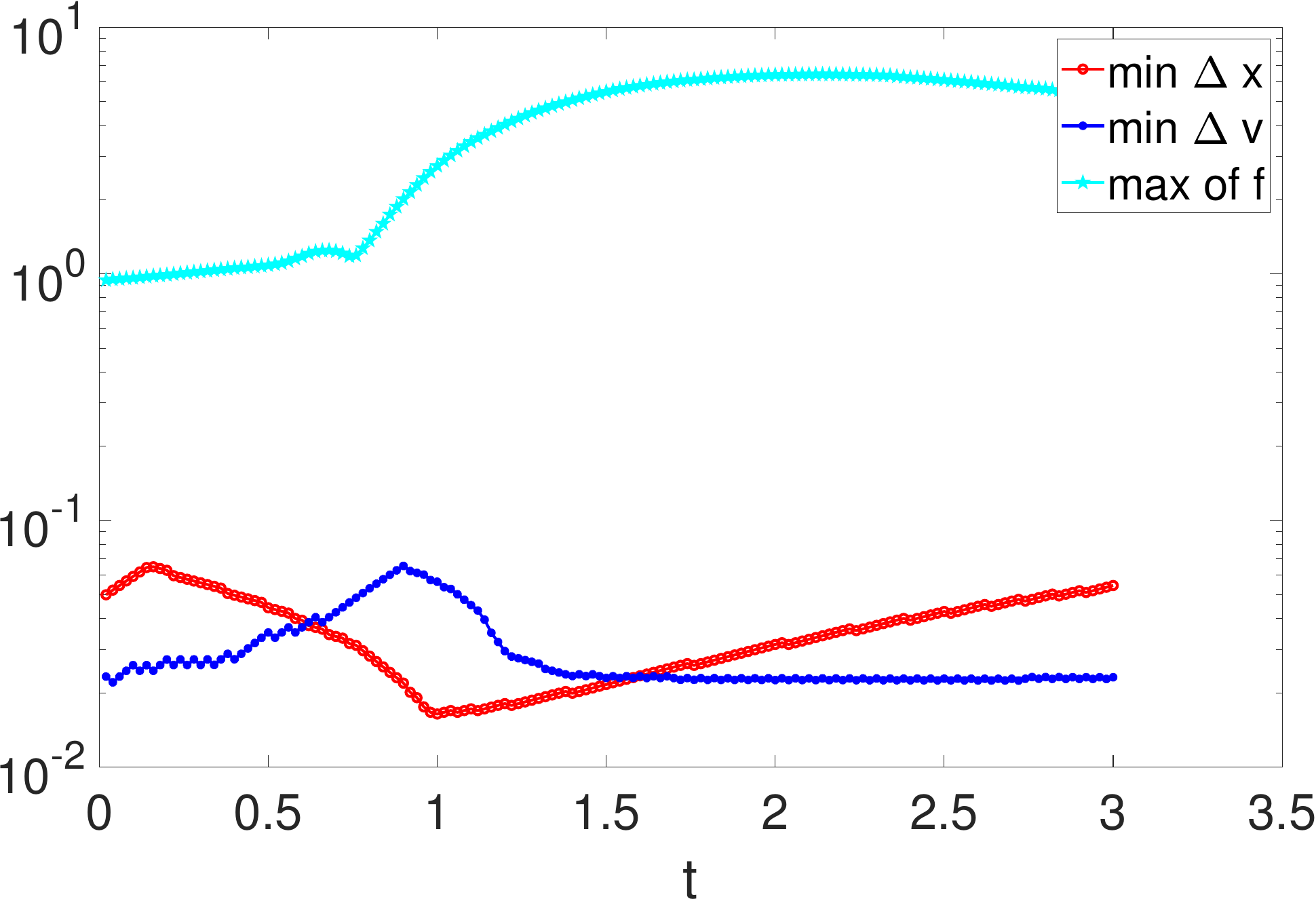}}
 \caption{Numerical solution of \eqref{eqn: gke} with kernel $W(v)=|v|^3$ for $N_x=121$, $N_v=121$, $\lambda = 4$. Top left: numerical solution $f(t=4,x,v)$. Top right: zoom-in plot of $f(t=4,x,v)$. Bottom left: contour plot of $f(t=4,x,v)$. Bottom right: record of minimum of $\Delta x$, $\Delta v$ and maximum of $f(t,x,v)$ at each time step.}
 \label{fig:gk3_ih_2_C2}
\end{figure}

Here, we present two representative numerical results for $\lambda=4$ and $\lambda=10$. Additional results for other values of $\lambda$ (specifically, $\lambda=2,6,8$) can be found in the Appendix~\ref{sec:sp_inhom_gamma_3}. For the case of $\lambda = 4$, the evolution of the solution unfolds in three distinct stages. At first, the two initial bumps move towards the origin. Subsequently, they converge at the origin, forming a concentration. Finally, due to the smallness of the concentration coefficient, the transport effect prevails over the concentration effect, leading to a separation of the two bumps after they pass through the origin. Such an evolution is illustrated in Figure~\ref{fig:gk3_ih_2_C2_2}. A more detailed plot at $t=3$ is presented in Figure~\ref{fig:gk3_ih_2_C2}, explicitly showing that the two bumps have bypassed each other. 
Most notably, the bottom right plot in Figure~\ref{fig:gk3_ih_2_C2} depicts the decreasing and then increasing changes in $\Delta x$ and $\Delta v$, revealing the concentration and separation motion of the two bumps. Additionally, the peak value of $f(t,x,v)$ undergoes a gradual increase and then decrease, indicating no blow-up. 

\begin{figure}[h!]
\centering
{\includegraphics[width=0.45\textwidth]{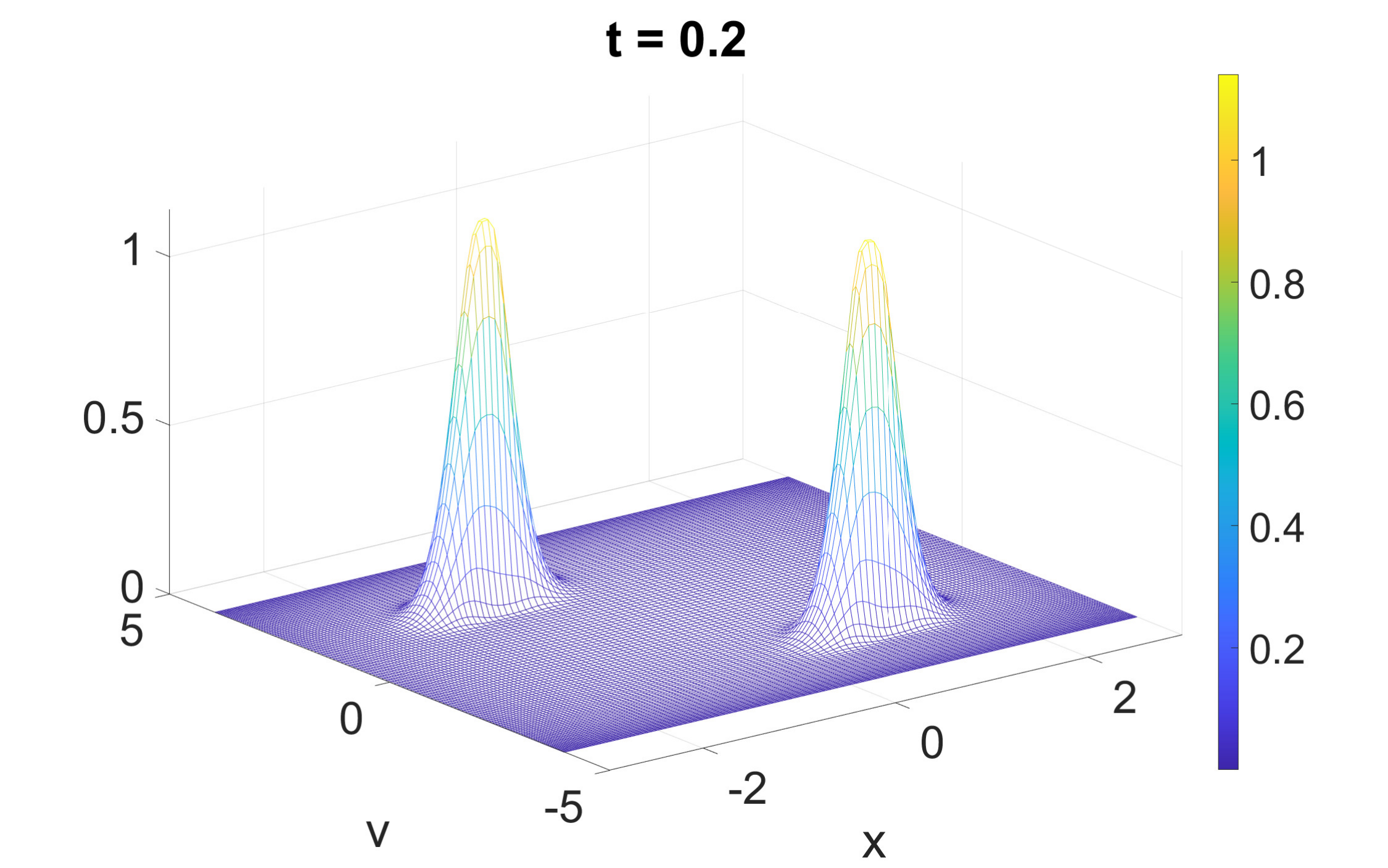}}
{\includegraphics[width=0.45\textwidth]{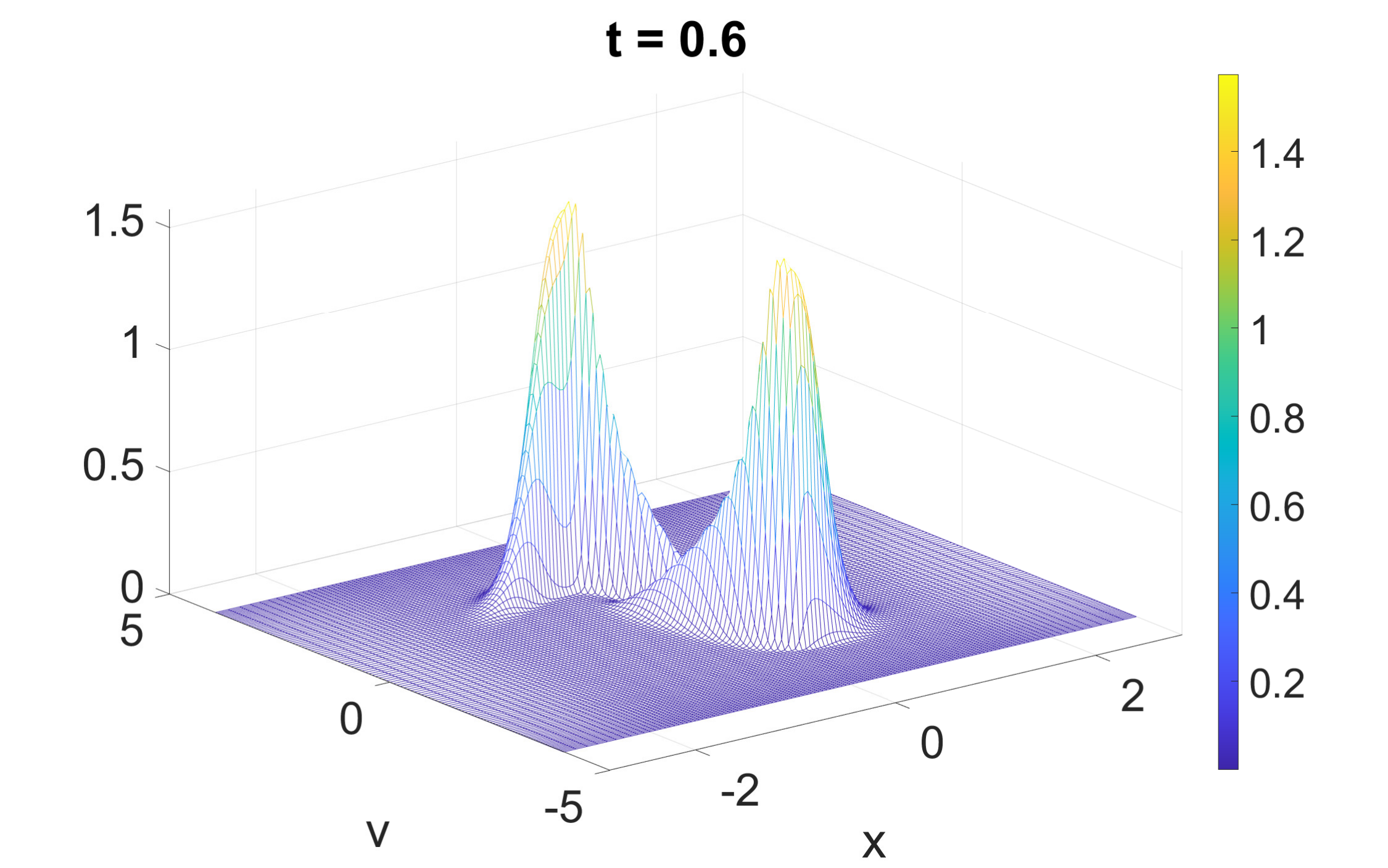}}
{\includegraphics[width=0.45\textwidth]{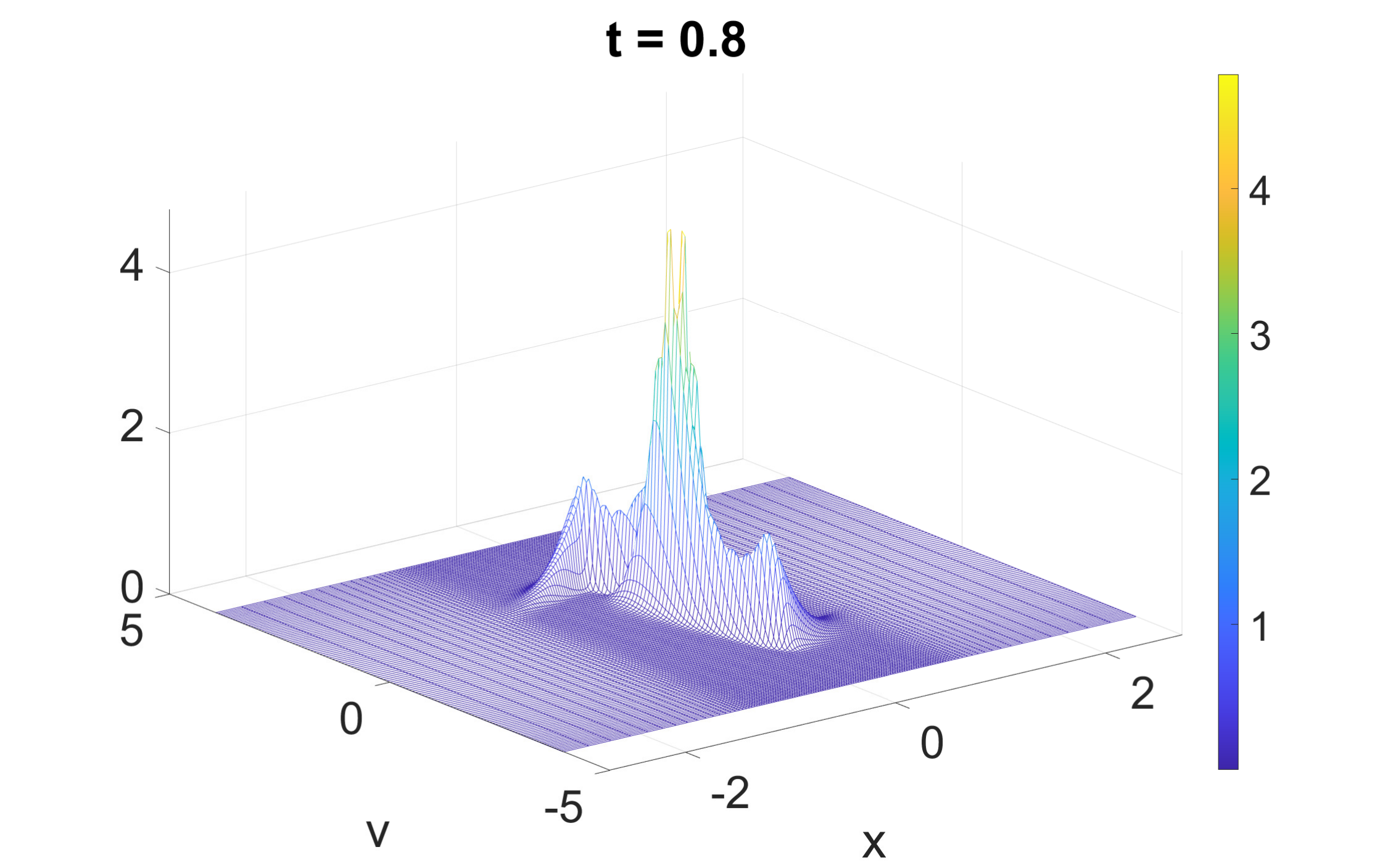}}
{\includegraphics[width=0.45\textwidth]{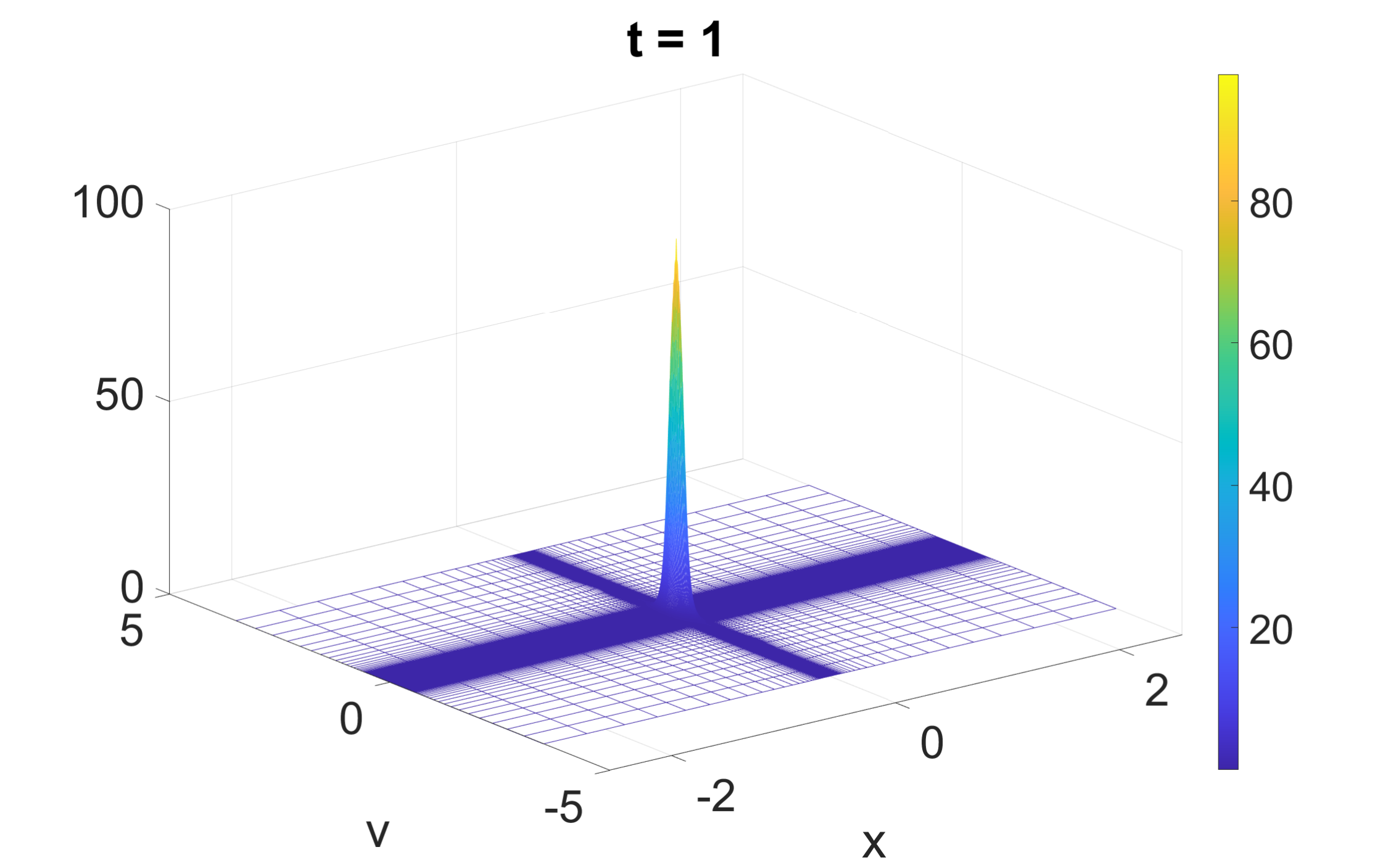}}
{\includegraphics[width=0.45\textwidth]{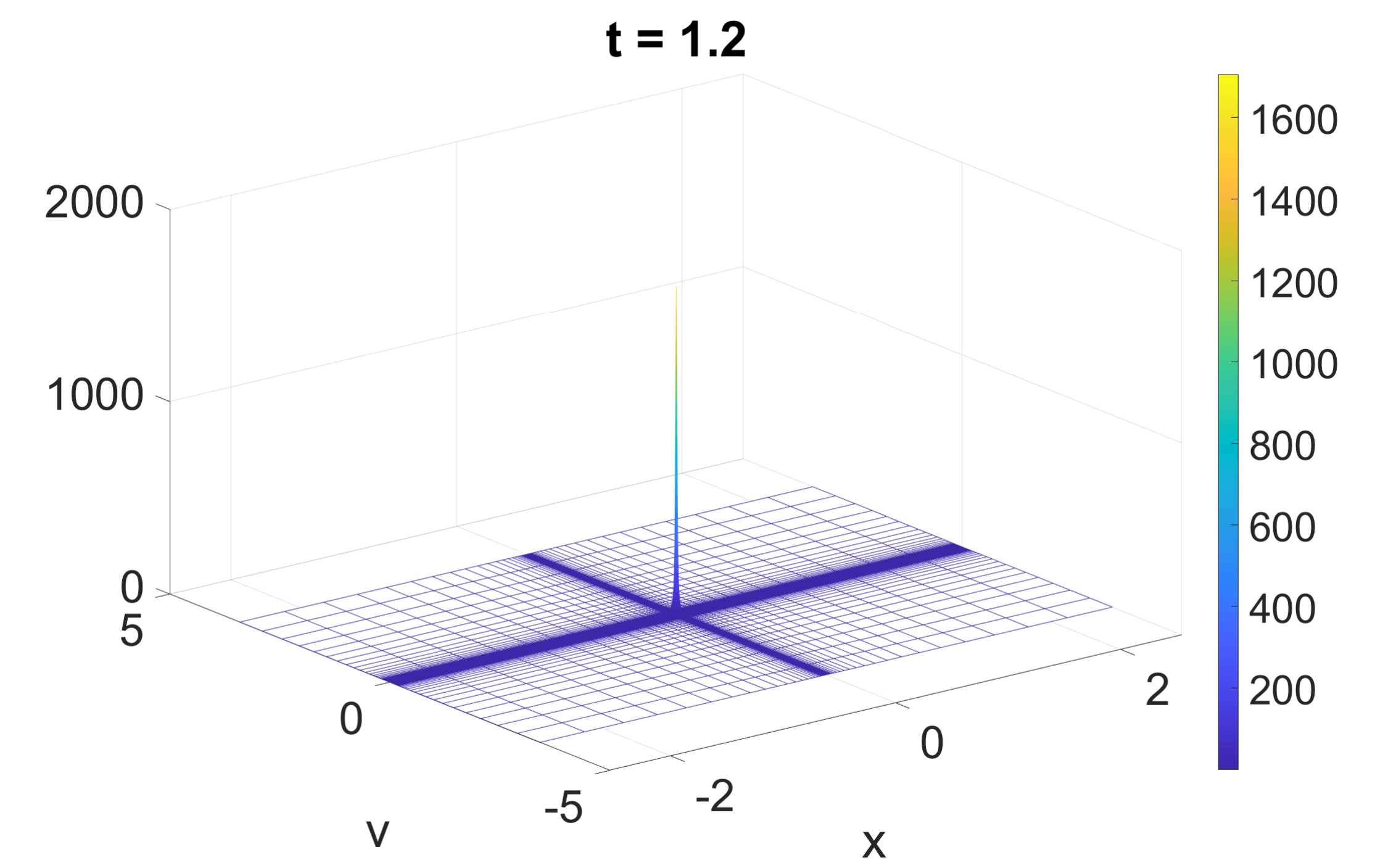}}
{\includegraphics[width=0.45\textwidth]{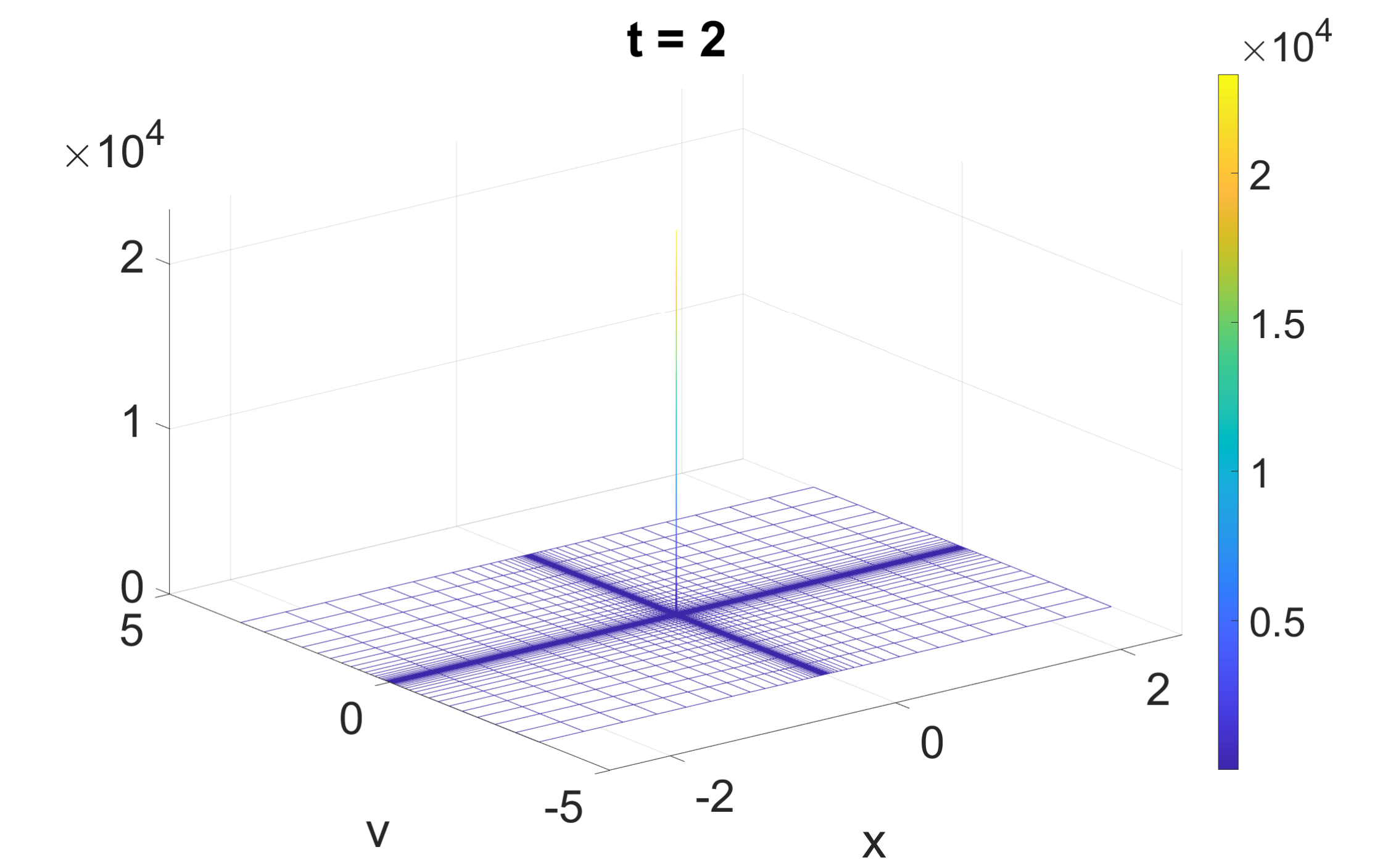}}
 \caption{Numerical solution of \eqref{eqn: gke} with kernel $W(v)=|v|^3$ for $N_x=121$, $N_v=121$, $\lambda = 10$. From top left to bottom right are numerical solution at $t=0.2, 0.6, 0.8, 1, 1.2, 2$.}
 \label{fig:gk3_ih_2_C5_2}
\end{figure}

For $\lambda=10$, the evolution of the solution also undergoes three stages. The first two stages are similar to those observed for $\lambda=4$, but with a more pronounced concentration effect. In this case, the solution rapidly forms a concentration that resembles a Dirac delta at the origin, as depicted in Figure~\ref{fig:gk3_ih_2_C5_2}.

\begin{figure}[h!]
\centering
{\includegraphics[width=0.45\textwidth]{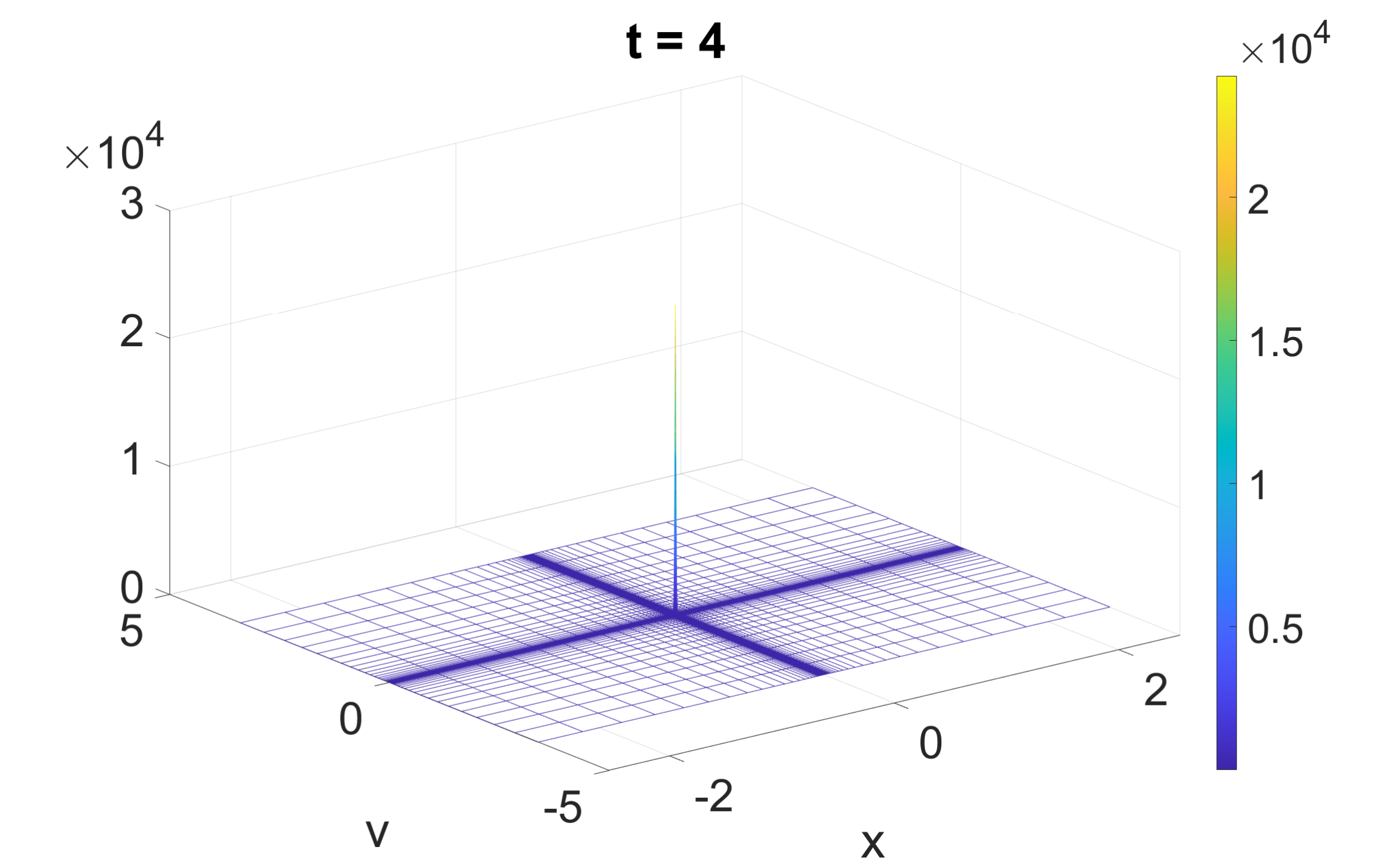}}
{\includegraphics[width=0.45\textwidth]{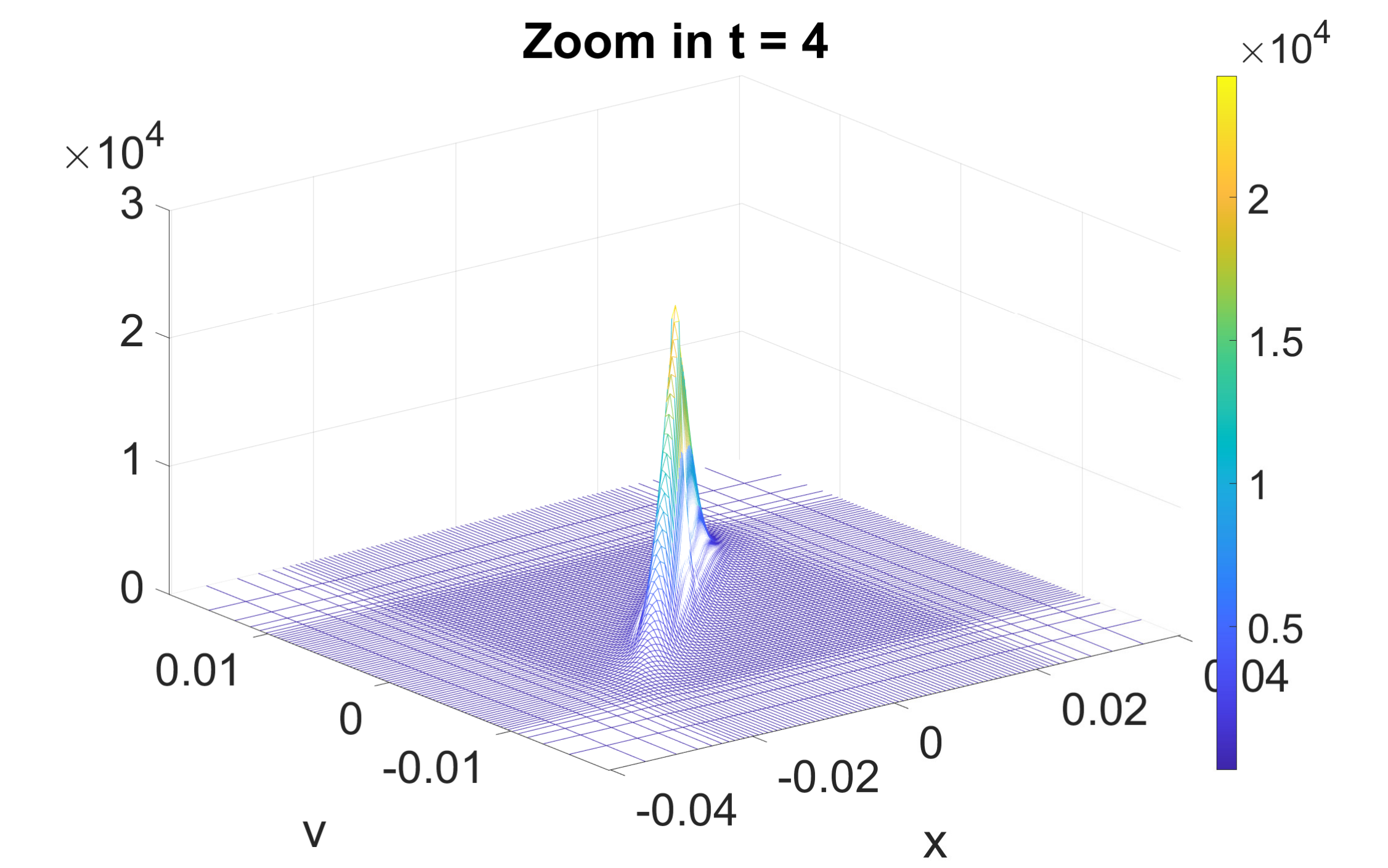}}
{\includegraphics[width=0.45\textwidth]{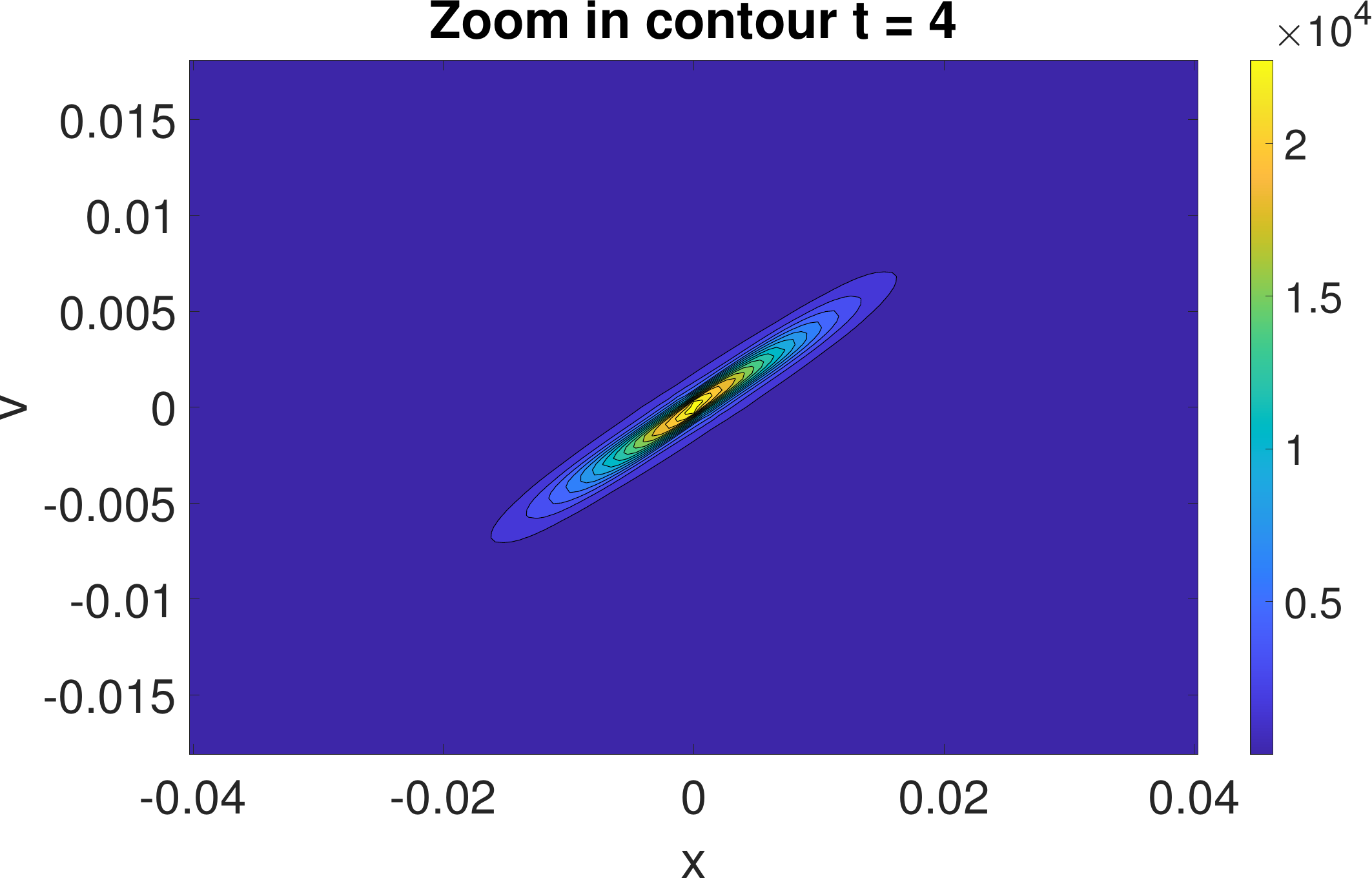}}
{\includegraphics[width=0.45\textwidth]{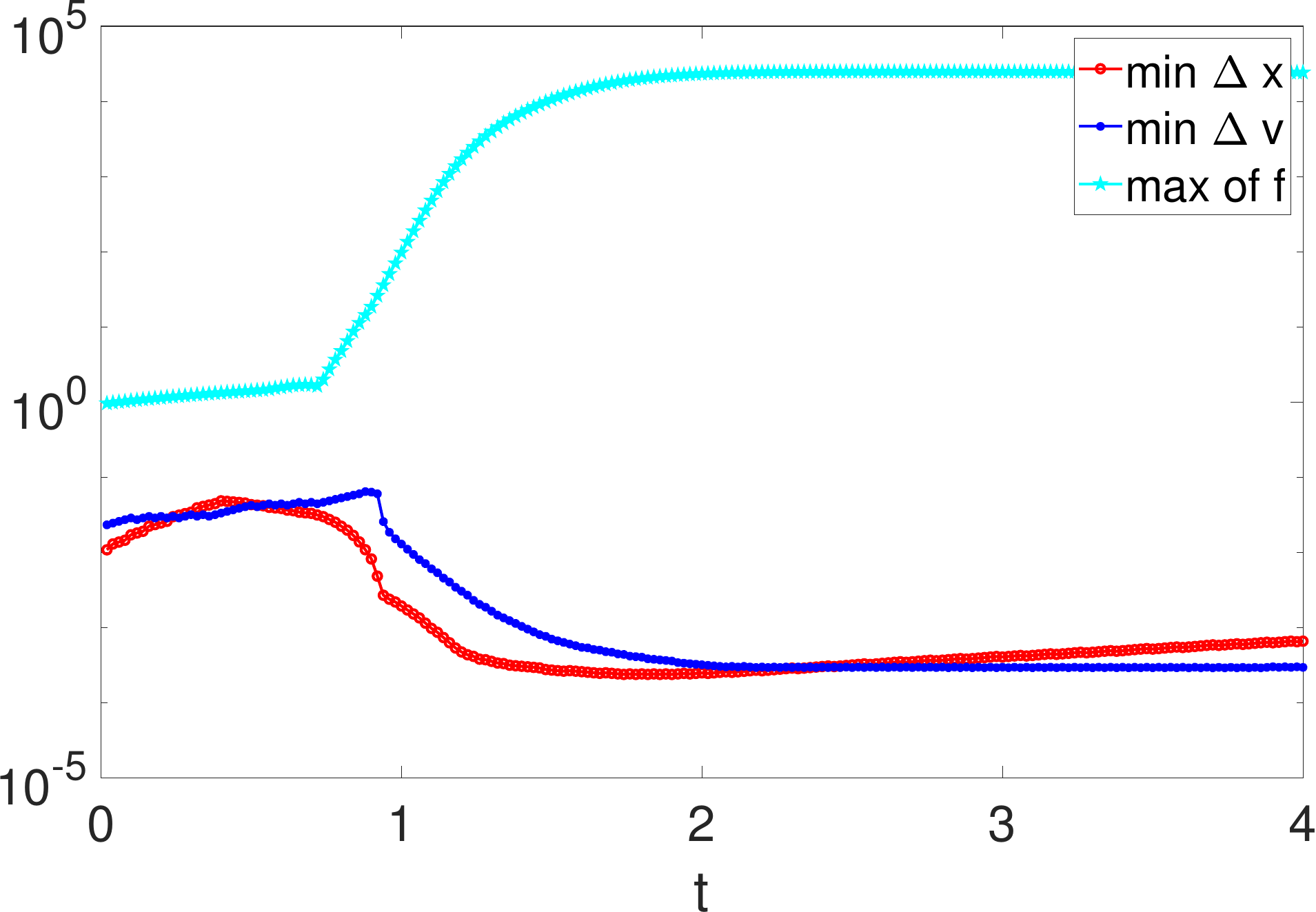}}
 \caption{Numerical solution of \eqref{eqn: gke} with kernel $W(v)=|v|^3$ for $N_x=121$, $N_v=121$, $\lambda = 10$. Top left: numerical solution $f(t=4,x,v)$. Top right: zoom-in plot of $f(t=4,x,v)$. Bottom left: contour plot of $f(t=4,x,v)$. Bottom right: record of minimum of $\Delta x$, $\Delta v$ and maximum of $f(t,x,v)$ at each time step.}
 \label{fig:gk3_ih_2_C5}
\end{figure}

While a stronger concentration effect might suggest a higher likelihood of a blow-up solution, we still observe that after an extended mixing time at the origin, the solution eventually spreads out along the $x$ direction due to the transport effect. This can be seen in the zoom-in plot and contour plot on the top right and bottom left of Fig.~\ref{fig:gk3_ih_2_C5}. This observation is further supported by the the bottom right plot in  Fig.~\ref{fig:gk3_ih_2_C5}, where the minimum of $\Delta x$ initially decreases and then increases, implying the concentration on $x$ at first and finally spreading out by the transport effect. 

We numerically verify the persistence of properties introduced in Proposition~\ref{prop:properties} over time for for $\gamma=3$ and $\lambda=4$ and $\lambda=10$, with the results plotted in the first two rows of Figure~\ref{fig:gk3_properties} for the settings of Figure~\ref{fig:gk3_ih_2_C2_2} and Figure~\ref{fig:gk3_ih_2_C5_2} respectively.
In Figure~\ref{fig:gk3_properties}, the total mass and momentum are preserved. The energy (second order moment), third order moment, fourth order moment and the entropy are decreasing over time. Here the entropy is defined by $U(f) = - f \log (f)$. 

\subsection{Spatially inhomogeneous case with $\gamma=3/2$}\label{sec:sp_inhom_gamma_1p5}
When $\gamma<2$, a blow-up solution can emerge within finite time in spatially homogeneous cases. When considering transport, examples with finite-time blow-up in velocity space can be readily constructed. However, it is nontrivial to devise an appropriate initial condition that leads to blow-up in both the $x$ and $v$ simultaneously or blow-up in the $x$ direction before the $v$ direction. For illustrative purposes, we provide examples of finite-time blow-up in the $v$ direction for $\gamma=3/2$.

\begin{figure}[h!]
\centering
{\includegraphics[width=0.45\textwidth]{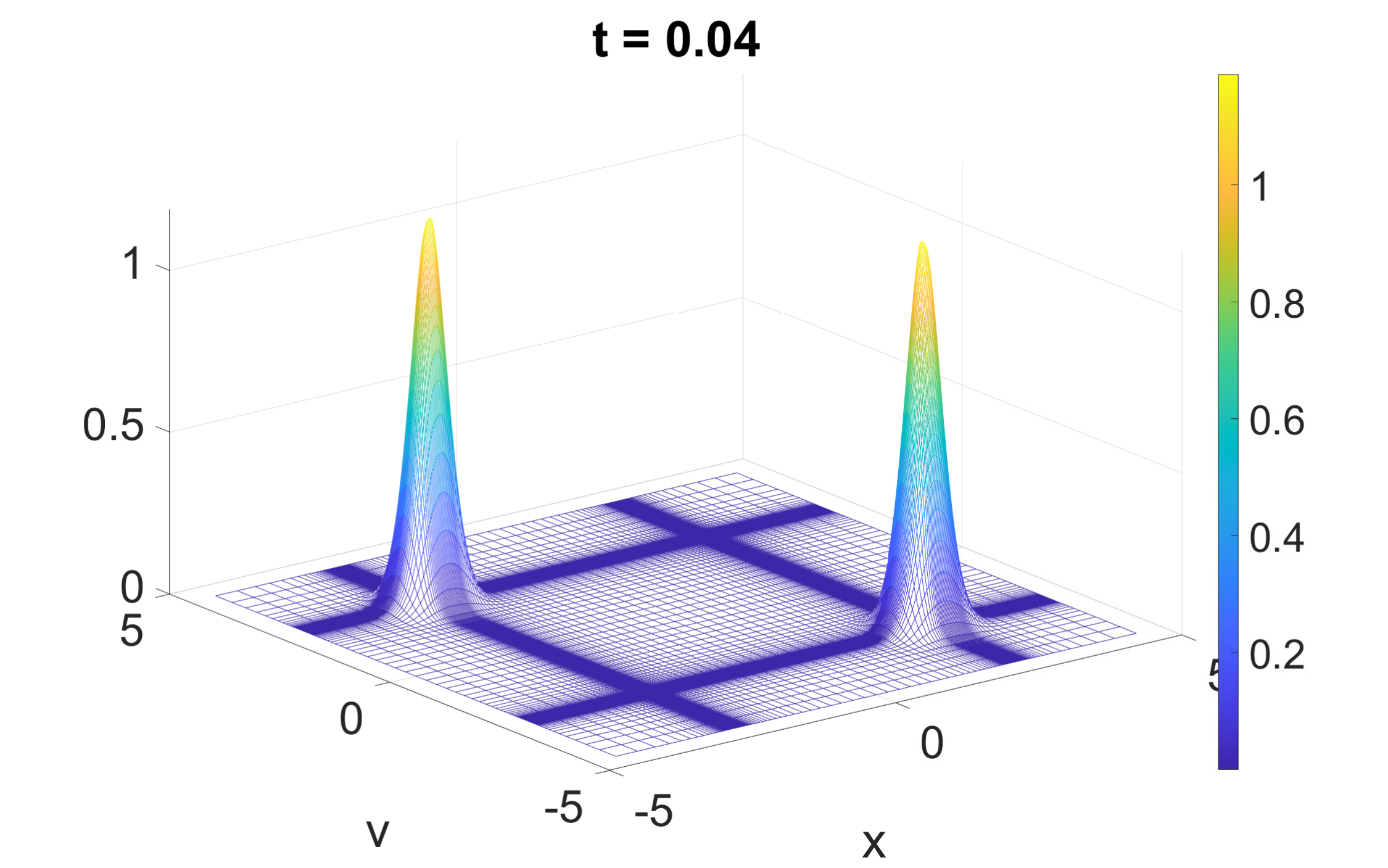}}
{\includegraphics[width=0.45\textwidth]{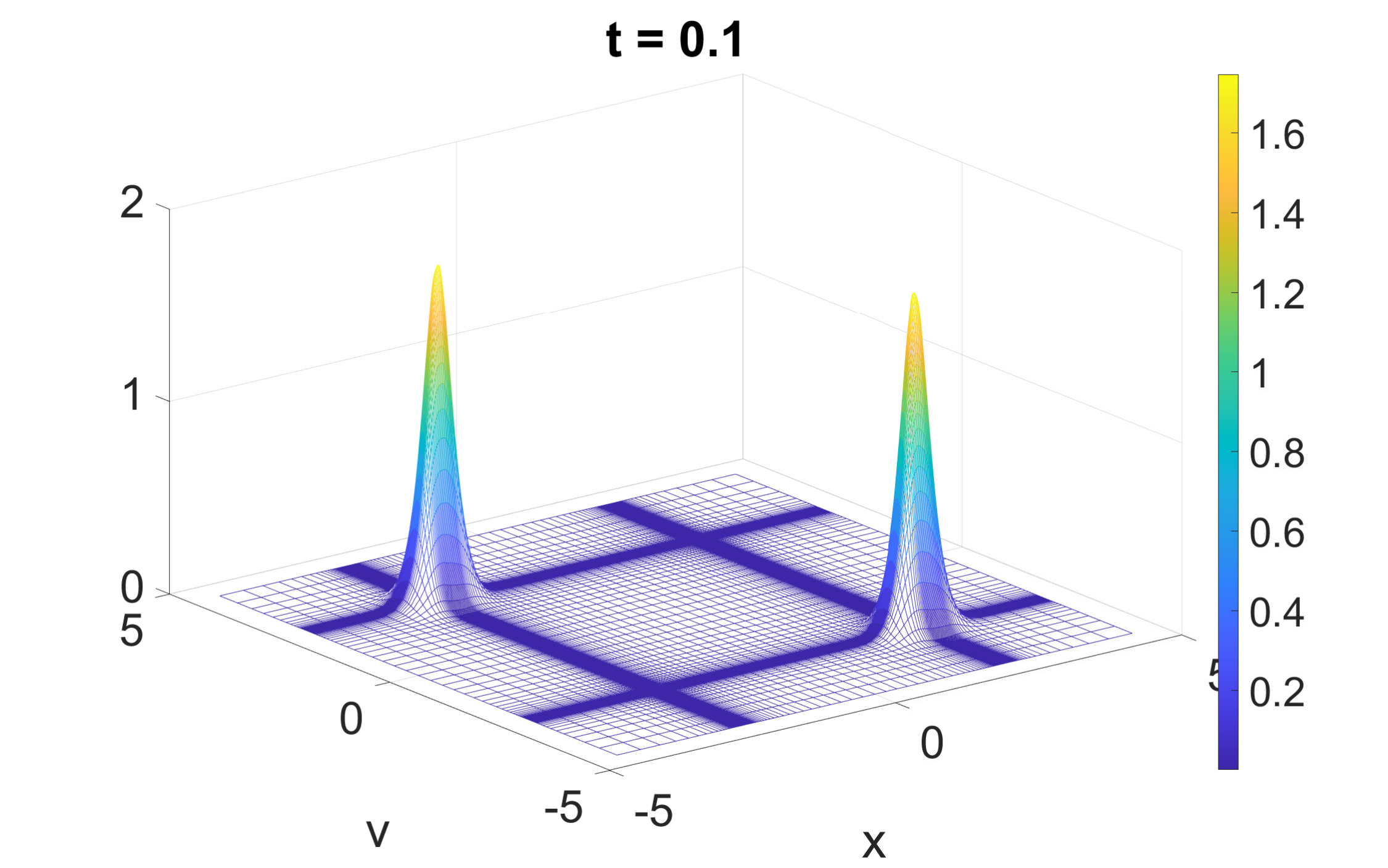}}

{\includegraphics[width=0.45\textwidth]{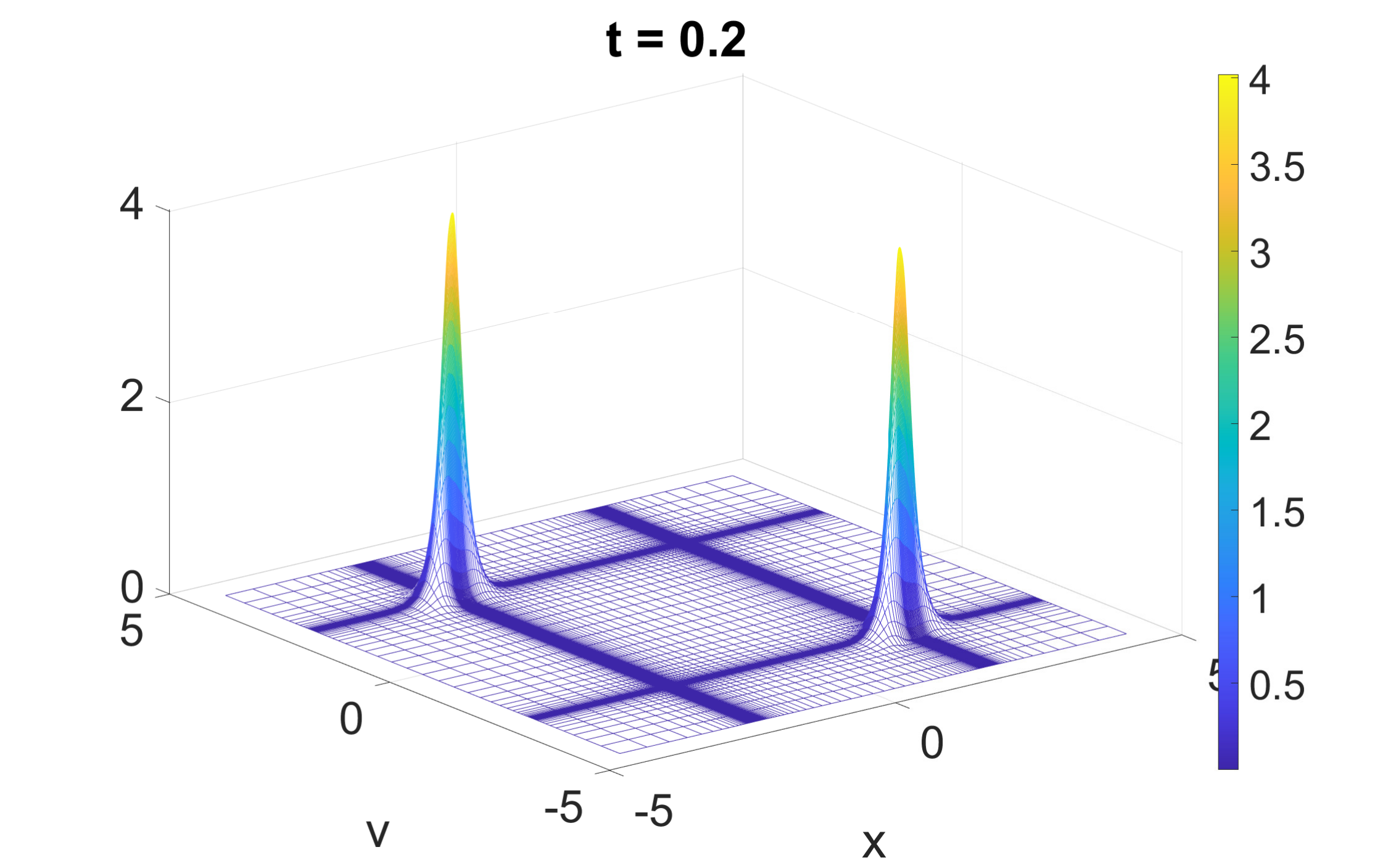}}
{\includegraphics[width=0.45\textwidth]{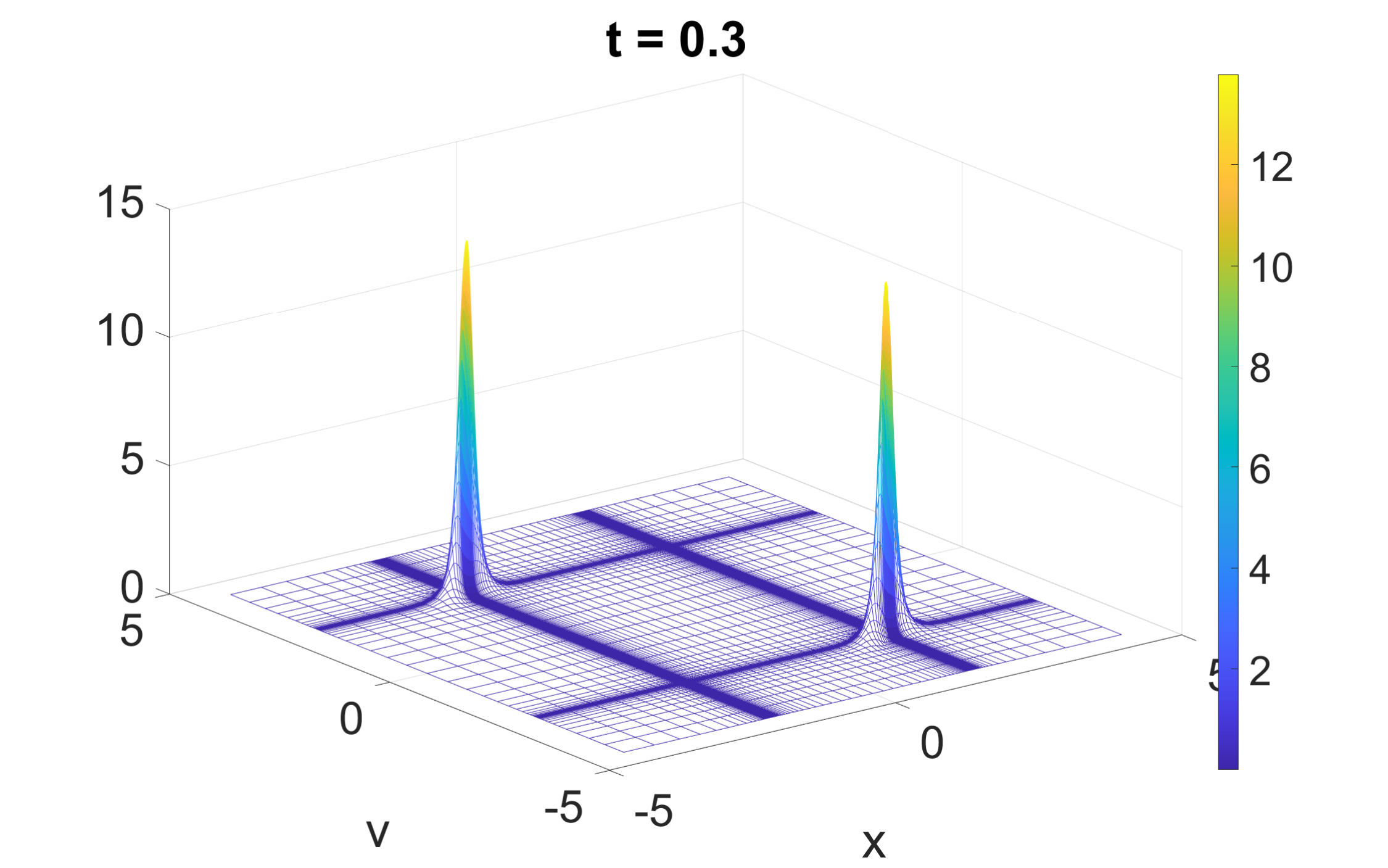}}

{\includegraphics[width=0.45\textwidth]{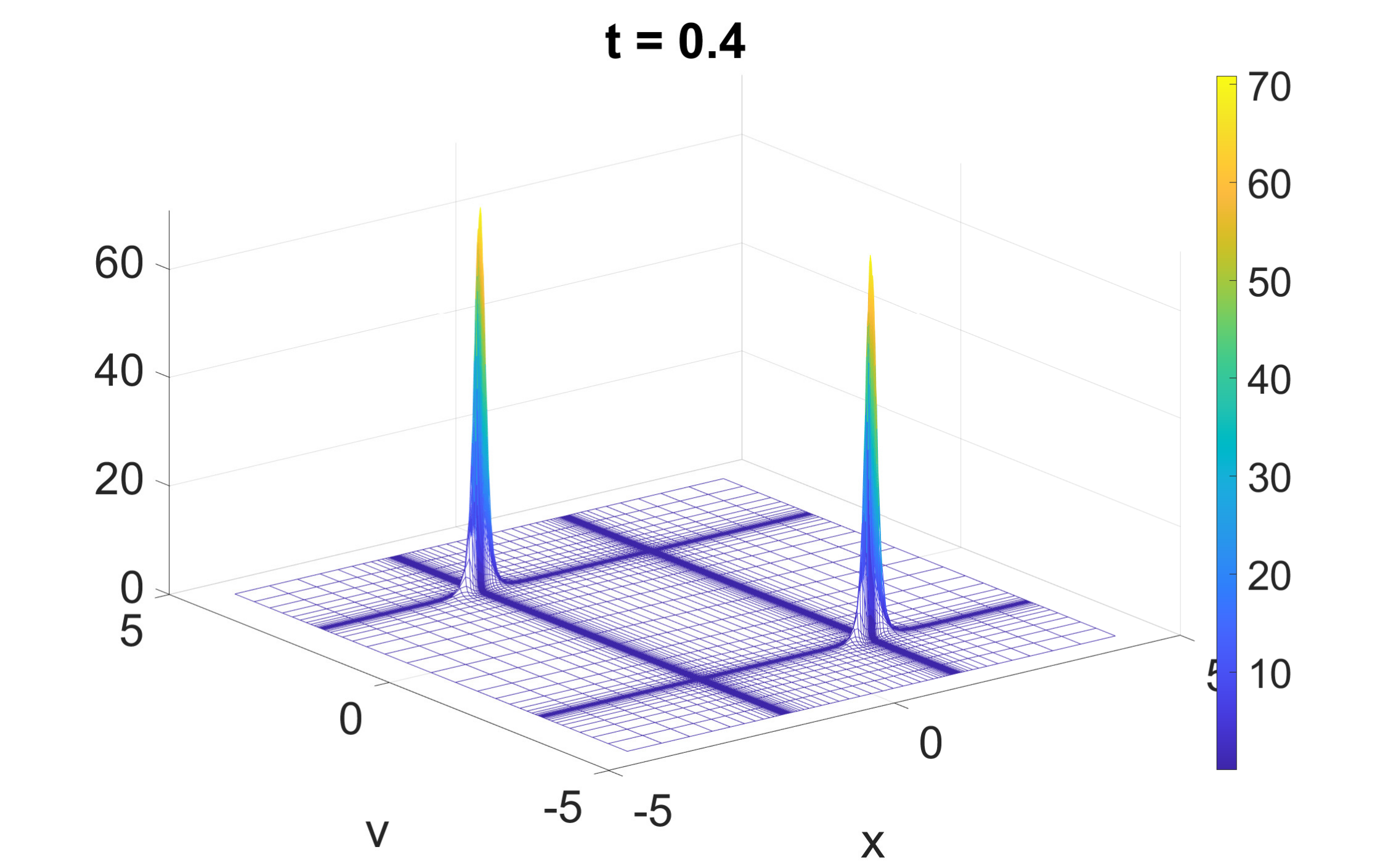}}
{\includegraphics[width=0.45\textwidth]{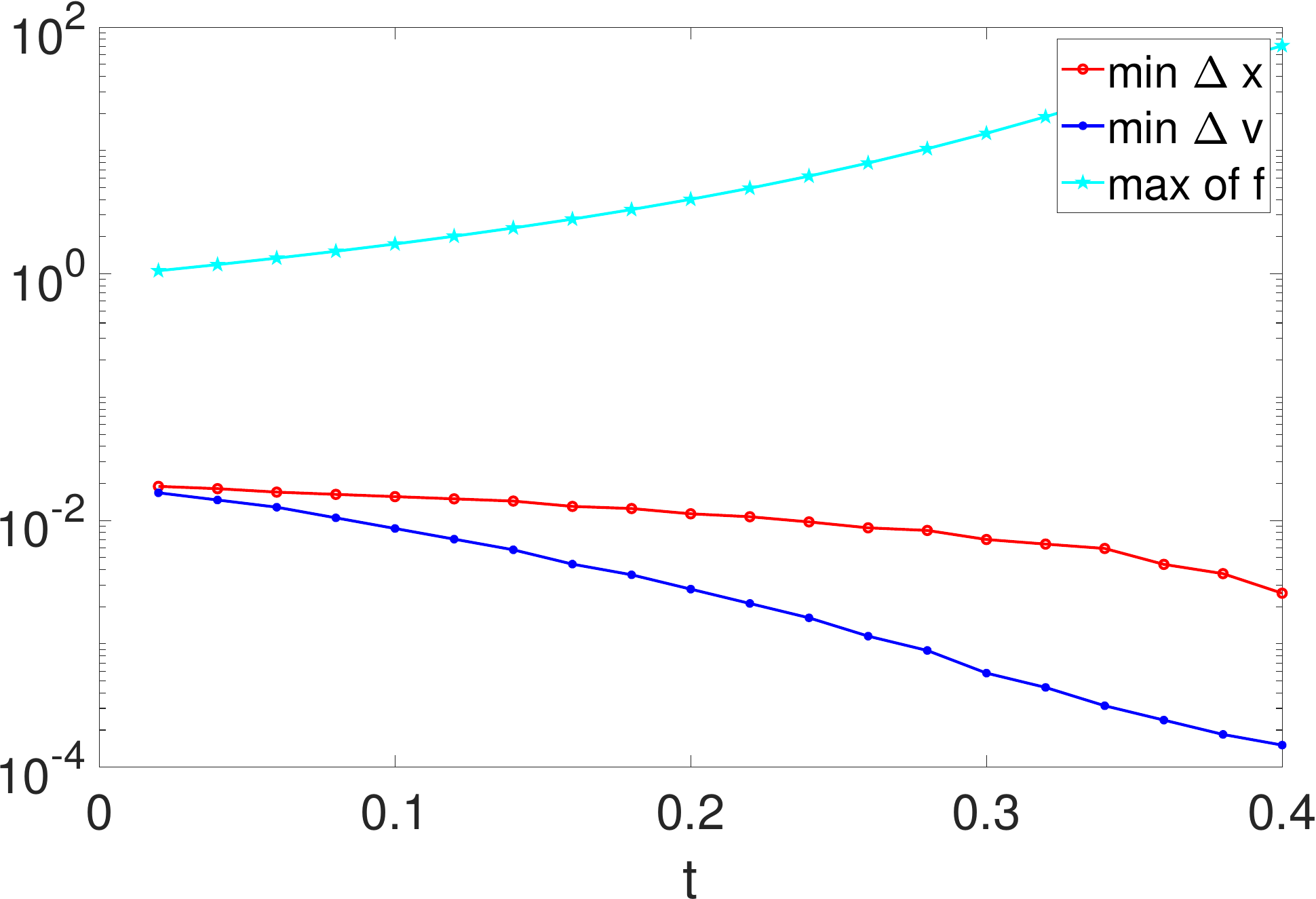}}

 \caption{Numerical solution of \eqref{eqn: gke} with kernel $W(v)=|v|^{3/2}$ for $N_x=121$, $N_v=121$, $a=b=6$, $c=d=2.5$ and $L_x=L_v=5$. Top left to bottom left are the numerical solution at $t = 0.04, 0.1, 0.2, 0.3, 0.4$ and the bottom right is the record of minimum of $\Delta x$, $\Delta v$ and maximum of $f(t,x,v)$ at each time step.}
 \label{fig:gk1p5_C5}
\end{figure}

\begin{figure}[h!]
\centering
{\includegraphics[width=0.45\textwidth]{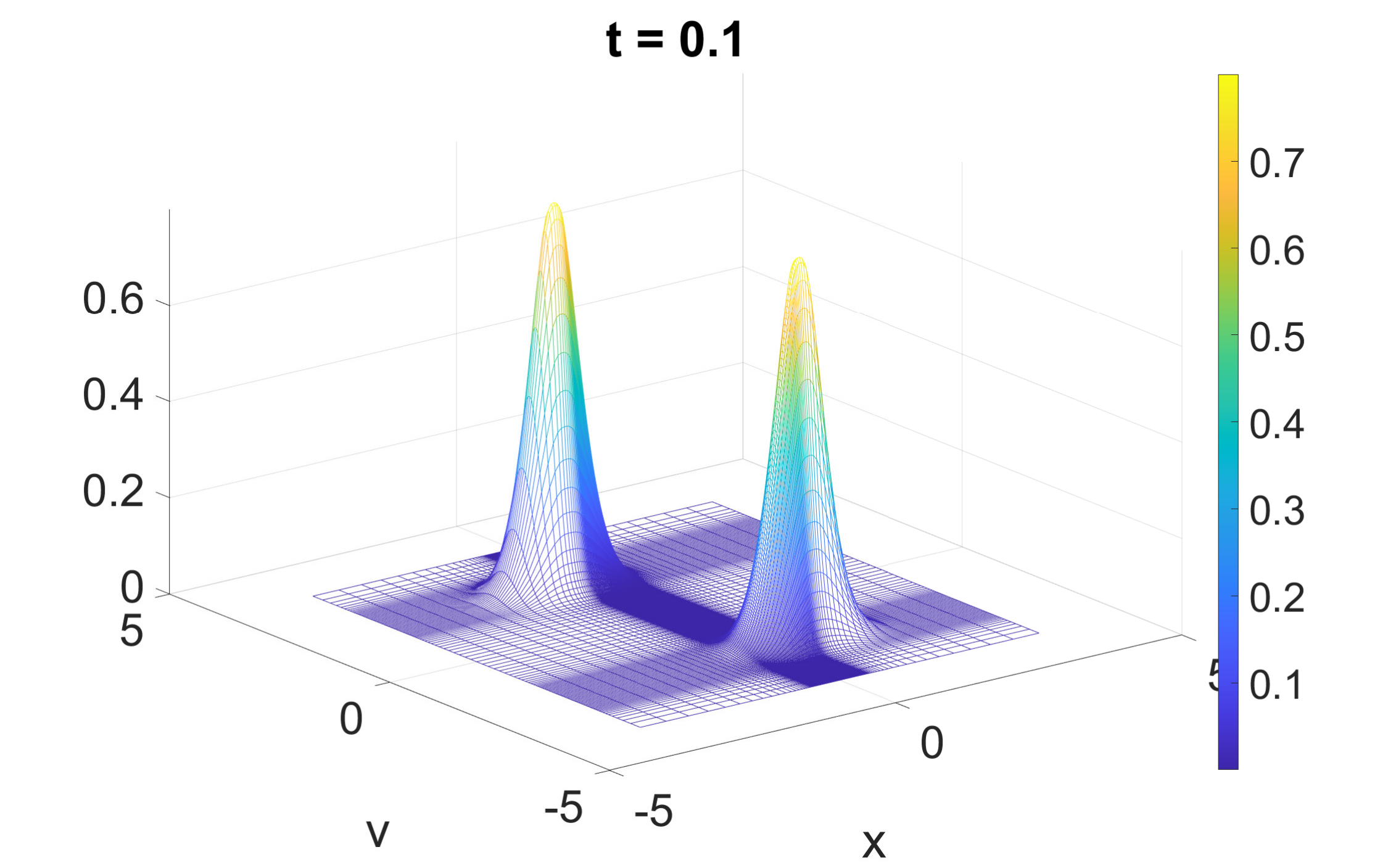}}
{\includegraphics[width=0.45\textwidth]{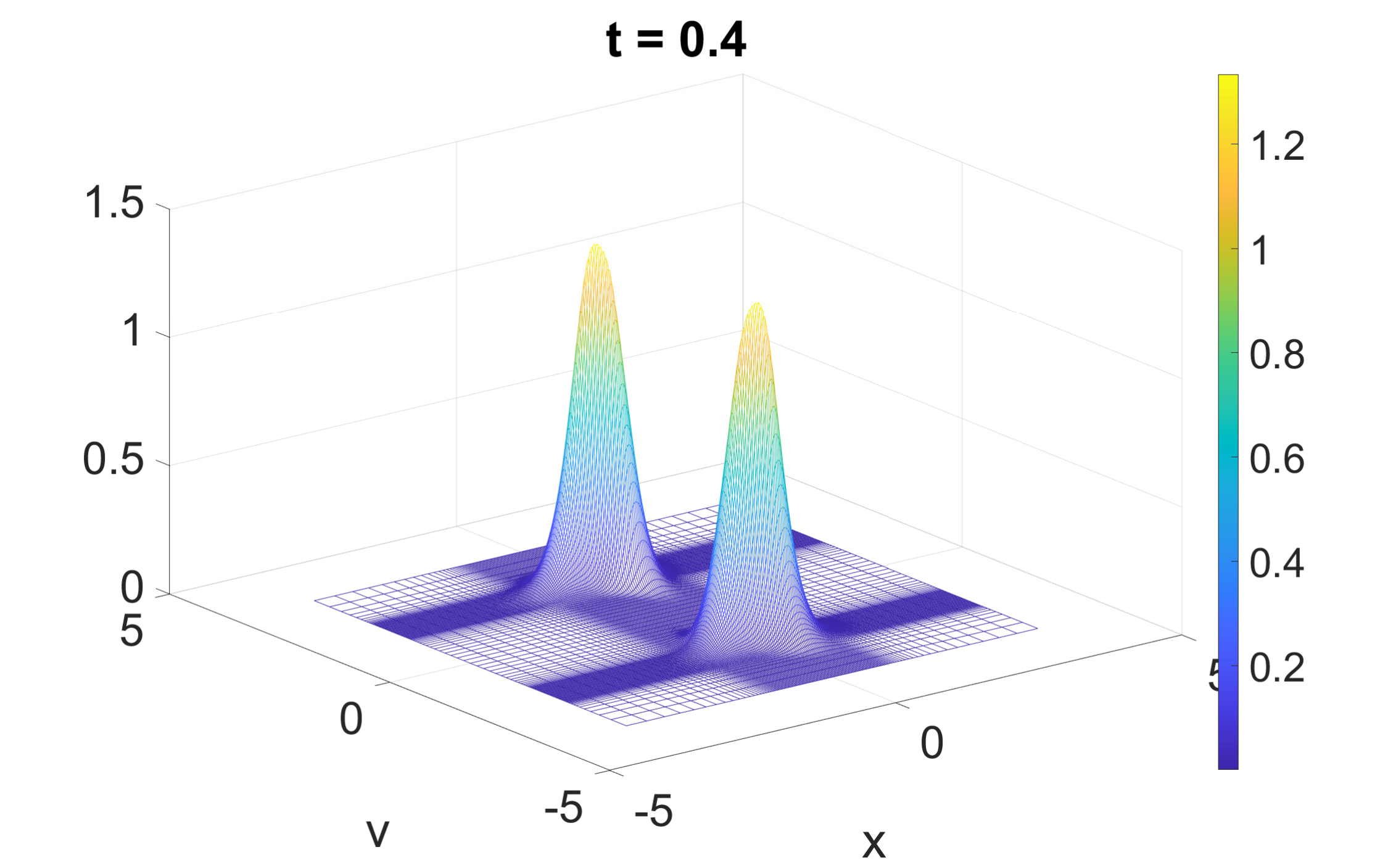}}
{\includegraphics[width=0.45\textwidth]{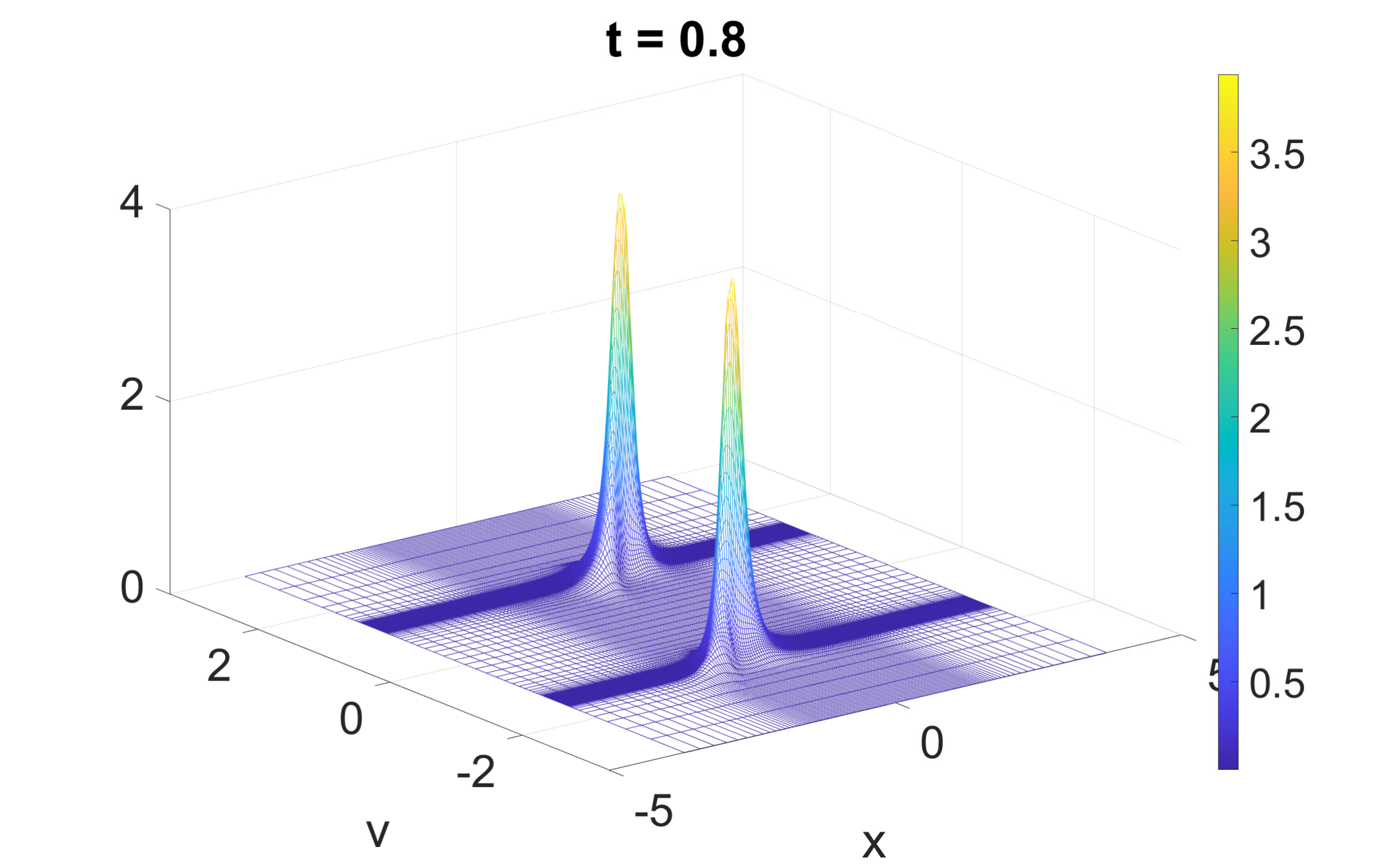}}
{\includegraphics[width=0.45\textwidth]{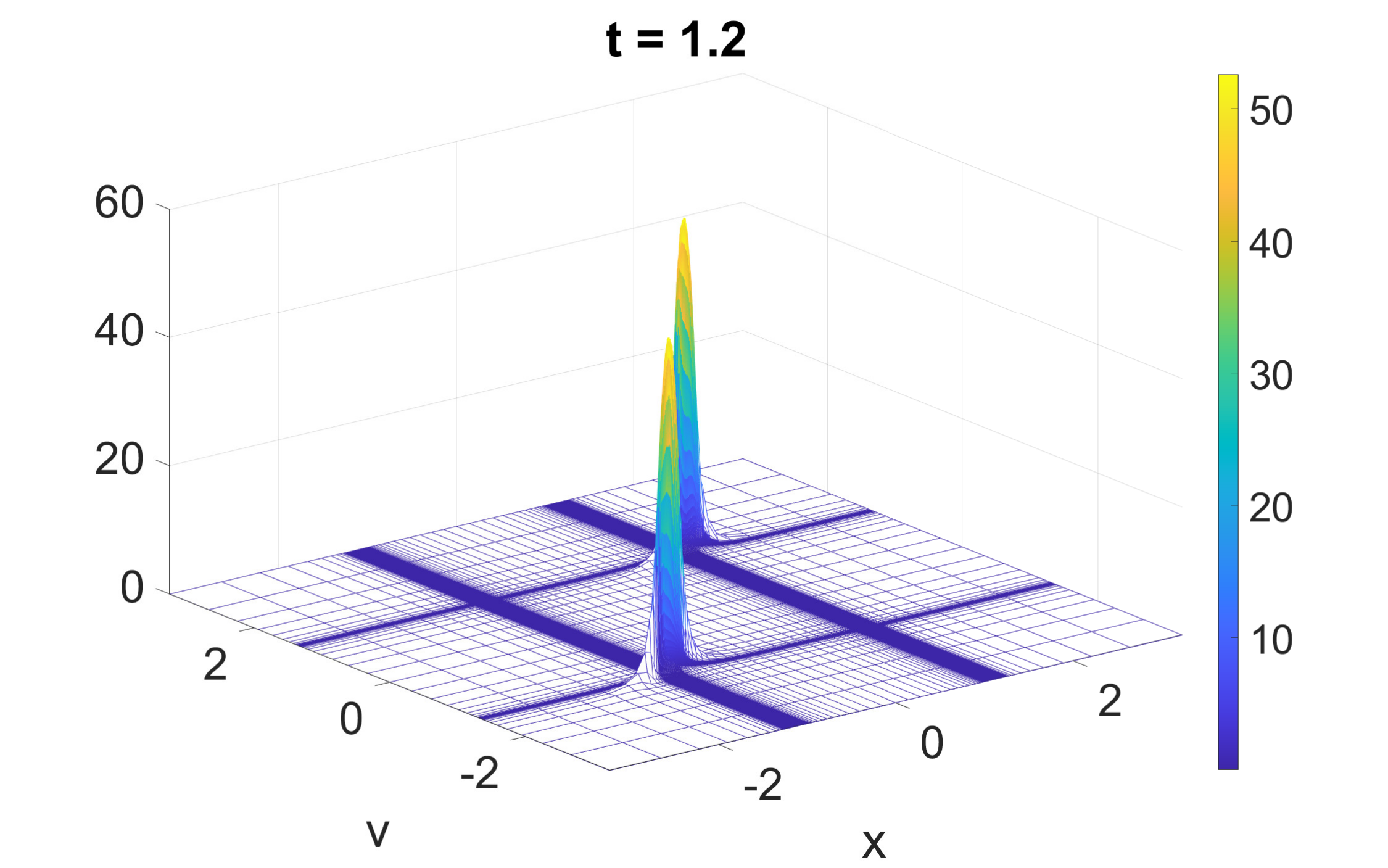}}
{\includegraphics[width=0.45\textwidth]{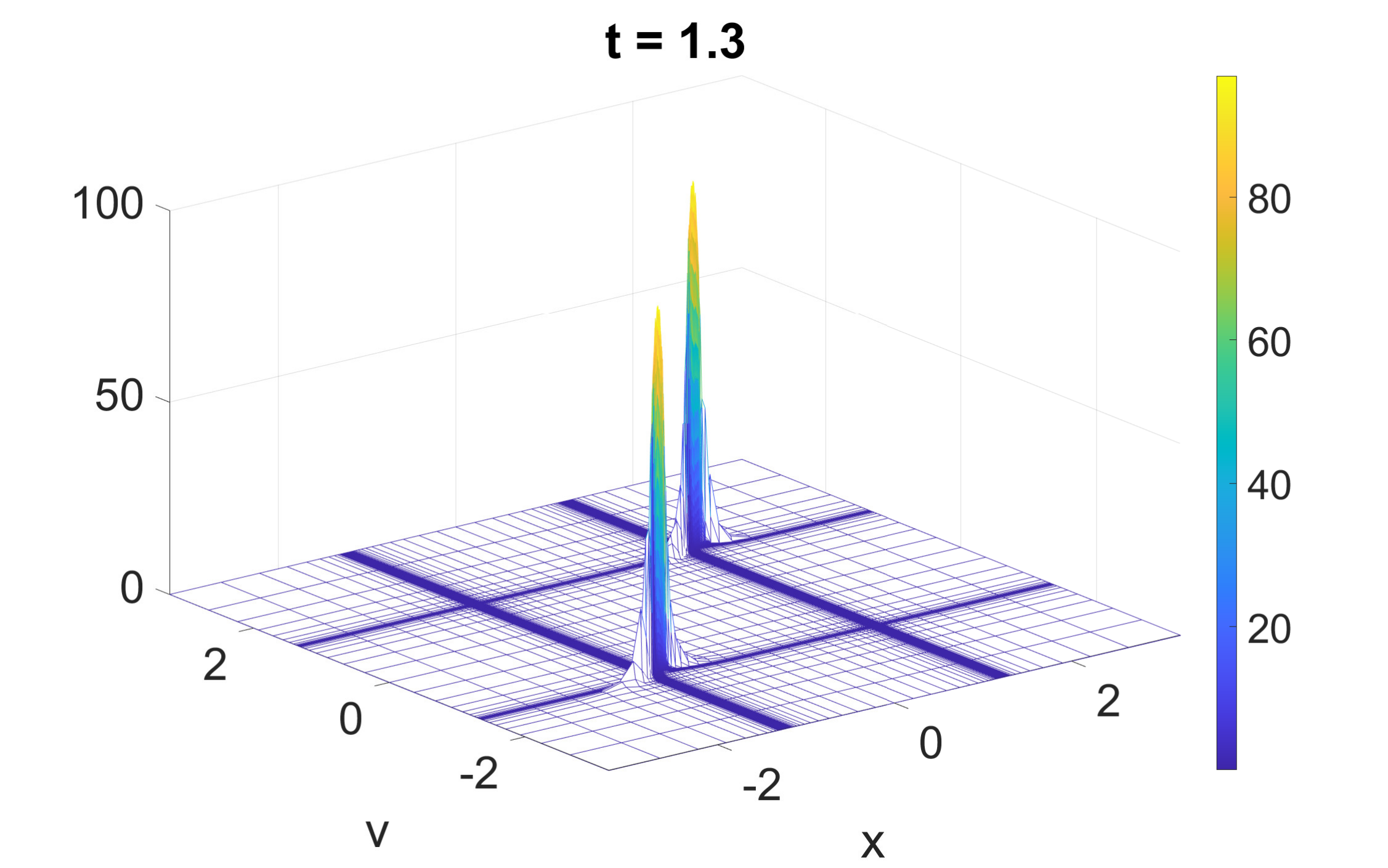}}
{\includegraphics[width=0.45\textwidth]{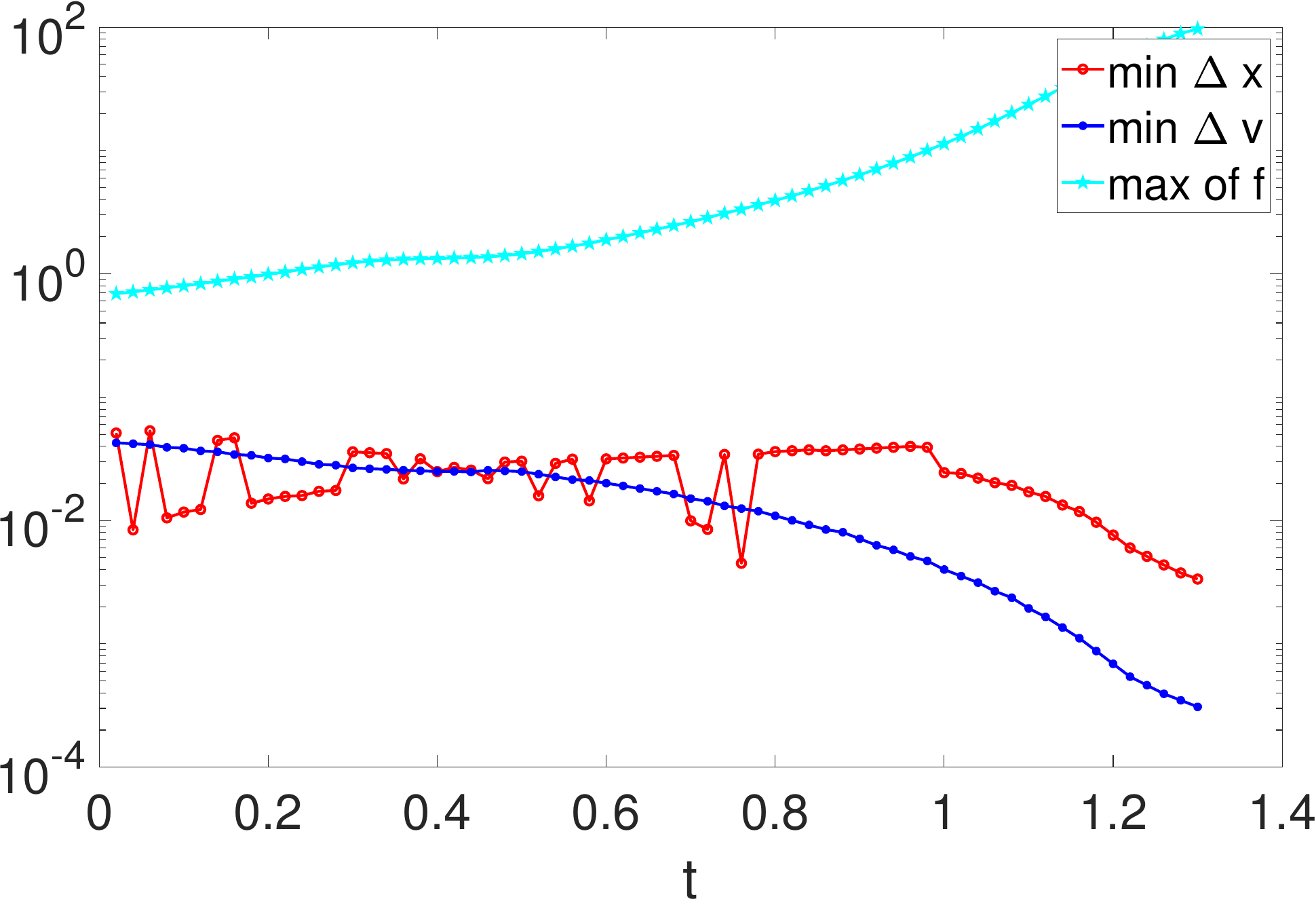}}
 \caption{Numerical solution of \eqref{eqn: gke} with kernel $W(v)=|v|^{3/2}$ for $N_x=121$, $N_v=121$, $a=6$, $b=3$, $c=d=2.5$ and $L_x=L_v=5$. Top left to bottom left are the numerical solution at $t = 0.1, 0.4, 0.8, 1.2 1.3$ and the bottom right is the record of minimum of $\Delta x$, $\Delta v$ and maximum of $f(t,x,v)$ at each time step.}
 \label{fig:gk1p5_C2p5}
\end{figure}

Typically, when two symmetric bumps initially locate at second quadrant and fourth quadrant respectively, the blow-up forms in $v$ direction before and after two bumps meets at $x=0$. More particularly, for initial condition $\eqref{eqn:f0}$ with $a=b=6$, $c=d=2.5$ and $\lambda = 10$, the blow-up arises on $v$ in finite time before two bumps meet at $x=0$, see result in Fig~\ref{fig:gk1p5_C5}. When $a=6$, $b=3$, $c=d=2.5$ and $\lambda=5$ are selected, the blow-up arises on $v$ in finite time after two bumps have passed the $x=0$, see Fig~\ref{fig:gk1p5_C2p5}. 

Finally, we conduct a numerical verification to ensure that the conservation or dissipation properties outlined in Proposition~\ref{prop:properties} persist over time up to the blow-up time. Conservations of mass and momentum, dissipation of moments and entropy are all well captured. The results are collected in the last two rows of Figure~\ref{fig:gk3_properties}, with the third row corresponding to the settings in Figure~\ref{fig:gk1p5_C5} and the bottom row corresponding to the settings in Figure~\ref{fig:gk1p5_C2p5}, respectively. 

\subsection{To Blow-up or not to blow-up: Numerical conclusions}

\begin{figure}[h!]
\centering
{\includegraphics[width=0.42\textwidth]{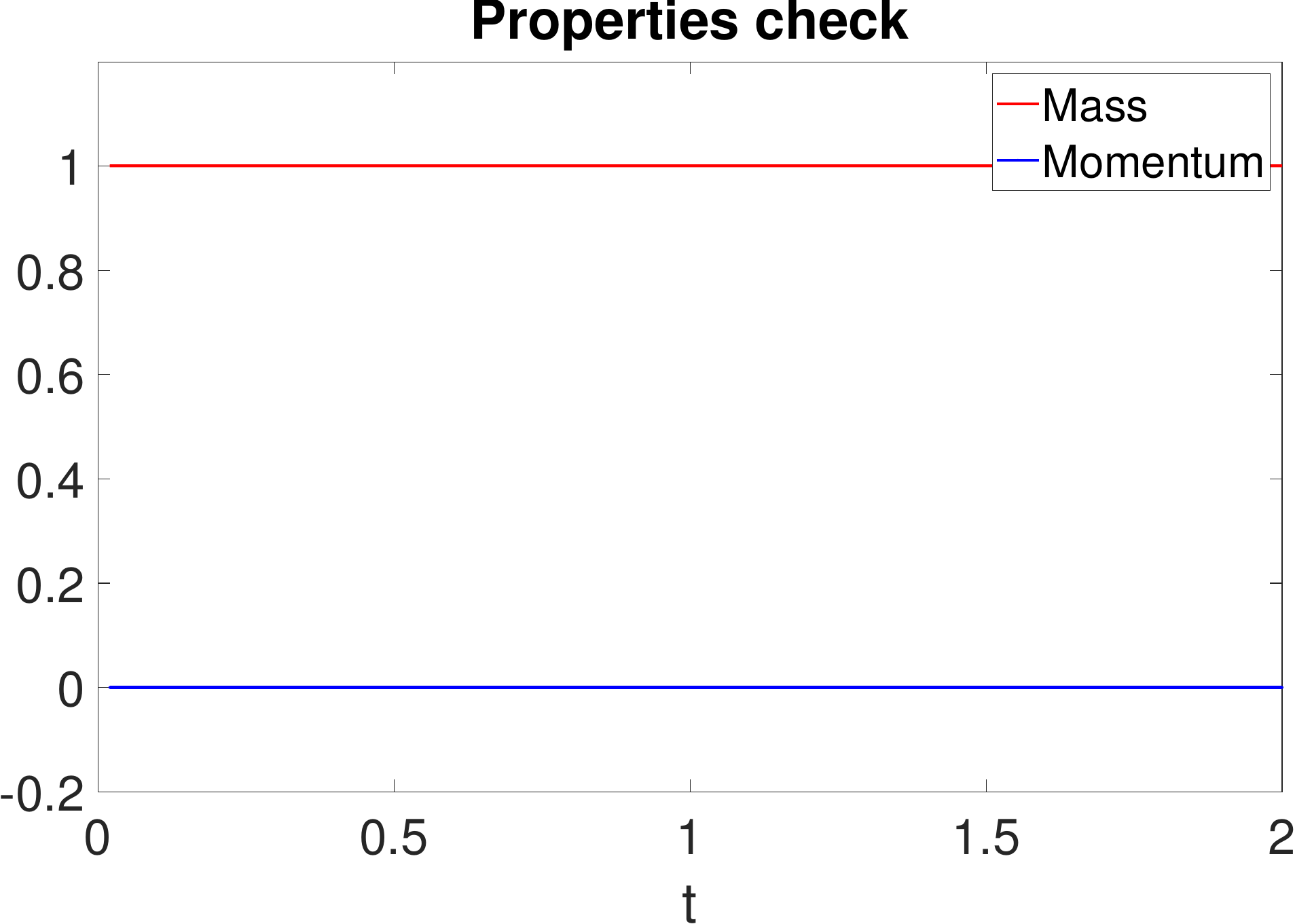}}
{\includegraphics[width=0.42\textwidth]{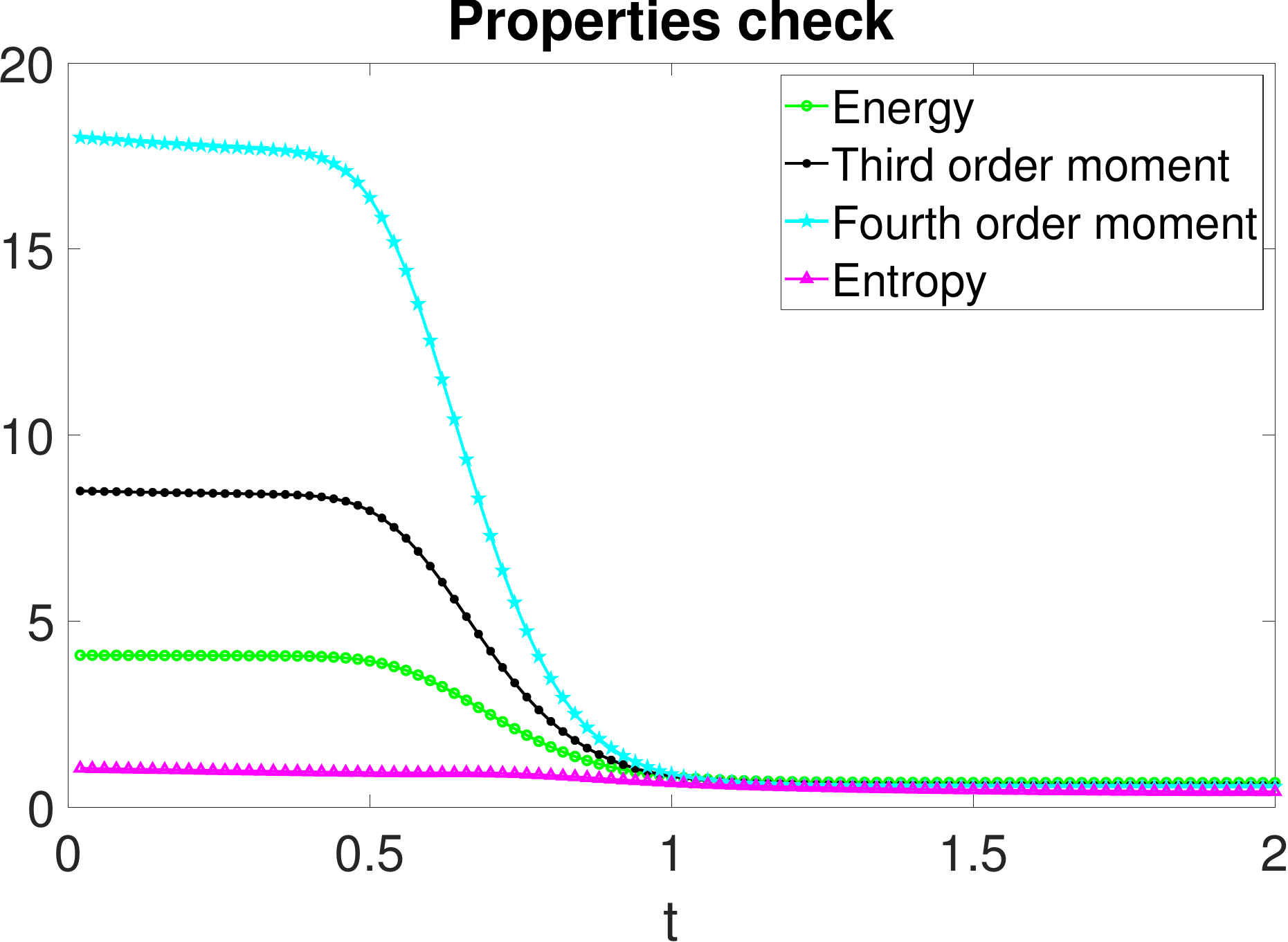}}
{\includegraphics[width=0.42\textwidth]{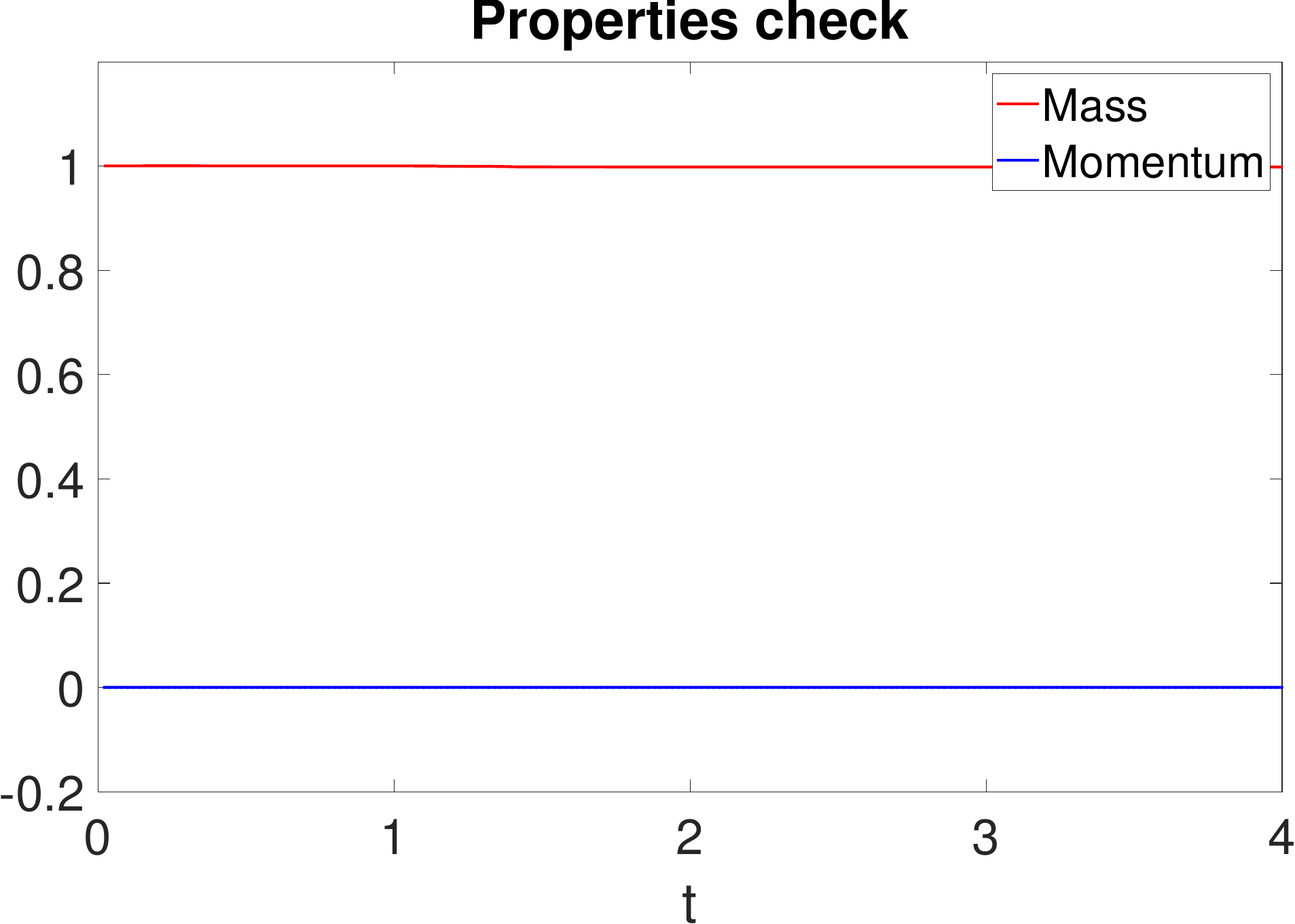}}
{\includegraphics[width=0.42\textwidth]{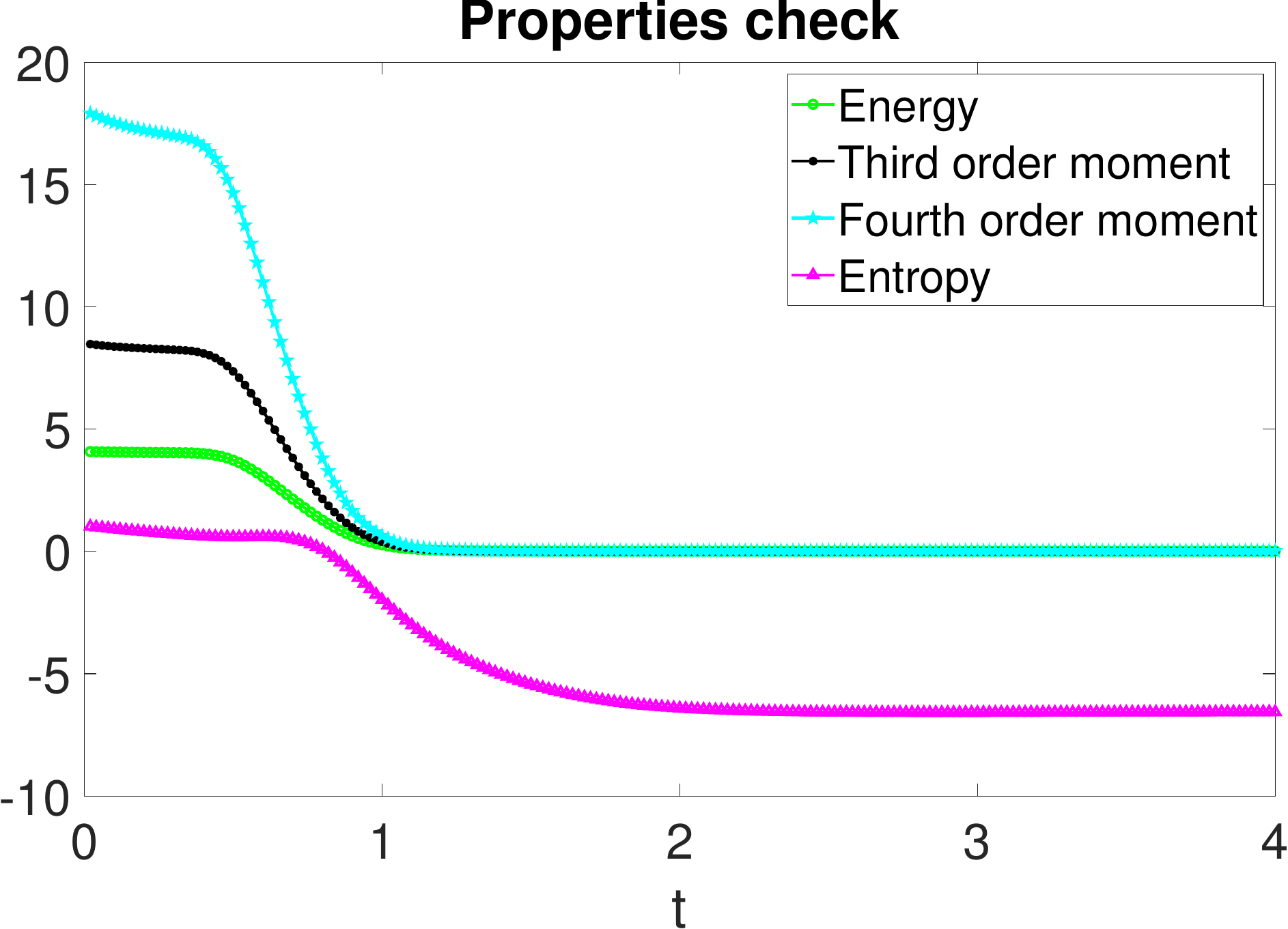}}
{\includegraphics[width=0.42\textwidth]{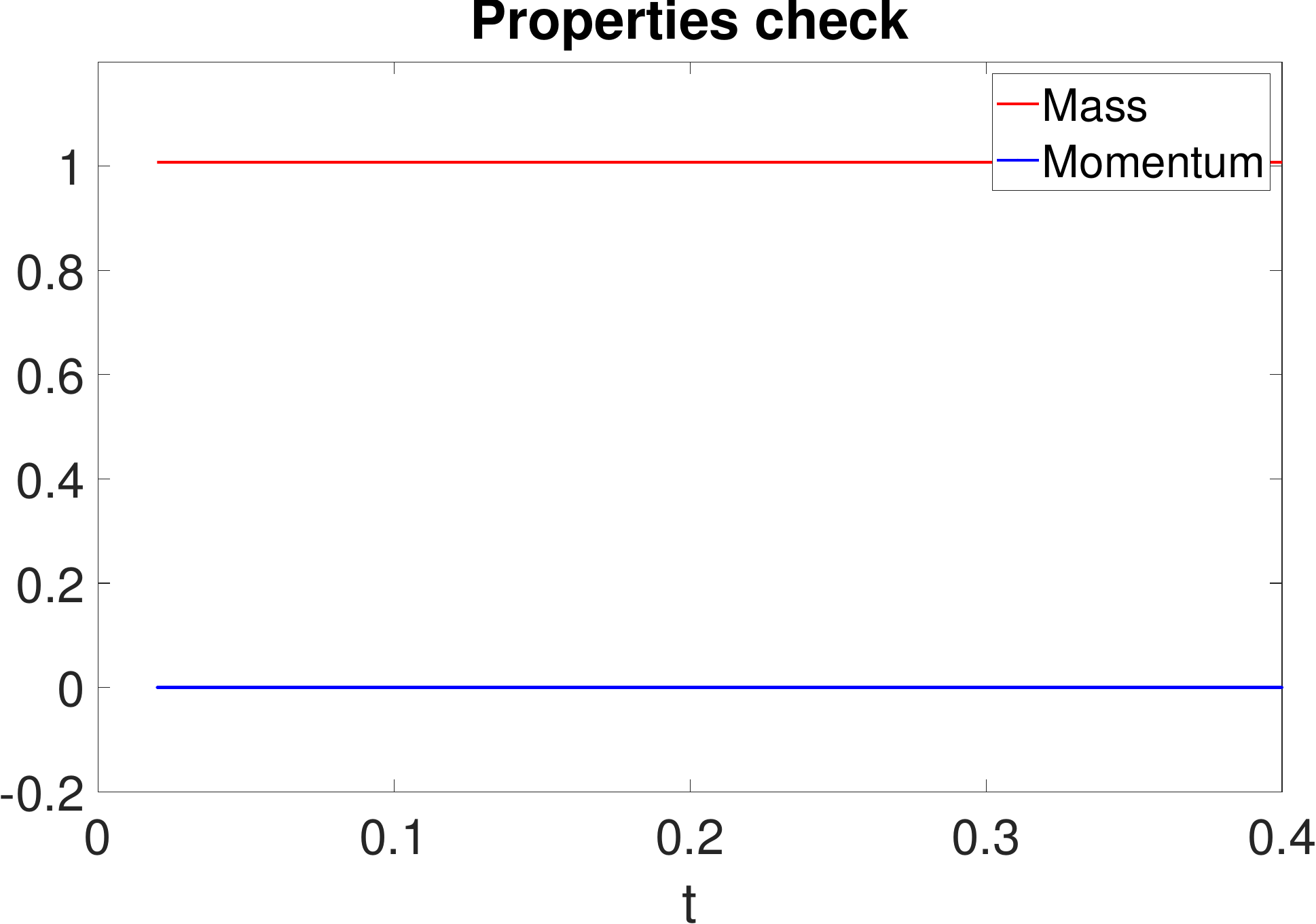}}
{\includegraphics[width=0.42\textwidth]{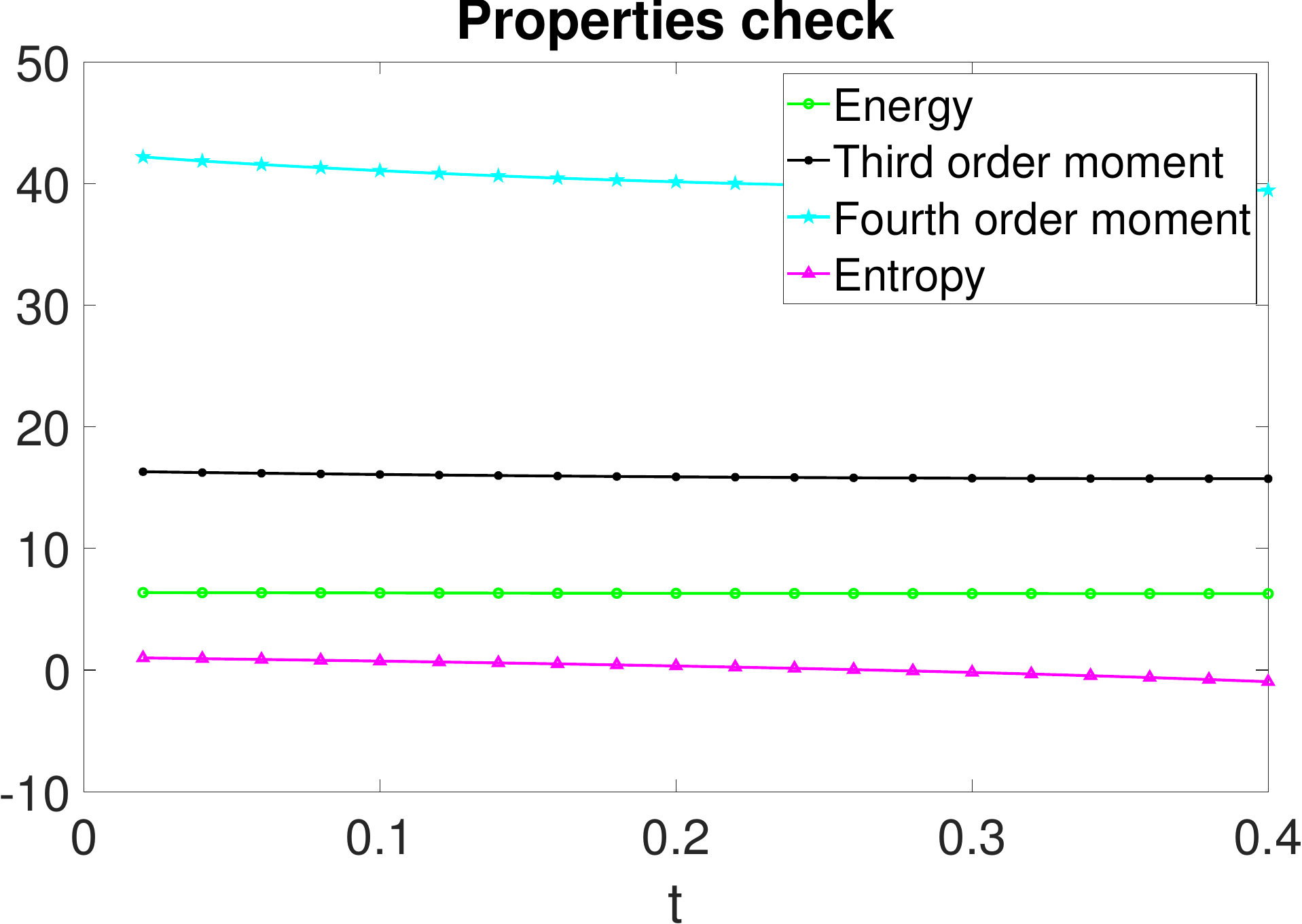}}
{\includegraphics[width=0.42\textwidth]{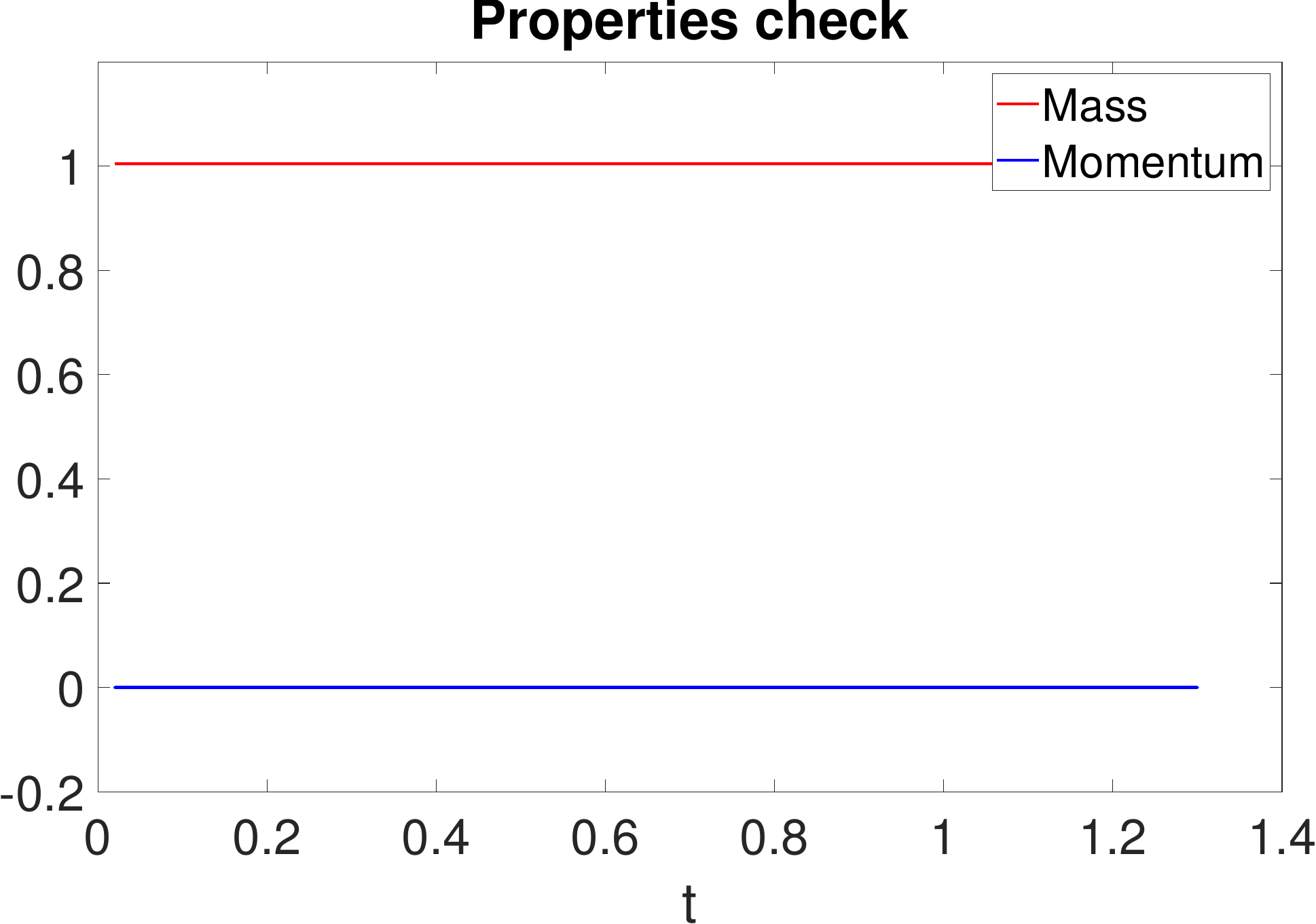}}
{\includegraphics[width=0.42\textwidth]{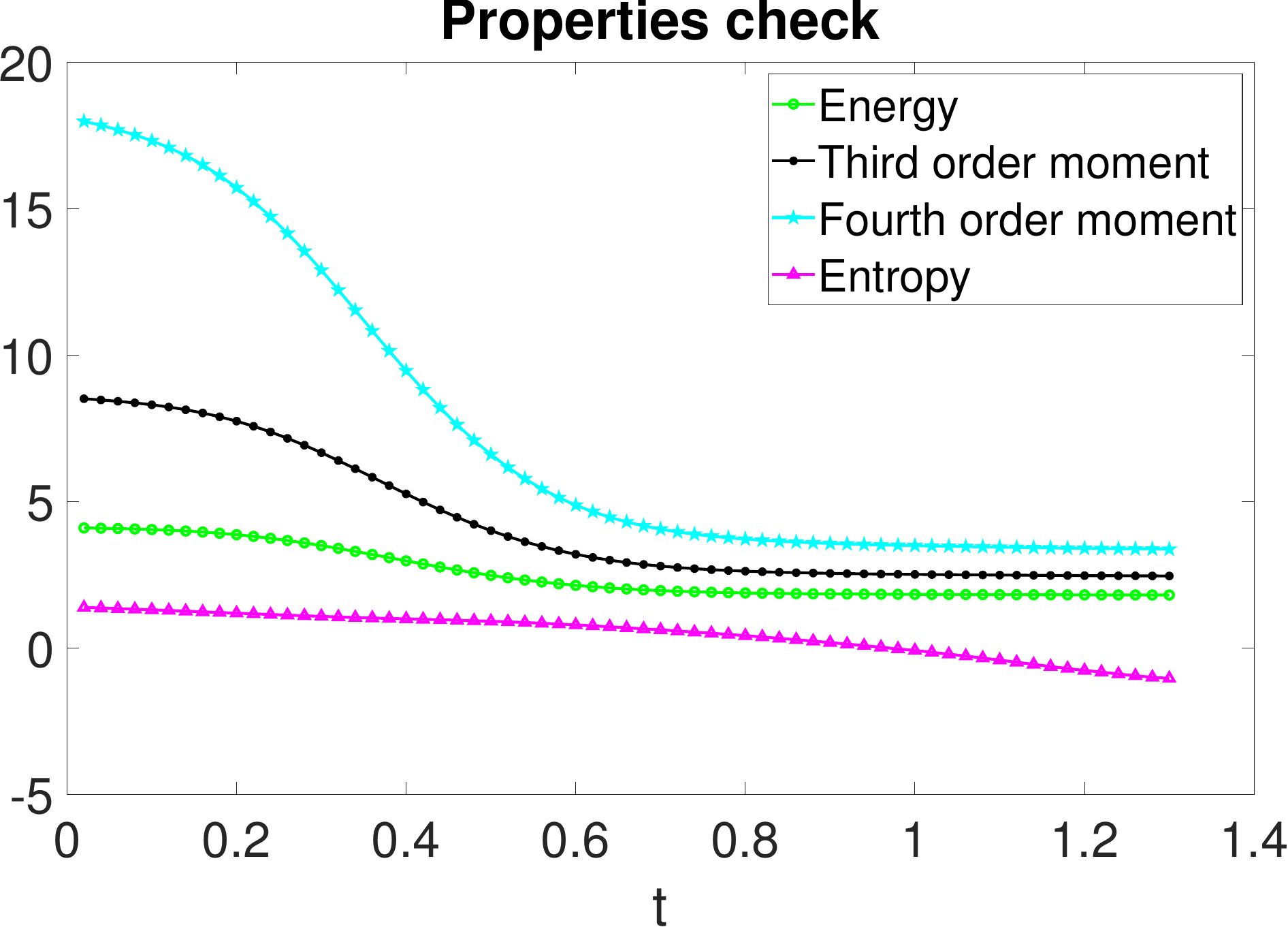}}
 \caption{Numerical conservation and dissipation properties. First and second row correspond to $\gamma=3$ and the numerical solutions in Figure~\ref{fig:gk3_ih_2_C2_2} and Figure~\ref{fig:gk3_ih_2_C5_2} respectively. Third and fourth row correspond to $\gamma=3/2$ and the numerical solutions in Figure~\ref{fig:gk1p5_C5} and  Figure~\ref{fig:gk1p5_C2p5} respectively.}
 \label{fig:gk3_properties}
\end{figure}

In the last two subsections, we have analysed, based on our time adaptive algorithm, a possible scenario of blow-up in the inhomogeneous problem \eqref{eqn: gke} with initial data \eqref{eqn:f0}. The numerical experiments firmly supports the evidence of no blow-up for $\gamma=3$ (see Figs. \ref{fig:gk3_ih_2_C2_2}--\ref{fig:gk3_ih_2_C5}) while numerical blow-up is detected for $\gamma=3/2$ (see Figs. \ref{fig:gk1p5_C5}--\ref{fig:gk1p5_C2p5}). Although theoretical results corroborating these numerical findings are lacking at present, we have conducted all possible theoretical checks on the numerical solutions to ensure they satisfy the correct conservations and dissipations, see Fig.~\ref{fig:gk3_properties}.

In this section, we have not shown numerical results for the case of $\gamma=2$. This case is quite particular since it is the borderline case between infinite-time blow-up and finite-time blow-up in the homogeneous case \cite{bertozzi2009blow}. Moreover, the homogeneous equation for $\gamma=2$ is equivalent to a linear equation when conservation of mass and momentum are taken into account for which it is easy to show that infinite time blow-up happens. Moreover, global-in-time measure solutions have been constructed for the inhomogeneous problem \eqref{eqn: gke} with $\gamma=2$ in \cite{AC16} allowing for possible concentrations of mass. Next section is devoted to show that such concentration of mass happens for infinite mass initial data. We will provide a self-similar infinite mass blowing-up solution. Moreover, we will give heuristic arguments and provide further numerical evidence that such a finite-time blow-up of solutions happens for finite-mass initial data for $\gamma=2$.

\section{Explicit solutions for $\gamma=2$}\label{sec: explicit}
In this section we give some explicit solutions to \eqref{eqn: gke} for $\gamma=2$ and its modifications, which provide intuition about whether \eqref{eqn: gke} has global wellposedness of $L^1$ solutions, see \cite{AC16} for global in time measure solutions.

\subsection{Infinite-mass self-similar solutions for $\gamma=2$}\label{sec: explicit1}

\begin{figure}[h!]
\centering
{\includegraphics[width=0.6\textwidth]{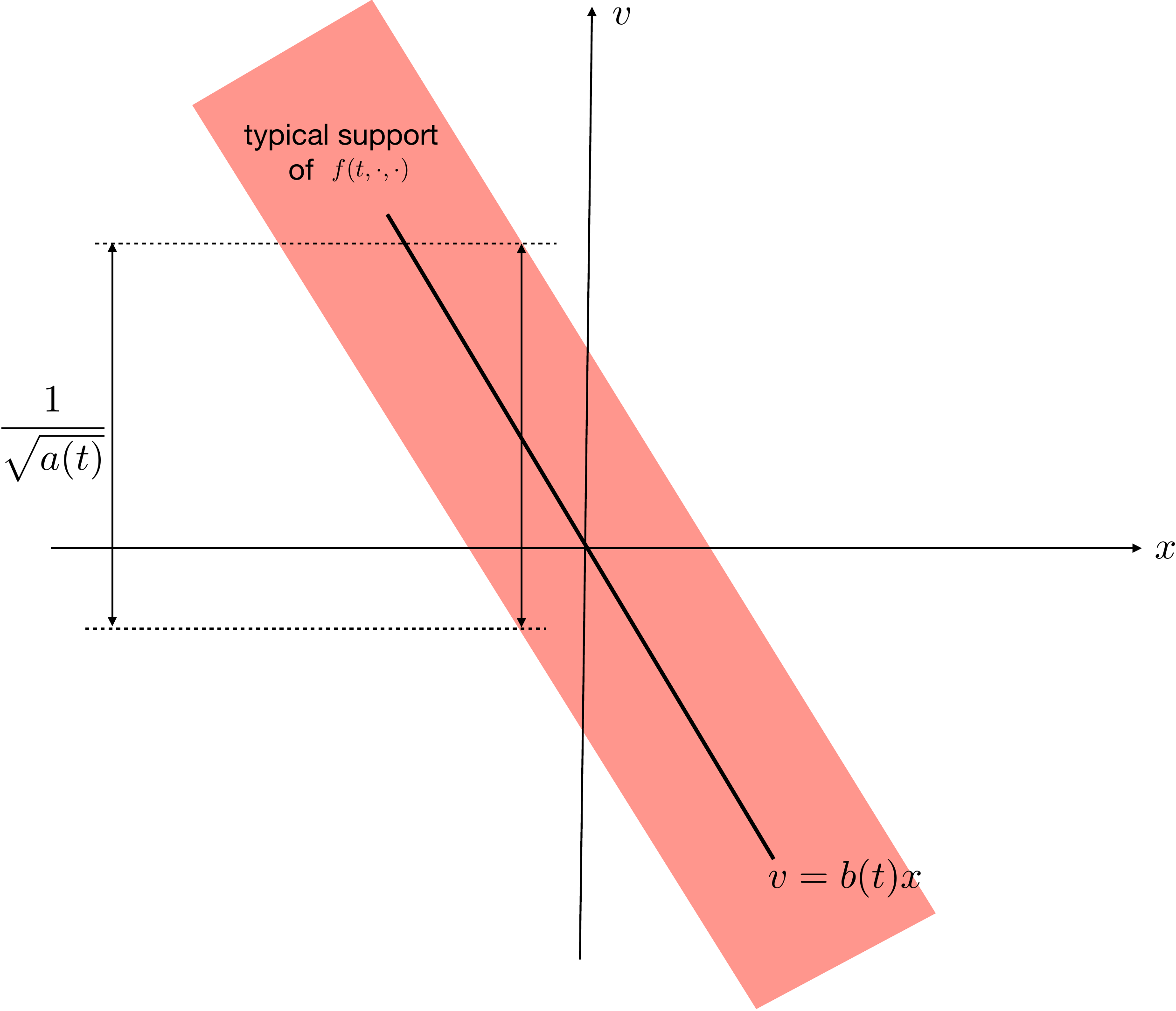}}
 \caption{Illustration of the infinite-mass solution to \eqref{eq} or \eqref{eqgamma} with initial data \eqref{ansatz}.}
 \label{fig:explicit1}
\end{figure}

We consider \eqref{eqn: gke} with $\gamma=2$. It can be rewritten as
\begin{equation}\label{eq}
\partial_t f + v\partial_x f = \frac{\lambda}{2} \rho \partial_v ((v-u)f),\quad \rho(t,x) = \int_\RR f(t,x,v)\rd{v},\quad \rho(t,x)u(t,x) = \int_\RR v f(t,x,v)\rd{v},
\end{equation}
using the macroscopic quantities $\rho$ and $u$. We will construct a class of infinite-mass self-similar solutions for \eqref{eq}.
Consider an ansatz of the form
\begin{equation}\label{ansatz}
f(t,x,v) = m(t)\phi\big(a(t)(v-b(t)x)^2\big)\,,
\end{equation}
where $m(t)>0,a(t)>0,b(t)$ are to be determined. Here $\phi$ is a smooth nonnegative function with $\int_\RR\phi(x^2)\rd{x}=1$. A typical choice is $\phi(x) = \frac{1}{\sqrt{2\pi}}e^{-x}$. It has density and bulk velocity
\begin{equation*}
\rho(t,x) = \rho(t) = \frac{m(t)}{\sqrt{a(t)}},\quad u(t,x) = b(t)x\,.
\end{equation*}
Since $\rho$ is constant in $x$, \eqref{ansatz} clearly has infinite total mass. See Figure \ref{fig:explicit1} for illustration.

Substituting \eqref{ansatz} into \eqref{eq} and dividing by $m$ on both sides, we get
\begin{equation*}
\frac{m'}{m}\phi - \Big( - a'(v-bx)^2 + 2a(v-bx)b'x + 2ab(v-bx)v\Big)\phi' =  \frac{\lambda m}{2\sqrt{a}}(\phi+2a(v-bx)^2\phi')\,,
\end{equation*}
i.e.,
\begin{equation*}
\Big(\underline{\frac{m'}{m}- \frac{\lambda m}{2\sqrt{a}}}\Big)\phi - (v-bx)\Big(\underline{ - a'(v-bx) + 2ab'x + 2abv + \lambda m \sqrt{a} (v-bx)}\Big)\phi'= 0\,.
\end{equation*}
Setting the underlined parts equal to zero (collecting $v$ terms and $x$ terms in the second part), we get
\begin{equation*}\left\{\begin{split}
& m' =  \frac{\lambda m^2}{2\sqrt{a}} \\
& a' = 2ab + \lambda m \sqrt{a} \\
& a'b + 2ab' - \lambda m \sqrt{a}b = 0\\
\end{split}\right.\end{equation*}
This shows that \eqref{ansatz} solves \eqref{eq} as long as $(m(t),a(t),b(t))$ solves the above ODE system.

\subsubsection{Solving the ODE system}
We substitute the second equation into the third and get
\begin{equation*}
\left\{\begin{split}
& m' =  \frac{\lambda m^2}{2\sqrt{a}} \\
& a' = 2ab + \lambda m \sqrt{a} \\
& b'=-b^2\\
\end{split}\right.
\end{equation*}
Using the density variable $\rho=m/\sqrt{a}$, we then have
\begin{equation}\label{ode1}\left\{\begin{split}
& \rho' = -b\rho \\
& m' = \frac{\lambda}{2} \rho m \\
& b'=-b^2\\
\end{split}\right.\end{equation}
Denote the initial condition of $(\rho,m,b)$ as $(\rho_0,m_0,b_0)$ with $\rho_0>0,\,m_0>0,\,b_0<0$. First notice that the $b$ equation has explicit solution
\begin{equation*}
b(t) = -\frac{1}{T-t},\quad T = \frac{1}{-b_0}>0\,,
\end{equation*}
which blows up to $-\infty$ at time $T$. Substituting into the $\rho$ equation, we get
\begin{equation}\label{rhot}
\rho(t) = \rho_0 \exp\Big(-\int_0^t b(s)\rd{s}\Big) = \frac{\rho_0 T}{T-t}\,.
\end{equation}
Then, substituting into the $m$ equation, we get
\begin{equation}\label{mt}
m(t) = m_0 \exp\Big(\frac{\lambda}{2}\int_0^t \rho(s)\rd{s}\Big)= m_0 \exp\Big(\frac{\lambda}{2}\rho_0 T\ln\frac{T}{T-t}\Big) = m_0 \Big(\frac{T}{T-t}\Big)^{\frac{\lambda}{2}\rho_0 T}\,.
\end{equation}

\subsubsection{Analysis of the blow-up behavior}\label{sec:412}
Although the solution to \eqref{ode1} always blows up at $T$, we need to analyze whether the corresponding blow-up behavior of the infinite-mass solution \eqref{ansatz} can be approximated by a finite-mass solution. This will be done by calculating the typical width of the velocity support for \eqref{ansatz}, which is
\begin{equation*}
\frac{1}{\sqrt{a(t)}} = \frac{\rho(t)}{m(t)} = \frac{\rho_0}{m_0} \Big(\frac{T}{T-t}\Big)^{1-\frac{\lambda}{2}\rho_0 T}\,.
\end{equation*}
Denote $\epsilon(t) = T-t$ as the time before blow-up. Then we get the relations
\begin{equation}\label{eqn:brhom}
b\sim -\epsilon^{-1},\quad \rho\sim \epsilon^{-1},\quad m\sim \epsilon^{-\frac{\lambda}{2}\rho_0 T},\quad \frac{1}{\sqrt{a}}\sim \epsilon^{\frac{\lambda}{2}\rho_0 T-1}\,,
\end{equation}
as $t\rightarrow T^-$ (where the notation $\sim$ means that the two quantities are related by a constant multiple). Take a box $[-\epsilon,\epsilon]\times [-1,1]$ in the $(x,v)$-plane (up to constant multiples of its size). Then the mean velocities $(x,b(t)x),\,x\in [-\epsilon,\epsilon]$ lies in the box, and the total mass $\iint_{x\in [-\epsilon,\epsilon]}f\rd{v}\rd{x}\sim 1$ since $\rho\sim \epsilon^{-1}$. We propose the threshold condition
\begin{equation}\label{thres}
\lambda\rho_0 T > 2\,.
\end{equation}
\begin{itemize}
    \item If \eqref{thres} holds, then the width of the velocity support $\frac{1}{\sqrt{a}}\rightarrow 0$ as $t\rightarrow T^-$. This means that most mass of $f$ with $x\in [-\epsilon,\epsilon]$ lies inside the box $[-\epsilon,\epsilon]\times [-1,1]$ with fixed height. It is likely that such a infinite-mass solution can be approximated by a finite-mass solution, leading to finite-time blow-up for a finite-mass solution.
    \item If \eqref{thres} fails with $\lambda\rho_0 T < 2$, then the velocity support $\frac{1}{\sqrt{a}}\rightarrow \infty$ as $t\rightarrow T^-$. Such an infinite-mass solution cannot be well-approximated by a finite-mass solution for $t$ close to $T$ because $\sup_{(x,v)\in \text{supp } f(t,\cdot,\cdot)}|v|$ for the latter cannot grow in time. Therefore we do not expect finite-time blow-up for a finite-mass solution if its initial data is close to \eqref{ansatz} in this case.
\end{itemize}

\subsection{A conditional blow-up result}

In this subsection we give a conditional blow-up result for finite-mass solutions to \eqref{eq} based on an analysis of the characteristic flow.

The characteristic flow $(X(t;x,v),V(t;x,v))$ of a $C^1$ solution to \eqref{eq} satisfies the ODE system
\begin{equation*}
    \dot{X} = V,\quad \dot{V} = \frac{\lambda}{2}\rho(t,X) (u(t,X,V)-V); \quad X(0;x,v)=x,\,V(0;x,v)=v.
\end{equation*}
For the infinite-mass solution \eqref{ansatz}, we have
\begin{equation*}
    \rho(t,X)=\rho(t),\quad u(t,X,V)=b(t)X,
\end{equation*}
and thus
\begin{equation*}
    \dot{X} = V,\quad \dot{V} = \frac{\lambda}{2}\rho(t) (b(t)X-V).
\end{equation*}
Along the characteristic originating from $(x,v)$, the quantity $\alpha=\alpha(x,v):=\sqrt{a}(V-b X)=\frac{m_0}{\rho_0}(v+\frac{x}{T})$ is conserved in time (by calculating its derivative explicitly). Therefore we get 
\begin{equation*}
    \dot{X}=bX+\frac{\alpha}{\sqrt{a}} = -\frac{1}{T-t}X +  \frac{\alpha\rho_0}{m_0} \Big(\frac{T}{T-t}\Big)^{1-\frac{\lambda}{2}\rho_0 T}
\end{equation*}
and thus its explicit solution
\begin{equation*}
    \frac{X(t;x,v)}{T-t} = \frac{x}{T} + \frac{\alpha \rho_0}{m_0}\frac{1}{1-\frac{\lambda}{2}\rho_0 T}\Big(\Big(\frac{T}{T-t}\Big)^{1-\frac{\lambda}{2}\rho_0 T}-1\Big)\,.
\end{equation*}

\begin{figure}[h!]
\centering
{\includegraphics[width=0.45\textwidth]{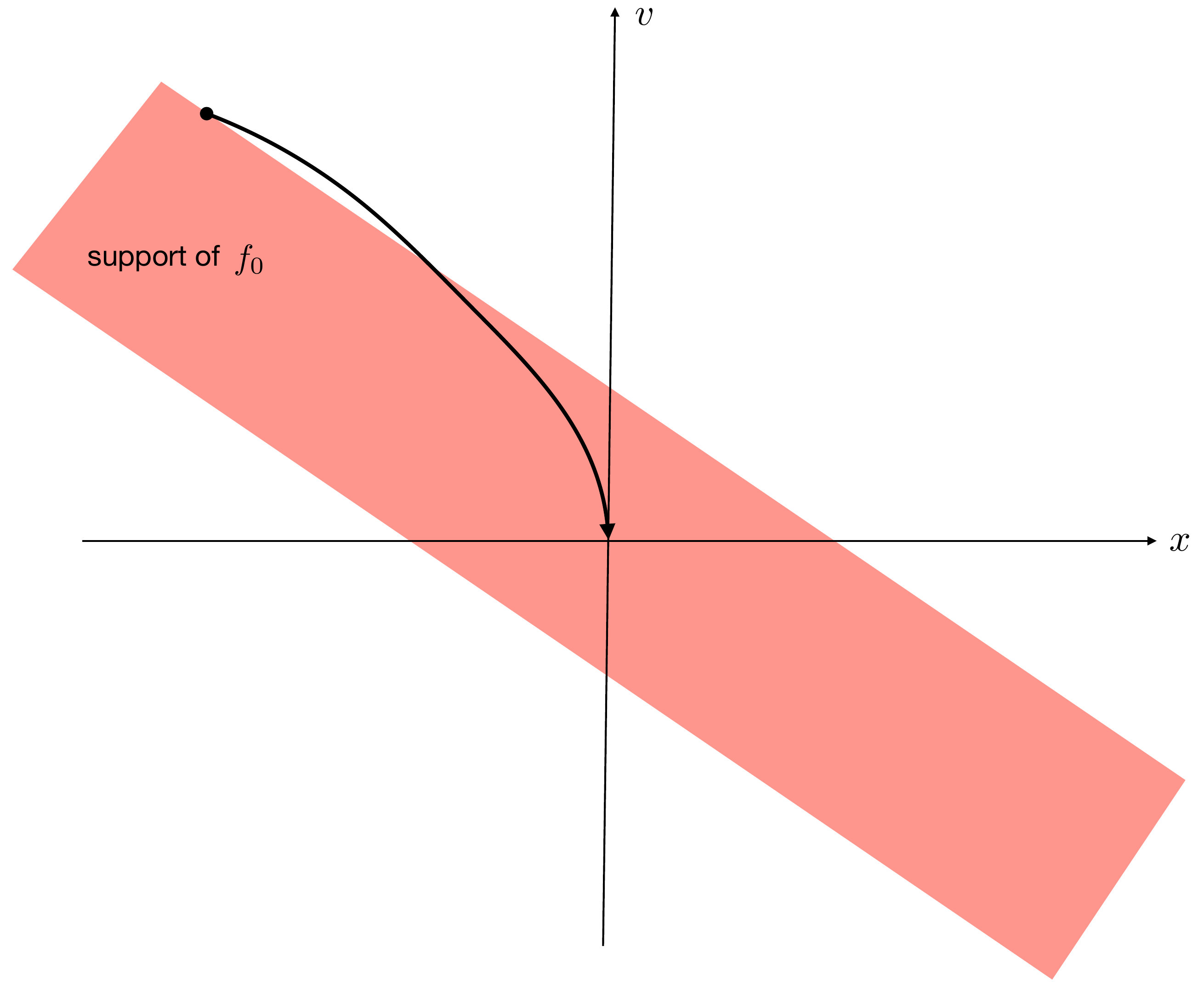}}
{\includegraphics[width=0.45\textwidth]{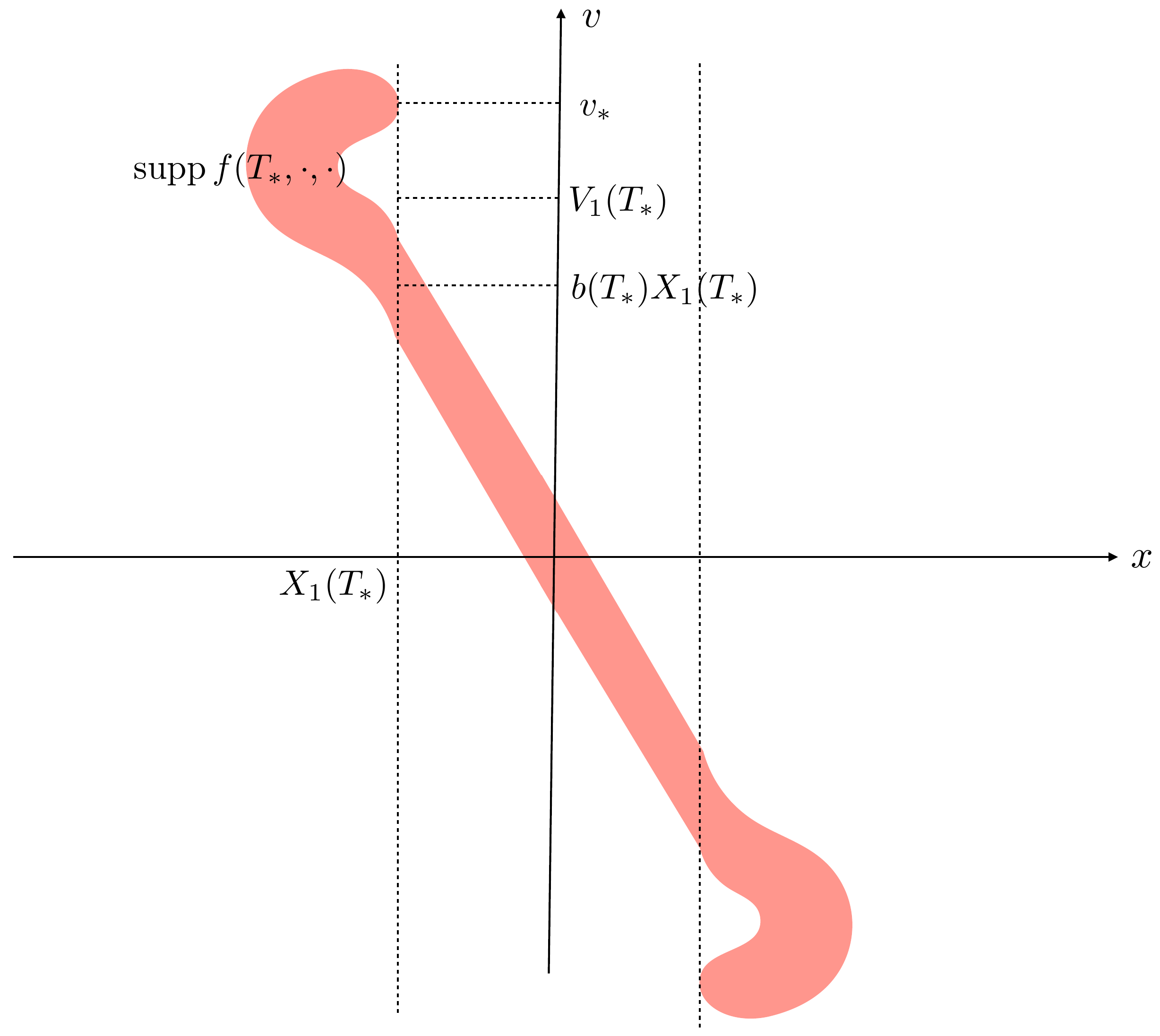}}
 \caption{Left: a particle trajectory for \eqref{eq} which satisfies $X(t;x,v)<0$ for any $t\in (0,T)$. Right: a possible scenario in which \eqref{condft} breaks down.}
 \label{fig:explicit2}
\end{figure}

Now we assume the super-critical case \eqref{thres}, and consider a characteristic starting from $(x,v)$ with $x<0$ and $\alpha(x,v)>0$ (i.e., $v>-\frac{x}{T}$). Then we have $1-\frac{\lambda}{2}\rho_0 T<0$ and thus 
\begin{equation*}
    \frac{X(t;x,v)}{T-t} \le \frac{x}{T} -(v+\frac{x}{T})\frac{1}{1-\frac{\lambda}{2}\rho_0 T}\,.
\end{equation*}
Therefore, for the nonempty set of velocities such that
\begin{equation*}
    -\frac{x}{T} < v \le -\frac{x}{T}\cdot\frac{\lambda}{2}\rho_0 T
\end{equation*}
then we have $X(t;x,v)<0$ for any $t\in (0,T)$, see Figure \ref{fig:explicit2} Left for illustration. If the profile $\phi$ in the ansatz is compactly supported, then by taking $x$ sufficiently negative it is always possible to choose $(x,v)\notin \textnormal{supp}\, f(0,\cdot,\cdot)$ satisfying 
\begin{equation}\label{xvcond}
    -\frac{x}{T}<v<-\frac{x}{T}\cdot\frac{\lambda}{2}\rho_0 T\,.
\end{equation}
Such a characteristic enables us to give a conditional finite-time blow-up result for finite-mass solutions.

\begin{theorem}\label{thm_condbu}
    Assume \eqref{thres} and denote $f(t,x,v)$ as the infinite-mass solution \eqref{ansatz} with $\phi$ compactly supported. Let $(x_1,v_1)\notin \textnormal{supp}\, f(0,\cdot,\cdot)$ satisfy $x_1<0$ and \eqref{xvcond}, and denote $(X_1(t),V_1(t))$ as its corresponding characteristic. 
    
    Let $\tilde{f}(t,x,v)$ be a $C^1$ solution to \eqref{eq}, whose initial data satisfies
    \begin{equation*}
        \tilde{f}(0,-x,-v) = \tilde{f}(0,x,v)
    \end{equation*}
    and
    \begin{equation*}
        \tilde{f}(0,x,v)=f(0,x,v),\quad \forall |x|\le |x_1|,\,\forall v.
    \end{equation*}
    Denote $\tilde{T}$ as the maximal existence time for the solution $\tilde{f}$.

    If there holds
    \begin{equation}\label{condft}
        (X_1(t),v)\notin \textnormal{supp}\,\tilde{f}(t,\cdot,\cdot),\quad \forall t<\max\{T,\tilde{T}\},\, v\ge V_1(t),
    \end{equation}
    then we have $\tilde{T} \ge T$, i.e., the solution $\tilde{f}$ has a finite-time blow-up.
\end{theorem}
The extra condition \eqref{condft} can be checked numerically as the following (see Figure \ref{fig:explicit2} Right for illustration). Suppose \eqref{condft} breaks down at time $T_*$, then the proof below shows that $\tilde{f}$ and $f$ agrees for any $t<T_*$ and $|x|<|X_1(t)|$. As a consequence, $(X_1(T_*),V_1(T_*))\notin \textnormal{supp}\,\tilde{f}(t,\cdot,\cdot)$ and $(X_1(T_*),b(T_*)X_1(T_*))\in \textnormal{supp}\,\tilde{f}(t,\cdot,\cdot)$. Since \eqref{condft} breaks down at time $T_*$, there exists some $v_*>V_1(t)$ such that $(X_1(T_*),v_*)\in \textnormal{supp}\,\tilde{f}(T_*,\cdot,\cdot)$. Since $V_1(T_*)$ is between $b(T_*)X_1(T_*)$ and $v_*$, we see that $\textnormal{supp}\,\tilde{f}(T_*,X_1(T_*),\cdot)$ has at least two connected components. Therefore, if one can check numerically that the level sets of $\tilde{f}(t,\cdot,\cdot)$ intersecting with any vertical line always have one connected component, this would indicate that \eqref{condft} holds.

\begin{proof}
Let $(X,V)$ denote the characteristic flow of $f$, and $(\tilde{X},\tilde{V})$ denote the characteristic flow of $\tilde{f}$. We claim that for any $(x,v)\in \textnormal{supp}\, \tilde{f}(0,\cdot,\cdot)$ with $x<x_1$, we have
\begin{equation*}
    \tilde{X}(t;x,v) < X_1(t),\quad \forall t<\max\{T,\tilde{T}\}\,.
\end{equation*}
Similar holds if $x>-x_1$ since $\tilde{f}$ is odd in $(x,v)$. 

The paragraphs before Theorem \ref{thm_condbu} show that $X_1(t)<0$ for any $t<T$. Therefore the claim would imply that 
\begin{equation*}
    \tilde{f}(t,x,v)=f(t,x,v),\quad \forall |x|\le |X_1(t)|,\,\forall v
\end{equation*}
by tracing back the characteristics in $\textnormal{supp}\,\tilde{f}(t,\cdot,\cdot)$. Therefore the finite-time blow-up of $f$ at time $T$ implies that of $\tilde{f}$ at $T$ provided that $\tilde{T}\ge T$, and the conclusion is obtained.

To prove the claim, we assume the contrary that $t_*>0$ is the earliest time such that $\tilde{X}(t_*;x,v) = X_1(t_*)$ is reached. Therefore we necessarily have $\tilde{V}(t_*;x,v)\ge V_1(t_*)$. Then, since $(x,v)\in \textnormal{supp}\, \tilde{f}(0,\cdot,\cdot)$, we may propagate to time $t_*$ and get $(\tilde{X}(t_*;x,v),\tilde{V}(t_*;x,v))\in \textnormal{supp}\, \tilde{f}(t_*,\cdot,\cdot)$. Combining with $\tilde{X}(t_*;x,v) = X_1(t_*)$ and $\tilde{V}(t_*;x,v)\ge V_1(t_*)$, we get a contradiction with \eqref{condft}.

\end{proof}

\subsection{Approximation of $\gamma>2$ solutions by $\gamma=2$ solutions}

We consider \eqref{eqn: gke} with $\gamma>2$, and provide a heuristic argument by approximating it with $\gamma=2$ solutions. In the spirit of the previous subsection, we assume that the typical width of the support of $f$ in the $v$-direction is $R_v>0$ which depends on $t$ but does not depend on $x$. Then we approximate \eqref{eqn: gke} by replacing  $|v-w|^{\gamma-2}$ with $R_v^{\gamma-2}$, and get
\begin{equation*}
\partial_t f + v\partial_x f = \frac{\lambda}{2} R_v^{\gamma-2} \rho\partial_v ((v-u)f)\,.
\end{equation*}
For the infinite-mass solution \eqref{ansatz}, we may take $R_v=\frac{1}{\sqrt{a}}$. This gives
\begin{equation}\label{eqgamma}
\partial_t f + v\partial_x f = \frac{\lambda}{2} \Big(\frac{1}{\sqrt{a}}\Big)^{\beta} \rho \partial_v ((v-u)f),\quad \beta:=\gamma-2>0
\end{equation}
which is analogous to \eqref{eq}, with an extra coefficient on the RHS. Notice that \eqref{ansatz} is still a solution to \eqref{eqgamma} if $(\rho,m,b)$ satisfies \eqref{ode1} with $\frac{\lambda}{2}$ replaced by $\frac{\lambda}{2} (\frac{1}{\sqrt{a}})^\beta$:
\begin{equation*} 
\left\{\begin{split}
& \rho' = -b\rho \\
& m' = \frac{\lambda}{2} \rho^{1+\beta} m^{1-\beta} \\
& b'=-b^2\\
\end{split}\right.
\end{equation*}
The $b,\rho$ solution is the same as before, given by \eqref{rhot} and \eqref{mt}. For the $m$ solution, we integrating the $m$ equation and get
\begin{equation*}
(m^\beta)' \sim \frac{1}{(T-t)^{1+\beta}},\quad m(t)\sim \frac{1}{T-t}\,.
\end{equation*}
Therefore we get the width of velocity support
\begin{equation*}
\frac{1}{\sqrt{a}}\sim 1
\end{equation*}
as $t\rightarrow T^-$.

\begin{figure}[h!]
\centering
{\includegraphics[width=0.5\textwidth]{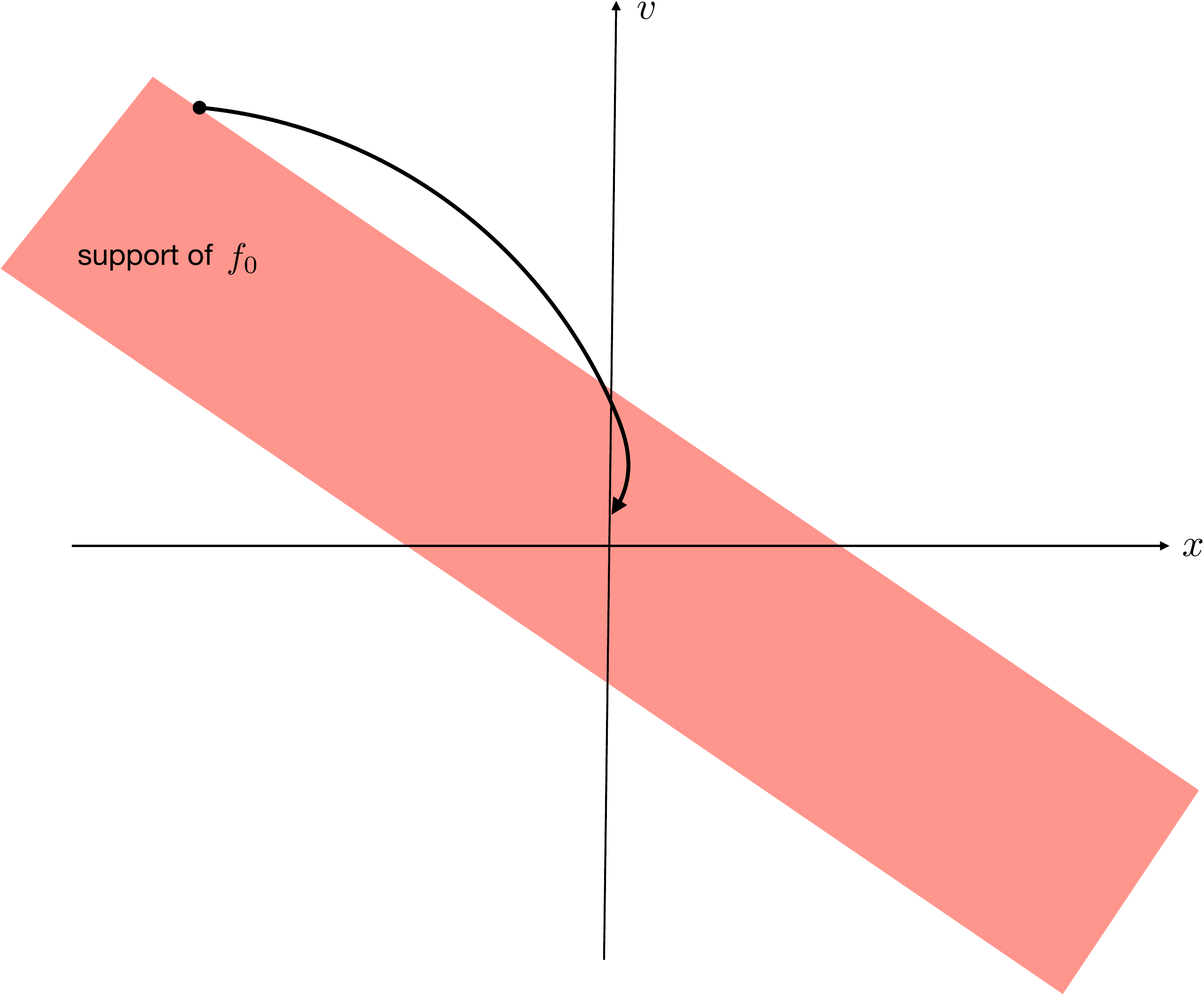}}
 \caption{A particle trajectory for \eqref{eqgamma} which passes $x=0$.}
 \label{fig:explicit3}
\end{figure}

To discuss the physical meaning of the last result, (following notations of the previous subsection) we notice that the total mass in the box $[-\epsilon,\epsilon]\times [-1,1]$ is bounded as $t\rightarrow T^-$. However, for any characteristic curve $(X(t),V(t))$ starting from $X_0<0,\,V_0=b_0X_0 + \frac{\alpha}{\sqrt{a_0}},\,\alpha\in\mathbb{R}$, we have
\begin{equation*}
    \dot{X}=V,\quad \dot{V}=-\frac{\lambda}{2} (\frac{1}{\sqrt{a}})^{\beta}\rho (V-b X)\,,
\end{equation*}
and therefore $\sqrt{a}(V-b X) = \alpha$ is conserved (by calculating its time derivative explicitly). This gives
\begin{equation*}
    \dot{X} = b X + \alpha \frac{\rho}{m} = -\frac{1}{T-t} X + \alpha C \,,
\end{equation*}
for some constant $C>0$ (since both $\rho,m$ behave like $1/(T-t)$). Then we integrate to get
\begin{equation*}
    X(t) = \frac{T-t}{T}\Big(X_0 + \alpha C T \ln\frac{T}{T-t}  \Big)\,.
\end{equation*}
Recall that $X_0<0$. If $\alpha>0$, then any characteristic $X(t)$ reaches 0 before time $T$. See Figure \ref{fig:explicit3} for illustration.

Therefore, in the infinite-mass solution, for any particle arbitrarily far from 0 initially, it always contribute to the mass in $[-\epsilon,\epsilon]\times [-1,1]$ when $t$ is close to $T$. Such contribution, for particles far enough, will not be possible in a finite-mass solution. Therefore we expect that this infinite-mass solution cannot be well-approximated by a finite-mass solution, i.e., a finite-mass solution of \eqref{eqgamma} with a similar initial data would not have a finite-time blow-up with the same mechanism.

Next we analyze what this might imply for the global wellposedness for \eqref{eqn: gke}. Suppose we take the initial data \eqref{ansatz} with $\phi$ supported in $[-1/2,1/2]$. Then the solution $f(t,x,v)$ to \eqref{eqgamma}, for fixed $t,x$, is always supported in $[bx-\frac{1}{2\sqrt{a}},bx+\frac{1}{2\sqrt{a}}]$ in the $v$-direction. For such a function $f$, we always have $|v-w|\le \frac{1}{\sqrt{a}}$ if $v\in \text{supp } f(t,x,\cdot)$. This suggests that the approximation from \eqref{eqn: gke} to \eqref{eqgamma} is making the aggregation stronger. Since \eqref{eqgamma} is not likely to have finite-time blow-up for finite-mass solutions, we deduce heuristically that the same is true for \eqref{eqn: gke} corroborating our numerical findings for $\gamma=3$ in Subsection \ref{sec:sp_inhom_gamma_3}.

\subsection{Numerics for spatially inhomogeneous case with $\gamma=2$}\label{sec:sp_inhom_gamma_2}
This section  dedicates to provide further numerical evidence to our analytical findings for the case $\gamma=2$.
As discussed previously in Section~\ref{sec:intro}, for the spatially homogeneous case, the solution is the critical case between finite-time and infinite time blow-up solution in $v$. As a result, it introduces the more intricate question of whether there will be blow-up in $x$ before the blow-up occurs in $v$, and also whether the shear in $x$ will accelerate the blow-up in $v$.

\begin{figure}[h!]
\centering
{\includegraphics[width=0.5\textwidth]{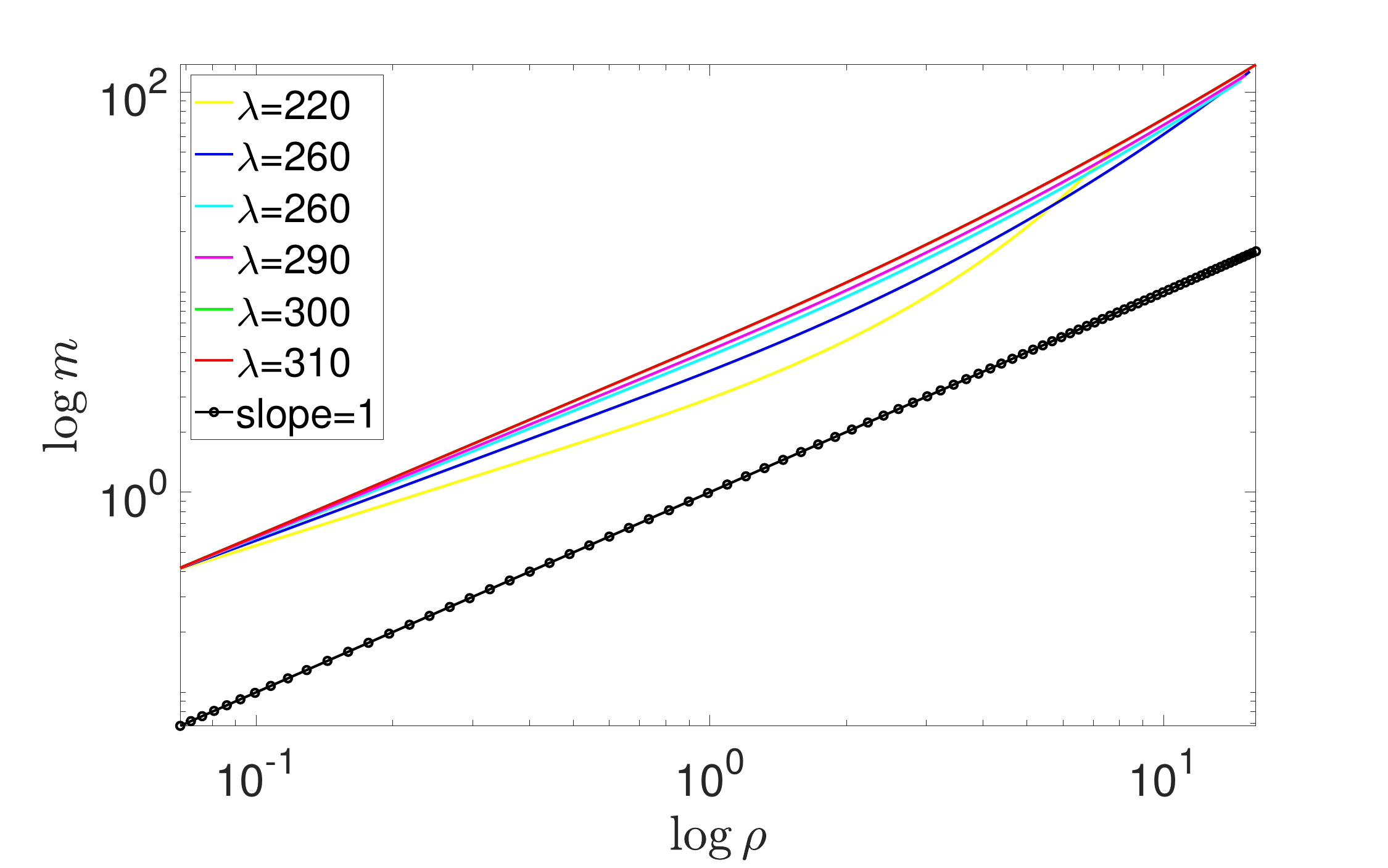}}
\caption{Numerical solution of \eqref{eqn: gke} with kernel $W(v)=|v|^2$ for $Nx=121$, $N_v=201$, $\delta = 0.5$, $\delta_0=0.5$ and adaptive time stepping strategy with $\epsilon=5e-6$. The figure shows that $\frac{\log \rho}{\log m}$ approaches to 1 as $\lambda \rightarrow 310$.}
\label{fig:loglog}
\end{figure}

To examine the solution behavior, we will adhere to the theoretical findings in this section and validate the conjecture in Section~\ref{sec:412}. Specifically, we consider the following initial condition:
\begin{equation}\label{eqn:ic_g2}
f_0(x,v) =
\begin{cases}
\frac{1}{\sqrt{2 \pi}} e^{-a(bx+v)^2}, ~ -x_1 \leq x \leq x_1; \\
\frac{1}{\sqrt{2 \pi}} e^{-a(bx+v)^2} e^{-1000(x-x_1)^2}, ~ \text{else}\,,
\end{cases}
\end{equation}
This initial condition is constructed to resemble the infinite mass case studied in Section~\ref{sec: explicit1}. Here $\rho_0 = \int f_0(0, v) \mathrm{d} v = \sqrt{\frac{1}{2a}}$. It is conjectured that when the relation \eqref{thres} is satisfied, with $a = 120$ and $b=10$, one expects a blow-up in the solution around $T = 0.1$ when $\lambda > \frac{2}{\rho_0 T} = 20\sqrt{240} \approx 310$. The blow-up is anticipated to behave similarly to the infinite mass scenario represented in \eqref{eqn:brhom}. Therefore, our first task is to verify that the asymptotic behavior of the solution indeed follows \eqref{eqn:brhom}. More particularly, as $\lambda$ approaches the blow-up threshold of 310 (i.e., $\frac{\lambda}{2} \rho_0 T \rightarrow 1$), we expect that $\rho \sim \epsilon^{-1}$ and $m \sim \epsilon^{-\frac{\lambda}{2} \rho_0 T}$. However, the exact $\epsilon(t)$ is unknown since the analytic blow-up time is also unknown. Therefore, we numerically check if $\frac{\log \rho}{\log m} \rightarrow 1$ instead, which is verified in Figure~\ref{fig:loglog}.

Next we numerically check the blow-up conditions stated in Theorem~\ref{thm_condbu}. Practically, this involves verifying that the level sets of the numerical solution intersect any vertical line with just one connected component, the numerical results are demonstrated in Figure~\ref{fig:gk2_bu}. Additional results with more values of $\lambda$ can be found in Appendix~\ref{sec:sup_g2}. 

\begin{figure}[h!]
\centering
{\includegraphics[width=0.45\textwidth]{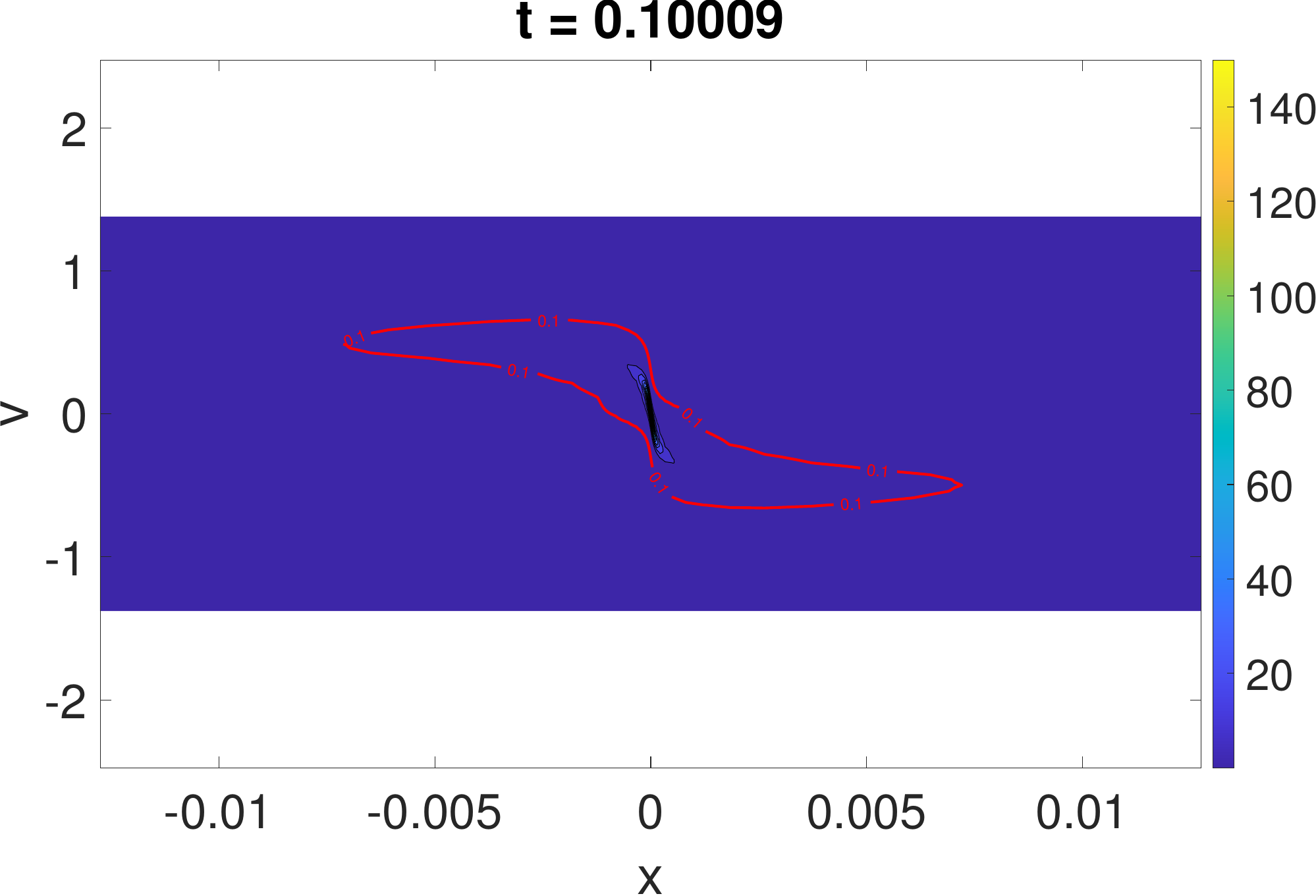}}
{\includegraphics[width=0.45\textwidth]{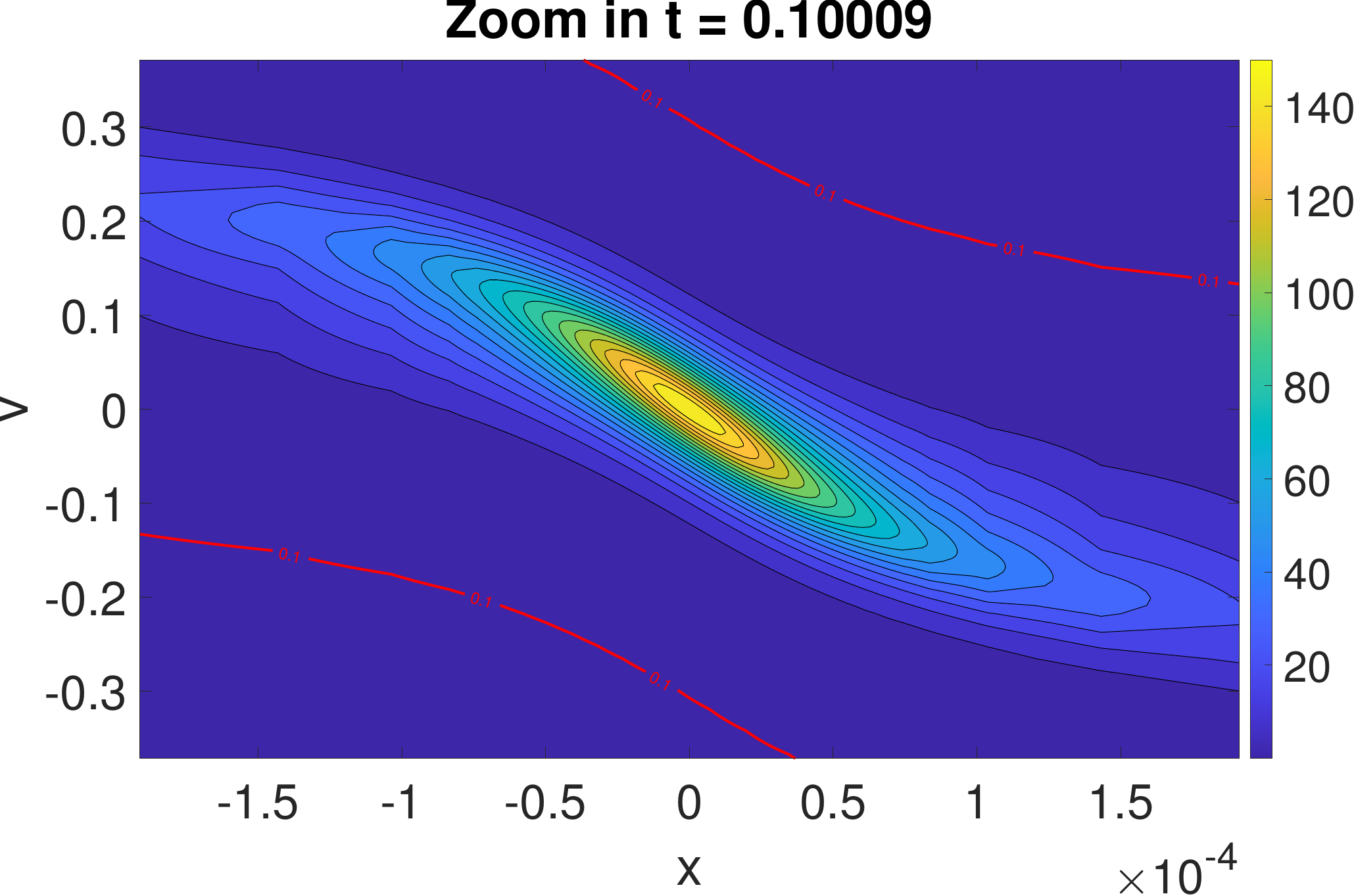}}

{\includegraphics[width=0.45\textwidth]{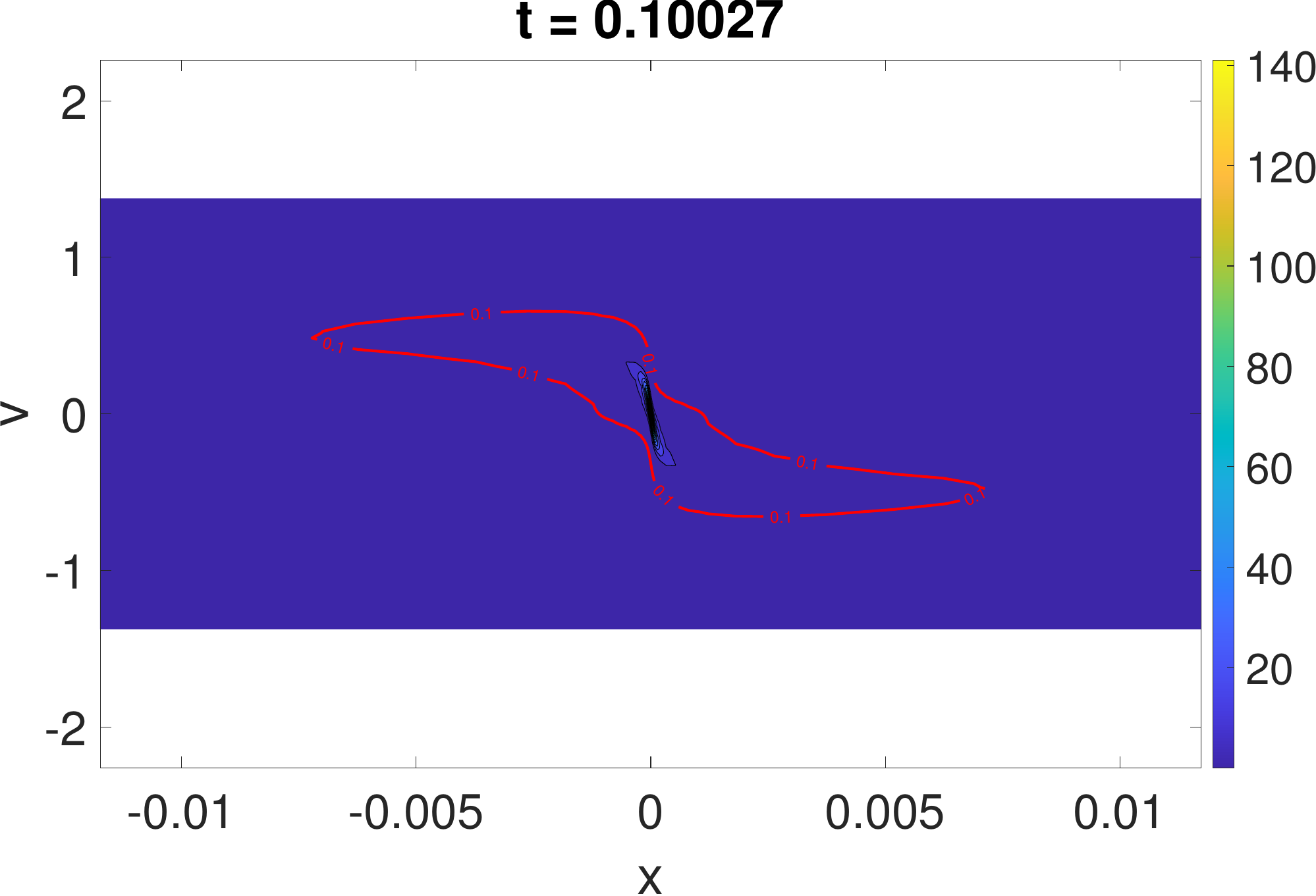}}
{\includegraphics[width=0.45\textwidth]{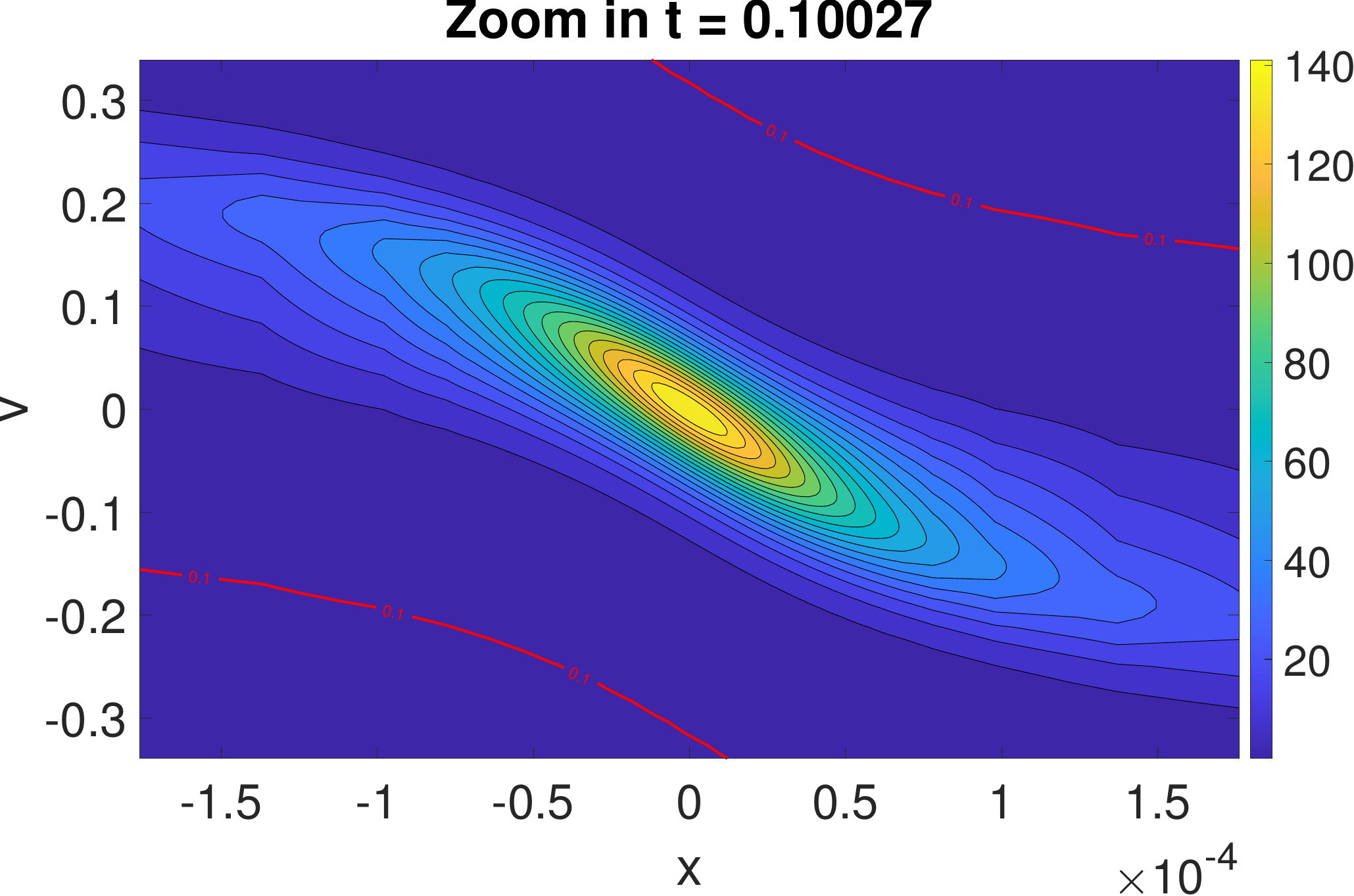}}

{\includegraphics[width=0.45\textwidth]{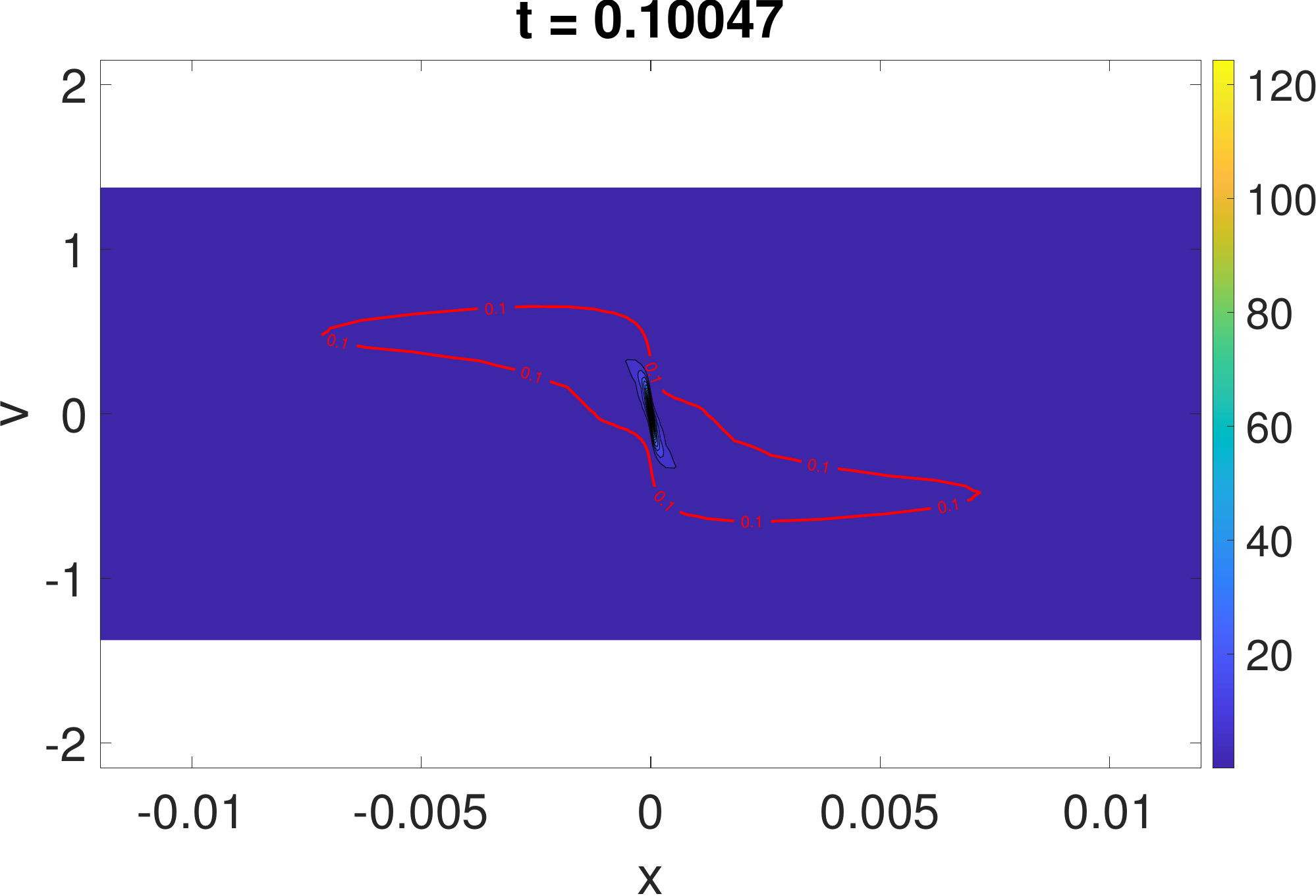}}
{\includegraphics[width=0.45\textwidth]{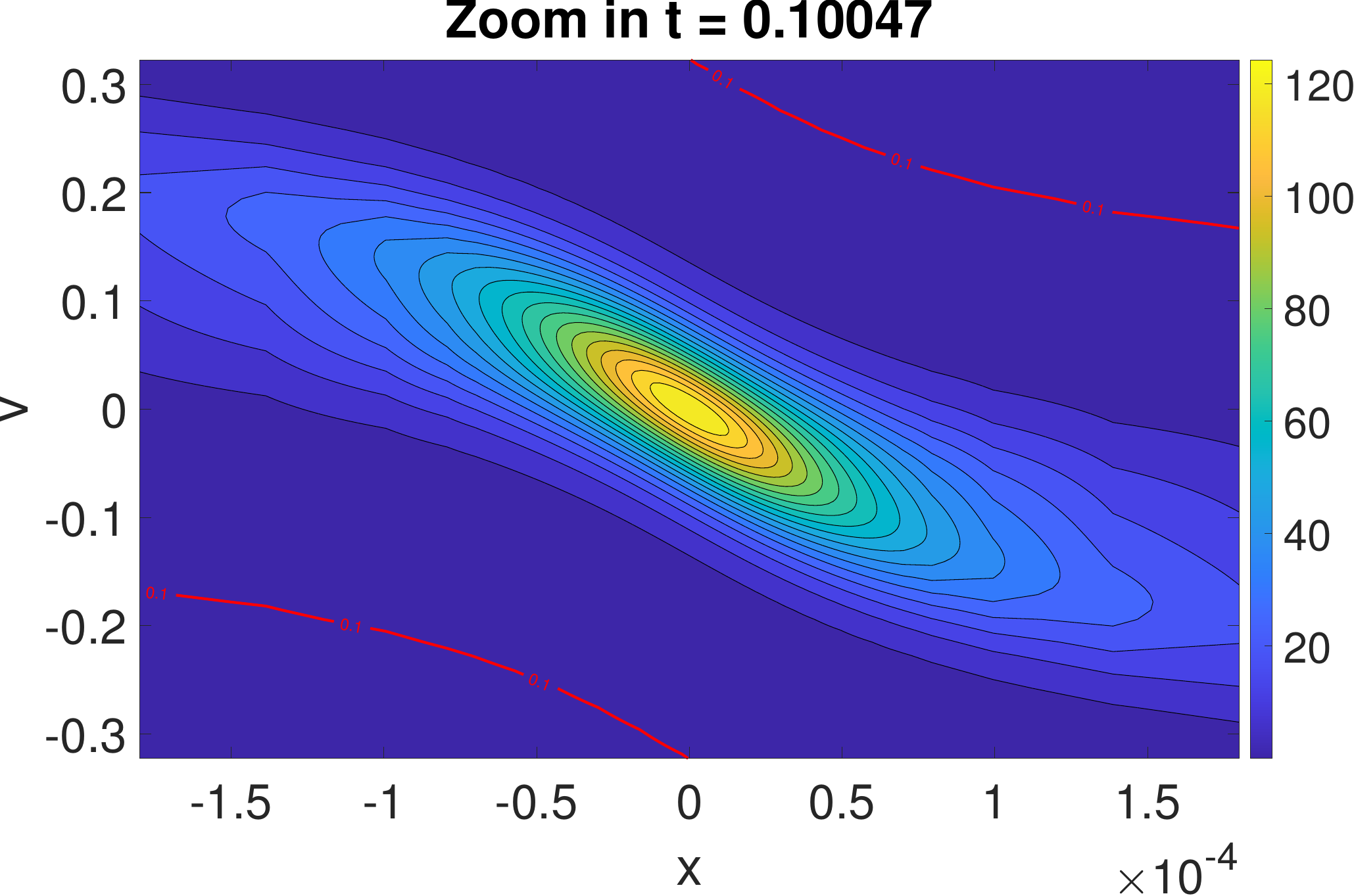}}

{\includegraphics[width=0.45\textwidth]{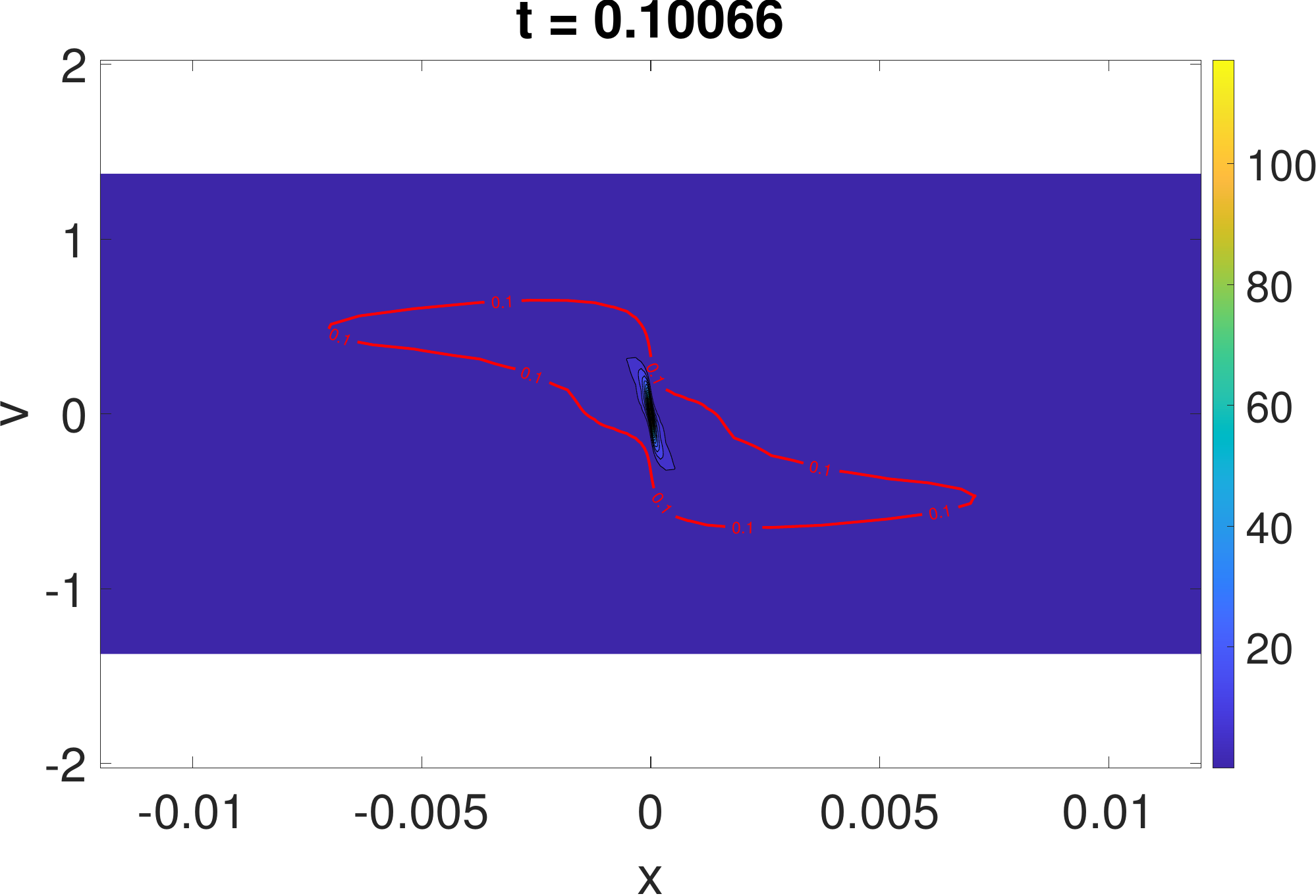}}
{\includegraphics[width=0.45\textwidth]{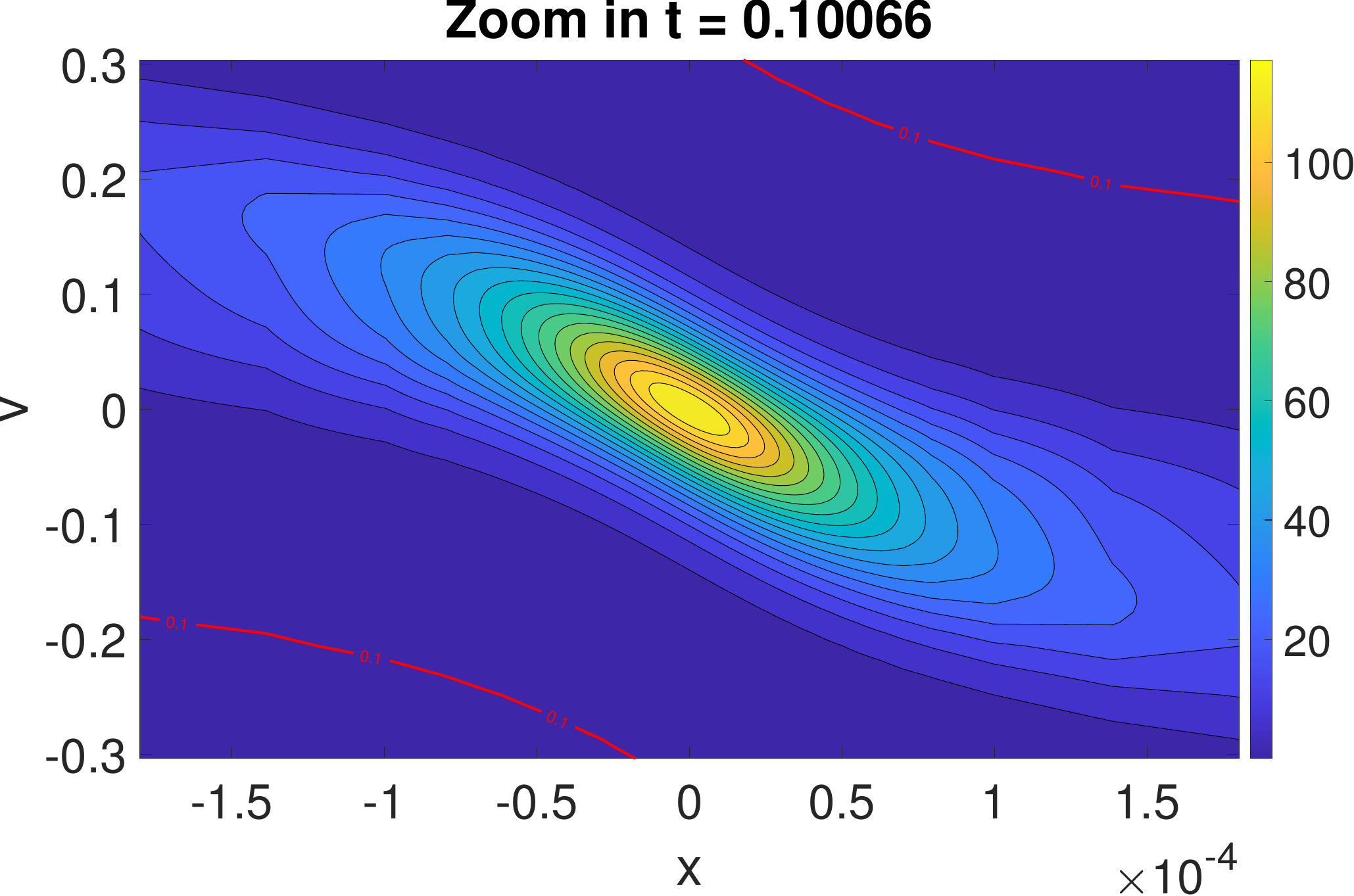}}
\caption{Numerical solution of \eqref{eqn: gke} with kernel $W(v)=|v|^2$ and initial condition~\eqref{eqn:ic_g2} for $\lambda = 310, 300, 290$ and $280$ with $N_x=121$, $N_v=201$, $\delta_0 = 0.5$ and $\delta = 0.5$ in both $x$ and $v$. From top to bottom the plots show numerical solutions and their zoom-in at the numerical blowup time.}\label{fig:gk2_bu}
\end{figure}

Furthermore, by experimenting with various aggregation constants $\lambda$, we are able to numerically validate both the blow-up criterion and the corresponding blow-up times $T=0.1$, see results in Table~\ref{tab:g2bu}. 

\

\begin{table}[h!]
\centering
 \begin{tabular}{|c|c|c|c|c|c|c|c|c|c|} 
 \hline
 $\lambda$ &220 & 260 & 280 & 290 & 300 & 304 & 306 & 308 & 310 \\
 \hline
 $T_b$ & 0.10349 &0.10154
 & 0.10066 
 & 0.10047 
 & 0.10027 
  & 0.10018 
  & 0.10016
  & 0.10012
 & 0.10009 \\
\hline
 \end{tabular}
 \caption{Numerical solution of \eqref{eqn: gke} with kernel $W(v)=|v|^2$ for $N_x=121$, $N_v=201$, $\delta = 0.5$, $\delta_0=0.5$ and adaptive time stepping strategy with $\epsilon=5e-6$. The table shows the numerical blow-up time $T_b$ for various $\lambda$.}\label{tab:g2bu}
\end{table}

\

\begin{figure}[h!]
\centering
{\includegraphics[width=0.45\textwidth]{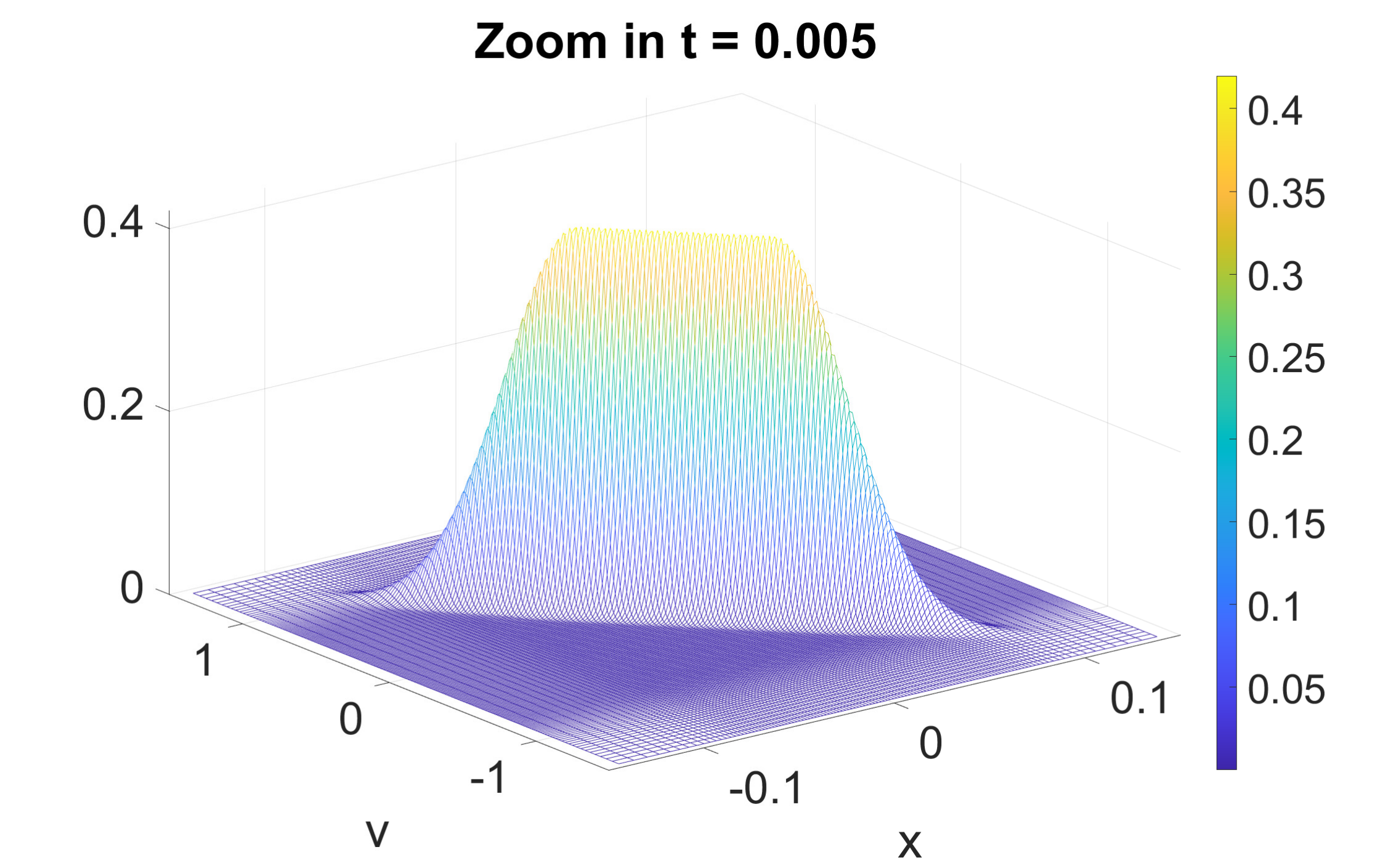}}
{\includegraphics[width=0.45\textwidth]{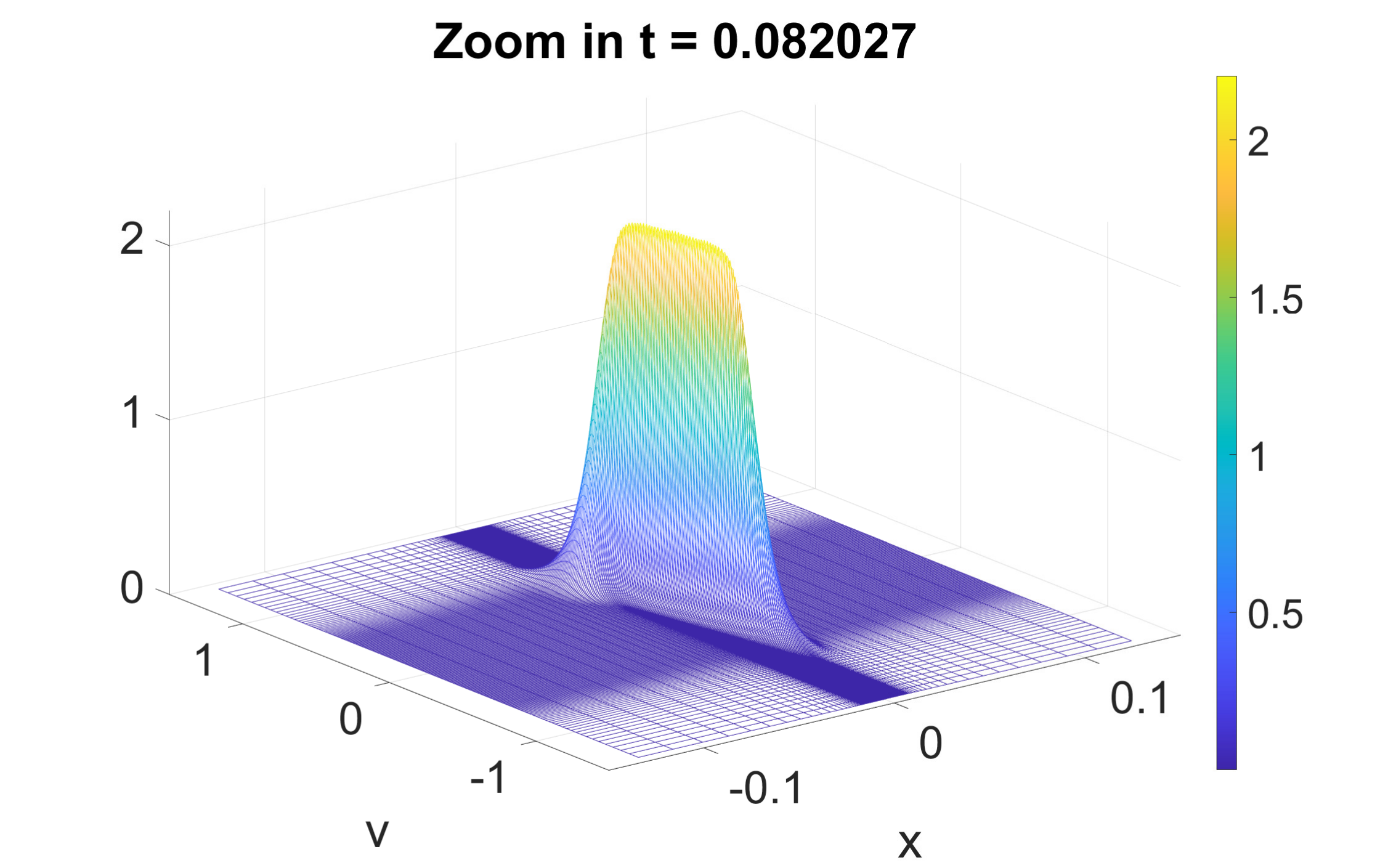}}

{\includegraphics[width=0.45\textwidth]{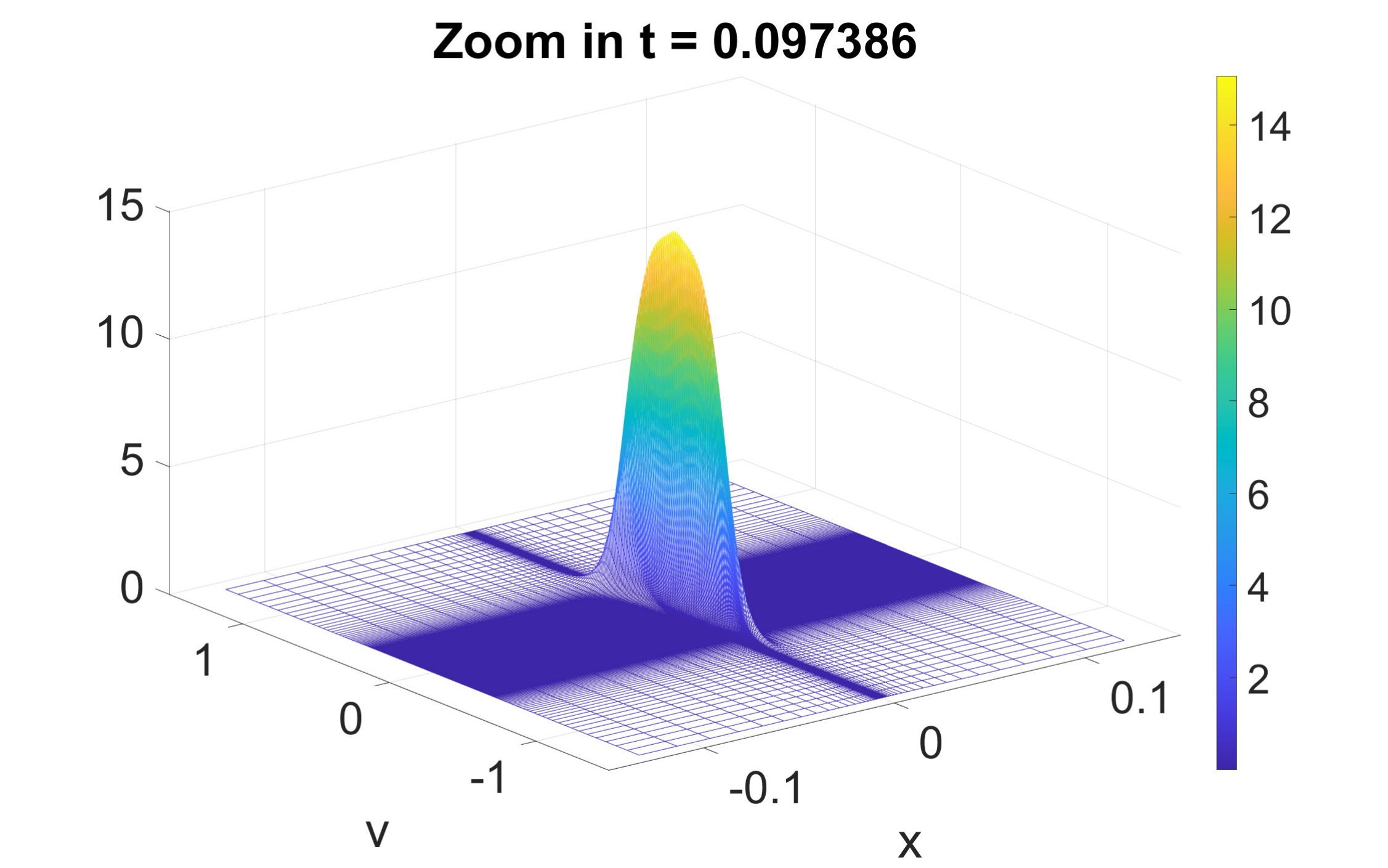}}
{\includegraphics[width=0.45\textwidth]{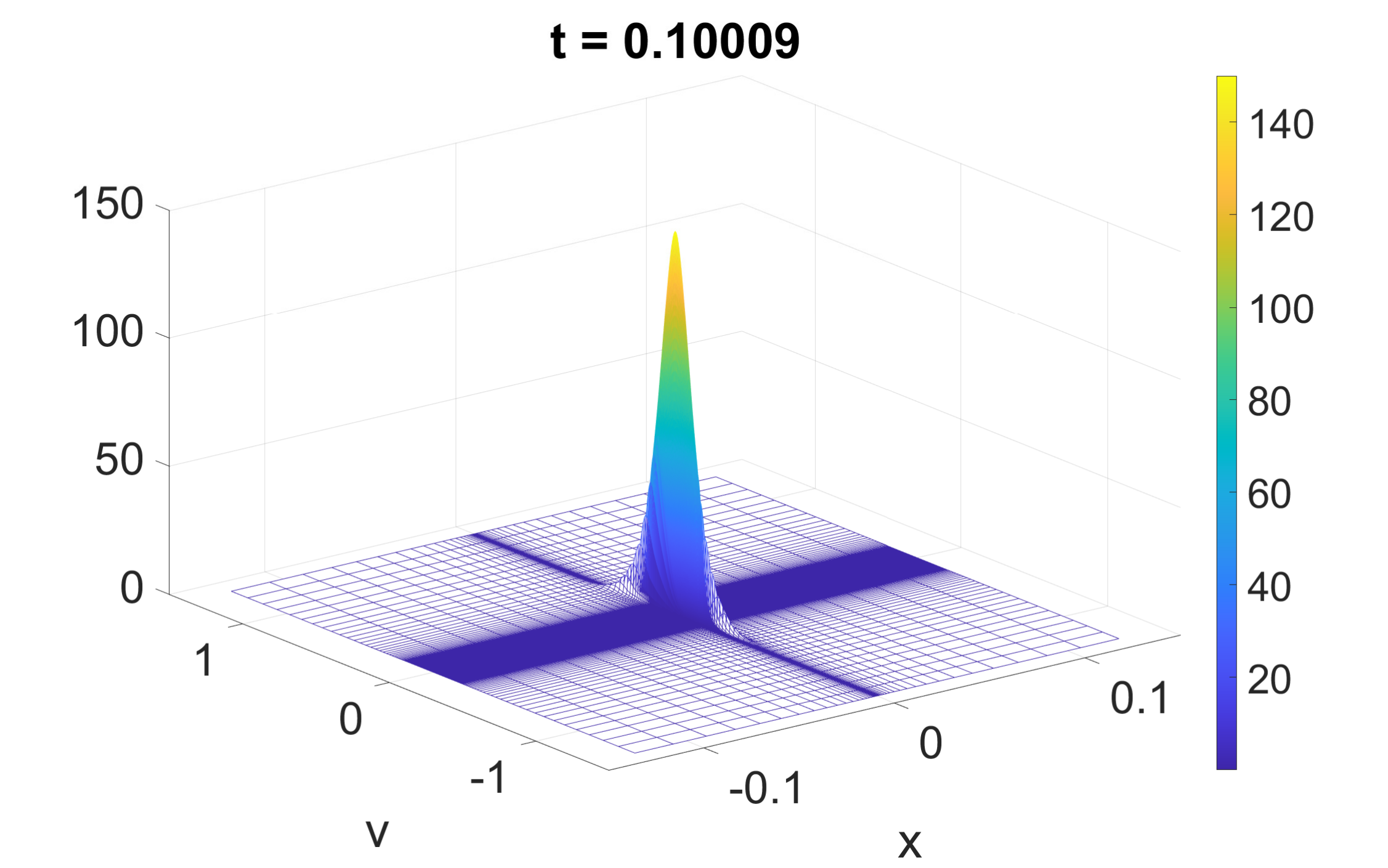}}

{\includegraphics[width=0.45\textwidth]{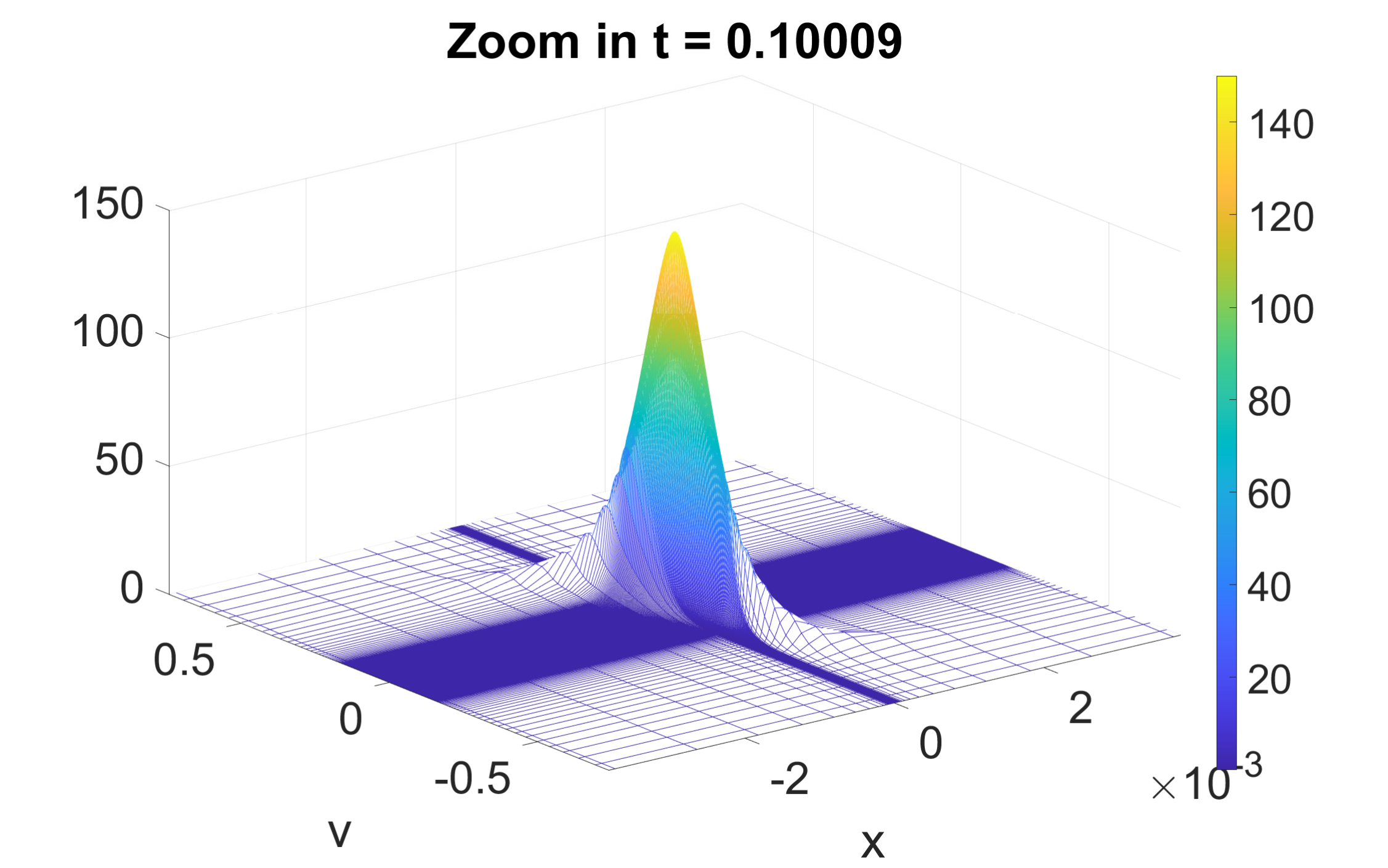}}
{\includegraphics[width=0.45\textwidth]{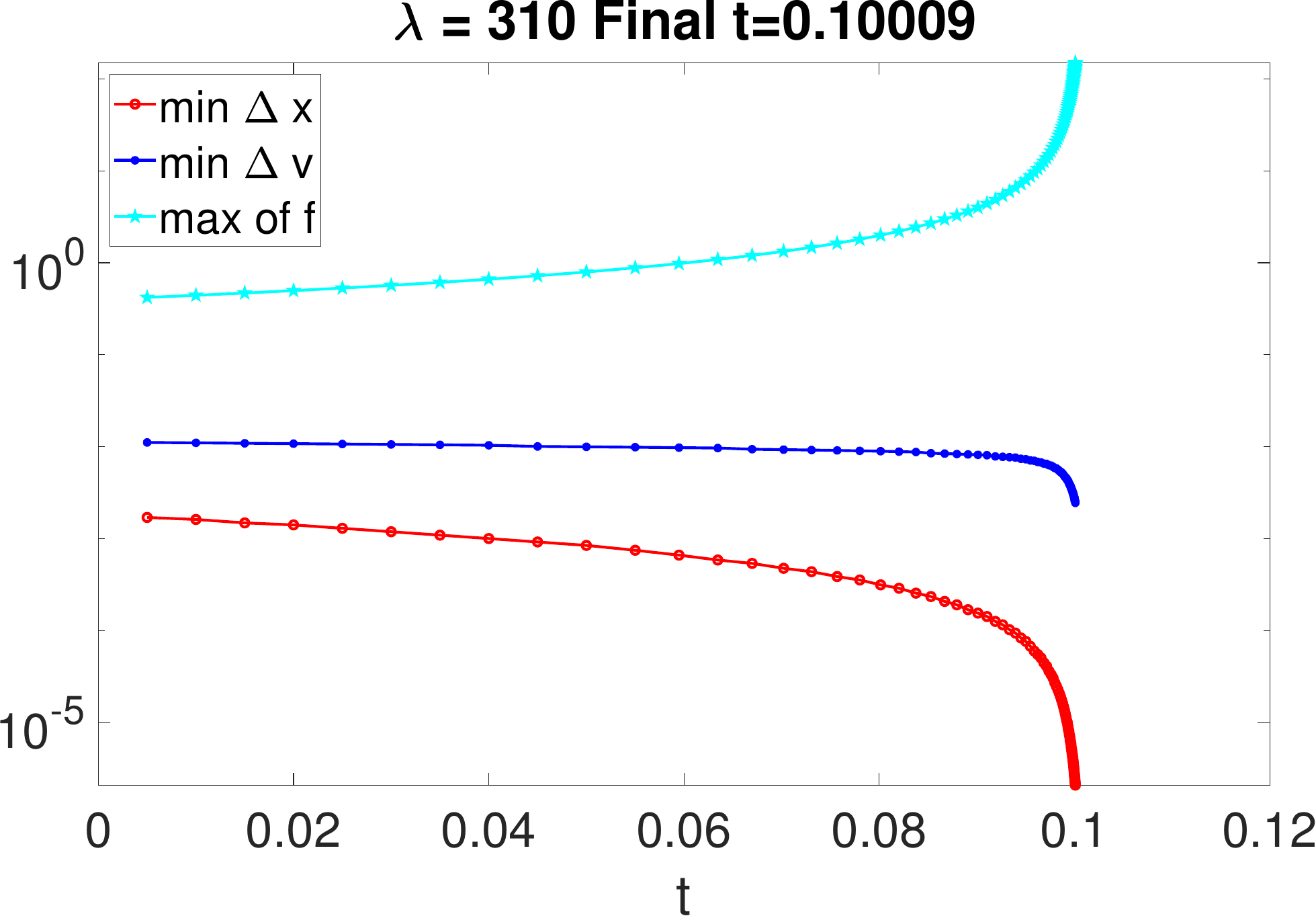}}

\caption{Numerical solution of \eqref{eqn: gke} with kernel $W(v)=|v|^2$ for $N_x=121$, $N_v=201$, $\lambda = 310$. Top two rows are time snapshots. Bottom left: the zoom-in plot at numerical blow-up time. Bottom right: record of minimum of $\Delta x$, $\Delta v$ and maximum of $f(t,x,v)$ at each time step.}\label{fig:gk2_evo_example2}
\end{figure}

In addition, we present the evolution for $\lambda=310$ in Figure~\ref{fig:gk2_evo_example2}. In contrast to the slow decay and the eventual increase of the minimum of $\Delta x$ in Figure~\ref{fig:gk3_ih_2_C5}, the minimum of $\Delta x$ decreases dramatically and rapidly falls below the blow-up criterion, which provides additional evidence of blow-up on $x$ direction.

\section{Conclusion}
In this paper, we investigate the well-posedness of the granular kinetic equation, a challenging problem characterized by the possibility of finite or infinite-time blow-up. To tackle this challenge, we propose a structure-preserving numerical scheme designed to mimic the underlying gradient flow of the associated energy functional. Consequently, our scheme upholds key physical properties, including conservation and entropy decay. Additionally, we integrate a mesh refinement mechanism that dynamically assigns a denser grid to regions of high concentration. This approach enhances resolution in concentrated regions, enabling the detection of potential blow-ups in the solution. To employ this numerical tool for investigating solution behavior, we meticulously choose a set of representative initial conditions for our experiments. Our numerical findings can be summarized in three aspects: 
\begin{itemize}
\item[1.] For $\gamma=3$, despite the initial formation of a pronounced concentration in the $x$ direction, the  transport effect eventually takes over and disperses the concentration. This observation prompts us to conjecture that for $\gamma>2$, the solution will not exhibit finite-time blow-up in either the $x$ or $v$ direction. 

\item[2.] In the critical case with $\gamma = 2$, through the construction of special solutions, we establish that there is always a finite-time blow-up for the infinite-mass solution. Subsequently, we perform a finite-mass approximation to the infinite-mass scenario, and numerical evidence suggests that the blow-up can occur in the $x$ direction within a finite time. 

\item[3.] For $\gamma<2$, we presented two examples, both demonstrating finite-time blow-up in the $v$ direction. It remains an open question whether a blow-up will also occur in the $x$ direction. While we believe that blow-up in both $x$ and $v$ directions can occur simultaneously, constructing a suitable initial condition for this scenario poses a significant challenge and will be a topic for future investigation.
\end{itemize}

\appendix
\section{Appendix}\label{sec:appendix}
Some Supplementary numerical examples are presented for completeness. 
\subsection{Supplementary results for spatially inhomogeneous cases with $\gamma=3$}\label{sec:sup_g3}
\begin{figure}[h!]
\centering
{\includegraphics[width=0.39\textwidth]{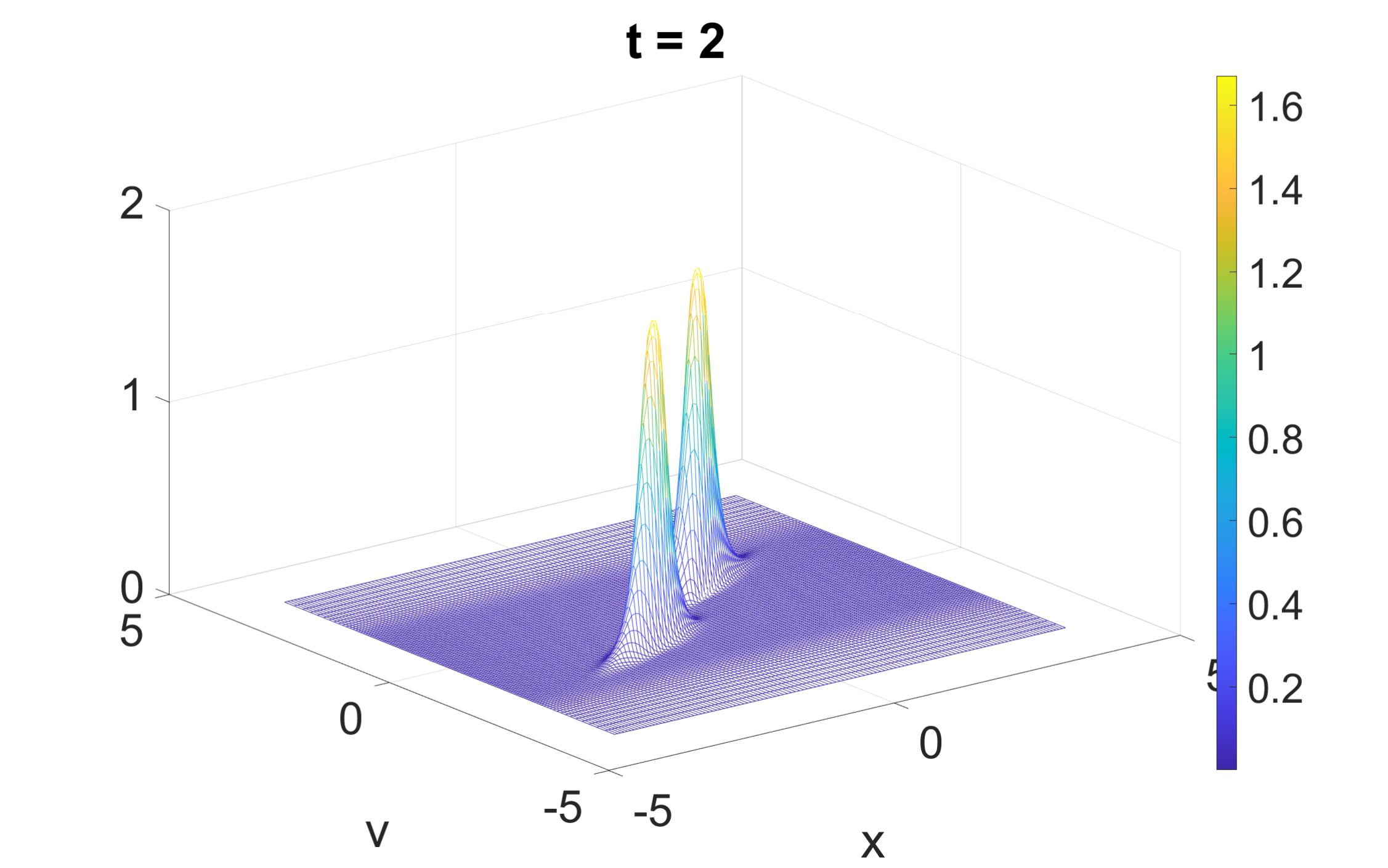}}
{\includegraphics[width=0.39\textwidth]{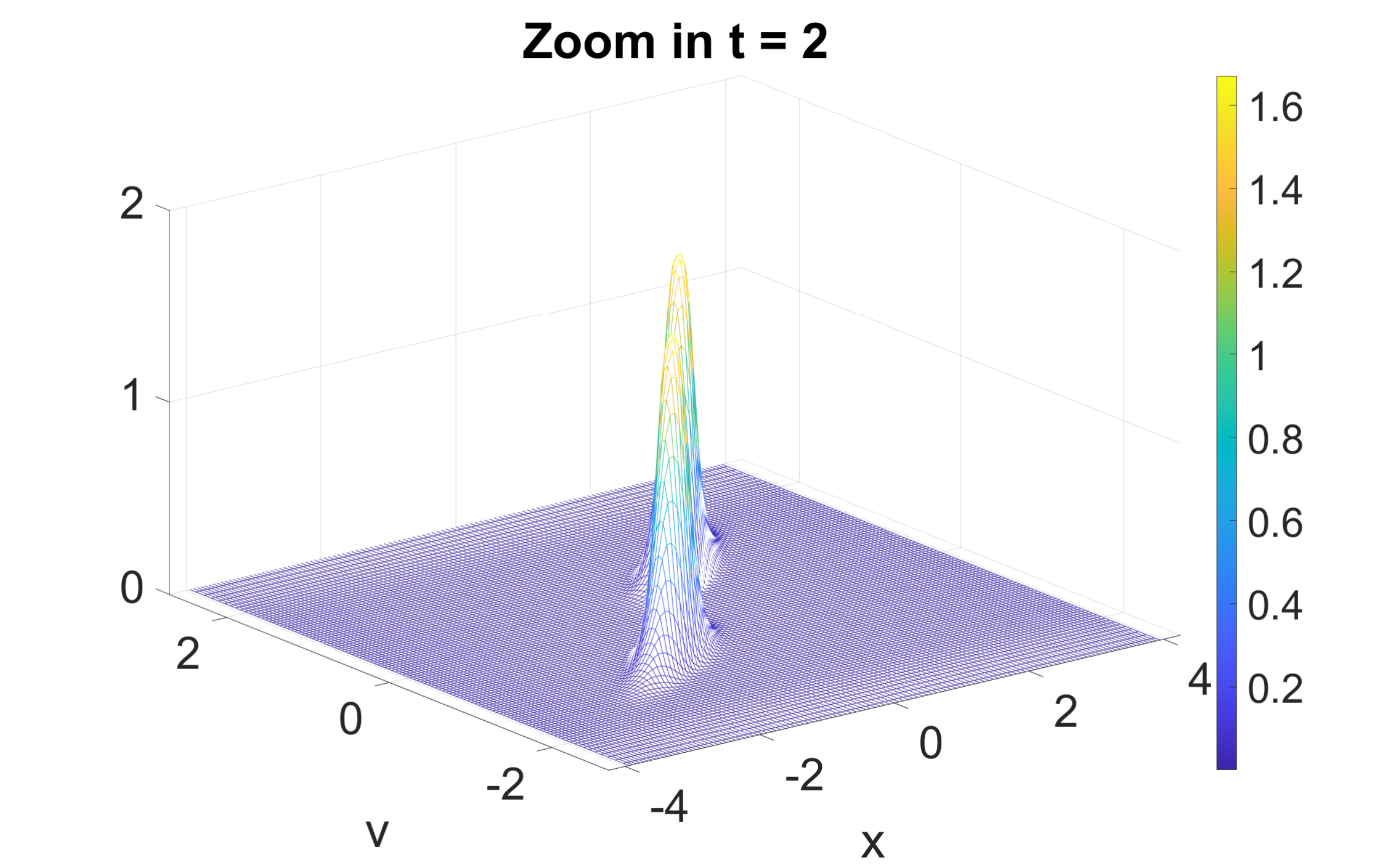}}
{\includegraphics[width=0.39\textwidth]{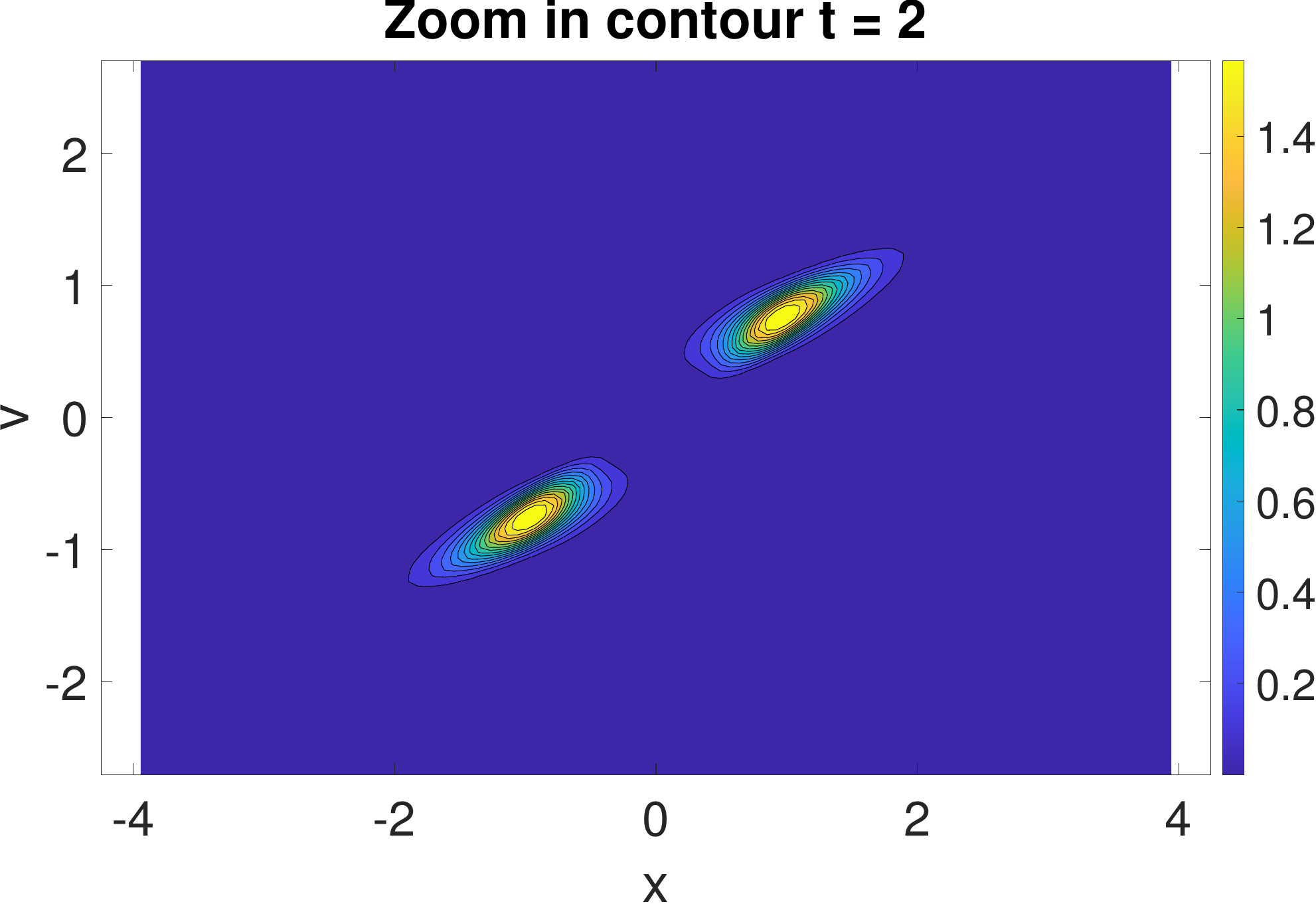}}
{\includegraphics[width=0.39\textwidth]{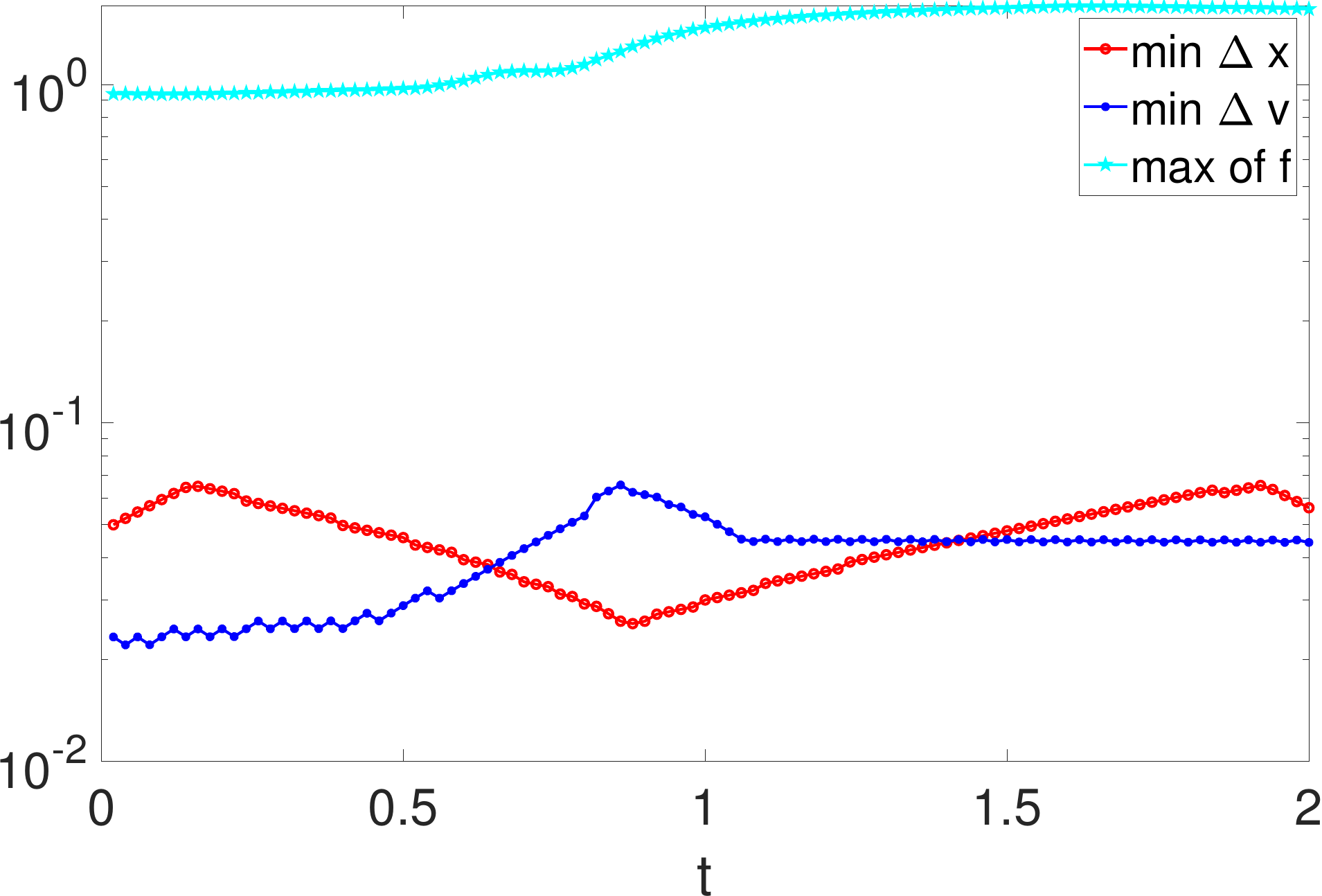}}
 \caption{Numerical solution of \eqref{eqn: gke} with kernel $W(v)=|v|^3$ for $N_x=121$, $N_v=121$, and $\lambda=2$. Top left is the $f(t=4,x,v)$, top right is the zoom in of $f(t=4,x,v)$, bottom left is the contour zoom in plot of $f(t=4,x,v)$ and bottom right is the record of minimum of $\Delta x$, $\Delta v$ and maximum of $f(t,x,v)$ at each time step.}
 \label{fig:gk3_ih_2_C1}
\end{figure}

\begin{figure}[h!]
\centering
{\includegraphics[width=0.39\textwidth]{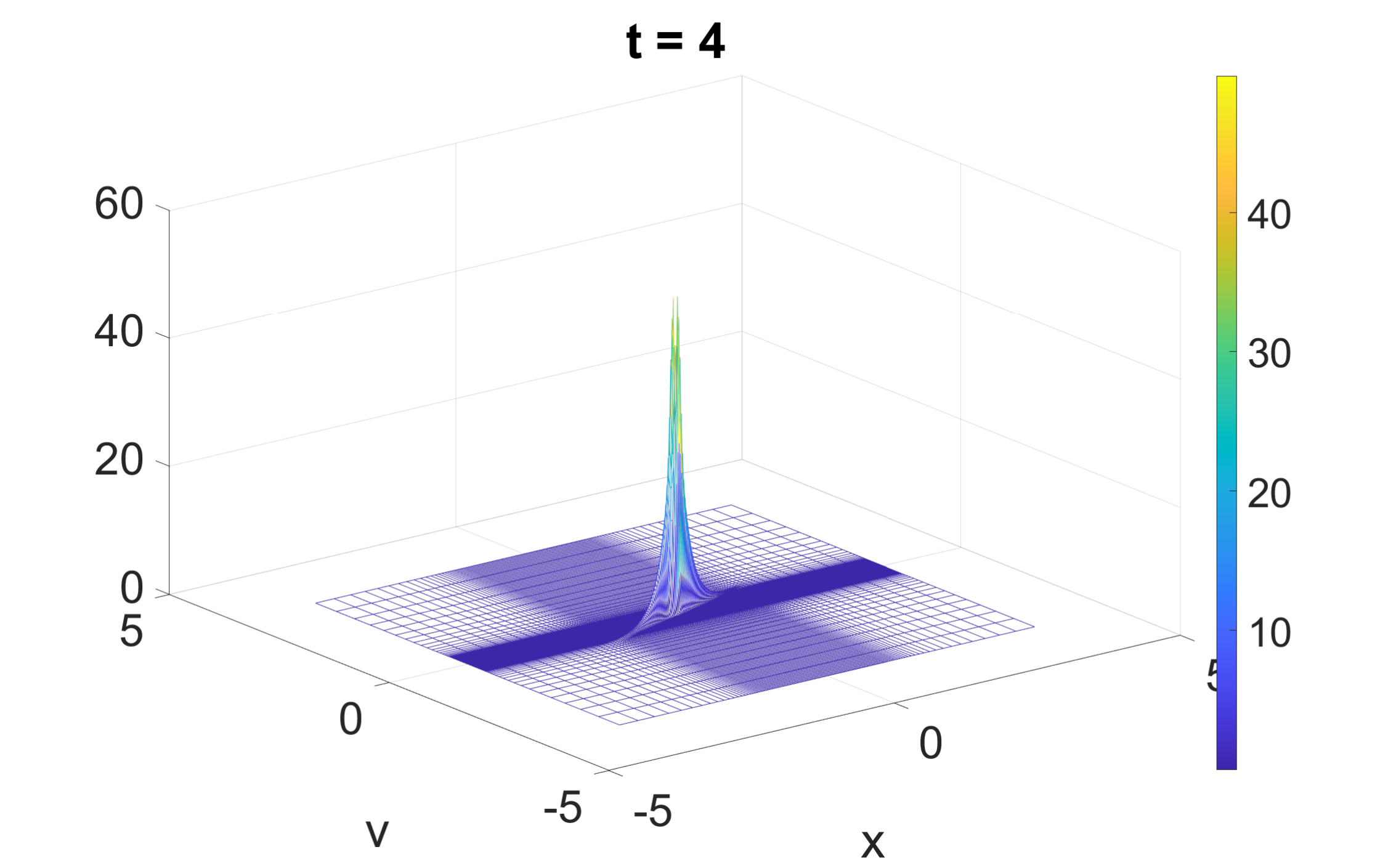}}
{\includegraphics[width=0.39\textwidth]{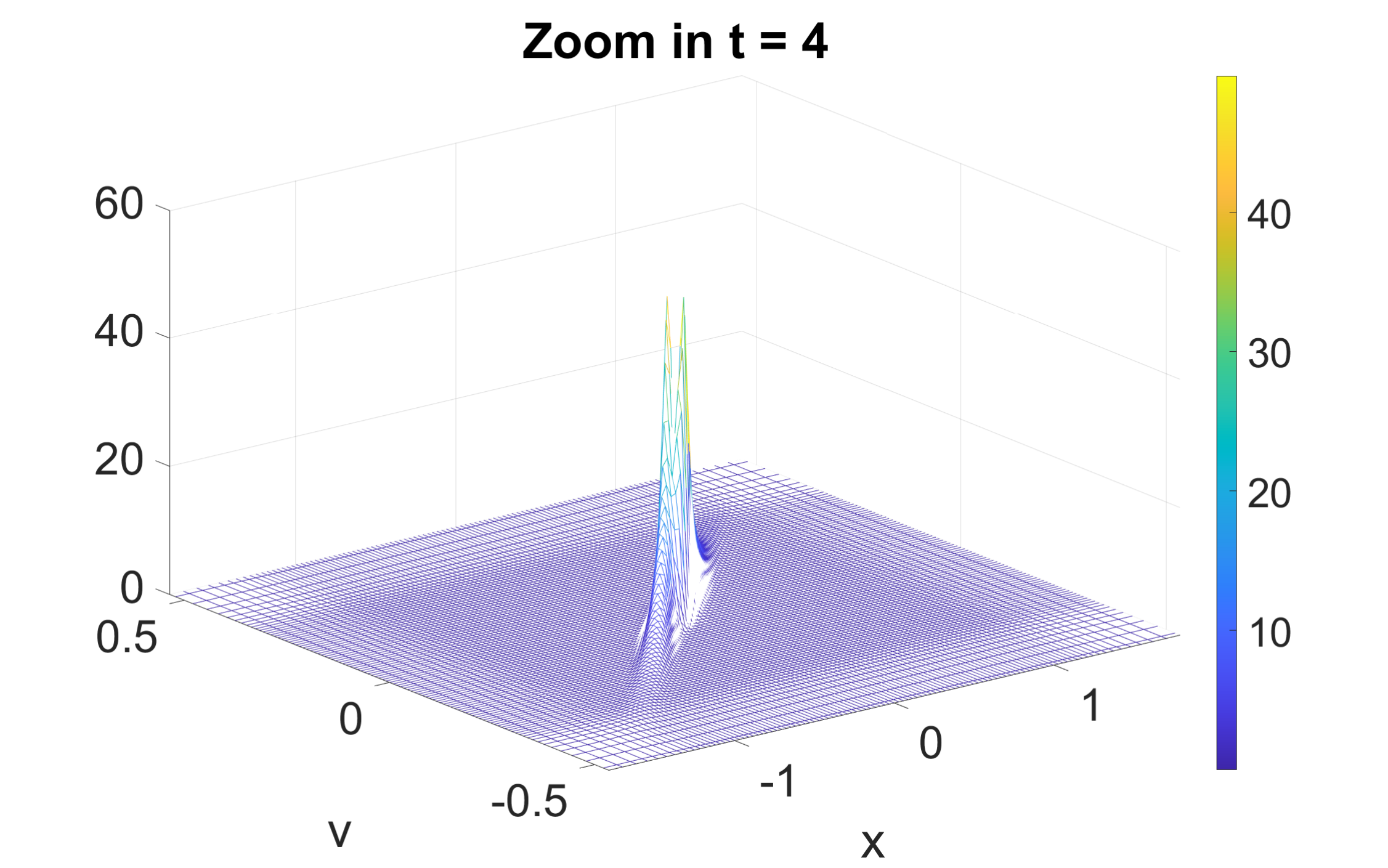}}
{\includegraphics[width=0.39\textwidth]{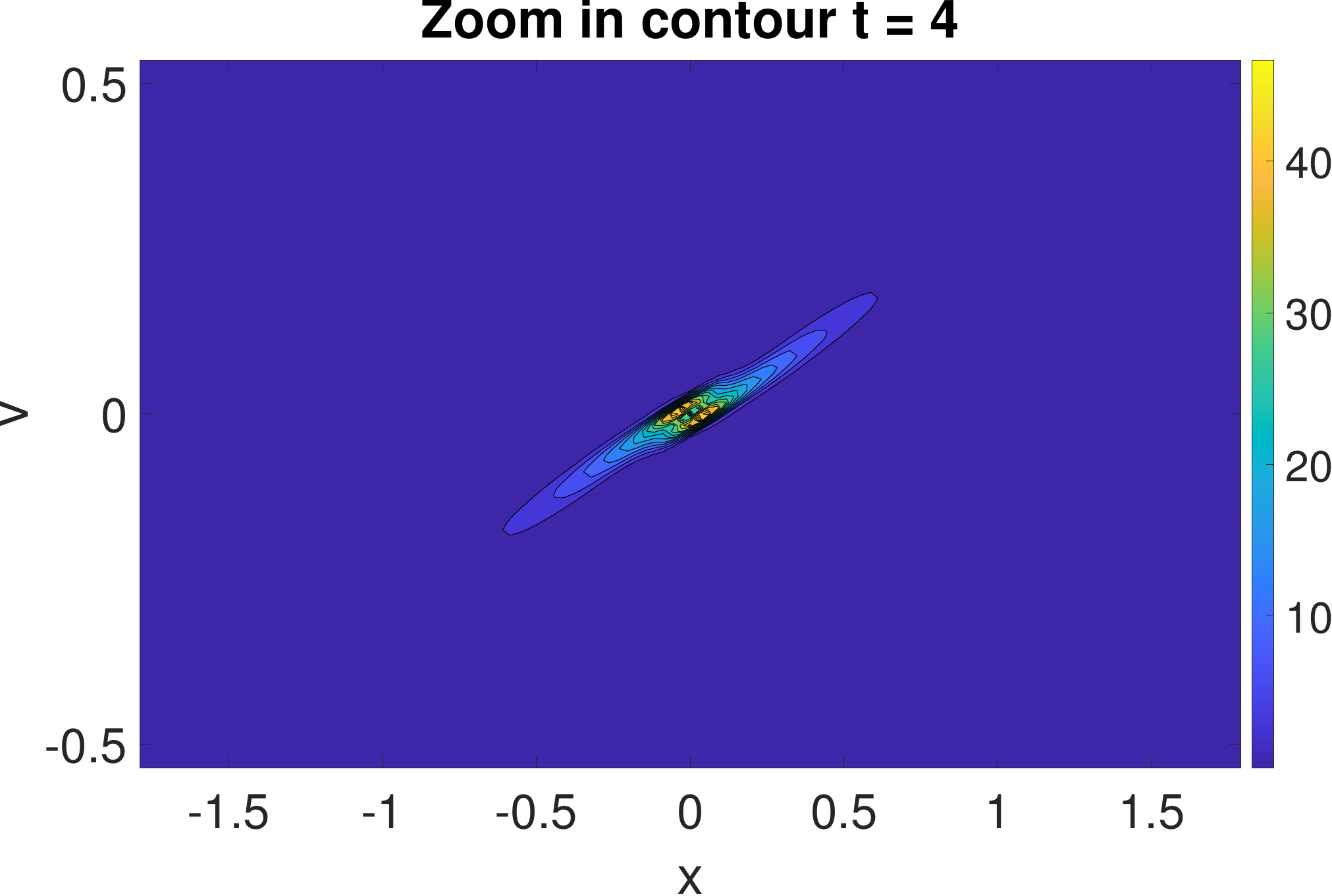}}
{\includegraphics[width=0.39\textwidth]{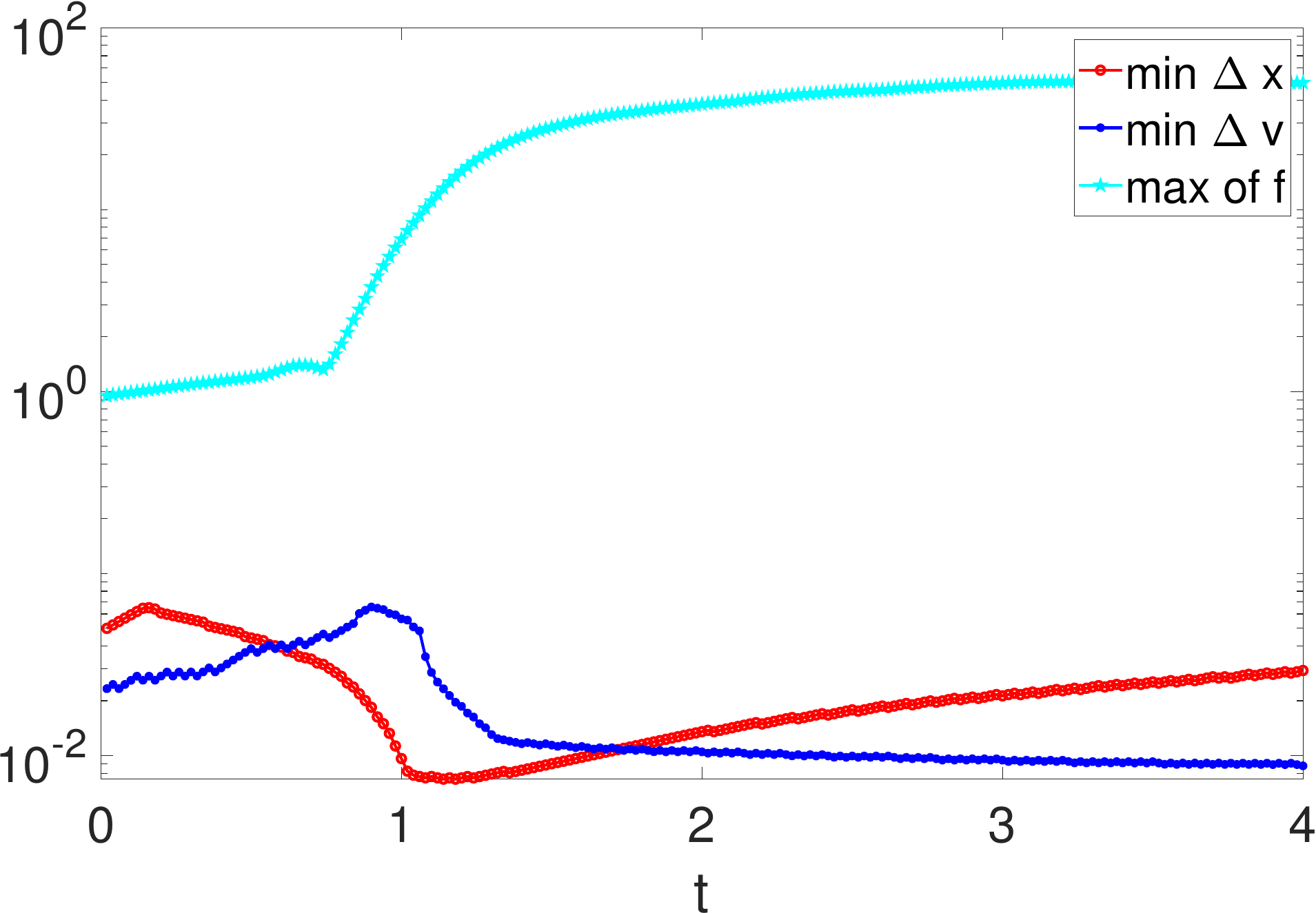}}
 \caption{Numerical solution of \eqref{eqn: gke} with kernel $W(v)=|v|^3$ for $N_x=121$, $N_v=121$, and $\lambda=4$. Top left is the $f(t=4,x,v)$, top right is the zoom in of $f(t=4,x,v)$, bottom left is the contour zoom in plot of $f(t=4,x,v)$ and bottom right is the record of minimum of $\Delta x$, $\Delta v$ and maximum of $f(t,x,v)$ at each time step.}
 \label{fig:gk3_ih_2_C3}
\end{figure}

\begin{figure}[h!]
\centering
{\includegraphics[width=0.39\textwidth]{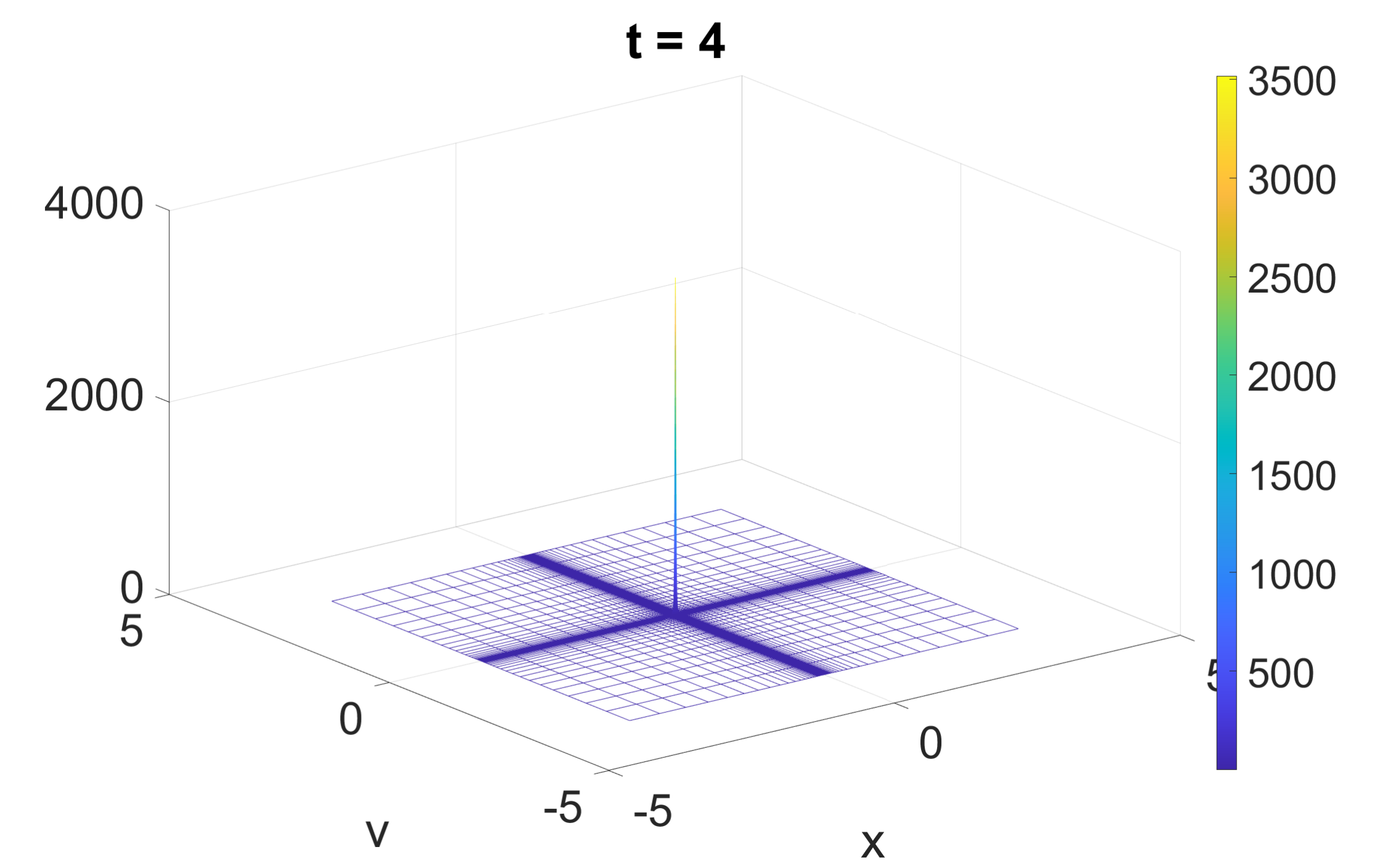}}
{\includegraphics[width=0.39\textwidth]{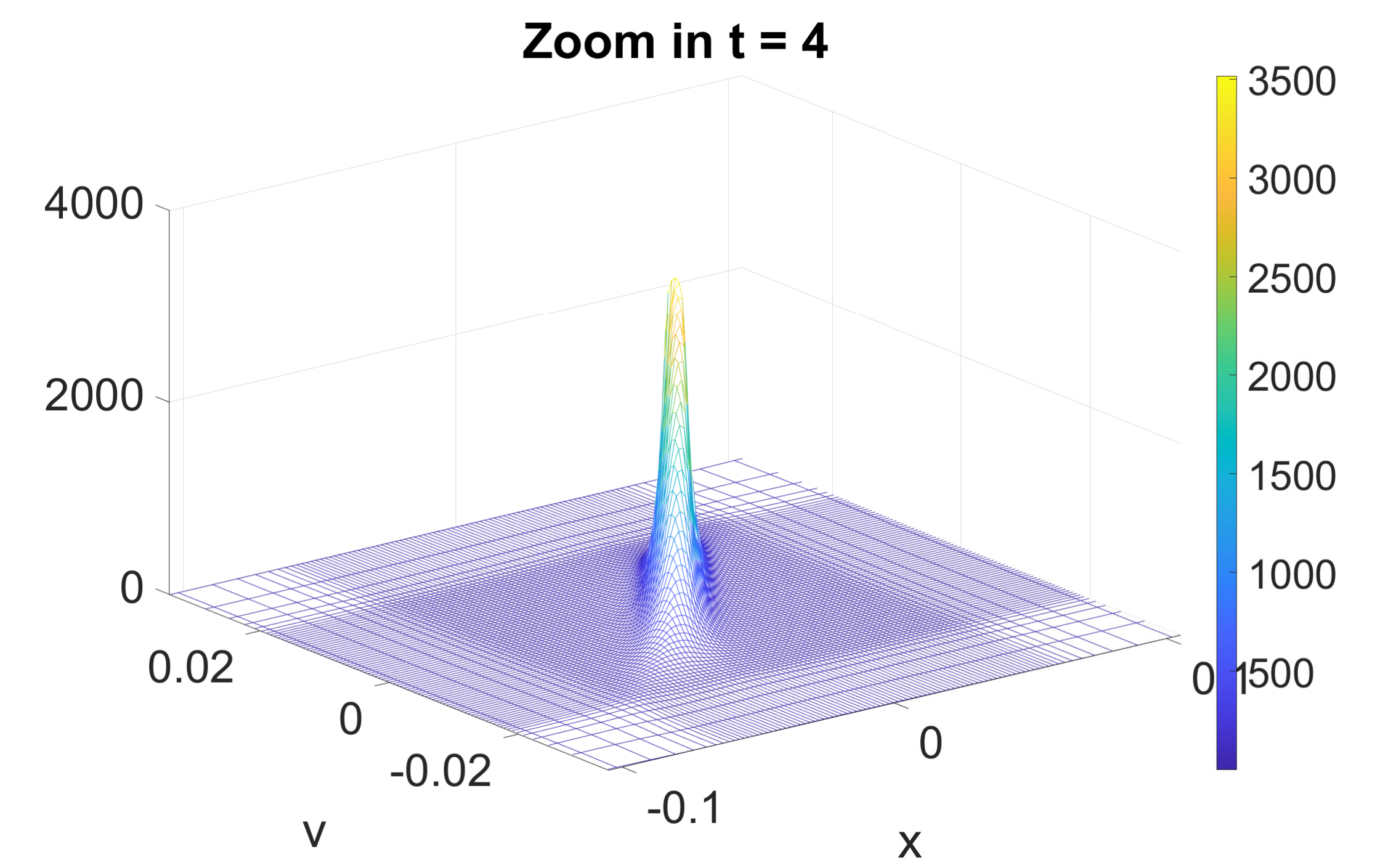}}
{\includegraphics[width=0.39\textwidth]{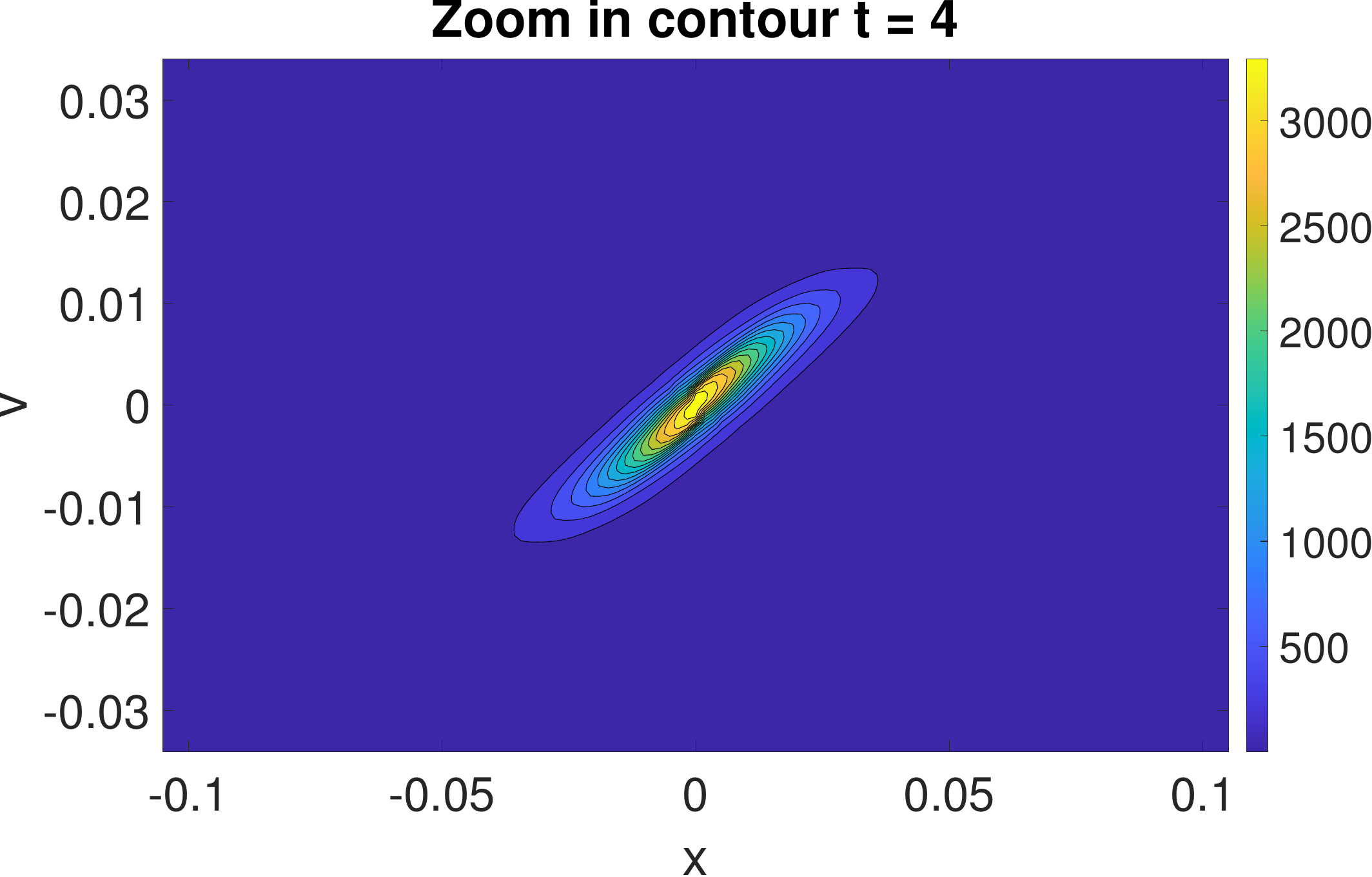}}
{\includegraphics[width=0.39\textwidth]{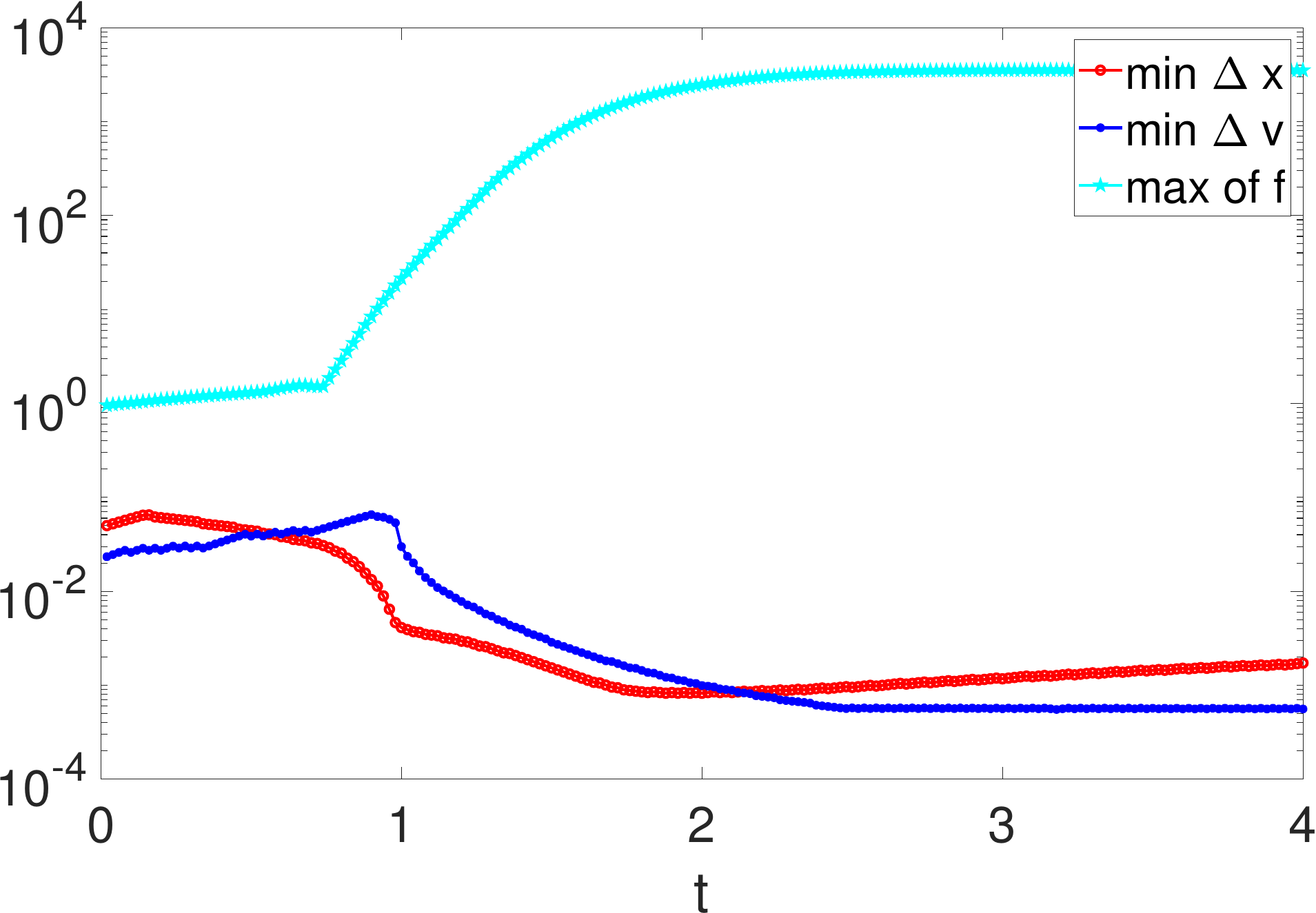}}
 \caption{Numerical solution of \eqref{eqn: gke} with kernel $W(v)=|v|^3$ for $N_x=121$, $N_v=121$, and $\lambda=6$. Top left is the $f(t=4,x,v)$, top right is the zoom in of $f(t=4,x,v)$, bottom left is the contour zoom in plot of $f(t=4,x,v)$ and bottom right is the record of minimum of $\Delta x$, $\Delta v$ and maximum of $f(t,x,v)$ at each time step.}
 \label{fig:gk3_ih_2_C4}
\end{figure}

Additionally, numerical results for $\lambda=2$, $\lambda=6$ and $\lambda=8$ are presented in Figure~\ref{fig:gk3_ih_2_C1}, Figure~\ref{fig:gk3_ih_2_C3} and Figure~\ref{fig:gk3_ih_2_C4} respectively. 
\newpage

\subsection{Supplementary results for spatially inhomogeneuous cases with $\gamma=2$}\label{sec:sup_g2}
\begin{figure}[h!]
\centering
{\includegraphics[width=0.41\textwidth]{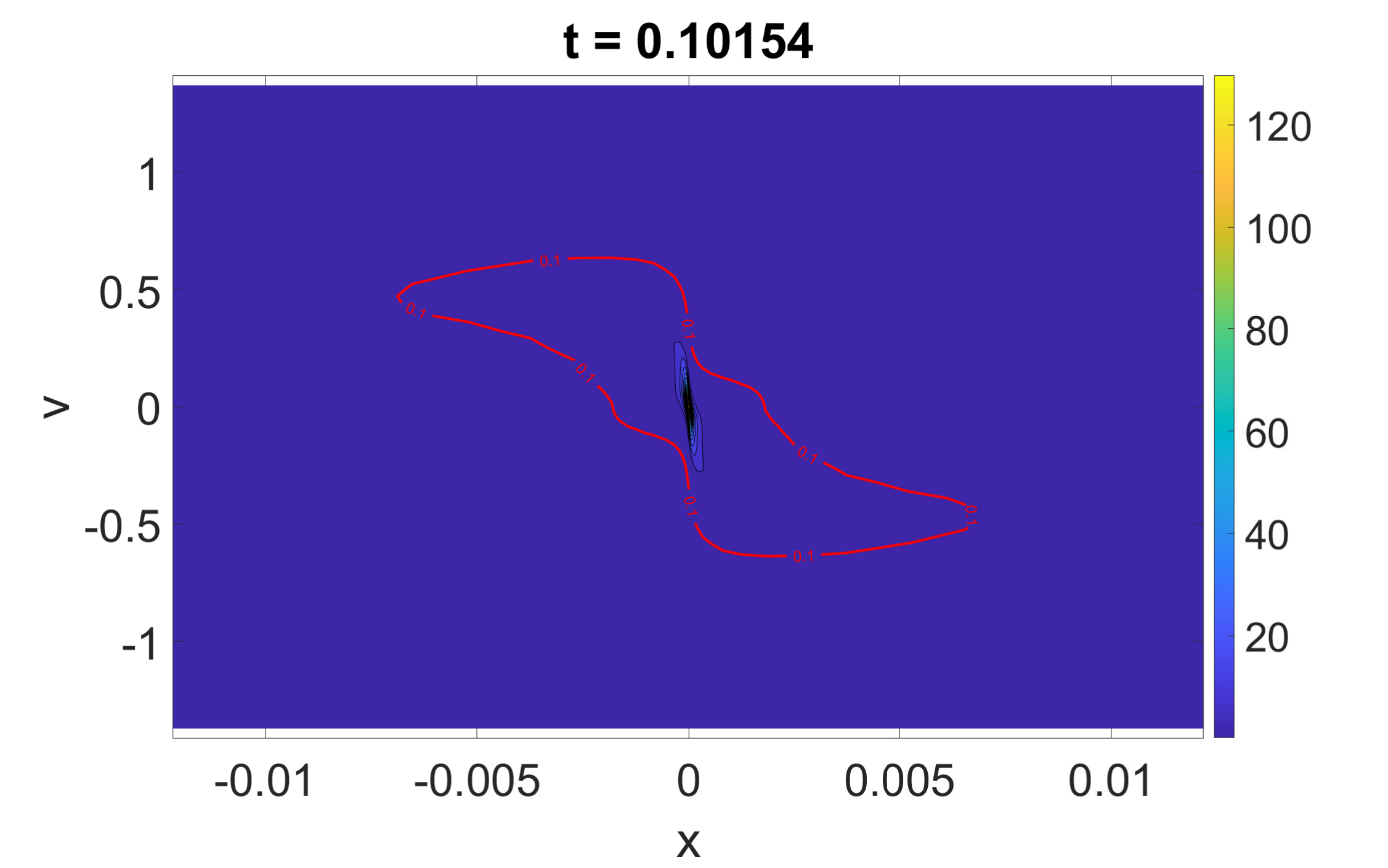}}
{\includegraphics[width=0.41\textwidth]{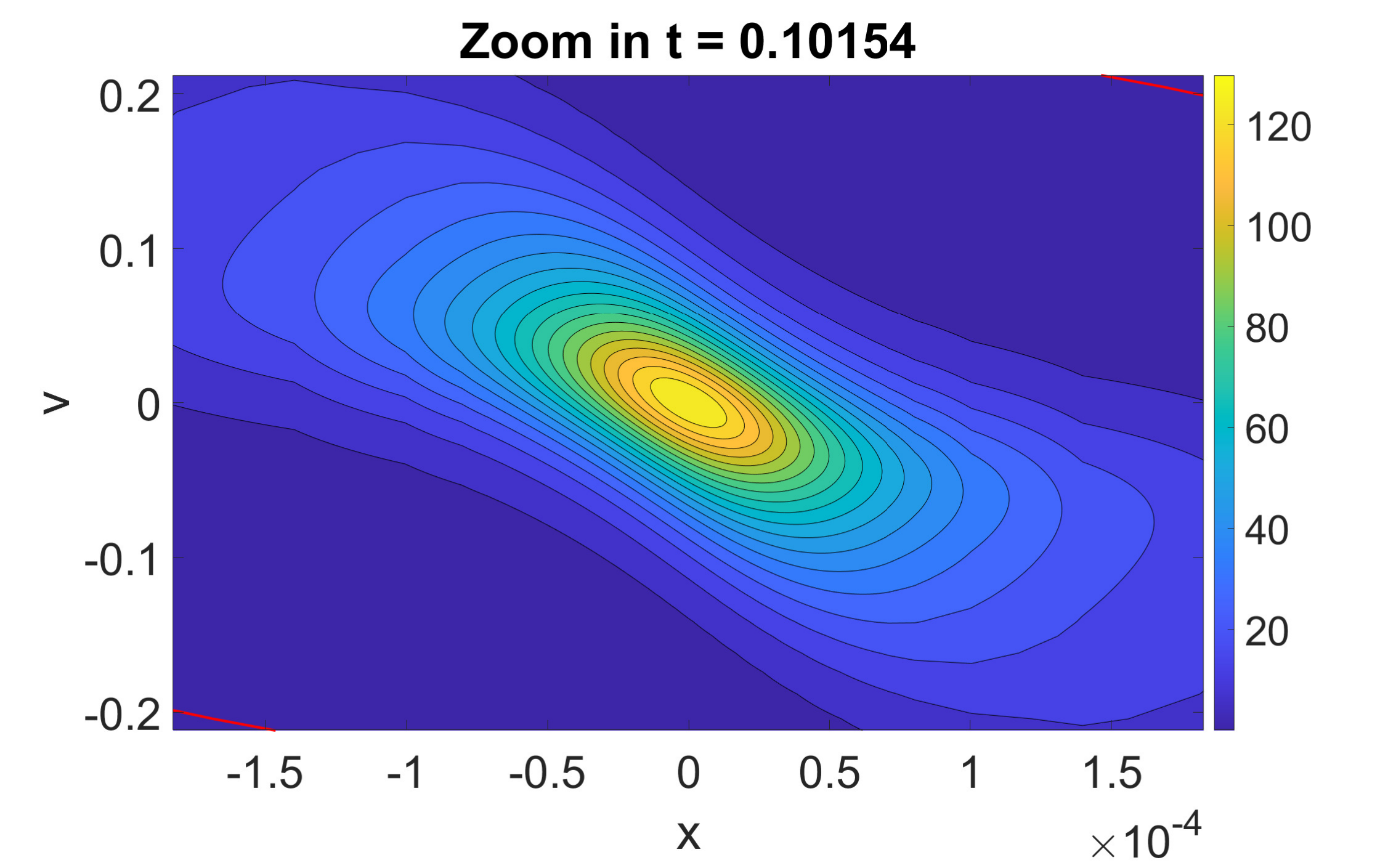}}

{\includegraphics[width=0.41\textwidth]{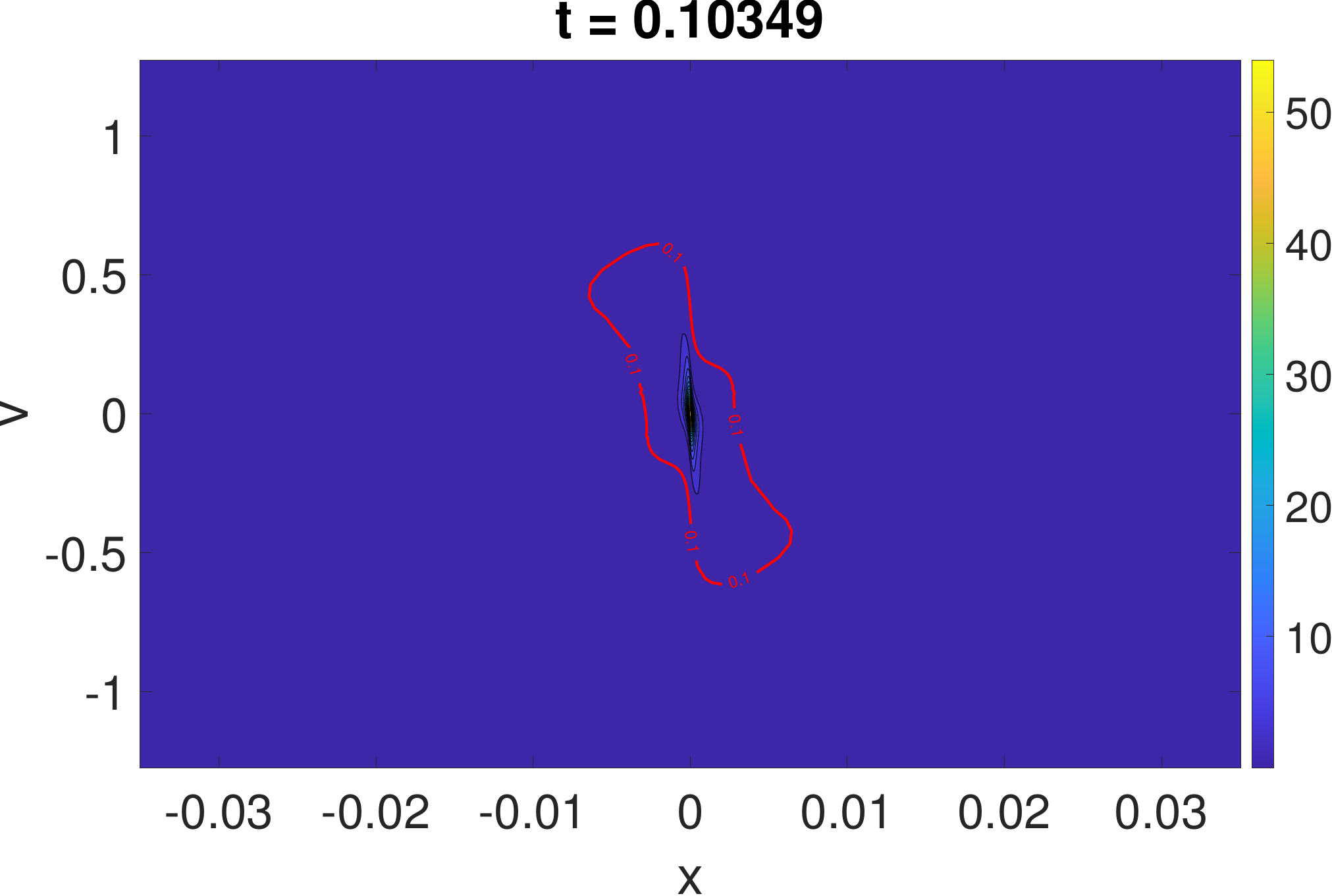}}
{\includegraphics[width=0.41\textwidth]{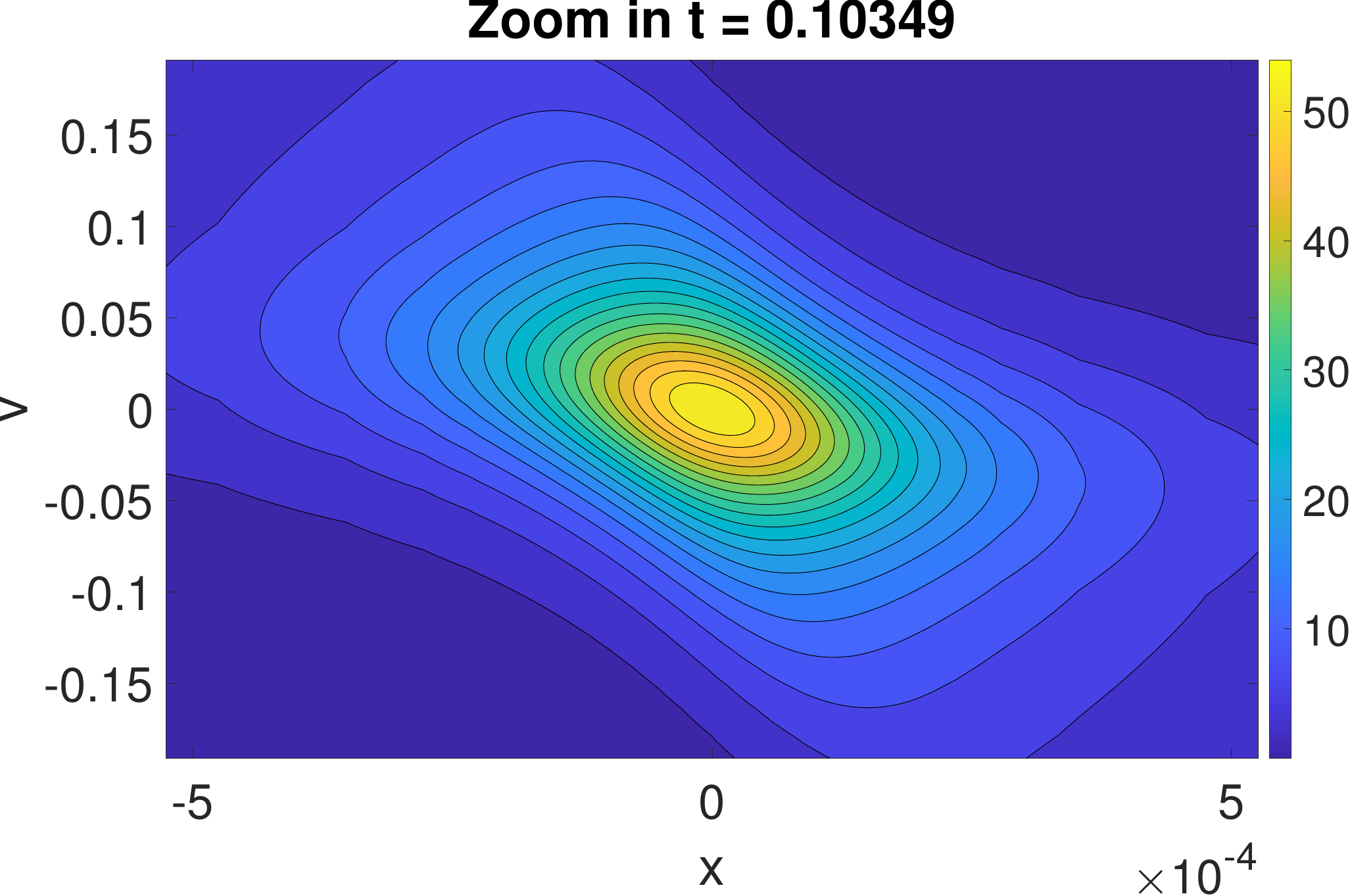}}

\caption{From top to bottom are numerical solutions of \eqref{eqn: gke} with kernel $W(v)=|v|^2$ and their zoom-in plots at the numerical blowup time for initial condition~\ref{eqn:ic_g2} with $\lambda = 260$ and $220$. Here we use $N_x=121$,  $N_v=201$, $\delta_0 = 0.5$ and $\delta = 0.5$ for both the $x$ and $v$.}\label{fig:gk2_bu_2}
\end{figure}

In this part, we present supplementary examples to numerically check the blow-up condition described in Theorem~\ref{thm_condbu}. In Figure~\ref{fig:gk2_bu_2}, we present supplementary results for $\gamma=2$ with initial condition \eqref{eqn:ic_g2} using different values of $\lambda$, and the snapshots of evolution with $\lambda=310$ are presented in Figure~\ref{fig:gk2_evo_example2}. In these two figures, the level sets at $f(t, x, v)=0.1$ intersect any vertical line with only one connected component, which verify the blow-up condition.

\begin{figure}[h!]
\centering
{\includegraphics[width=0.41\textwidth]{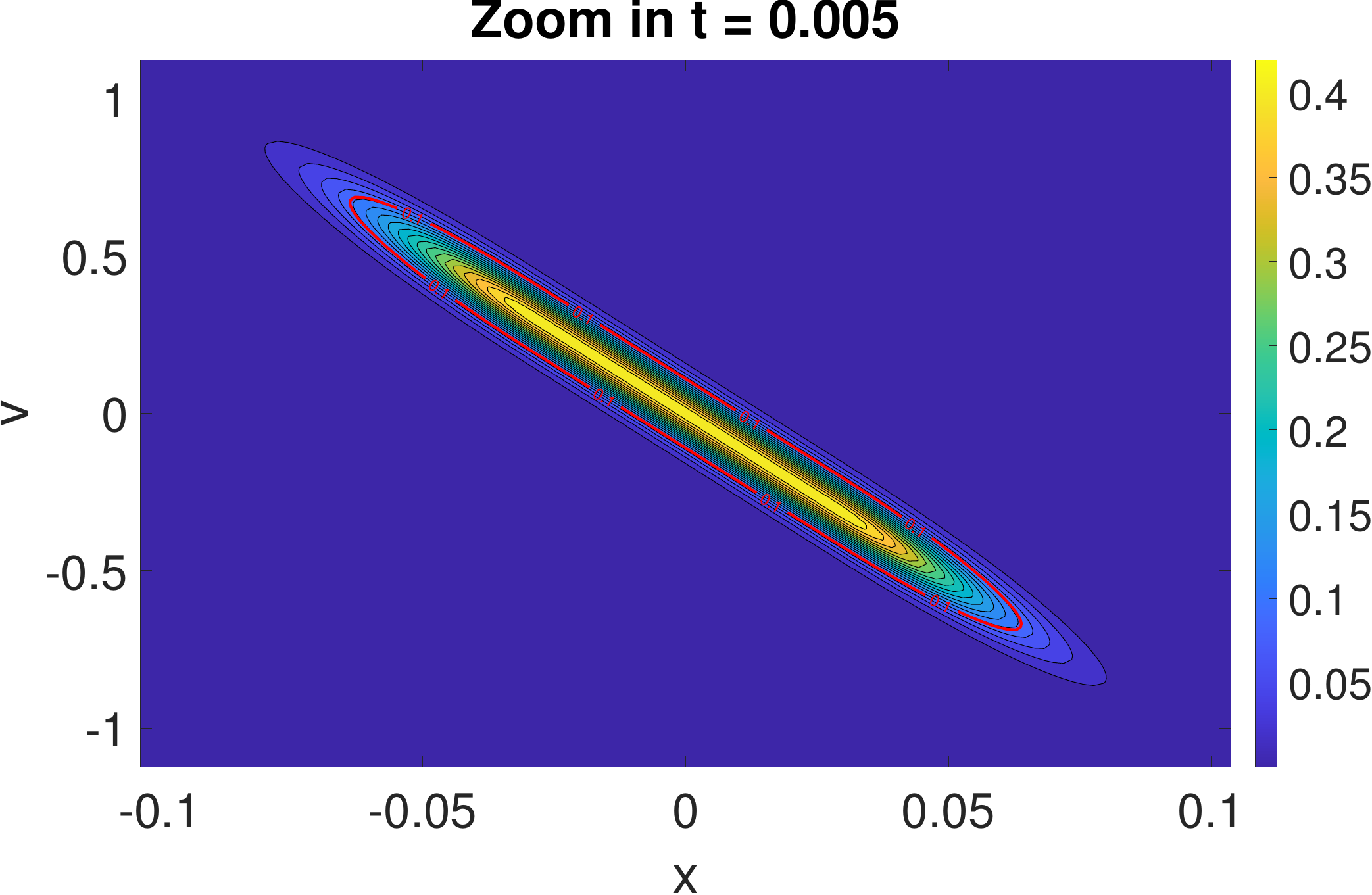}}
{\includegraphics[width=0.41\textwidth]{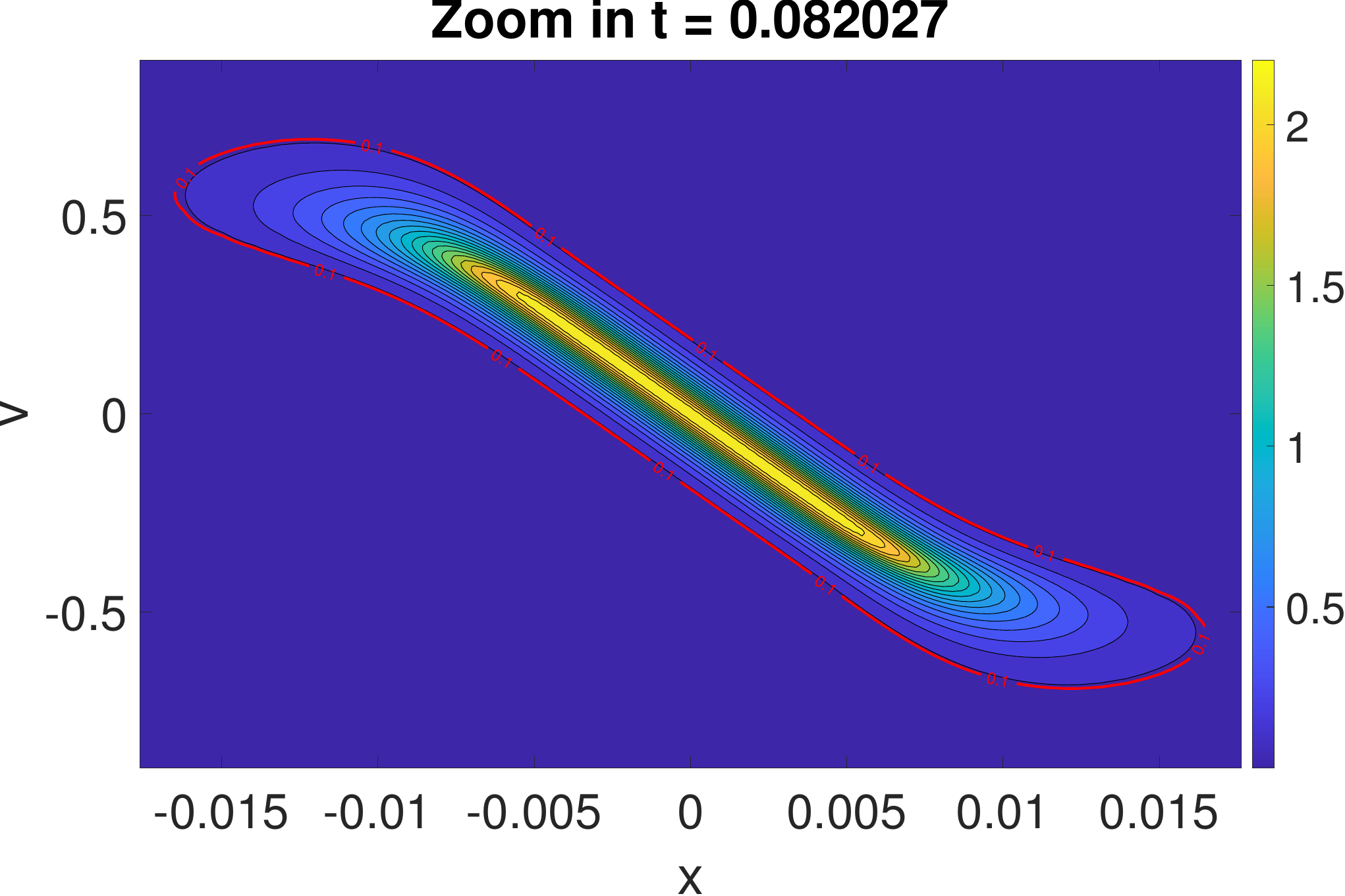}}

{\includegraphics[width=0.41\textwidth]{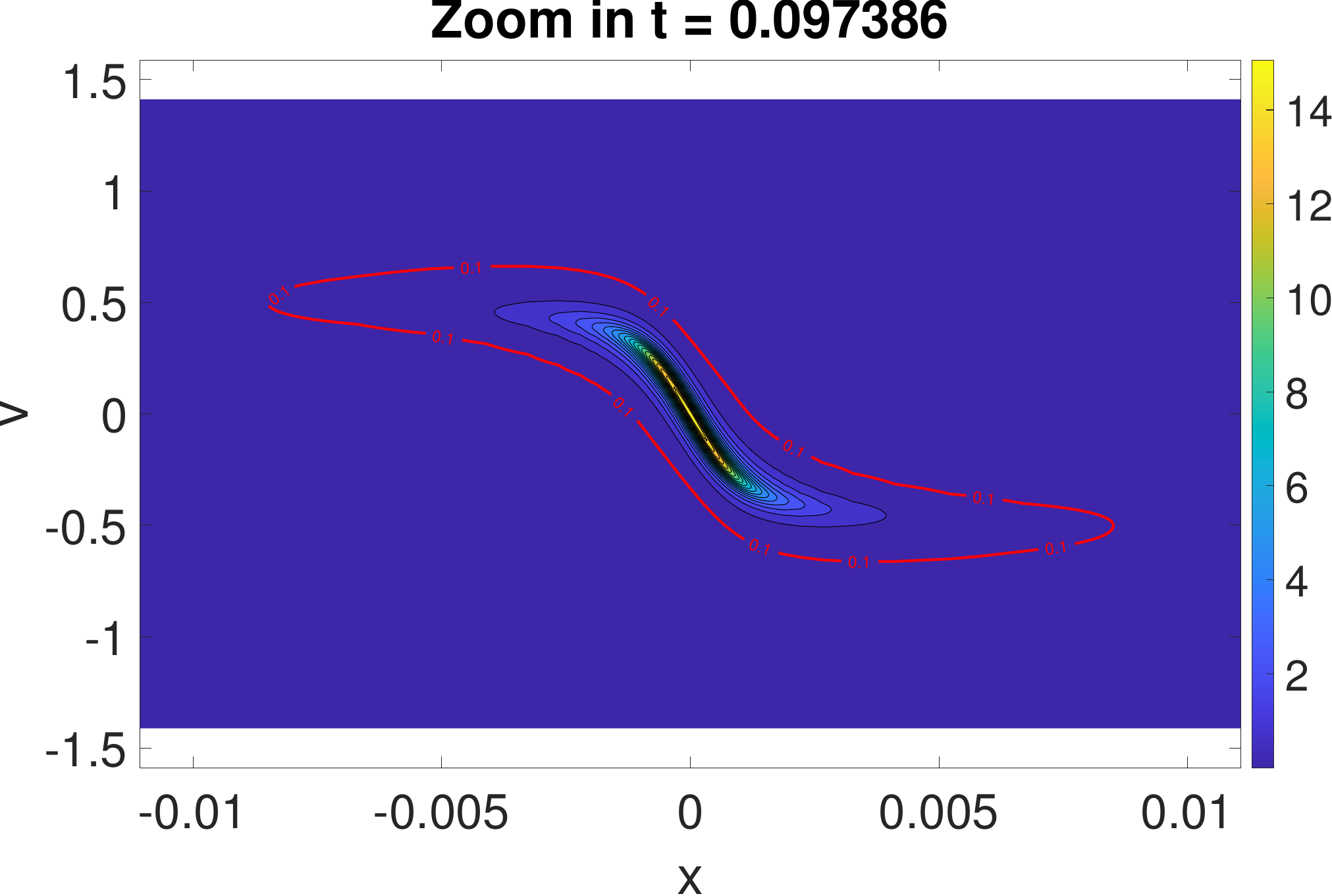}}
{\includegraphics[width=0.41\textwidth]{GK_fig4_contour_CC_155_IC_22-eps-converted-to}}
\caption{Snapshots of the contour plots for the evolution of the numerical solution to \eqref{eqn: gke} with kernel $W(v)=|v|^2$ for $N_x=121$, $N_v=201$, $\lambda = 310$.}\label{fig:gk2_evo_example}
\end{figure}

\section*{Acknowledgements}
JAC and RS were supported by the Advanced Grant Nonlocal-CPD (Nonlocal PDEs for Complex Particle Dynamics: Phase Transitions, Patterns and Synchronization) of the European Research Council Executive Agency (ERC) under the European Union’s Horizon 2020 research and innovation programme (grant agreement No. 883363). JAC was also partially supported by the Engineering and Physical Sciences Research Council (EPSRC) under grants EP/T022132/1 and EP/V051121/1. LW and WX are partially supported by NSF grant DMS-1846854. We also acknowledge the generous support from the Simons Foundation for the authors to participate in the program ``Frontiers in kinetic theory: connecting microscopic to macroscopic scale" held at the Isaac Newton Institute (UK) in Spring 2022, where the part of the work was finished.  

\newpage

\bibliographystyle{abbrv}
\bibliography{reference}

\begin{thebibliography}{10}

\bibitem{agueh2016local}
M.~Agueh.
\newblock Local existence of weak solutions to kinetic models of granular
  media.
\newblock {\em Archive for Rational Mechanics and Analysis}, 221(2):917--959,
  2016.

\bibitem{AC16}
M.~Agueh and G.~Carlier.
\newblock Generalized solutions of a kinetic granular media equation by a
  gradient flow approach.
\newblock {\em Calc. Var. Partial Differential Equations}, 55(2):Art. 37, 26,
  2016.

\bibitem{ACI15}
M.~Agueh, G.~Carlier, and R.~Illner.
\newblock Remarks on a class of kinetic models of granular media: asymptotics
  and entropy bounds.
\newblock {\em Kinet. Relat. Models}, 8(2):201--214, 2015.

\bibitem{benedetto1998non}
D.~Benedetto, E.~Caglioti, J.~A. Carrillo, and M.~Pulvirenti.
\newblock A non-maxwellian steady distribution for one-dimensional granular
  media.
\newblock {\em Journal of statistical physics}, 91(5):979--990, 1998.

\bibitem{benedetto1997kinetic}
D.~Benedetto, E.~Caglioti, and M.~Pulvirenti.
\newblock A kinetic equation for granular media.
\newblock {\em ESAIM: Mathematical Modelling and Numerical Analysis},
  31(5):615--641, 1997.

\bibitem{benedetto1997one}
D.~Benedetto, E.~Caglioti, and M.~Pulvirenti.
\newblock A one dimensional boltzmann equation with inelastic collisions.
\newblock {\em Rendiconti del Seminario Matematico e Fisico di Milano},
  67(1):169--179, 1997.

\bibitem{bertozzi2009blow}
A.~L. Bertozzi, J.~A. Carrillo, and T.~Laurent.
\newblock Blow-up in multidimensional aggregation equations with mildly
  singular interaction kernels.
\newblock {\em Nonlinearity}, 22(3):683, 2009.

\bibitem{BCDP15}
G.~A. Bonaschi, J.~A. Carrillo, M.~Di~Francesco, and M.~A. Peletier.
\newblock Equivalence of gradient flows and entropy solutions for singular
  nonlocal interaction equations in 1{D}.
\newblock {\em ESAIM Control Optim. Calc. Var.}, 21(2):414--441, 2015.

\bibitem{Bougieetal04}
J.~Bougie, J.~Kreft, J.~B. Swift, and H.~Swinney.
\newblock Onset of patterns in an oscillated granular layer: continuum and
  molecular dynamics simulations.
\newblock {\em Phys. Rev. E}, 71:021301, 2005.

\bibitem{BMSS02}
J.~Bougie, S.~J. Moon, J.~Swift, and H.~Swinney.
\newblock Shocks in vertically oscillated granular layers.
\newblock {\em Phys. Rev. E}, 66:051301, 2002.

\bibitem{brilliantov:2004}
N.~V. Brilliantov and T.~P\"{o}schel.
\newblock {\em {Kinetic {T}heory of {G}ranular {G}ases}}.
\newblock Oxford University Press, USA, 2004.

\bibitem{BSSP04}
N.~V. Brilliantov, C.~Salue\~na, T.~Schwager, and T.~P\"oschel.
\newblock Transient structures in a granular gas.
\newblock {\em Phys. Rev. Lett.}, 93:134301, 2004.

\bibitem{caglioti2002homogeneous}
E.~Caglioti and C.~Villani.
\newblock Homogeneous cooling states are not always good approximations to
  granular flows.
\newblock {\em Archive for rational mechanics and analysis}, 163(4):329--343,
  2002.

\bibitem{CHMRReview}
J.~A. Carrillo, J.~Hu, Z.~Ma, and T.~Rey.
\newblock Recent development in kinetic theory of granular materials: analysis
  and numerical methods.
\newblock In {\em Trails in kinetic theory---foundational aspects and numerical
  methods}, volume~25 of {\em SEMA SIMAI Springer Ser.}, pages 1--36. Springer,
  Cham, [2021] \copyright 2021.

\bibitem{carrillo2003kinetic}
J.~A. Carrillo, R.~J. McCann, and C.~Villani.
\newblock Kinetic equilibration rates for granular media and related equations:
  entropy dissipation and mass transportation estimates.
\newblock {\em Revista Matematica Iberoamericana}, 19(3):971--1018, 2003.

\bibitem{carrillo:2008granular}
J.~A. Carrillo, T.~Po\"{e}schel, and C.~Salue\~{n}a.
\newblock {Granular hydrodynamics and pattern formation in vertically
  oscillated granular disk layers}.
\newblock {\em J. Fluid Mech.}, 597:119--144, 2008.

\bibitem{DLK95}
Y.~Du, H.~Li, and L.~P. Kadanoff.
\newblock Breakdown of hydrodynamics in a one-dimensional system of inelastic
  particles.
\newblock {\em Phys. Rev. Lett.}, 74:1268--1271, Feb 1995.

\bibitem{gamba2004boltzmann}
I.~M. Gamba, V.~Panferov, and C.~Villani.
\newblock On the boltzmann equation for diffusively excited granular media.
\newblock {\em Communications in Mathematical Physics}, 246(3):503--541, 2004.

\bibitem{Goldhirsch:2003}
I.~Goldhirsch.
\newblock Rapid granular flows.
\newblock {\em Ann. Rev. Fluid Mech.}, 35:267, 2003.

\bibitem{HM03}
S.~A. Hill and G.~F. Mazenko.
\newblock Granular clustering in a hydrodynamic simulation.
\newblock {\em Phys. Rev. E}, 67:061302, 2003.

\bibitem{huang2012asymptotics}
Y.~Huang and A.~Bertozzi.
\newblock Asymptotics of blowup solutions for the aggregation equation.
\newblock {\em Discrete \& Continuous Dynamical Systems-B}, 17(4):1309, 2012.

\bibitem{JR85}
J.~Jenkins and M.~W. Richman.
\newblock Grad's $13$-moment system for a dense gas of inelastic spheres.
\newblock {\em Arch. Rational Mech. Anal.}, 87:355, 1985.

\bibitem{JR852}
J.~Jenkins and M.~W. Richman.
\newblock Kinetic theory for plane flows of a dense gas of identical, rough,
  inelastic, circular disks.
\newblock {\em Phys. Fluids}, 28:3485, 1985.

\bibitem{Li:2004}
H.~Li and G.~Toscani.
\newblock {Long-Time Asymptotics of Kinetic Models of Granular Flows}.
\newblock {\em Arch. Ration. Mech. Anal.}, 172(3):407--428, May 2004.

\bibitem{li2020fisher}
W.~Li, J.~Lu, and L.~Wang.
\newblock Fisher information regularization schemes for wasserstein gradient
  flows.
\newblock {\em Journal of Computational Physics}, 416:109449, 2020.

\bibitem{luo2014potentially}
G.~Luo and T.~Y. Hou.
\newblock Potentially singular solutions of the 3d axisymmetric euler
  equations.
\newblock {\em Proceedings of the National Academy of Sciences},
  111(36):12968--12973, 2014.

\bibitem{luo2014toward}
G.~Luo and T.~Y. Hou.
\newblock Toward the finite-time blowup of the 3d axisymmetric euler equations:
  a numerical investigation.
\newblock {\em Multiscale Modeling \& Simulation}, 12(4):1722--1776, 2014.

\bibitem{mcnamara1992inelastic}
S.~McNamara and W.~Young.
\newblock Inelastic collapse and clumping in a one-dimensional granular medium.
\newblock {\em Physics of Fluids A: Fluid Dynamics}, 4(3):496--504, 1992.

\bibitem{MY93}
S.~McNamara and W.~R. Young.
\newblock {Kinetics of a one-dimensional granular medium in the quasielastic
  limit}.
\newblock {\em Physics of Fluids A: Fluid Dynamics}, 5(1):34--45, 01 1993.

\bibitem{MeloUmbanhowarSwinney:1995}
F.~Melo, P.~Umbanhowar, and H.~L. Swinney.
\newblock Hexagons, kinks, and disorder in oscillated granular layers.
\newblock {\em Phys. Rev. Lett.}, 75:3838--3841, 1995.

\bibitem{RBSS02}
E.~Rericha, C.~Bizon, M.~Shattuck, and H.~Swinney.
\newblock Shocks in supersonic sand.
\newblock {\em Phys. Rev. Lett.}, 88:1, 2002.

\bibitem{toscani2000one}
G.~Toscani.
\newblock One-dimensional kinetic models of granular flows.
\newblock {\em ESAIM: Mathematical Modelling and Numerical Analysis},
  34(6):1277--1291, 2000.

\bibitem{TsimringAranson:1997}
L.~S. Tsimring and I.~S. Aranson.
\newblock Localized and cellular patterns in a vibrated granular layer.
\newblock {\em Phys. Rev. Lett.}, 79:213--216, 1997.

\bibitem{MeloUmbanhowarSwinney:1996}
P.~B. Umbanhowar, F.~Melo, and H.~L. Swinney.
\newblock Localized excitations in a vertically vibrated granular layer.
\newblock {\em Nature}, 382:793--796, 1996.

\end{thebibliography}

\end{document}